\documentclass{amsart}
\usepackage{a4wide,amssymb,color,hyperref}

\theoremstyle{plain}

\theoremstyle{definition}

\theoremstyle{remark}

\newtheorem{example}{Example}[section]

\title{Point-transitive Steiner systems S(2,6,111/121/126), S(2,7,169/175)}
\author{Ivan Hetman}
\date{April 2025}

\begin{document}

\begin{abstract} In this paper new Steiner systems $S(2,6,111)$, $S(2,6,121)$, $S(2,6,126)$, $S(2,7,169)$, $S(2,7,175)$ and possibly others with point-transitive (commutative except $S(2,6,111)$ case) automorphism groups are introduced.
\end{abstract}
\maketitle

\section{Introduction}
In \cite{Het} results of an algorithm for finding cyclic difference families (CDF) were presented. Since that time lots of generalizations were made. Some were made by author himself, some with Taras Banakh and Alex Ravsky elaborated more general alrogithm which will be published later. This paper is a collection of results produced by one of the first generalizations of algorithm for CDF search.

In \cite{HoCD} there is collection of known results which are quite sporadic. It seems that lots of authors found one design for admissible $v,k$ and stopped any attempts to find more. More systematic approach is present on Vedran Krcadinac website. I assume that the author included results that are related to his work. Still, I was not able to find other sources for Steiner designs, so I'm assuming that both Handbook and Krcadinac website contain most recent data regarding Steiner systems.

There is also website by Dávid Mezőfi, Gábor Nagy who developed paramodification method. This method allows obtaining new Steiner designs from old ones. They have large collection of $S(2,4,28)$ and $S(2,5,65)$ which are unitals. Some designs on Krcadinac website are obtained by paramodifications as well.

In future paper we will give detailed description of algorithm which produces examples published below. But this paper is one of placeholder to publish new Steiner systems. Here I publish results that were obtained by generalization of CDF algorithm to commutative group and one result obtained in first naive attempt to generalize CDF algorithm to non-commutative groups.

\section{Important properties}
Firstly I need to describe fingerprint. Due to the fact that it's complicated to integrate nauty with Java, using Practical graph isomorphism \cite{McKay} paper, I implemented naive copy of nauty in Java. It works well on Steiner systems with $\frac{v(v-1)}{k(k-1)}\le 1000$, but then running time grows significantly. Assuming that I need to check whether generated designs are new, a design fingerprint was implemented. This fingerprint has name {\bf hyperbolic frequency} because basic motivation that led to this construction was search of "hyperbolic" designs \cite{MOQ}.

The hyperbolic frequency is calculated as follows. We go through all non-collinear (not on same block) triples $oxy$ and point $p$ on the line (block) $xy$ which is distinct from $x$ and $y$ (so we check all such quadruples). Then we take number of points $u$ on the line $oy$ such that line $pu$ doesn't intersect $ox$. Then we sum frequencies and obtain something like $\{1=31968, 2=534132, 3=2203128, 4=2358972\}$. "Hyperbolic" in this case means that frequencies don't contain $0$ and $1$ as keys. Interesting property is that frequency with $0$ key only means that design corresponds to projective space. $1$ key only means that we have either affine space or some subset of Steiner triple systems. Single key $x\ge 2$ corresponds to classical unitals. For more details see mathoverflow answers \cite{MOQ} where lots of known and two supposedly unknown earlier "hyperbolic" designs are presented.

Hyperbolic frequency fingerprint works very well for Steiner designs with $k\ge 5$. If we don't face affine/projective design, when we have same hyperbolic freequency, then almost sure that two designs are isomorphic which then is checked by nauty reimplementation. When frequencies are different, obviously that designs are non-isomorphic and additional check is not needed.

Designs will be presented in form of difference family. To know more about the properties of difference families we refer a reader to VI.16 chapter of \cite{HoCD}. Also one more important property of difference family is the automorphism multiplier. If we act with automorphism group of generator group on difference family, it is the number of automorphisms that fix this difference family. For instance, the known design $S(2,7,175)$ has multiplier $24$ and that means that its automorphism group is a multiple of $24*175=4200$ and that is confirmed in \cite{Krc}.

\section{Commutative groups results}
Results will come in form multiplier-fingerprint-difference family. All this results are enumerations, so they give exhaustive list of difference families that generate non-isomorphic Steiner systems. Known (where I know reference) are marked bold and citation is added in the end. If you know where this designs were introduced earlier or consider this results so trivial that they were not worth mentioning - please drop me an email and I will add a correct citation.

Let's reconfirm (supposedly?) known small results for $k\in\{3,4\}$. I will not mark them bold because very probably they were obtained by someone else. Cyclic groups are avoided because there is lots of enumerations of them like \cite{BaiTop} including also our joint future paper \cite{HetBaiTop}.

\begin{example} $\mathbb Z_5 \times \mathbb Z_5$, $v=25$, $k=3$.
\begin{enumerate}
    \item 3 \{0=1200, 1=12000\} [[(0, 0), (0, 1), (1, 0)], [(0, 0), (0, 2), (2, 0)], [(0, 0), (1, 1), (2, 4)], [(0, 0), (1, 2), (3, 3)]]
    \item 3 \{0=1200, 1=12000\} [[(0, 0), (0, 1), (1, 0)], [(0, 0), (0, 2), (2, 1)], [(0, 0), (1, 1), (2, 3)], [(0, 0), (1, 3), (3, 0)]]
    \item 3 \{1=13200\} [[(0, 0), (0, 1), (1, 0)], [(0, 0), (0, 2), (2, 1)], [(0, 0), (1, 1), (2, 3)], [(0, 0), (1, 3), (3, 3)]]
\end{enumerate}
\end{example}

Two difference families with same fingerprint are non-isomorphic.

\begin{example} There is 3862 multiplier-nonisomorphic designs for $\mathbb Z_7 \times \mathbb Z_7$, $v=49$, $k=3$. Their multiplier distribution is $\{1=3442, 2=100, 3=286, 6=24, 9=5, 18=4, 63=1\}$. Few interesting examples with big multiplier (and therefore large automorphism group) are given below.
\begin{enumerate}
    \item 63 \{0=9408, 1=98784\} [[(0, 0), (0, 1), (0, 3)], [(0, 0), (1, 0), (3, 0)], [(0, 0), (1, 1), (3, 3)], [(0, 0), (1, 2), (3, 6)], [(0, 0), (1, 3), (3, 2)], [(0, 0), (1, 4), (3, 5)], [(0, 0), (1, 5), (3, 1)], [(0, 0), (1, 6), (3, 4)]]
    \item 18 \{0=9408, 1=98784\} [[(0, 0), (0, 1), (0, 3)], [(0, 0), (1, 0), (2, 3)], [(0, 0), (1, 1), (3, 3)], [(0, 0), (1, 2), (5, 3)], [(0, 0), (1, 4), (3, 5)], [(0, 0), (1, 5), (4, 0)], [(0, 0), (1, 6), (5, 2)], [(0, 0), (2, 0), (4, 6)]]
    \item 18 \{0=9408, 1=98784\} [[(0, 0), (0, 1), (0, 3)], [(0, 0), (1, 0), (2, 3)], [(0, 0), (1, 1), (3, 3)], [(0, 0), (1, 2), (5, 3)], [(0, 0), (1, 4), (5, 6)], [(0, 0), (1, 5), (4, 0)], [(0, 0), (1, 6), (3, 4)], [(0, 0), (2, 0), (4, 6)]]
    \item 18 \{0=44688, 1=63504\} [[(0, 0), (0, 1), (0, 3)], [(0, 0), (1, 0), (2, 1)], [(0, 0), (1, 2), (3, 6)], [(0, 0), (1, 3), (3, 2)], [(0, 0), (1, 4), (4, 0)], [(0, 0), (1, 5), (5, 4)], [(0, 0), (1, 6), (5, 2)], [(0, 0), (2, 0), (4, 2)]]
    \item 18 \{0=44688, 1=63504\} [[(0, 0), (0, 1), (0, 3)], [(0, 0), (1, 0), (2, 1)], [(0, 0), (1, 2), (5, 3)], [(0, 0), (1, 3), (3, 2)], [(0, 0), (1, 4), (4, 0)], [(0, 0), (1, 5), (5, 4)], [(0, 0), (1, 6), (3, 4)], [(0, 0), (2, 0), (4, 2)]]
    \item 9 \{0=9408, 1=98784\} [[(0, 0), (0, 1), (0, 3)], [(0, 0), (1, 0), (3, 0)], [(0, 0), (1, 1), (3, 3)], [(0, 0), (1, 2), (3, 6)], [(0, 0), (1, 3), (3, 2)], [(0, 0), (1, 4), (3, 5)], [(0, 0), (1, 5), (3, 1)], [(0, 0), (1, 6), (5, 2)]]
    \item 9 \{0=9408, 1=98784\} [[(0, 0), (0, 1), (0, 3)], [(0, 0), (1, 0), (3, 0)], [(0, 0), (1, 1), (3, 5)], [(0, 0), (1, 2), (3, 3)], [(0, 0), (1, 3), (3, 1)], [(0, 0), (1, 4), (3, 6)], [(0, 0), (1, 5), (3, 4)], [(0, 0), (1, 6), (3, 2)]]
    \item 9 \{0=9408, 1=98784\} [[(0, 0), (0, 1), (0, 3)], [(0, 0), (1, 0), (3, 0)], [(0, 0), (1, 1), (3, 5)], [(0, 0), (1, 2), (3, 3)], [(0, 0), (1, 3), (5, 2)], [(0, 0), (1, 4), (3, 6)], [(0, 0), (1, 5), (5, 1)], [(0, 0), (1, 6), (5, 4)]]
    \item 9 \{0=9408, 1=98784\} [[(0, 0), (0, 1), (0, 3)], [(0, 0), (1, 0), (2, 3)], [(0, 0), (1, 1), (3, 3)], [(0, 0), (1, 2), (5, 3)], [(0, 0), (1, 4), (3, 5)], [(0, 0), (1, 5), (4, 0)], [(0, 0), (1, 6), (3, 4)], [(0, 0), (2, 0), (4, 6)]]
    \item 9 \{0=44688, 1=63504\} [[(0, 0), (0, 1), (0, 3)], [(0, 0), (1, 0), (2, 1)], [(0, 0), (1, 2), (3, 6)], [(0, 0), (1, 3), (3, 2)], [(0, 0), (1, 4), (4, 0)], [(0, 0), (1, 5), (5, 4)], [(0, 0), (1, 6), (3, 4)], [(0, 0), (2, 0), (4, 2)]]
\end{enumerate}
\end{example}

I didn't check for nauty isomorphism of Steiner systems above, so some of them could be isomorphic (which is proved to be possible).

\begin{example} $\mathbb Z_5 \times \mathbb Z_5$, $v=25$, $k=4$.
\begin{enumerate}
    \item 6 \{0=900, 1=12000, 2=12300\} [[(0, 0), (0, 1), (1, 0), (2, 2)], [(0, 0), (0, 2), (1, 3), (3, 2)]]
\end{enumerate}
\end{example}

\begin{example} $\mathbb Z_7 \times \mathbb Z_7$, $v=49$, $k=4$.
\begin{enumerate}
    \item 3 \{0=2646, 1=40572, 2=168462\} [[(0, 0), (0, 1), (1, 0), (3, 6)], [(0, 0), (0, 2), (2, 3), (3, 4)], [(0, 0), (0, 3), (1, 5), (5, 5)], [(0, 0), (1, 3), (2, 0), (4, 4)]]
    \item 3 \{0=1764, 1=36456, 2=173460\} [[(0, 0), (0, 1), (1, 0), (3, 6)], [(0, 0), (0, 2), (2, 3), (3, 4)], [(0, 0), (0, 3), (1, 5), (5, 5)], [(0, 0), (1, 3), (4, 6), (6, 3)]]
    \item 3 \{0=1764, 1=37632, 2=172284\} [[(0, 0), (0, 1), (1, 0), (2, 2)], [(0, 0), (0, 2), (1, 3), (4, 1)], [(0, 0), (0, 3), (2, 0), (5, 4)], [(0, 0), (1, 4), (3, 2), (4, 0)]]
    \item 1 \{0=1764, 1=39984, 2=169932\} [[(0, 0), (0, 1), (1, 0), (2, 2)], [(0, 0), (0, 2), (1, 3), (4, 1)], [(0, 0), (0, 3), (2, 0), (5, 4)], [(0, 0), (1, 4), (4, 4), (5, 2)]]
    \item 3 \{0=882, 1=34692, 2=176106\} [[(0, 0), (0, 1), (1, 0), (2, 2)], [(0, 0), (0, 2), (1, 3), (4, 1)], [(0, 0), (0, 3), (2, 6), (5, 3)], [(0, 0), (1, 4), (3, 2), (4, 0)]]
    \item 1 \{0=1176, 1=37632, 2=172872\} [[(0, 0), (0, 1), (1, 0), (2, 2)], [(0, 0), (0, 2), (1, 3), (4, 1)], [(0, 0), (0, 3), (2, 6), (5, 3)], [(0, 0), (1, 4), (4, 4), (5, 2)]]
    \item 1 \{0=2058, 1=44100, 2=165522\} [[(0, 0), (0, 1), (1, 0), (2, 2)], [(0, 0), (0, 2), (1, 3), (5, 4)], [(0, 0), (0, 3), (3, 0), (4, 5)], [(0, 0), (1, 4), (3, 1), (5, 1)]]
    \item 1 \{0=882, 1=37044, 2=173754\} [[(0, 0), (0, 1), (1, 0), (2, 2)], [(0, 0), (0, 2), (1, 3), (5, 4)], [(0, 0), (0, 3), (3, 0), (4, 5)], [(0, 0), (1, 4), (3, 3), (5, 3)]]
    \item 3 \{0=10584, 1=76440, 2=124656\} [[(0, 0), (0, 1), (1, 0), (2, 2)], [(0, 0), (0, 2), (1, 3), (5, 4)], [(0, 0), (0, 3), (3, 5), (4, 3)], [(0, 0), (1, 4), (3, 3), (5, 3)]]
    \item 3 \{0=6174, 1=55860, 2=149646\} [[(0, 0), (0, 1), (1, 0), (2, 2)], [(0, 0), (0, 2), (1, 5), (4, 4)], [(0, 0), (0, 3), (2, 6), (4, 3)], [(0, 0), (1, 1), (3, 1), (6, 3)]]
    \item 1 \{0=5292, 1=48216, 2=158172\} [[(0, 0), (0, 1), (1, 0), (2, 2)], [(0, 0), (0, 2), (1, 5), (5, 2)], [(0, 0), (0, 3), (3, 2), (4, 6)], [(0, 0), (1, 1), (3, 0), (5, 3)]]
    \item 1 \{0=2352, 1=42336, 2=166992\} [[(0, 0), (0, 1), (1, 0), (2, 2)], [(0, 0), (0, 2), (1, 5), (5, 2)], [(0, 0), (0, 3), (3, 4), (4, 1)], [(0, 0), (1, 1), (3, 5), (5, 1)]]
\end{enumerate}
\end{example}

\begin{example} $\mathbb Z_3 \times \mathbb Z_3 \times \mathbb Z_5$, $v=45$, $k=5$.
\begin{enumerate}
    \item {\bf 8 \{0=3600, 1=32400, 2=125280, 3=76320\} [[(0, 0, 0), (0, 1, 0), (0, 2, 1), (1, 0, 2), (2, 1, 2)], [(0, 0, 0), (0, 1, 2), (1, 0, 1), (2, 0, 0), (2, 2, 2)], [(0, 0, 0), (0, 0, 1), (0, 0, 2), (0, 0, 3), (0, 0, 4)]]} \cite{HoCD}
\end{enumerate}
\end{example}

Now let's start with commutative groups of order $81$ and $k=5$. Fortunately, for $k\ge 5$ all fingerprints are different that mean that we have non-isomorphic Steiner systems. From here on, I will mark known results with bold.

\begin{example} $\mathbb Z_3 \times \mathbb Z_3 \times \mathbb Z_3 \times \mathbb Z_3$, $v=81$, $k=5$.
\begin{enumerate}
    \item 5 \{0=3240, 1=24300, 2=427680, 3=1022220\} [[(0, 0, 0, 0), (0, 0, 0, 1), (0, 0, 1, 0), (0, 1, 0, 0), (1, 0, 0, 0)], [(0, 0, 0, 0), (0, 0, 1, 1), (0, 1, 1, 2), (1, 2, 1, 0), (2, 0, 2, 0)], [(0, 0, 0, 0), (0, 1, 1, 0), (1, 1, 1, 1), (1, 2, 2, 2), (2, 0, 2, 2)], [(0, 0, 0, 0), (0, 1, 2, 1), (1, 0, 1, 2), (1, 2, 2, 0), (2, 1, 0, 2)]]
\end{enumerate}
\end{example}

\begin{example} $\mathbb Z_3 \times \mathbb Z_{27}$, $v=81$, $k=5$.
\begin{enumerate}
    \item 1 \{1=42768, 2=476280, 3=958392\} [[(0, 0), (0, 1), (0, 3), (0, 7), (1, 8)], [(0, 0), (0, 5), (1, 15), (1, 26), (2, 8)], [(0, 0), (0, 8), (0, 17), (1, 3), (2, 21)], [(0, 0), (0, 12), (1, 2), (1, 16), (2, 0)]]
    \item 1 \{0=648, 1=45198, 2=461700, 3=969894\} [[(0, 0), (0, 1), (0, 3), (0, 7), (1, 8)], [(0, 0), (0, 5), (1, 15), (1, 26), (2, 8)], [(0, 0), (0, 8), (0, 17), (1, 3), (2, 21)], [(0, 0), (0, 12), (1, 12), (2, 10), (2, 23)]]
    \item 1 \{1=50058, 2=479196, 3=948186\} [[(0, 0), (0, 1), (0, 3), (0, 7), (1, 8)], [(0, 0), (0, 5), (1, 15), (1, 26), (2, 8)], [(0, 0), (0, 8), (0, 18), (1, 14), (2, 5)], [(0, 0), (0, 12), (1, 2), (1, 16), (2, 0)]]
    \item 1 \{1=34506, 2=471420, 3=971514\} [[(0, 0), (0, 1), (0, 3), (0, 7), (1, 8)], [(0, 0), (0, 5), (1, 15), (1, 26), (2, 8)], [(0, 0), (0, 8), (0, 18), (1, 14), (2, 5)], [(0, 0), (0, 12), (1, 12), (2, 10), (2, 23)]]
    \item 1 \{0=1296, 1=36450, 2=463644, 3=976050\} [[(0, 0), (0, 1), (0, 3), (0, 7), (1, 8)], [(0, 0), (0, 5), (1, 24), (2, 6), (2, 17)], [(0, 0), (0, 8), (0, 17), (1, 3), (2, 21)], [(0, 0), (0, 12), (1, 2), (1, 16), (2, 0)]]
    \item 1 \{1=37908, 2=452952, 3=986580\} [[(0, 0), (0, 1), (0, 3), (0, 7), (1, 8)], [(0, 0), (0, 5), (1, 24), (2, 6), (2, 17)], [(0, 0), (0, 8), (0, 17), (1, 3), (2, 21)], [(0, 0), (0, 12), (1, 12), (2, 10), (2, 23)]]
    \item 1 \{0=648, 1=47628, 2=474336, 3=954828\} [[(0, 0), (0, 1), (0, 3), (0, 7), (1, 8)], [(0, 0), (0, 5), (1, 24), (2, 6), (2, 17)], [(0, 0), (0, 8), (0, 18), (1, 14), (2, 5)], [(0, 0), (0, 12), (1, 2), (1, 16), (2, 0)]]
    \item 1 \{0=648, 1=34020, 2=456840, 3=985932\} [[(0, 0), (0, 1), (0, 3), (0, 7), (1, 8)], [(0, 0), (0, 5), (1, 24), (2, 6), (2, 17)], [(0, 0), (0, 8), (0, 18), (1, 14), (2, 5)], [(0, 0), (0, 12), (1, 12), (2, 10), (2, 23)]]
    \item 1 \{0=1296, 1=48114, 2=469476, 3=958554\} [[(0, 0), (0, 1), (0, 3), (0, 10), (1, 2)], [(0, 0), (0, 4), (1, 0), (2, 15), (2, 20)], [(0, 0), (0, 6), (0, 14), (1, 9), (2, 23)], [(0, 0), (0, 11), (1, 5), (1, 17), (2, 3)]]
    \item 1 \{0=648, 1=52974, 2=477252, 3=946566\} [[(0, 0), (0, 1), (0, 3), (0, 10), (1, 2)], [(0, 0), (0, 4), (1, 0), (2, 15), (2, 20)], [(0, 0), (0, 6), (0, 14), (1, 9), (2, 23)], [(0, 0), (0, 11), (1, 8), (2, 6), (2, 21)]]
    \item 1 \{0=648, 1=41796, 2=476280, 3=958716\} [[(0, 0), (0, 1), (0, 3), (0, 10), (1, 2)], [(0, 0), (0, 4), (1, 0), (2, 15), (2, 20)], [(0, 0), (0, 6), (0, 19), (1, 10), (2, 24)], [(0, 0), (0, 11), (1, 5), (1, 17), (2, 3)]]
    \item 1 \{0=1944, 1=41310, 2=479196, 3=954990\} [[(0, 0), (0, 1), (0, 3), (0, 10), (1, 2)], [(0, 0), (0, 4), (1, 0), (2, 15), (2, 20)], [(0, 0), (0, 6), (0, 19), (1, 10), (2, 24)], [(0, 0), (0, 11), (1, 8), (2, 6), (2, 21)]]
    \item 1 \{0=1296, 1=50544, 2=487944, 3=937656\} [[(0, 0), (0, 1), (0, 3), (0, 10), (1, 2)], [(0, 0), (0, 4), (1, 9), (2, 7), (2, 19)], [(0, 0), (0, 5), (1, 0), (1, 11), (2, 14)], [(0, 0), (0, 6), (0, 14), (1, 21), (2, 10)]]
    \item 1 \{1=55404, 2=484056, 3=937980\} [[(0, 0), (0, 1), (0, 3), (0, 10), (1, 2)], [(0, 0), (0, 4), (1, 9), (2, 7), (2, 19)], [(0, 0), (0, 5), (1, 0), (1, 11), (2, 14)], [(0, 0), (0, 6), (0, 19), (1, 23), (2, 12)]]
    \item 1 \{1=47142, 2=481140, 3=949158\} [[(0, 0), (0, 1), (0, 3), (0, 10), (1, 2)], [(0, 0), (0, 4), (1, 9), (2, 7), (2, 19)], [(0, 0), (0, 5), (1, 18), (2, 5), (2, 21)], [(0, 0), (0, 6), (0, 14), (1, 21), (2, 10)]]
    \item 1 \{1=34020, 2=466560, 3=976860\} [[(0, 0), (0, 1), (0, 3), (0, 10), (1, 2)], [(0, 0), (0, 4), (1, 9), (2, 7), (2, 19)], [(0, 0), (0, 5), (1, 18), (2, 5), (2, 21)], [(0, 0), (0, 6), (0, 19), (1, 23), (2, 12)]]
    \item 1 \{0=1296, 1=44226, 2=463644, 3=968274\} [[(0, 0), (0, 1), (0, 3), (0, 10), (1, 2)], [(0, 0), (0, 4), (1, 11), (1, 16), (2, 4)], [(0, 0), (0, 6), (0, 14), (1, 9), (2, 23)], [(0, 0), (0, 11), (1, 5), (1, 17), (2, 3)]]
    \item 1 \{1=36936, 2=478224, 3=962280\} [[(0, 0), (0, 1), (0, 3), (0, 10), (1, 2)], [(0, 0), (0, 4), (1, 11), (1, 16), (2, 4)], [(0, 0), (0, 6), (0, 14), (1, 9), (2, 23)], [(0, 0), (0, 11), (1, 8), (2, 6), (2, 21)]]
    \item 1 \{1=54432, 2=480168, 3=942840\} [[(0, 0), (0, 1), (0, 3), (0, 10), (1, 2)], [(0, 0), (0, 4), (1, 11), (1, 16), (2, 4)], [(0, 0), (0, 6), (0, 19), (1, 10), (2, 24)], [(0, 0), (0, 11), (1, 5), (1, 17), (2, 3)]]
    \item 1 \{1=58320, 2=478224, 3=940896\} [[(0, 0), (0, 1), (0, 3), (0, 10), (1, 2)], [(0, 0), (0, 4), (1, 11), (1, 16), (2, 4)], [(0, 0), (0, 6), (0, 19), (1, 10), (2, 24)], [(0, 0), (0, 11), (1, 8), (2, 6), (2, 21)]]
    \item 1 \{0=1944, 1=45684, 2=472392, 3=957420\} [[(0, 0), (0, 1), (0, 3), (0, 10), (1, 2)], [(0, 0), (0, 4), (1, 12), (1, 24), (2, 22)], [(0, 0), (0, 5), (1, 0), (1, 11), (2, 14)], [(0, 0), (0, 6), (0, 14), (1, 21), (2, 10)]]
    \item 1 \{0=648, 1=48600, 2=468504, 3=959688\} [[(0, 0), (0, 1), (0, 3), (0, 10), (1, 2)], [(0, 0), (0, 4), (1, 12), (1, 24), (2, 22)], [(0, 0), (0, 5), (1, 0), (1, 11), (2, 14)], [(0, 0), (0, 6), (0, 19), (1, 23), (2, 12)]]
    \item 1 \{1=62208, 2=491832, 3=923400\} [[(0, 0), (0, 1), (0, 3), (0, 10), (1, 2)], [(0, 0), (0, 4), (1, 12), (1, 24), (2, 22)], [(0, 0), (0, 5), (1, 18), (2, 5), (2, 21)], [(0, 0), (0, 6), (0, 14), (1, 21), (2, 10)]]
    \item 1 \{0=648, 1=56862, 2=502524, 3=917406\} [[(0, 0), (0, 1), (0, 3), (0, 10), (1, 2)], [(0, 0), (0, 4), (1, 12), (1, 24), (2, 22)], [(0, 0), (0, 5), (1, 18), (2, 5), (2, 21)], [(0, 0), (0, 6), (0, 19), (1, 23), (2, 12)]]
    \item 1 \{0=1944, 1=52488, 2=505440, 3=917568\} [[(0, 0), (0, 1), (0, 3), (0, 11), (1, 0)], [(0, 0), (0, 4), (1, 6), (2, 17), (2, 23)], [(0, 0), (0, 5), (0, 12), (1, 25), (2, 20)], [(0, 0), (0, 9), (1, 3), (2, 4), (2, 18)]]
    \item 1 \{0=2592, 1=61722, 2=481140, 3=931986\} [[(0, 0), (0, 1), (0, 3), (0, 11), (1, 0)], [(0, 0), (0, 4), (1, 6), (2, 17), (2, 23)], [(0, 0), (0, 5), (0, 12), (1, 25), (2, 20)], [(0, 0), (0, 9), (1, 5), (1, 18), (2, 6)]]
    \item 1 \{0=1296, 1=53460, 2=497664, 3=925020\} [[(0, 0), (0, 1), (0, 3), (0, 11), (1, 0)], [(0, 0), (0, 4), (1, 6), (2, 17), (2, 23)], [(0, 0), (0, 5), (0, 20), (1, 12), (2, 7)], [(0, 0), (0, 9), (1, 3), (2, 4), (2, 18)]]
    \item 1 \{1=49086, 2=492804, 3=935550\} [[(0, 0), (0, 1), (0, 3), (0, 11), (1, 0)], [(0, 0), (0, 4), (1, 6), (2, 17), (2, 23)], [(0, 0), (0, 5), (0, 20), (1, 12), (2, 7)], [(0, 0), (0, 9), (1, 5), (1, 18), (2, 6)]]
    \item 1 \{0=1296, 1=46656, 2=480168, 3=949320\} [[(0, 0), (0, 1), (0, 3), (0, 11), (1, 0)], [(0, 0), (0, 4), (1, 8), (1, 14), (2, 25)], [(0, 0), (0, 5), (0, 12), (1, 25), (2, 20)], [(0, 0), (0, 9), (1, 3), (2, 4), (2, 18)]]
    \item 1 \{0=1944, 1=64638, 2=500580, 3=910278\} [[(0, 0), (0, 1), (0, 3), (0, 11), (1, 0)], [(0, 0), (0, 4), (1, 8), (1, 14), (2, 25)], [(0, 0), (0, 5), (0, 12), (1, 25), (2, 20)], [(0, 0), (0, 9), (1, 5), (1, 18), (2, 6)]]
    \item 1 \{0=2592, 1=60264, 2=487944, 3=926640\} [[(0, 0), (0, 1), (0, 3), (0, 11), (1, 0)], [(0, 0), (0, 4), (1, 8), (1, 14), (2, 25)], [(0, 0), (0, 5), (0, 20), (1, 12), (2, 7)], [(0, 0), (0, 9), (1, 3), (2, 4), (2, 18)]]
    \item 1 \{0=648, 1=66096, 2=503496, 3=907200\} [[(0, 0), (0, 1), (0, 3), (0, 11), (1, 0)], [(0, 0), (0, 4), (1, 8), (1, 14), (2, 25)], [(0, 0), (0, 5), (0, 20), (1, 12), (2, 7)], [(0, 0), (0, 9), (1, 5), (1, 18), (2, 6)]]
    \item 1 \{0=648, 1=61722, 2=490860, 3=924210\} [[(0, 0), (0, 1), (0, 3), (0, 11), (1, 1)], [(0, 0), (0, 4), (1, 9), (2, 12), (2, 17)], [(0, 0), (0, 6), (0, 13), (1, 2), (2, 20)], [(0, 0), (0, 9), (1, 4), (2, 3), (2, 15)]]
    \item 1 \{0=1296, 1=65124, 2=505440, 3=905580\} [[(0, 0), (0, 1), (0, 3), (0, 11), (1, 1)], [(0, 0), (0, 4), (1, 9), (2, 12), (2, 17)], [(0, 0), (0, 6), (0, 13), (1, 2), (2, 20)], [(0, 0), (0, 9), (1, 6), (1, 21), (2, 5)]]
    \item 1 \{0=1944, 1=52974, 2=479196, 3=943326\} [[(0, 0), (0, 1), (0, 3), (0, 11), (1, 1)], [(0, 0), (0, 4), (1, 9), (2, 12), (2, 17)], [(0, 0), (0, 6), (0, 20), (1, 13), (2, 4)], [(0, 0), (0, 9), (1, 4), (2, 3), (2, 15)]]
    \item 1 \{0=1296, 1=51030, 2=502524, 3=922590\} [[(0, 0), (0, 1), (0, 3), (0, 11), (1, 1)], [(0, 0), (0, 4), (1, 9), (2, 12), (2, 17)], [(0, 0), (0, 6), (0, 20), (1, 13), (2, 4)], [(0, 0), (0, 9), (1, 6), (1, 21), (2, 5)]]
    \item 1 \{0=2592, 1=52488, 2=487944, 3=934416\} [[(0, 0), (0, 1), (0, 3), (0, 11), (1, 1)], [(0, 0), (0, 4), (1, 14), (1, 19), (2, 22)], [(0, 0), (0, 6), (0, 13), (1, 2), (2, 20)], [(0, 0), (0, 9), (1, 4), (2, 3), (2, 15)]]
    \item 1 \{0=648, 1=61722, 2=498636, 3=916434\} [[(0, 0), (0, 1), (0, 3), (0, 11), (1, 1)], [(0, 0), (0, 4), (1, 14), (1, 19), (2, 22)], [(0, 0), (0, 6), (0, 13), (1, 2), (2, 20)], [(0, 0), (0, 9), (1, 6), (1, 21), (2, 5)]]
    \item 1 \{0=2592, 1=59292, 2=493776, 3=921780\} [[(0, 0), (0, 1), (0, 3), (0, 11), (1, 1)], [(0, 0), (0, 4), (1, 14), (1, 19), (2, 22)], [(0, 0), (0, 6), (0, 20), (1, 13), (2, 4)], [(0, 0), (0, 9), (1, 4), (2, 3), (2, 15)]]
    \item 1 \{1=50544, 2=495720, 3=931176\} [[(0, 0), (0, 1), (0, 3), (0, 11), (1, 1)], [(0, 0), (0, 4), (1, 14), (1, 19), (2, 22)], [(0, 0), (0, 6), (0, 20), (1, 13), (2, 4)], [(0, 0), (0, 9), (1, 6), (1, 21), (2, 5)]]
    \item 1 \{0=648, 1=64638, 2=504468, 3=907686\} [[(0, 0), (0, 1), (0, 3), (0, 12), (1, 1)], [(0, 0), (0, 4), (1, 13), (2, 5), (2, 10)], [(0, 0), (0, 6), (0, 13), (1, 18), (2, 25)], [(0, 0), (0, 8), (1, 4), (1, 14), (2, 24)]]
    \item 1 \{0=648, 1=56862, 2=492804, 3=927126\} [[(0, 0), (0, 1), (0, 3), (0, 12), (1, 1)], [(0, 0), (0, 4), (1, 13), (2, 5), (2, 10)], [(0, 0), (0, 6), (0, 13), (1, 18), (2, 25)], [(0, 0), (0, 8), (1, 11), (2, 4), (2, 21)]]
    \item 1 \{1=60750, 2=457812, 3=958878\} [[(0, 0), (0, 1), (0, 3), (0, 12), (1, 1)], [(0, 0), (0, 4), (1, 13), (2, 5), (2, 10)], [(0, 0), (0, 6), (0, 20), (1, 8), (2, 15)], [(0, 0), (0, 8), (1, 4), (1, 14), (2, 24)]]
    \item 1 \{0=648, 1=55404, 2=501552, 3=919836\} [[(0, 0), (0, 1), (0, 3), (0, 12), (1, 1)], [(0, 0), (0, 4), (1, 13), (2, 5), (2, 10)], [(0, 0), (0, 6), (0, 20), (1, 8), (2, 15)], [(0, 0), (0, 8), (1, 11), (2, 4), (2, 21)]]
    \item 1 \{0=648, 1=52002, 2=463644, 3=961146\} [[(0, 0), (0, 1), (0, 3), (0, 12), (1, 1)], [(0, 0), (0, 4), (1, 14), (1, 24), (2, 8)], [(0, 0), (0, 5), (0, 19), (1, 22), (2, 1)], [(0, 0), (0, 6), (1, 8), (1, 15), (2, 20)]]
    \item 1 \{1=46656, 2=486000, 3=944784\} [[(0, 0), (0, 1), (0, 3), (0, 12), (1, 1)], [(0, 0), (0, 4), (1, 14), (1, 24), (2, 8)], [(0, 0), (0, 5), (0, 19), (1, 22), (2, 1)], [(0, 0), (0, 6), (1, 13), (2, 18), (2, 25)]]
    \item 1 \{0=648, 1=56376, 2=482112, 3=938304\} [[(0, 0), (0, 1), (0, 3), (0, 12), (1, 1)], [(0, 0), (0, 4), (1, 21), (1, 26), (2, 18)], [(0, 0), (0, 6), (0, 13), (1, 18), (2, 25)], [(0, 0), (0, 8), (1, 11), (2, 4), (2, 21)]]
    \item 1 \{0=648, 1=72900, 2=491832, 3=912060\} [[(0, 0), (0, 1), (0, 3), (0, 12), (1, 1)], [(0, 0), (0, 4), (1, 21), (1, 26), (2, 18)], [(0, 0), (0, 6), (0, 20), (1, 8), (2, 15)], [(0, 0), (0, 8), (1, 11), (2, 4), (2, 21)]]
    \item 1 \{0=648, 1=42282, 2=496692, 3=937818\} [[(0, 0), (0, 1), (0, 3), (0, 12), (1, 5)], [(0, 0), (0, 4), (0, 10), (1, 16), (2, 14)], [(0, 0), (0, 5), (1, 0), (1, 14), (2, 24)], [(0, 0), (0, 7), (1, 1), (2, 8), (2, 16)]]
    \item 1 \{0=648, 1=43254, 2=479196, 3=954342\} [[(0, 0), (0, 1), (0, 3), (0, 12), (1, 5)], [(0, 0), (0, 4), (0, 10), (1, 16), (2, 14)], [(0, 0), (0, 5), (1, 0), (1, 14), (2, 24)], [(0, 0), (0, 7), (1, 18), (1, 26), (2, 6)]]
    \item 1 \{1=51516, 2=474336, 3=951588\} [[(0, 0), (0, 1), (0, 3), (0, 12), (1, 5)], [(0, 0), (0, 4), (0, 10), (1, 16), (2, 14)], [(0, 0), (0, 5), (1, 8), (2, 5), (2, 18)], [(0, 0), (0, 7), (1, 1), (2, 8), (2, 16)]]
    \item 1 \{0=2592, 1=73872, 2=495720, 3=905256\} [[(0, 0), (0, 1), (0, 3), (0, 12), (1, 5)], [(0, 0), (0, 4), (0, 10), (1, 16), (2, 14)], [(0, 0), (0, 5), (1, 8), (2, 5), (2, 18)], [(0, 0), (0, 7), (1, 18), (1, 26), (2, 6)]]
    \item 1 \{0=648, 1=77760, 2=495720, 3=903312\} [[(0, 0), (0, 1), (0, 3), (0, 12), (1, 5)], [(0, 0), (0, 4), (0, 21), (1, 17), (2, 15)], [(0, 0), (0, 5), (1, 0), (1, 14), (2, 24)], [(0, 0), (0, 7), (1, 1), (2, 8), (2, 16)]]
    \item 1 \{0=1296, 1=47628, 2=480168, 3=948348\} [[(0, 0), (0, 1), (0, 3), (0, 12), (1, 5)], [(0, 0), (0, 4), (0, 21), (1, 17), (2, 15)], [(0, 0), (0, 5), (1, 0), (1, 14), (2, 24)], [(0, 0), (0, 7), (1, 18), (1, 26), (2, 6)]]
    \item 1 \{0=1296, 1=74844, 2=489888, 3=911412\} [[(0, 0), (0, 1), (0, 3), (0, 12), (1, 5)], [(0, 0), (0, 4), (0, 21), (1, 17), (2, 15)], [(0, 0), (0, 5), (1, 8), (2, 5), (2, 18)], [(0, 0), (0, 7), (1, 1), (2, 8), (2, 16)]]
    \item 1 \{0=1944, 1=64638, 2=502524, 3=908334\} [[(0, 0), (0, 1), (0, 3), (0, 12), (1, 5)], [(0, 0), (0, 4), (0, 21), (1, 17), (2, 15)], [(0, 0), (0, 5), (1, 8), (2, 5), (2, 18)], [(0, 0), (0, 7), (1, 18), (1, 26), (2, 6)]]
    \item 1 \{0=1296, 1=51030, 2=502524, 3=922590\} [[(0, 0), (0, 1), (0, 3), (1, 3), (1, 8)], [(0, 0), (0, 4), (0, 10), (1, 23), (2, 17)], [(0, 0), (0, 7), (0, 15), (2, 9), (2, 18)], [(0, 0), (0, 11), (1, 12), (1, 26), (2, 16)]]
    \item 1 \{1=39366, 2=492804, 3=945270\} [[(0, 0), (0, 1), (0, 3), (1, 3), (1, 8)], [(0, 0), (0, 4), (0, 10), (1, 23), (2, 17)], [(0, 0), (0, 7), (0, 15), (2, 9), (2, 18)], [(0, 0), (0, 11), (1, 22), (2, 12), (2, 26)]]
    \item 1 \{1=34992, 2=486000, 3=956448\} [[(0, 0), (0, 1), (0, 3), (1, 3), (1, 8)], [(0, 0), (0, 4), (0, 10), (1, 23), (2, 17)], [(0, 0), (0, 7), (0, 19), (1, 16), (1, 25)], [(0, 0), (0, 11), (1, 12), (1, 26), (2, 16)]]
    \item 1 \{1=46170, 2=483084, 3=948186\} [[(0, 0), (0, 1), (0, 3), (1, 3), (1, 8)], [(0, 0), (0, 4), (0, 10), (1, 23), (2, 17)], [(0, 0), (0, 7), (0, 19), (1, 16), (1, 25)], [(0, 0), (0, 11), (1, 22), (2, 12), (2, 26)]]
    \item 1 \{1=65610, 2=483084, 3=928746\} [[(0, 0), (0, 1), (0, 3), (1, 3), (1, 8)], [(0, 0), (0, 4), (0, 21), (1, 14), (2, 8)], [(0, 0), (0, 7), (0, 15), (2, 9), (2, 18)], [(0, 0), (0, 11), (1, 12), (1, 26), (2, 16)]]
    \item 1 \{0=648, 1=48114, 2=488916, 3=939762\} [[(0, 0), (0, 1), (0, 3), (1, 3), (1, 8)], [(0, 0), (0, 4), (0, 21), (1, 14), (2, 8)], [(0, 0), (0, 7), (0, 15), (2, 9), (2, 18)], [(0, 0), (0, 11), (1, 22), (2, 12), (2, 26)]]
    \item 1 \{1=48114, 2=502524, 3=926802\} [[(0, 0), (0, 1), (0, 3), (1, 3), (1, 8)], [(0, 0), (0, 4), (0, 21), (1, 14), (2, 8)], [(0, 0), (0, 7), (0, 19), (1, 16), (1, 25)], [(0, 0), (0, 11), (1, 12), (1, 26), (2, 16)]]
    \item 1 \{0=648, 1=46170, 2=492804, 3=937818\} [[(0, 0), (0, 1), (0, 3), (1, 3), (1, 8)], [(0, 0), (0, 4), (0, 21), (1, 14), (2, 8)], [(0, 0), (0, 7), (0, 19), (1, 16), (1, 25)], [(0, 0), (0, 11), (1, 22), (2, 12), (2, 26)]]
    \item 1 \{0=1296, 1=89910, 2=500580, 3=885654\} [[(0, 0), (0, 1), (0, 4), (0, 9), (1, 1)], [(0, 0), (0, 2), (1, 9), (1, 16), (2, 12)], [(0, 0), (0, 6), (0, 16), (1, 18), (2, 23)], [(0, 0), (0, 12), (1, 6), (2, 1), (2, 14)]]
    \item 1 \{0=648, 1=59292, 2=486000, 3=931500\} [[(0, 0), (0, 1), (0, 4), (0, 9), (1, 1)], [(0, 0), (0, 2), (1, 9), (1, 16), (2, 12)], [(0, 0), (0, 6), (0, 16), (1, 18), (2, 23)], [(0, 0), (0, 12), (1, 11), (1, 25), (2, 6)]]
    \item 1 \{0=1296, 1=76788, 2=513216, 3=886140\} [[(0, 0), (0, 1), (0, 4), (0, 9), (1, 1)], [(0, 0), (0, 2), (1, 9), (1, 16), (2, 12)], [(0, 0), (0, 6), (0, 17), (1, 10), (2, 15)], [(0, 0), (0, 12), (1, 6), (2, 1), (2, 14)]]
    \item 1 \{0=3240, 1=62694, 2=500580, 3=910926\} [[(0, 0), (0, 1), (0, 4), (0, 9), (1, 1)], [(0, 0), (0, 2), (1, 9), (1, 16), (2, 12)], [(0, 0), (0, 6), (0, 17), (1, 10), (2, 15)], [(0, 0), (0, 12), (1, 11), (1, 25), (2, 6)]]
    \item 1 \{0=3888, 1=83592, 2=497664, 3=892296\} [[(0, 0), (0, 1), (0, 4), (0, 9), (1, 1)], [(0, 0), (0, 2), (1, 17), (2, 13), (2, 20)], [(0, 0), (0, 6), (0, 16), (1, 18), (2, 23)], [(0, 0), (0, 12), (1, 6), (2, 1), (2, 14)]]
    \item 1 \{0=1296, 1=80190, 2=483084, 3=912870\} [[(0, 0), (0, 1), (0, 4), (0, 9), (1, 1)], [(0, 0), (0, 2), (1, 17), (2, 13), (2, 20)], [(0, 0), (0, 6), (0, 16), (1, 18), (2, 23)], [(0, 0), (0, 12), (1, 11), (1, 25), (2, 6)]]
    \item 1 \{0=648, 1=59292, 2=505440, 3=912060\} [[(0, 0), (0, 1), (0, 4), (0, 9), (1, 1)], [(0, 0), (0, 2), (1, 17), (2, 13), (2, 20)], [(0, 0), (0, 6), (0, 17), (1, 10), (2, 15)], [(0, 0), (0, 12), (1, 6), (2, 1), (2, 14)]]
    \item 1 \{0=1296, 1=62208, 2=526824, 3=887112\} [[(0, 0), (0, 1), (0, 4), (0, 9), (1, 1)], [(0, 0), (0, 2), (1, 17), (2, 13), (2, 20)], [(0, 0), (0, 6), (0, 17), (1, 10), (2, 15)], [(0, 0), (0, 12), (1, 11), (1, 25), (2, 6)]]
\end{enumerate}
\end{example}

\begin{example} $\mathbb Z_9 \times \mathbb Z_9$, $v=81$, $k=5$.
\begin{enumerate}
    \item 1 \{0=648, 1=43740, 2=466560, 3=966492\} [[(0, 0), (0, 1), (0, 3), (1, 0), (4, 1)], [(0, 0), (0, 4), (2, 2), (6, 1), (7, 5)], [(0, 0), (1, 1), (2, 3), (3, 6), (7, 8)], [(0, 0), (1, 5), (3, 5), (4, 3), (6, 0)]]
    \item 1 \{0=1296, 1=35964, 2=474336, 3=965844\} [[(0, 0), (0, 1), (0, 3), (1, 0), (4, 1)], [(0, 0), (0, 4), (2, 2), (6, 1), (7, 5)], [(0, 0), (1, 1), (2, 3), (3, 6), (7, 8)], [(0, 0), (1, 5), (4, 5), (6, 2), (7, 0)]]
    \item 1 \{0=1296, 1=63180, 2=484056, 3=928908\} [[(0, 0), (0, 1), (0, 3), (1, 0), (6, 2)], [(0, 0), (0, 4), (2, 1), (4, 1), (6, 4)], [(0, 0), (1, 1), (2, 4), (6, 6), (8, 2)], [(0, 0), (1, 2), (2, 7), (6, 3), (7, 7)]]
    \item 1 \{0=1296, 1=66582, 2=483084, 3=926478\} [[(0, 0), (0, 1), (0, 3), (1, 0), (6, 2)], [(0, 0), (0, 4), (2, 1), (4, 1), (6, 4)], [(0, 0), (1, 1), (2, 4), (6, 6), (8, 2)], [(0, 0), (1, 2), (3, 4), (4, 8), (8, 4)]]
    \item 1 \{0=648, 1=62694, 2=500580, 3=913518\} [[(0, 0), (0, 1), (0, 3), (1, 0), (6, 2)], [(0, 0), (0, 4), (2, 1), (4, 4), (6, 0)], [(0, 0), (1, 1), (2, 8), (6, 4), (7, 7)], [(0, 0), (1, 2), (2, 7), (3, 2), (7, 5)]]
    \item 1 \{0=1944, 1=71928, 2=505440, 3=898128\} [[(0, 0), (0, 1), (0, 3), (1, 0), (6, 2)], [(0, 0), (0, 4), (2, 1), (4, 4), (6, 0)], [(0, 0), (1, 1), (2, 8), (6, 4), (7, 7)], [(0, 0), (1, 2), (3, 6), (7, 0), (8, 4)]]
    \item 3 \{1=49572, 2=491832, 3=936036\} [[(0, 0), (0, 1), (0, 3), (1, 0), (6, 2)], [(0, 0), (0, 4), (2, 1), (4, 4), (6, 0)], [(0, 0), (1, 1), (3, 3), (4, 6), (8, 2)], [(0, 0), (1, 2), (3, 6), (7, 0), (8, 4)]]
    \item 1 \{0=1296, 1=64152, 2=482112, 3=929880\} [[(0, 0), (0, 1), (0, 3), (1, 0), (6, 2)], [(0, 0), (0, 4), (2, 2), (4, 3), (5, 8)], [(0, 0), (1, 1), (2, 4), (6, 6), (8, 2)], [(0, 0), (1, 2), (2, 6), (4, 6), (7, 6)]]
    \item 1 \{1=55890, 2=481140, 3=940410\} [[(0, 0), (0, 1), (0, 3), (1, 0), (6, 2)], [(0, 0), (0, 4), (2, 2), (4, 3), (5, 8)], [(0, 0), (1, 1), (2, 4), (6, 6), (8, 2)], [(0, 0), (1, 2), (3, 5), (6, 5), (8, 5)]]
    \item 1 \{0=1296, 1=57348, 2=511272, 3=907524\} [[(0, 0), (0, 1), (0, 3), (1, 0), (6, 2)], [(0, 0), (0, 4), (2, 2), (4, 3), (5, 8)], [(0, 0), (1, 1), (2, 8), (4, 4), (8, 6)], [(0, 0), (1, 2), (2, 6), (4, 6), (7, 6)]]
    \item 1 \{0=648, 1=50058, 2=508356, 3=918378\} [[(0, 0), (0, 1), (0, 3), (1, 0), (6, 2)], [(0, 0), (0, 4), (2, 2), (4, 3), (5, 8)], [(0, 0), (1, 1), (2, 8), (4, 4), (8, 6)], [(0, 0), (1, 2), (3, 5), (6, 5), (8, 5)]]
    \item 3 \{1=33534, 2=434484, 3=1009422\} [[(0, 0), (0, 1), (0, 3), (1, 0), (6, 2)], [(0, 0), (0, 4), (3, 4), (5, 0), (7, 3)], [(0, 0), (1, 1), (3, 3), (4, 6), (8, 2)], [(0, 0), (1, 2), (3, 6), (7, 0), (8, 4)]]
    \item 3 \{0=1944, 1=36450, 2=434484, 3=1004562\} [[(0, 0), (0, 1), (1, 0), (1, 3), (3, 8)], [(0, 0), (0, 2), (3, 3), (4, 0), (5, 7)], [(0, 0), (0, 4), (2, 1), (6, 7), (7, 2)], [(0, 0), (1, 1), (3, 4), (5, 4), (6, 0)]]
    \item 3 \{1=34992, 2=425736, 3=1016712\} [[(0, 0), (0, 1), (1, 0), (1, 3), (3, 8)], [(0, 0), (0, 2), (3, 3), (4, 0), (5, 7)], [(0, 0), (0, 4), (2, 1), (6, 7), (7, 2)], [(0, 0), (1, 1), (4, 1), (5, 6), (7, 6)]]
\end{enumerate}
\end{example}

\begin{example} There are no difference families for group $\mathbb Z_3 \times \mathbb Z_3 \times \mathbb Z_9$.
\end{example}

\begin{example} $\mathbb Z_{11} \times \mathbb Z_{11}$, $v=121$, $k=6$.
\begin{enumerate}
    \item 1 \{1=45496, 2=645414, 3=2774772, 4=3213518\} [[(0, 0), (0, 1), (0, 3), (1, 0), (1, 4), (4, 2)], [(0, 0), (0, 5), (2, 0), (4, 9), (6, 0), (8, 5)], [(0, 0), (1, 2), (3, 1), (4, 6), (6, 9), (9, 3)], [(0, 0), (1, 6), (2, 4), (6, 7), (7, 3), (9, 10)]]
    {\bf \item 1 \{1=51304, 2=633798, 3=2798004, 4=3196094\} [[(0, 0), (0, 1), (0, 3), (1, 0), (1, 4), (4, 2)], [(0, 0), (0, 5), (2, 0), (4, 9), (6, 0), (8, 5)], [(0, 0), (1, 2), (3, 1), (4, 6), (6, 9), (9, 3)], [(0, 0), (1, 6), (3, 7), (5, 3), (6, 10), (10, 2)]] } \cite{HoCD}
    \item 1 \{1=40656, 2=669372, 3=2726856, 4=3242316\} [[(0, 0), (0, 1), (0, 3), (1, 0), (3, 1), (5, 5)], [(0, 0), (0, 4), (1, 9), (4, 8), (6, 4), (9, 6)], [(0, 0), (0, 5), (1, 6), (2, 8), (4, 3), (7, 9)], [(0, 0), (1, 3), (2, 10), (4, 10), (5, 3), (9, 9)]]
    \item 1 \{1=41624, 2=702042, 3=2661516, 4=3274018\} [[(0, 0), (0, 1), (0, 3), (1, 0), (3, 1), (5, 5)], [(0, 0), (0, 4), (1, 9), (4, 8), (6, 4), (9, 6)], [(0, 0), (0, 5), (1, 6), (2, 8), (4, 3), (7, 9)], [(0, 0), (1, 3), (3, 5), (7, 0), (8, 4), (10, 4)]]
    \item 1 \{1=42592, 2=655578, 3=2713788, 4=3267242\} [[(0, 0), (0, 1), (0, 3), (1, 0), (3, 1), (5, 5)], [(0, 0), (0, 4), (1, 9), (4, 8), (6, 4), (9, 6)], [(0, 0), (0, 5), (4, 7), (7, 2), (9, 8), (10, 10)], [(0, 0), (1, 3), (2, 10), (4, 10), (5, 3), (9, 9)]]
    \item 1 \{0=1210, 1=35816, 2=667920, 3=2707496, 4=3266758\} [[(0, 0), (0, 1), (0, 3), (1, 0), (3, 1), (5, 5)], [(0, 0), (0, 4), (1, 9), (4, 8), (6, 4), (9, 6)], [(0, 0), (0, 5), (4, 7), (7, 2), (9, 8), (10, 10)], [(0, 0), (1, 3), (3, 5), (7, 0), (8, 4), (10, 4)]]
    \item 3 \{1=43560, 2=542322, 3=2766060, 4=3327258\} [[(0, 0), (0, 1), (0, 3), (1, 0), (3, 6), (4, 8)], [(0, 0), (0, 4), (2, 1), (3, 4), (5, 2), (7, 5)], [(0, 0), (0, 5), (1, 9), (2, 10), (5, 1), (6, 8)], [(0, 0), (1, 5), (2, 0), (4, 4), (6, 0), (9, 9)]]
    \item 3 \{0=3630, 1=37752, 2=596772, 3=2738472, 4=3302574\} [[(0, 0), (0, 1), (0, 3), (1, 0), (3, 6), (4, 8)], [(0, 0), (0, 4), (2, 1), (3, 4), (5, 2), (7, 5)], [(0, 0), (0, 5), (1, 9), (2, 10), (5, 1), (6, 8)], [(0, 0), (1, 5), (3, 7), (6, 5), (8, 1), (10, 5)]]
    \item 1 \{1=47432, 2=691878, 3=2702172, 4=3237718\} [[(0, 0), (0, 1), (0, 3), (1, 0), (3, 6), (4, 8)], [(0, 0), (0, 4), (3, 2), (4, 6), (5, 0), (6, 6)], [(0, 0), (0, 5), (2, 1), (4, 4), (5, 2), (9, 3)], [(0, 0), (1, 1), (2, 4), (3, 0), (6, 7), (8, 1)]]
    \item 3 \{1=52272, 2=616374, 3=2673132, 4=3337422\} [[(0, 0), (0, 1), (0, 3), (1, 0), (3, 6), (4, 8)], [(0, 0), (0, 4), (3, 2), (4, 6), (5, 0), (6, 6)], [(0, 0), (0, 5), (2, 1), (4, 4), (5, 2), (9, 3)], [(0, 0), (1, 1), (4, 0), (6, 5), (9, 1), (10, 8)]]
    \item 1 \{0=3630, 1=51304, 2=668646, 3=2722500, 4=3233120\} [[(0, 0), (0, 1), (0, 3), (1, 0), (3, 6), (4, 8)], [(0, 0), (0, 4), (3, 2), (4, 6), (5, 0), (6, 6)], [(0, 0), (0, 5), (2, 2), (6, 3), (7, 1), (9, 4)], [(0, 0), (1, 1), (2, 4), (3, 0), (6, 7), (8, 1)]]
    \item 3 \{1=75504, 2=598950, 3=2713788, 4=3290958\} [[(0, 0), (0, 1), (0, 3), (1, 0), (3, 6), (4, 8)], [(0, 0), (0, 4), (5, 9), (6, 4), (7, 9), (8, 2)], [(0, 0), (0, 5), (2, 2), (6, 3), (7, 1), (9, 4)], [(0, 0), (1, 1), (2, 4), (3, 0), (6, 7), (8, 1)]]
    \item 3 \{1=46464, 2=596772, 3=2694912, 4=3341052\} [[(0, 0), (0, 1), (0, 3), (1, 0), (5, 10), (7, 5)], [(0, 0), (0, 4), (2, 3), (3, 8), (7, 0), (8, 1)], [(0, 0), (0, 5), (1, 7), (2, 2), (3, 0), (9, 0)], [(0, 0), (1, 3), (3, 1), (4, 5), (6, 6), (8, 2)]]
    \item 3 \{1=43560, 2=485694, 3=2716692, 4=3433254\} [[(0, 0), (0, 1), (0, 3), (1, 0), (5, 10), (7, 5)], [(0, 0), (0, 4), (2, 3), (3, 8), (7, 0), (8, 1)], [(0, 0), (0, 5), (1, 7), (2, 2), (3, 0), (9, 0)], [(0, 0), (1, 3), (4, 1), (6, 8), (8, 9), (9, 2)]]
    \item 3 \{0=7260, 1=84216, 2=714384, 3=2654256, 4=3219084\} [[(0, 0), (0, 1), (1, 0), (1, 3), (5, 3), (8, 10)], [(0, 0), (0, 2), (2, 8), (4, 9), (6, 5), (8, 8)], [(0, 0), (0, 4), (1, 8), (2, 4), (7, 9), (10, 6)], [(0, 0), (0, 5), (2, 10), (3, 0), (5, 9), (6, 4)]]
    \item 3 \{0=7260, 1=66792, 2=703494, 3=2690556, 4=3211098\} [[(0, 0), (0, 1), (1, 0), (1, 3), (5, 3), (8, 10)], [(0, 0), (0, 2), (2, 8), (4, 9), (6, 5), (8, 8)], [(0, 0), (0, 4), (1, 8), (2, 4), (7, 9), (10, 6)], [(0, 0), (0, 5), (5, 1), (6, 7), (8, 5), (9, 6)]]
    \item 3 \{1=31944, 2=572814, 3=2763156, 4=3311286\} [[(0, 0), (0, 1), (0, 3), (1, 0), (2, 7), (4, 10)], [(0, 0), (0, 4), (1, 6), (4, 4), (6, 4), (9, 10)], [(0, 0), (0, 5), (1, 3), (5, 8), (7, 10), (8, 0)], [(0, 0), (1, 4), (3, 1), (6, 9), (7, 3), (9, 2)]]
    \item 1 \{1=28072, 2=562650, 3=2780580, 4=3307898\} [[(0, 0), (0, 1), (0, 3), (1, 0), (2, 7), (4, 10)], [(0, 0), (0, 4), (1, 6), (4, 4), (6, 4), (9, 10)], [(0, 0), (0, 5), (1, 3), (5, 8), (7, 10), (8, 0)], [(0, 0), (1, 4), (3, 2), (5, 1), (6, 6), (9, 3)]]
    \item 3 \{1=34848, 2=625086, 3=2783484, 4=3235782\} [[(0, 0), (0, 1), (0, 3), (1, 0), (2, 7), (4, 10)], [(0, 0), (0, 4), (1, 6), (4, 4), (6, 4), (9, 10)], [(0, 0), (0, 5), (3, 5), (4, 6), (6, 8), (10, 2)], [(0, 0), (1, 4), (3, 1), (6, 9), (7, 3), (9, 2)]]
    \item 1 \{1=23232, 2=561198, 3=2774772, 4=3319998\} [[(0, 0), (0, 1), (0, 3), (1, 0), (2, 7), (4, 10)], [(0, 0), (0, 4), (1, 6), (4, 4), (6, 4), (9, 10)], [(0, 0), (0, 5), (3, 5), (4, 6), (6, 8), (10, 2)], [(0, 0), (1, 4), (3, 2), (5, 1), (6, 6), (9, 3)]]
\end{enumerate}
\end{example}

\begin{example} $\mathbb Z_{2} \times \mathbb Z_{3} \times \mathbb Z_{3} \times \mathbb Z_{7}$, $v=126$, $k=6$.
\begin{enumerate}
    \item 1 \{1=37296, 2=659988, 3=3019464, 4=3843252\} [[(0, 0, 0, 0), (0, 0, 0, 1), (0, 0, 1, 2), (0, 1, 0, 3), (0, 2, 1, 4), (1, 2, 2, 0)], [(0, 0, 0, 0), (0, 0, 0, 2), (1, 0, 1, 1), (1, 0, 2, 5), (1, 2, 0, 5), (1, 2, 1, 1)], [(0, 0, 0, 0), (0, 0, 0, 3), (0, 1, 2, 5), (1, 1, 2, 3), (1, 2, 0, 2), (1, 2, 2, 4)], [(0, 0, 0, 0), (0, 0, 1, 6), (0, 1, 1, 5), (0, 2, 2, 0), (1, 0, 1, 5), (1, 1, 1, 2)], [(0, 0, 0, 0), (0, 0, 1, 0), (0, 0, 2, 0), (1, 0, 0, 0), (1, 0, 1, 0), (1, 0, 2, 0)]]
    \item 1 \{1=48384, 2=650916, 3=3019464, 4=3841236\} [[(0, 0, 0, 0), (0, 0, 0, 1), (0, 0, 1, 2), (0, 1, 0, 3), (0, 2, 1, 4), (1, 2, 2, 0)], [(0, 0, 0, 0), (0, 0, 0, 2), (1, 0, 1, 1), (1, 0, 2, 5), (1, 2, 0, 5), (1, 2, 1, 1)], [(0, 0, 0, 0), (0, 0, 0, 3), (0, 1, 2, 5), (1, 1, 2, 3), (1, 2, 0, 2), (1, 2, 2, 4)], [(0, 0, 0, 0), (0, 0, 1, 6), (0, 1, 2, 6), (0, 2, 0, 1), (1, 0, 0, 1), (1, 2, 0, 4)], [(0, 0, 0, 0), (0, 0, 1, 0), (0, 0, 2, 0), (1, 0, 0, 0), (1, 0, 1, 0), (1, 0, 2, 0)]]
    \item 1 \{1=35280, 2=620676, 3=2995272, 4=3908772\} [[(0, 0, 0, 0), (0, 0, 0, 1), (0, 0, 1, 2), (0, 1, 0, 3), (0, 2, 1, 4), (1, 2, 2, 0)], [(0, 0, 0, 0), (0, 0, 0, 2), (1, 0, 1, 1), (1, 0, 2, 5), (1, 2, 0, 5), (1, 2, 1, 1)], [(0, 0, 0, 0), (0, 0, 0, 3), (0, 2, 1, 5), (1, 1, 0, 1), (1, 1, 1, 6), (1, 2, 1, 0)], [(0, 0, 0, 0), (0, 0, 1, 6), (0, 1, 1, 5), (0, 2, 2, 0), (1, 0, 1, 5), (1, 1, 1, 2)], [(0, 0, 0, 0), (0, 0, 1, 0), (0, 0, 2, 0), (1, 0, 0, 0), (1, 0, 1, 0), (1, 0, 2, 0)]]
    \item 1 \{1=57456, 2=681912, 3=3033072, 4=3787560\} [[(0, 0, 0, 0), (0, 0, 0, 1), (0, 0, 1, 2), (0, 1, 0, 3), (0, 2, 1, 4), (1, 2, 2, 0)], [(0, 0, 0, 0), (0, 0, 0, 2), (1, 0, 1, 1), (1, 0, 2, 5), (1, 2, 0, 5), (1, 2, 1, 1)], [(0, 0, 0, 0), (0, 0, 0, 3), (0, 2, 1, 5), (1, 1, 0, 1), (1, 1, 1, 6), (1, 2, 1, 0)], [(0, 0, 0, 0), (0, 0, 1, 6), (0, 1, 2, 6), (0, 2, 0, 1), (1, 0, 0, 1), (1, 2, 0, 4)], [(0, 0, 0, 0), (0, 0, 1, 0), (0, 0, 2, 0), (1, 0, 0, 0), (1, 0, 1, 0), (1, 0, 2, 0)]]
    \item 1 \{1=44352, 2=725004, 3=3028536, 4=3762108\} [[(0, 0, 0, 0), (0, 0, 0, 1), (0, 0, 1, 2), (0, 1, 0, 3), (0, 2, 1, 4), (1, 2, 2, 0)], [(0, 0, 0, 0), (0, 0, 0, 2), (1, 0, 1, 4), (1, 0, 2, 1), (1, 1, 0, 4), (1, 1, 2, 1)], [(0, 0, 0, 0), (0, 0, 0, 3), (0, 1, 2, 5), (1, 1, 2, 3), (1, 2, 0, 2), (1, 2, 2, 4)], [(0, 0, 0, 0), (0, 0, 1, 6), (0, 1, 1, 5), (0, 2, 2, 0), (1, 0, 1, 5), (1, 1, 1, 2)], [(0, 0, 0, 0), (0, 0, 1, 0), (0, 0, 2, 0), (1, 0, 0, 0), (1, 0, 1, 0), (1, 0, 2, 0)]]
    \item 1 \{1=33264, 2=662256, 3=3005856, 4=3858624\} [[(0, 0, 0, 0), (0, 0, 0, 1), (0, 0, 1, 2), (0, 1, 0, 3), (0, 2, 1, 4), (1, 2, 2, 0)], [(0, 0, 0, 0), (0, 0, 0, 2), (1, 0, 1, 4), (1, 0, 2, 1), (1, 1, 0, 4), (1, 1, 2, 1)], [(0, 0, 0, 0), (0, 0, 0, 3), (0, 1, 2, 5), (1, 1, 2, 3), (1, 2, 0, 2), (1, 2, 2, 4)], [(0, 0, 0, 0), (0, 0, 1, 6), (0, 1, 2, 6), (0, 2, 0, 1), (1, 0, 0, 1), (1, 2, 0, 4)], [(0, 0, 0, 0), (0, 0, 1, 0), (0, 0, 2, 0), (1, 0, 0, 0), (1, 0, 1, 0), (1, 0, 2, 0)]]
    \item 1 \{1=41328, 2=639576, 3=3002832, 4=3876264\} [[(0, 0, 0, 0), (0, 0, 0, 1), (0, 0, 1, 2), (0, 1, 0, 3), (0, 2, 1, 4), (1, 2, 2, 0)], [(0, 0, 0, 0), (0, 0, 0, 2), (1, 0, 1, 4), (1, 0, 2, 1), (1, 1, 0, 4), (1, 1, 2, 1)], [(0, 0, 0, 0), (0, 0, 0, 3), (0, 2, 1, 5), (1, 1, 0, 1), (1, 1, 1, 6), (1, 2, 1, 0)], [(0, 0, 0, 0), (0, 0, 1, 6), (0, 1, 1, 5), (0, 2, 2, 0), (1, 0, 1, 5), (1, 1, 1, 2)], [(0, 0, 0, 0), (0, 0, 1, 0), (0, 0, 2, 0), (1, 0, 0, 0), (1, 0, 1, 0), (1, 0, 2, 0)]]
    \item 1 \{0=1260, 1=43344, 2=656208, 3=3046176, 4=3813012\} [[(0, 0, 0, 0), (0, 0, 0, 1), (0, 0, 1, 2), (0, 1, 0, 3), (0, 2, 1, 4), (1, 2, 2, 0)], [(0, 0, 0, 0), (0, 0, 0, 2), (1, 0, 1, 4), (1, 0, 2, 1), (1, 1, 0, 4), (1, 1, 2, 1)], [(0, 0, 0, 0), (0, 0, 0, 3), (0, 2, 1, 5), (1, 1, 0, 1), (1, 1, 1, 6), (1, 2, 1, 0)], [(0, 0, 0, 0), (0, 0, 1, 6), (0, 1, 2, 6), (0, 2, 0, 1), (1, 0, 0, 1), (1, 2, 0, 4)], [(0, 0, 0, 0), (0, 0, 1, 0), (0, 0, 2, 0), (1, 0, 0, 0), (1, 0, 1, 0), (1, 0, 2, 0)]]
\end{enumerate}
\end{example}

\begin{example} $\mathbb Z_{13} \times \mathbb Z_{13}$, $v=169$, $k=7$.
\begin{enumerate}
    \item 9 \{0=12168, 1=60840, 2=328536, 3=3598686, 4=10263708, 5=8733582\} [[(0, 0), (0, 1), (0, 3), (0, 9), (1, 0), (3, 0), (9, 0)], [(0, 0), (1, 1), (2, 7), (3, 5), (6, 12), (8, 8), (10, 10)], [(0, 0), (1, 2), (2, 10), (5, 8), (7, 11), (9, 3), (12, 4)], [(0, 0), (1, 3), (2, 8), (3, 2), (6, 10), (8, 11), (10, 4)]]
    \item 9 \{0=18252, 1=30420, 2=450216, 3=3486132, 4=10239372, 5=8773128\} [[(0, 0), (0, 1), (0, 3), (0, 9), (1, 0), (3, 0), (9, 0)], [(0, 0), (1, 1), (2, 7), (3, 5), (6, 12), (8, 8), (10, 10)], [(0, 0), (1, 2), (2, 10), (5, 8), (7, 11), (9, 3), (12, 4)], [(0, 0), (1, 3), (4, 12), (6, 5), (8, 6), (11, 1), (12, 8)]]
    {\bf \item 3 \{1=25350, 2=612456, 3=3887676, 4=10026432, 5=8445606\} [[(0, 0), (0, 1), (1, 0), (1, 2), (3, 5), (7, 12), (8, 7)], [(0, 0), (0, 3), (3, 1), (4, 4), (5, 0), (5, 8), (9, 5)], [(0, 0), (0, 4), (1, 11), (2, 4), (4, 0), (10, 5), (11, 2)], [(0, 0), (0, 6), (2, 1), (3, 6), (6, 0), (8, 10), (10, 3)]]} \cite{HoCD}
    \item 1 \{1=21970, 2=452920, 3=3835962, 4=10106876, 5=8579792\} [[(0, 0), (0, 1), (0, 3), (0, 9), (1, 0), (3, 7), (9, 6)], [(0, 0), (1, 1), (2, 3), (5, 0), (6, 3), (7, 8), (8, 1)], [(0, 0), (1, 7), (3, 3), (5, 3), (8, 12), (10, 1), (11, 12)], [(0, 0), (1, 8), (2, 4), (5, 6), (7, 12), (9, 11), (11, 3)]]
    \item 1 \{1=15210, 2=478608, 3=3756870, 4=10101468, 5=8645364\} [[(0, 0), (0, 1), (0, 3), (0, 9), (1, 0), (3, 7), (9, 6)], [(0, 0), (1, 1), (2, 3), (5, 0), (6, 3), (7, 8), (8, 1)], [(0, 0), (1, 7), (3, 3), (5, 3), (8, 12), (10, 1), (11, 12)], [(0, 0), (1, 8), (3, 5), (5, 10), (7, 9), (9, 2), (12, 4)]]
    \item 3 \{1=25350, 2=539448, 3=3686904, 4=9945312, 5=8800506\} [[(0, 0), (0, 1), (0, 3), (0, 9), (1, 0), (3, 7), (9, 6)], [(0, 0), (1, 1), (2, 3), (5, 0), (6, 3), (7, 8), (8, 1)], [(0, 0), (1, 7), (3, 8), (4, 6), (6, 8), (9, 4), (11, 4)], [(0, 0), (1, 8), (3, 5), (5, 10), (7, 9), (9, 2), (12, 4)]]
    \item 3 \{1=40560, 2=470496, 3=3677778, 4=10101468, 5=8707218\} [[(0, 0), (0, 1), (0, 3), (0, 9), (1, 0), (3, 7), (9, 6)], [(0, 0), (1, 1), (6, 0), (7, 6), (8, 11), (9, 1), (12, 11)], [(0, 0), (1, 7), (3, 3), (5, 3), (8, 12), (10, 1), (11, 12)], [(0, 0), (1, 8), (2, 4), (5, 6), (7, 12), (9, 11), (11, 3)]]
\end{enumerate}
\end{example}

\begin{example} $\mathbb Z_{5} \times \mathbb Z_{5} \times \mathbb Z_{7}$, $v=175$, $k=7$.
\begin{enumerate}
    {\bf \item 24 \{2=495600, 3=3313800, 4=11306400, 5=10462200\} [[(0, 0, 0), (0, 1, 0), (0, 3, 1), (1, 0, 4), (1, 2, 6), (4, 1, 4), (4, 4, 6)], [(0, 0, 0), (0, 1, 1), (0, 2, 0), (2, 2, 4), (2, 3, 6), (3, 0, 4), (3, 4, 6)], [(0, 0, 0), (0, 1, 3), (1, 1, 3), (1, 2, 0), (2, 0, 2), (3, 1, 4), (4, 2, 2)], [(0, 0, 0), (0, 2, 3), (1, 0, 2), (2, 4, 2), (3, 2, 3), (3, 4, 0), (4, 2, 4)], [(0, 0, 0), (0, 0, 1), (0, 0, 2), (0, 0, 3), (0, 0, 4), (0, 0, 5), (0, 0, 6)]]} \cite{HoCD} \cite{Krc}
    {\bf \item 6 \{1=78750, 2=659400, 3=4164300, 4=10777200, 5=9898350\} [[(0, 0, 0), (0, 1, 0), (0, 3, 1), (1, 0, 4), (1, 2, 6), (4, 1, 4), (4, 4, 6)], [(0, 0, 0), (0, 1, 1), (0, 2, 0), (2, 2, 4), (2, 3, 6), (3, 0, 4), (3, 4, 6)], [(0, 0, 0), (0, 1, 3), (1, 1, 3), (1, 2, 0), (2, 0, 2), (3, 1, 4), (4, 2, 2)], [(0, 0, 0), (0, 2, 3), (1, 0, 6), (2, 0, 0), (2, 3, 3), (3, 3, 1), (4, 2, 1)], [(0, 0, 0), (0, 0, 1), (0, 0, 2), (0, 0, 3), (0, 0, 4), (0, 0, 5), (0, 0, 6)]]} \cite{Krc}
    \item 8 \{1=133000, 2=828800, 3=4237800, 4=10698800, 5=9679600\} [[(0, 0, 0), (0, 1, 0), (0, 3, 1), (1, 0, 4), (1, 2, 6), (4, 1, 4), (4, 4, 6)], [(0, 0, 0), (0, 1, 1), (0, 2, 0), (2, 2, 4), (2, 3, 6), (3, 0, 4), (3, 4, 6)], [(0, 0, 0), (0, 1, 3), (1, 4, 1), (2, 0, 6), (3, 1, 1), (4, 0, 0), (4, 4, 3)], [(0, 0, 0), (0, 2, 3), (1, 0, 6), (2, 0, 0), (2, 3, 3), (3, 3, 1), (4, 2, 1)], [(0, 0, 0), (0, 0, 1), (0, 0, 2), (0, 0, 3), (0, 0, 4), (0, 0, 5), (0, 0, 6)]]
\end{enumerate}
\end{example}
It's kind of strange that first two designs were found by different authors and third was missing for some reason.

\section{Enumeration of S(2,6,111) point-transitive systems}

Before implementing effective generalized algorithm for non-commutative groups, other algorithm was used to generate results below. Later this results were confirmed by generalized algorithm which will be published later. It is obvious that group $\mathbb Z_{37} \times \mathbb Z_3$ can't produce any, so group $\mathbb Z_{37} \rtimes \mathbb Z_3$ is only one option. One of such designs was constructed by W.~H.~Mills \cite{Mills}. Now I'm finishing his work and claiming that there are $30$ non-isomorphic Steiner systems with automorphism group $\mathbb Z_{37} \rtimes \mathbb Z_3$ (the order of their automorphism group is exactly 111). Below is raw data with all lines of this designs. The first design was obtained by Mills, others seems to be unknown.

\begin{enumerate}
\item 1 \{1=31968, 2=534132, 3=2203128, 4=2358972\} [[8, 39, 40, 41, 42, 70], [0, 1, 2, 3, 31, 80], [30, 31, 32, 33, 61, 110], [6, 7, 8, 9, 37, 86], [23, 54, 55, 56, 57, 85], [15, 16, 17, 18, 46, 95], [27, 28, 29, 30, 58, 107], [24, 25, 26, 27, 55, 104], [26, 57, 58, 59, 60, 88], [29, 60, 61, 62, 63, 91], [25, 74, 105, 106, 107, 108], [41, 72, 73, 74, 75, 103], [7, 56, 87, 88, 89, 90], [10, 59, 90, 91, 92, 93], [13, 62, 93, 94, 95, 96], [22, 71, 102, 103, 104, 105], [9, 10, 11, 12, 40, 89], [44, 75, 76, 77, 78, 106], [18, 19, 20, 21, 49, 98], [11, 42, 43, 44, 45, 73], [3, 4, 5, 6, 34, 83], [16, 65, 96, 97, 98, 99], [19, 68, 99, 100, 101, 102], [17, 48, 49, 50, 51, 79], [47, 78, 79, 80, 81, 109], [20, 51, 52, 53, 54, 82], [12, 13, 14, 15, 43, 92], [0, 28, 77, 108, 109, 110], [4, 53, 84, 85, 86, 87], [1, 50, 81, 82, 83, 84], [32, 63, 64, 65, 66, 94], [2, 33, 34, 35, 36, 64], [21, 22, 23, 24, 52, 101], [38, 69, 70, 71, 72, 100], [14, 45, 46, 47, 48, 76], [35, 66, 67, 68, 69, 97], [5, 36, 37, 38, 39, 67], [5, 26, 28, 43, 69, 102], [18, 22, 25, 30, 50, 89], [31, 40, 54, 62, 92, 108], [26, 66, 70, 73, 78, 98], [1, 27, 60, 74, 95, 97], [5, 21, 55, 64, 78, 86], [29, 45, 79, 88, 102, 110], [13, 39, 72, 86, 107, 109], [10, 19, 33, 41, 71, 87], [6, 14, 44, 60, 94, 103], [17, 19, 34, 60, 93, 107], [11, 27, 61, 70, 84, 92], [20, 22, 37, 63, 96, 110], [1, 10, 24, 32, 62, 78], [6, 40, 49, 63, 71, 101], [2, 41, 81, 85, 88, 93], [21, 35, 56, 58, 73, 99], [25, 34, 48, 56, 86, 102], [7, 33, 66, 80, 101, 103], [9, 43, 52, 66, 74, 104], [14, 53, 93, 97, 100, 105], [17, 33, 67, 76, 90, 98], [15, 19, 22, 27, 47, 86], [20, 60, 64, 67, 72, 92], [2, 42, 46, 49, 54, 74], [30, 34, 37, 42, 62, 101], [35, 75, 79, 82, 87, 107], [7, 16, 30, 38, 68, 84], [0, 14, 35, 37, 52, 78], [6, 10, 13, 18, 38, 77], [11, 51, 55, 58, 63, 83], [0, 33, 47, 68, 70, 85], [24, 57, 71, 92, 94, 109], [14, 30, 64, 73, 87, 95], [8, 47, 87, 91, 94, 99], [2, 18, 52, 61, 75, 83], [3, 7, 10, 15, 35, 74], [8, 48, 52, 55, 60, 80], [3, 36, 50, 71, 73, 88], [28, 37, 51, 59, 89, 105], [0, 34, 43, 57, 65, 95], [22, 31, 45, 53, 83, 99], [0, 8, 38, 54, 88, 97], [11, 32, 34, 49, 75, 108], [12, 46, 55, 69, 77, 107], [17, 56, 96, 100, 103, 108], [16, 25, 39, 47, 77, 93], [21, 54, 68, 89, 91, 106], [13, 22, 36, 44, 74, 90], [3, 23, 62, 102, 106, 109], [2, 32, 48, 82, 91, 105], [11, 50, 90, 94, 97, 102], [6, 20, 41, 43, 58, 84], [5, 35, 51, 85, 94, 108], [15, 48, 62, 83, 85, 100], [23, 63, 67, 70, 75, 95], [12, 26, 47, 49, 64, 90], [9, 17, 47, 63, 97, 106], [15, 49, 58, 72, 80, 110], [5, 45, 49, 52, 57, 77], [0, 20, 59, 99, 103, 106], [7, 21, 29, 59, 75, 109], [11, 13, 28, 54, 87, 101], [38, 78, 82, 85, 90, 110], [3, 37, 46, 60, 68, 98], [4, 13, 27, 35, 65, 81], [24, 28, 31, 36, 56, 95], [15, 29, 50, 52, 67, 93], [0, 4, 7, 12, 32, 71], [5, 44, 84, 88, 91, 96], [19, 28, 42, 50, 80, 96], [3, 17, 38, 40, 55, 81], [12, 20, 50, 66, 100, 109], [5, 7, 22, 48, 81, 95], [39, 43, 46, 51, 71, 110], [9, 42, 56, 77, 79, 94], [27, 31, 34, 39, 59, 98], [24, 38, 59, 61, 76, 102], [1, 6, 26, 65, 105, 109], [9, 13, 16, 21, 41, 80], [14, 54, 58, 61, 66, 86], [27, 41, 62, 64, 79, 105], [29, 69, 73, 76, 81, 101], [4, 30, 63, 77, 98, 100], [8, 29, 31, 46, 72, 105], [23, 39, 73, 82, 96, 104], [9, 23, 44, 46, 61, 87], [3, 11, 41, 57, 91, 100], [2, 23, 25, 40, 66, 99], [12, 45, 59, 80, 82, 97], [8, 24, 58, 67, 81, 89], [8, 10, 25, 51, 84, 98], [30, 44, 65, 67, 82, 108], [4, 18, 26, 56, 72, 106], [2, 4, 19, 45, 78, 92], [10, 36, 69, 83, 104, 106], [17, 57, 61, 64, 69, 89], [12, 16, 19, 24, 44, 83], [1, 16, 42, 75, 89, 110], [36, 40, 43, 48, 68, 107], [20, 36, 70, 79, 93, 101], [6, 39, 53, 74, 76, 91], [33, 37, 40, 45, 65, 104], [14, 16, 31, 57, 90, 104], [26, 42, 76, 85, 99, 107], [32, 72, 76, 79, 84, 104], [1, 15, 23, 53, 69, 103], [18, 51, 65, 86, 88, 103], [18, 32, 53, 55, 70, 96], [21, 25, 28, 33, 53, 92], [1, 4, 9, 29, 68, 108], [15, 51, 56, 70, 97, 104], [24, 29, 43, 70, 77, 99], [4, 11, 33, 69, 74, 88], [33, 38, 52, 79, 86, 108], [6, 42, 47, 61, 88, 95], [4, 31, 38, 60, 96, 101], [3, 39, 44, 58, 85, 92], [0, 5, 19, 46, 53, 75], [25, 32, 54, 90, 95, 109], [2, 16, 43, 50, 72, 108], [21, 57, 62, 76, 103, 110], [1, 8, 30, 66, 71, 85], [27, 32, 46, 73, 80, 102], [21, 26, 40, 67, 74, 96], [15, 20, 34, 61, 68, 90], [5, 27, 63, 68, 82, 109], [16, 23, 45, 81, 86, 100], [7, 14, 36, 72, 77, 91], [2, 24, 60, 65, 79, 106], [13, 20, 42, 78, 83, 97], [7, 34, 41, 63, 99, 104], [9, 45, 50, 64, 91, 98], [3, 8, 22, 49, 56, 78], [0, 36, 41, 55, 82, 89], [9, 14, 28, 55, 62, 84], [18, 23, 37, 64, 71, 93], [19, 26, 48, 84, 89, 103], [10, 37, 44, 66, 102, 107], [12, 48, 53, 67, 94, 101], [18, 54, 59, 73, 100, 107], [6, 11, 25, 52, 59, 81], [13, 40, 47, 69, 105, 110], [22, 29, 51, 87, 92, 106], [1, 28, 35, 57, 93, 98], [12, 17, 31, 58, 65, 87], [10, 17, 39, 75, 80, 94], [30, 35, 49, 76, 83, 105], [11, 22, 60, 66, 77, 82], [17, 28, 66, 72, 83, 88], [33, 39, 50, 55, 95, 106], [3, 14, 19, 59, 70, 108], [24, 30, 41, 46, 86, 97], [6, 12, 23, 28, 68, 79], [1, 39, 45, 56, 61, 101], [12, 18, 29, 34, 74, 85], [26, 37, 75, 81, 92, 97], [4, 42, 48, 59, 64, 104], [15, 21, 32, 37, 77, 88], [27, 33, 44, 49, 89, 100], [7, 45, 51, 62, 67, 107], [2, 13, 51, 57, 68, 73], [8, 19, 57, 63, 74, 79], [4, 44, 55, 93, 99, 110], [23, 34, 72, 78, 89, 94], [38, 49, 87, 93, 104, 109], [29, 40, 78, 84, 95, 100], [18, 24, 35, 40, 80, 91], [14, 25, 63, 69, 80, 85], [3, 9, 20, 25, 65, 76], [30, 36, 47, 52, 92, 103], [35, 46, 84, 90, 101, 106], [20, 31, 69, 75, 86, 91], [9, 15, 26, 31, 71, 82], [5, 10, 50, 61, 99, 105], [36, 42, 53, 58, 98, 109], [10, 48, 54, 65, 70, 110], [0, 6, 17, 22, 62, 73], [8, 13, 53, 64, 102, 108], [32, 43, 81, 87, 98, 103], [21, 27, 38, 43, 83, 94], [1, 41, 52, 90, 96, 107], [2, 7, 47, 58, 96, 102], [0, 11, 16, 56, 67, 105], [5, 16, 54, 60, 71, 76], [14, 26, 32, 41, 51, 101], [14, 38, 50, 56, 65, 75], [2, 12, 62, 86, 98, 104], [13, 34, 50, 58, 70, 76], [0, 10, 21, 45, 63, 72], [6, 45, 55, 66, 90, 108], [1, 22, 38, 46, 58, 64], [0, 24, 42, 51, 90, 100], [25, 46, 62, 70, 82, 88], [41, 65, 77, 83, 92, 102], [1, 7, 55, 76, 92, 100], [0, 18, 27, 66, 76, 87], [2, 26, 38, 44, 53, 63], [7, 18, 42, 60, 69, 108], [0, 9, 48, 58, 69, 93], [31, 52, 68, 76, 88, 94], [11, 23, 29, 38, 48, 98], [22, 43, 59, 67, 79, 85], [15, 33, 42, 81, 91, 102], [26, 50, 62, 68, 77, 87], [12, 22, 33, 57, 75, 84], [4, 16, 22, 70, 91, 107], [4, 10, 58, 79, 95, 103], [21, 31, 42, 66, 84, 93], [3, 21, 30, 69, 79, 90], [4, 15, 39, 57, 66, 105], [21, 39, 48, 87, 97, 108], [18, 36, 45, 84, 94, 105], [19, 40, 56, 64, 76, 82], [15, 25, 36, 60, 78, 87], [34, 55, 71, 79, 91, 97], [18, 28, 39, 63, 81, 90], [1, 49, 70, 86, 94, 106], [7, 13, 61, 82, 98, 106], [3, 12, 51, 61, 72, 96], [4, 52, 73, 89, 97, 109], [6, 15, 54, 64, 75, 99], [23, 35, 41, 50, 60, 110], [4, 20, 28, 40, 46, 94], [17, 29, 35, 44, 54, 104], [46, 67, 83, 91, 103, 109], [11, 35, 47, 53, 62, 72], [2, 14, 20, 29, 39, 89], [12, 30, 39, 78, 88, 99], [8, 20, 26, 35, 45, 95], [47, 71, 83, 89, 98, 108], [5, 11, 20, 30, 80, 104], [10, 31, 47, 55, 67, 73], [3, 53, 77, 89, 95, 104], [23, 47, 59, 65, 74, 84], [2, 11, 21, 71, 95, 107], [16, 37, 53, 61, 73, 79], [20, 32, 38, 47, 57, 107], [40, 61, 77, 85, 97, 103], [2, 8, 17, 27, 77, 101], [7, 28, 44, 52, 64, 70], [8, 14, 23, 33, 83, 107], [9, 18, 57, 67, 78, 102], [35, 59, 71, 77, 86, 96], [9, 19, 30, 54, 72, 81], [5, 17, 23, 32, 42, 92], [8, 32, 44, 50, 59, 69], [6, 56, 80, 92, 98, 107], [4, 25, 41, 49, 61, 67], [38, 62, 74, 80, 89, 99], [9, 27, 36, 75, 85, 96], [44, 68, 80, 86, 95, 105], [6, 30, 48, 57, 96, 106], [28, 49, 65, 73, 85, 91], [8, 16, 28, 34, 82, 103], [3, 13, 24, 48, 66, 75], [29, 53, 65, 71, 80, 90], [5, 29, 41, 47, 56, 66], [8, 18, 68, 92, 104, 110], [5, 14, 24, 74, 98, 110], [36, 46, 57, 81, 99, 108], [20, 44, 56, 62, 71, 81], [10, 26, 34, 46, 52, 100], [37, 58, 74, 82, 94, 100], [5, 13, 25, 31, 79, 100], [6, 16, 27, 51, 69, 78], [27, 37, 48, 72, 90, 99], [11, 19, 31, 37, 85, 106], [7, 19, 25, 73, 94, 110], [6, 24, 33, 72, 82, 93], [0, 50, 74, 86, 92, 101], [19, 35, 43, 55, 61, 109], [33, 43, 54, 78, 96, 105], [3, 42, 52, 63, 87, 105], [5, 15, 65, 89, 101, 107], [11, 17, 26, 36, 86, 110], [17, 41, 53, 59, 68, 78], [43, 64, 80, 88, 100, 106], [10, 16, 64, 85, 101, 109], [24, 34, 45, 69, 87, 96], [15, 24, 63, 73, 84, 108], [13, 29, 37, 49, 55, 103], [1, 17, 25, 37, 43, 91], [12, 21, 60, 70, 81, 105], [9, 59, 83, 95, 101, 110], [1, 12, 36, 54, 63, 102], [0, 39, 49, 60, 84, 102], [16, 32, 40, 52, 58, 106], [14, 22, 34, 40, 88, 109], [32, 56, 68, 74, 83, 93], [30, 40, 51, 75, 93, 102], [3, 27, 45, 54, 93, 103], [7, 23, 31, 43, 49, 97], [1, 13, 19, 67, 88, 104], [9, 33, 51, 60, 99, 109], [2, 10, 22, 28, 76, 97], [43, 63, 76, 86, 89, 93], [1, 11, 14, 18, 79, 99], [40, 60, 73, 83, 86, 90], [19, 39, 52, 62, 65, 69], [7, 17, 20, 24, 85, 105], [15, 28, 38, 41, 45, 106], [1, 21, 34, 44, 47, 51], [22, 42, 55, 65, 68, 72], [34, 54, 67, 77, 80, 84], [52, 72, 85, 95, 98, 102], [2, 6, 67, 87, 100, 110], [7, 27, 40, 50, 53, 57], [10, 30, 43, 53, 56, 60], [5, 8, 12, 73, 93, 106], [9, 22, 32, 35, 39, 100], [49, 69, 82, 92, 95, 99], [3, 64, 84, 97, 107, 110], [10, 20, 23, 27, 88, 108], [31, 51, 64, 74, 77, 81], [46, 66, 79, 89, 92, 96], [8, 11, 15, 76, 96, 109], [25, 45, 58, 68, 71, 75], [4, 14, 17, 21, 82, 102], [4, 24, 37, 47, 50, 54], [18, 31, 41, 44, 48, 109], [2, 5, 9, 70, 90, 103], [28, 48, 61, 71, 74, 78], [58, 78, 91, 101, 104, 108], [0, 61, 81, 94, 104, 107], [0, 13, 23, 26, 30, 91], [55, 75, 88, 98, 101, 105], [3, 16, 26, 29, 33, 94], [16, 36, 49, 59, 62, 66], [37, 57, 70, 80, 83, 87], [6, 19, 29, 32, 36, 97], [13, 33, 46, 56, 59, 63], [12, 25, 35, 38, 42, 103], [22, 26, 61, 80, 93, 108], [6, 21, 46, 50, 85, 104], [10, 29, 42, 57, 82, 86], [19, 38, 51, 66, 91, 95], [13, 17, 52, 71, 84, 99], [7, 26, 39, 54, 79, 83], [16, 20, 55, 74, 87, 102], [3, 28, 32, 67, 86, 99], [10, 14, 49, 68, 81, 96], [5, 18, 33, 58, 62, 97], [8, 21, 36, 61, 65, 100], [12, 37, 41, 76, 95, 108], [34, 53, 66, 81, 106, 110], [2, 37, 56, 69, 84, 109], [4, 8, 43, 62, 75, 90], [2, 15, 30, 55, 59, 94], [12, 27, 52, 56, 91, 110], [4, 23, 36, 51, 76, 80], [19, 23, 58, 77, 90, 105], [16, 35, 48, 63, 88, 92], [31, 50, 63, 78, 103, 107], [13, 32, 45, 60, 85, 89], [14, 27, 42, 67, 71, 106], [17, 30, 45, 70, 74, 109], [3, 18, 43, 47, 82, 101], [7, 11, 46, 65, 78, 93], [1, 20, 33, 48, 73, 77], [11, 24, 39, 64, 68, 103], [28, 47, 60, 75, 100, 104], [1, 5, 40, 59, 72, 87], [9, 24, 49, 53, 88, 107], [25, 44, 57, 72, 97, 101], [22, 41, 54, 69, 94, 98], [0, 25, 29, 64, 83, 96], [0, 15, 40, 44, 79, 98], [6, 31, 35, 70, 89, 102], [9, 34, 38, 73, 92, 105]]
\item 1 \{0=1110, 1=28416, 2=474858, 3=2181372, 4=2442444\} [[8, 39, 40, 41, 42, 70], [0, 1, 2, 3, 31, 80], [30, 31, 32, 33, 61, 110], [6, 7, 8, 9, 37, 86], [23, 54, 55, 56, 57, 85], [15, 16, 17, 18, 46, 95], [27, 28, 29, 30, 58, 107], [24, 25, 26, 27, 55, 104], [26, 57, 58, 59, 60, 88], [29, 60, 61, 62, 63, 91], [25, 74, 105, 106, 107, 108], [41, 72, 73, 74, 75, 103], [7, 56, 87, 88, 89, 90], [10, 59, 90, 91, 92, 93], [13, 62, 93, 94, 95, 96], [22, 71, 102, 103, 104, 105], [9, 10, 11, 12, 40, 89], [44, 75, 76, 77, 78, 106], [18, 19, 20, 21, 49, 98], [11, 42, 43, 44, 45, 73], [3, 4, 5, 6, 34, 83], [16, 65, 96, 97, 98, 99], [19, 68, 99, 100, 101, 102], [17, 48, 49, 50, 51, 79], [47, 78, 79, 80, 81, 109], [20, 51, 52, 53, 54, 82], [12, 13, 14, 15, 43, 92], [0, 28, 77, 108, 109, 110], [4, 53, 84, 85, 86, 87], [1, 50, 81, 82, 83, 84], [32, 63, 64, 65, 66, 94], [2, 33, 34, 35, 36, 64], [21, 22, 23, 24, 52, 101], [38, 69, 70, 71, 72, 100], [14, 45, 46, 47, 48, 76], [35, 66, 67, 68, 69, 97], [5, 36, 37, 38, 39, 67], [6, 27, 39, 57, 62, 79], [0, 5, 22, 60, 81, 93], [12, 17, 34, 72, 93, 105], [55, 65, 78, 86, 89, 95], [9, 13, 22, 25, 80, 91], [13, 51, 72, 84, 102, 107], [1, 10, 13, 68, 79, 108], [6, 11, 28, 66, 87, 99], [3, 7, 16, 19, 74, 85], [6, 18, 36, 41, 58, 96], [12, 24, 42, 47, 64, 102], [0, 18, 23, 40, 78, 99], [18, 39, 51, 69, 74, 91], [1, 11, 24, 32, 35, 41], [2, 5, 11, 82, 92, 105], [1, 56, 67, 96, 100, 109], [3, 21, 26, 43, 81, 102], [15, 36, 48, 66, 71, 88], [7, 36, 40, 49, 52, 107], [70, 80, 93, 101, 104, 110], [3, 11, 14, 20, 91, 101], [53, 64, 93, 97, 106, 109], [9, 14, 31, 69, 90, 102], [0, 12, 30, 35, 52, 90], [12, 33, 45, 63, 68, 85], [43, 53, 66, 74, 77, 83], [10, 39, 43, 52, 55, 110], [20, 31, 60, 64, 73, 76], [2, 73, 83, 96, 104, 107], [3, 15, 33, 38, 55, 93], [6, 10, 19, 22, 77, 88], [40, 50, 63, 71, 74, 80], [50, 61, 90, 94, 103, 106], [37, 47, 60, 68, 71, 77], [17, 28, 57, 61, 70, 73], [8, 21, 29, 32, 38, 109], [44, 55, 84, 88, 97, 100], [22, 32, 45, 53, 56, 62], [7, 10, 65, 76, 105, 109], [6, 14, 17, 23, 94, 104], [15, 27, 45, 50, 67, 105], [64, 74, 87, 95, 98, 104], [7, 17, 30, 38, 41, 47], [34, 44, 57, 65, 68, 74], [36, 57, 69, 87, 92, 109], [2, 15, 23, 26, 32, 103], [15, 20, 37, 75, 96, 108], [9, 27, 32, 49, 87, 108], [4, 14, 27, 35, 38, 44], [21, 42, 54, 72, 77, 94], [7, 45, 66, 78, 96, 101], [28, 38, 51, 59, 62, 68], [3, 24, 36, 54, 59, 76], [18, 22, 31, 34, 89, 100], [61, 71, 84, 92, 95, 101], [4, 7, 62, 73, 102, 106], [0, 4, 13, 16, 71, 82], [18, 30, 48, 53, 70, 108], [30, 51, 63, 81, 86, 103], [21, 25, 34, 37, 92, 103], [2, 8, 79, 89, 102, 110], [1, 39, 60, 72, 90, 95], [31, 41, 54, 62, 65, 71], [15, 19, 28, 31, 86, 97], [49, 59, 72, 80, 83, 89], [5, 8, 14, 85, 95, 108], [4, 33, 37, 46, 49, 104], [10, 20, 33, 41, 44, 50], [1, 30, 34, 43, 46, 101], [19, 29, 42, 50, 53, 59], [6, 24, 29, 46, 84, 105], [16, 54, 75, 87, 105, 110], [24, 28, 37, 40, 95, 106], [12, 20, 23, 29, 100, 110], [25, 35, 48, 56, 59, 65], [46, 56, 69, 77, 80, 86], [3, 8, 25, 63, 84, 96], [9, 21, 39, 44, 61, 99], [2, 13, 42, 46, 55, 58], [35, 46, 75, 79, 88, 91], [4, 42, 63, 75, 93, 98], [1, 4, 59, 70, 99, 103], [5, 16, 45, 49, 58, 61], [32, 43, 72, 76, 85, 88], [38, 49, 78, 82, 91, 94], [13, 23, 36, 44, 47, 53], [27, 31, 40, 43, 98, 109], [2, 19, 57, 78, 90, 108], [8, 19, 48, 52, 61, 64], [23, 34, 63, 67, 76, 79], [29, 40, 69, 73, 82, 85], [0, 21, 33, 51, 56, 73], [9, 30, 42, 60, 65, 82], [24, 45, 57, 75, 80, 97], [10, 48, 69, 81, 99, 104], [14, 25, 54, 58, 67, 70], [41, 52, 81, 85, 94, 97], [33, 54, 66, 84, 89, 106], [58, 68, 81, 89, 92, 98], [0, 8, 11, 17, 88, 98], [47, 58, 87, 91, 100, 103], [26, 37, 66, 70, 79, 82], [5, 18, 26, 29, 35, 106], [12, 16, 25, 28, 83, 94], [5, 76, 86, 99, 107, 110], [16, 26, 39, 47, 50, 56], [67, 77, 90, 98, 101, 107], [9, 17, 20, 26, 97, 107], [11, 22, 51, 55, 64, 67], [27, 48, 60, 78, 83, 100], [52, 62, 75, 83, 86, 92], [5, 13, 21, 64, 70, 88], [3, 17, 42, 66, 103, 108], [8, 33, 57, 94, 99, 105], [21, 58, 63, 69, 83, 108], [20, 35, 39, 86, 104, 109], [22, 28, 46, 74, 82, 90], [3, 46, 52, 70, 98, 106], [14, 22, 30, 73, 79, 97], [1, 9, 52, 58, 76, 104], [2, 17, 21, 68, 86, 91], [19, 24, 30, 44, 69, 93], [32, 50, 55, 77, 92, 96], [26, 34, 42, 85, 91, 109], [6, 49, 55, 73, 101, 109], [35, 53, 58, 80, 95, 99], [31, 37, 55, 83, 91, 99], [20, 28, 36, 79, 85, 103], [5, 30, 54, 91, 96, 102], [3, 9, 23, 48, 72, 109], [0, 47, 65, 70, 92, 107], [8, 26, 31, 53, 68, 72], [0, 6, 20, 45, 69, 106], [4, 12, 55, 61, 79, 107], [12, 49, 54, 60, 74, 99], [19, 25, 43, 71, 79, 87], [18, 55, 60, 66, 80, 105], [16, 21, 27, 41, 66, 90], [34, 40, 58, 86, 94, 102], [9, 33, 70, 75, 81, 95], [16, 22, 40, 68, 76, 84], [26, 44, 49, 71, 86, 90], [29, 47, 52, 74, 89, 93], [5, 23, 28, 50, 65, 69], [1, 23, 38, 42, 89, 107], [9, 46, 51, 57, 71, 96], [8, 23, 27, 74, 92, 97], [1, 6, 12, 26, 51, 75], [3, 40, 45, 51, 65, 90], [7, 13, 31, 59, 67, 75], [18, 42, 79, 84, 90, 104], [17, 25, 33, 76, 82, 100], [14, 18, 65, 83, 88, 110], [16, 44, 52, 60, 103, 109], [7, 15, 58, 64, 82, 110], [0, 43, 49, 67, 95, 103], [13, 19, 37, 65, 73, 81], [28, 33, 39, 53, 78, 102], [0, 37, 42, 48, 62, 87], [4, 9, 15, 29, 54, 78], [2, 6, 53, 71, 76, 98], [10, 16, 34, 62, 70, 78], [6, 43, 48, 54, 68, 93], [0, 14, 39, 63, 100, 105], [1, 19, 47, 55, 63, 106], [13, 18, 24, 38, 63, 87], [1, 7, 25, 53, 61, 69], [20, 38, 43, 65, 80, 84], [11, 16, 38, 53, 57, 104], [6, 30, 67, 72, 78, 92], [4, 22, 50, 58, 66, 109], [14, 29, 33, 80, 98, 103], [13, 41, 49, 57, 100, 106], [8, 13, 35, 50, 54, 101], [2, 7, 29, 44, 48, 95], [25, 31, 49, 77, 85, 93], [37, 43, 61, 89, 97, 105], [41, 59, 64, 86, 101, 105], [0, 24, 61, 66, 72, 86], [21, 45, 82, 87, 93, 107], [12, 36, 73, 78, 84, 98], [34, 39, 45, 59, 84, 108], [7, 35, 43, 51, 94, 100], [4, 10, 28, 56, 64, 72], [31, 36, 42, 56, 81, 105], [3, 50, 68, 73, 95, 110], [4, 26, 41, 45, 92, 110], [17, 32, 36, 83, 101, 106], [10, 15, 21, 35, 60, 84], [11, 29, 34, 56, 71, 75], [2, 27, 51, 88, 93, 99], [2, 20, 25, 47, 62, 66], [23, 41, 46, 68, 83, 87], [5, 9, 56, 74, 79, 101], [14, 19, 41, 56, 60, 107], [17, 22, 44, 59, 63, 110], [25, 30, 36, 50, 75, 99], [23, 31, 39, 82, 88, 106], [8, 12, 59, 77, 82, 104], [14, 32, 37, 59, 74, 78], [44, 62, 67, 89, 104, 108], [10, 38, 46, 54, 97, 103], [5, 10, 32, 47, 51, 98], [38, 56, 61, 83, 98, 102], [4, 32, 40, 48, 91, 97], [22, 27, 33, 47, 72, 96], [28, 34, 52, 80, 88, 96], [11, 26, 30, 77, 95, 100], [5, 20, 24, 71, 89, 94], [40, 46, 64, 92, 100, 108], [8, 16, 24, 67, 73, 91], [2, 10, 18, 61, 67, 85], [3, 27, 64, 69, 75, 89], [1, 29, 37, 45, 88, 94], [7, 12, 18, 32, 57, 81], [11, 15, 62, 80, 85, 107], [15, 39, 76, 81, 87, 101], [11, 36, 60, 97, 102, 108], [11, 19, 27, 70, 76, 94], [17, 35, 40, 62, 77, 81], [24, 48, 85, 90, 96, 110], [15, 52, 57, 63, 77, 102], [25, 41, 88, 95, 102, 109], [16, 23, 30, 37, 64, 80], [1, 28, 44, 91, 98, 105], [18, 33, 59, 71, 107, 109], [1, 14, 40, 57, 66, 93], [6, 13, 40, 56, 103, 110], [23, 25, 45, 60, 86, 98], [19, 26, 33, 40, 67, 83], [0, 27, 46, 59, 85, 102], [1, 8, 15, 22, 49, 65], [5, 31, 48, 57, 84, 103], [40, 47, 54, 61, 88, 104], [15, 24, 51, 70, 83, 109], [1, 21, 36, 62, 74, 110], [0, 15, 41, 53, 89, 91], [7, 23, 70, 77, 84, 91], [13, 26, 52, 69, 78, 105], [2, 9, 16, 43, 59, 106], [2, 38, 40, 60, 75, 101], [46, 53, 60, 67, 94, 110], [5, 7, 27, 42, 68, 80], [8, 20, 56, 58, 78, 93], [6, 33, 52, 65, 91, 108], [25, 42, 51, 78, 97, 110], [29, 31, 51, 66, 92, 104], [5, 17, 53, 55, 75, 90], [12, 38, 50, 86, 88, 108], [11, 23, 59, 61, 81, 96], [10, 26, 73, 80, 87, 94], [2, 49, 56, 63, 70, 97], [8, 34, 51, 60, 87, 106], [15, 34, 47, 73, 90, 99], [5, 52, 59, 66, 73, 100], [3, 10, 37, 53, 100, 107], [19, 36, 45, 72, 91, 104], [4, 20, 67, 74, 81, 88], [8, 10, 30, 45, 71, 83], [4, 21, 30, 57, 76, 89], [11, 37, 54, 63, 90, 109], [16, 29, 55, 72, 81, 108], [0, 9, 36, 55, 68, 94], [13, 30, 39, 66, 85, 98], [25, 32, 39, 46, 73, 89], [22, 38, 85, 92, 99, 106], [6, 21, 47, 59, 95, 97], [31, 38, 45, 52, 79, 95], [20, 22, 42, 57, 83, 95], [7, 24, 33, 60, 79, 92], [3, 22, 35, 61, 78, 87], [6, 15, 42, 61, 74, 100], [12, 31, 44, 70, 87, 96], [12, 27, 53, 65, 101, 103], [1, 18, 27, 54, 73, 86], [28, 35, 42, 49, 76, 92], [3, 30, 49, 62, 88, 105], [2, 28, 45, 54, 81, 100], [9, 35, 47, 83, 85, 105], [11, 13, 33, 48, 74, 86], [35, 37, 57, 72, 98, 110], [22, 39, 48, 75, 94, 107], [5, 12, 19, 46, 62, 109], [5, 41, 43, 63, 78, 104], [0, 26, 38, 74, 76, 96], [0, 7, 34, 50, 97, 104], [7, 14, 21, 28, 55, 71], [9, 24, 50, 62, 98, 100], [4, 11, 18, 25, 52, 68], [2, 4, 24, 39, 65, 77], [14, 16, 36, 51, 77, 89], [16, 32, 79, 86, 93, 100], [17, 29, 65, 67, 87, 102], [10, 27, 36, 63, 82, 95], [17, 19, 39, 54, 80, 92], [3, 12, 39, 58, 71, 97], [23, 35, 71, 73, 93, 108], [10, 23, 49, 66, 75, 102], [32, 34, 54, 69, 95, 107], [7, 20, 46, 63, 72, 99], [11, 47, 49, 69, 84, 110], [43, 50, 57, 64, 91, 107], [6, 25, 38, 64, 81, 90], [20, 32, 68, 70, 90, 105], [3, 29, 41, 77, 79, 99], [21, 40, 53, 79, 96, 105], [11, 58, 65, 72, 79, 106], [9, 18, 45, 64, 77, 103], [3, 18, 44, 56, 92, 94], [19, 35, 82, 89, 96, 103], [12, 21, 48, 67, 80, 106], [13, 29, 76, 83, 90, 97], [13, 20, 27, 34, 61, 77], [0, 19, 32, 58, 75, 84], [24, 43, 56, 82, 99, 108], [8, 44, 46, 66, 81, 107], [2, 14, 50, 52, 72, 87], [37, 44, 51, 58, 85, 101], [9, 28, 41, 67, 84, 93], [14, 61, 68, 75, 82, 109], [14, 26, 62, 64, 84, 99], [1, 17, 64, 71, 78, 85], [15, 30, 56, 68, 104, 106], [4, 31, 47, 94, 101, 108], [16, 33, 42, 69, 88, 101], [18, 37, 50, 76, 93, 102], [8, 55, 62, 69, 76, 103], [34, 41, 48, 55, 82, 98], [6, 32, 44, 80, 82, 102], [4, 17, 43, 60, 69, 96], [22, 29, 36, 43, 70, 86], [26, 28, 48, 63, 89, 101], [10, 17, 24, 31, 58, 74], [0, 10, 25, 29, 57, 101], [10, 14, 42, 86, 96, 106], [32, 42, 52, 67, 71, 99], [3, 47, 57, 67, 82, 86], [14, 24, 34, 49, 53, 81], [38, 48, 58, 73, 77, 105], [17, 27, 37, 52, 56, 84], [9, 53, 63, 73, 88, 92], [1, 5, 33, 77, 87, 97], [11, 21, 31, 46, 50, 78], [0, 44, 54, 64, 79, 83], [2, 30, 74, 84, 94, 109], [5, 15, 25, 40, 44, 72], [4, 19, 23, 51, 95, 105], [20, 30, 40, 55, 59, 87], [26, 36, 46, 61, 65, 93], [41, 51, 61, 76, 80, 108], [35, 45, 55, 70, 74, 102], [27, 71, 81, 91, 106, 110], [7, 22, 26, 54, 98, 108], [12, 56, 66, 76, 91, 95], [7, 11, 39, 83, 93, 103], [8, 18, 28, 43, 47, 75], [1, 16, 20, 48, 92, 102], [4, 8, 36, 80, 90, 100], [24, 68, 78, 88, 103, 107], [21, 65, 75, 85, 100, 104], [3, 13, 28, 32, 60, 104], [6, 50, 60, 70, 85, 89], [18, 62, 72, 82, 97, 101], [6, 16, 31, 35, 63, 107], [29, 39, 49, 64, 68, 96], [15, 59, 69, 79, 94, 98], [13, 17, 45, 89, 99, 109], [9, 19, 34, 38, 66, 110], [2, 12, 22, 37, 41, 69], [23, 33, 43, 58, 62, 90]]
\item 1 \{0=2220, 1=39072, 2=519480, 3=2240424, 4=2327004\} [[8, 39, 40, 41, 42, 70], [0, 1, 2, 3, 31, 80], [30, 31, 32, 33, 61, 110], [6, 7, 8, 9, 37, 86], [23, 54, 55, 56, 57, 85], [15, 16, 17, 18, 46, 95], [27, 28, 29, 30, 58, 107], [24, 25, 26, 27, 55, 104], [26, 57, 58, 59, 60, 88], [29, 60, 61, 62, 63, 91], [25, 74, 105, 106, 107, 108], [41, 72, 73, 74, 75, 103], [7, 56, 87, 88, 89, 90], [10, 59, 90, 91, 92, 93], [13, 62, 93, 94, 95, 96], [22, 71, 102, 103, 104, 105], [9, 10, 11, 12, 40, 89], [44, 75, 76, 77, 78, 106], [18, 19, 20, 21, 49, 98], [11, 42, 43, 44, 45, 73], [3, 4, 5, 6, 34, 83], [16, 65, 96, 97, 98, 99], [19, 68, 99, 100, 101, 102], [17, 48, 49, 50, 51, 79], [47, 78, 79, 80, 81, 109], [20, 51, 52, 53, 54, 82], [12, 13, 14, 15, 43, 92], [0, 28, 77, 108, 109, 110], [4, 53, 84, 85, 86, 87], [1, 50, 81, 82, 83, 84], [32, 63, 64, 65, 66, 94], [2, 33, 34, 35, 36, 64], [21, 22, 23, 24, 52, 101], [38, 69, 70, 71, 72, 100], [14, 45, 46, 47, 48, 76], [35, 66, 67, 68, 69, 97], [5, 36, 37, 38, 39, 67], [1, 4, 28, 47, 69, 104], [47, 51, 66, 70, 77, 95], [22, 41, 63, 98, 106, 109], [2, 20, 83, 87, 102, 106], [17, 80, 84, 99, 103, 110], [53, 57, 72, 76, 83, 101], [0, 35, 43, 46, 70, 89], [4, 21, 43, 54, 74, 90], [0, 20, 36, 61, 78, 100], [6, 28, 39, 59, 75, 100], [2, 10, 13, 37, 56, 78], [1, 25, 44, 66, 101, 109], [17, 39, 74, 82, 85, 109], [13, 24, 44, 60, 85, 102], [23, 27, 42, 46, 53, 71], [19, 36, 58, 69, 89, 105], [21, 56, 64, 67, 91, 110], [10, 21, 41, 57, 82, 99], [11, 15, 30, 34, 41, 59], [13, 32, 54, 89, 97, 100], [13, 30, 52, 63, 83, 99], [9, 34, 51, 73, 84, 104], [4, 23, 45, 80, 88, 91], [5, 68, 72, 87, 91, 98], [12, 34, 45, 65, 81, 106], [3, 23, 39, 64, 81, 103], [3, 38, 46, 49, 73, 92], [20, 24, 39, 43, 50, 68], [14, 36, 71, 79, 82, 106], [1, 12, 32, 48, 73, 90], [7, 18, 38, 54, 79, 96], [17, 33, 58, 75, 97, 108], [2, 6, 21, 25, 32, 50], [9, 31, 42, 62, 78, 103], [7, 24, 46, 57, 77, 93], [41, 45, 60, 64, 71, 89], [14, 30, 55, 72, 94, 105], [4, 7, 31, 50, 72, 107], [2, 65, 69, 84, 88, 95], [9, 13, 20, 38, 101, 105], [16, 35, 57, 92, 100, 103], [14, 18, 33, 37, 44, 62], [26, 34, 37, 61, 80, 102], [7, 10, 34, 53, 75, 110], [12, 16, 23, 41, 104, 108], [8, 71, 75, 90, 94, 101], [5, 23, 86, 90, 105, 109], [7, 26, 48, 83, 91, 94], [3, 7, 14, 32, 95, 99], [15, 50, 58, 61, 85, 104], [17, 25, 28, 52, 71, 93], [12, 47, 55, 58, 82, 101], [5, 9, 24, 28, 35, 53], [26, 30, 45, 49, 56, 74], [22, 39, 61, 72, 92, 108], [56, 60, 75, 79, 86, 104], [6, 10, 17, 35, 98, 102], [11, 19, 22, 46, 65, 87], [8, 12, 27, 31, 38, 56], [5, 13, 16, 40, 59, 81], [2, 24, 59, 67, 70, 94], [9, 29, 45, 70, 87, 109], [1, 18, 40, 51, 71, 87], [11, 33, 68, 76, 79, 103], [15, 37, 48, 68, 84, 109], [6, 26, 42, 67, 84, 106], [16, 33, 55, 66, 86, 102], [0, 22, 33, 53, 69, 94], [17, 21, 36, 40, 47, 65], [18, 53, 61, 64, 88, 107], [5, 21, 46, 63, 85, 96], [32, 40, 43, 67, 86, 108], [32, 36, 51, 55, 62, 80], [29, 37, 40, 64, 83, 105], [3, 25, 36, 56, 72, 97], [4, 15, 35, 51, 76, 93], [6, 41, 49, 52, 76, 95], [19, 30, 50, 66, 91, 108], [20, 28, 31, 55, 74, 96], [9, 44, 52, 55, 79, 98], [2, 18, 43, 60, 82, 93], [6, 31, 48, 70, 81, 101], [3, 28, 45, 67, 78, 98], [11, 74, 78, 93, 97, 104], [23, 31, 34, 58, 77, 99], [59, 63, 78, 82, 89, 107], [0, 4, 11, 29, 92, 96], [38, 42, 57, 61, 68, 86], [15, 40, 57, 79, 90, 110], [11, 27, 52, 69, 91, 102], [8, 30, 65, 73, 76, 100], [14, 22, 25, 49, 68, 90], [50, 54, 69, 73, 80, 98], [14, 77, 81, 96, 100, 107], [0, 25, 42, 64, 75, 95], [5, 27, 62, 70, 73, 97], [12, 37, 54, 76, 87, 107], [62, 66, 81, 85, 92, 110], [8, 16, 19, 43, 62, 84], [0, 15, 19, 26, 44, 107], [44, 48, 63, 67, 74, 92], [16, 27, 47, 63, 88, 105], [1, 20, 42, 77, 85, 88], [29, 33, 48, 52, 59, 77], [3, 18, 22, 29, 47, 110], [35, 39, 54, 58, 65, 83], [8, 24, 49, 66, 88, 99], [19, 38, 60, 95, 103, 106], [10, 29, 51, 86, 94, 97], [10, 27, 49, 60, 80, 96], [1, 8, 26, 89, 93, 108], [5, 8, 32, 79, 102, 107], [6, 11, 20, 23, 47, 94], [6, 24, 62, 79, 87, 108], [21, 37, 58, 70, 92, 106], [17, 34, 42, 63, 72, 90], [10, 22, 44, 58, 84, 100], [1, 27, 43, 64, 76, 98], [9, 47, 64, 72, 93, 102], [3, 24, 33, 51, 89, 106], [19, 31, 53, 67, 93, 109], [5, 22, 30, 51, 60, 78], [2, 16, 42, 58, 79, 91], [11, 25, 51, 67, 88, 100], [16, 24, 45, 54, 72, 110], [18, 23, 32, 35, 59, 106], [2, 5, 29, 76, 99, 104], [2, 19, 27, 48, 57, 75], [28, 51, 56, 65, 68, 92], [0, 38, 55, 63, 84, 93], [4, 26, 40, 66, 82, 103], [10, 33, 38, 47, 50, 74], [26, 43, 51, 72, 81, 99], [7, 19, 41, 55, 81, 97], [18, 34, 55, 67, 89, 103], [20, 67, 90, 95, 104, 107], [8, 22, 48, 64, 85, 97], [25, 48, 53, 62, 65, 89], [13, 21, 42, 51, 69, 107], [8, 25, 33, 54, 63, 81], [0, 16, 37, 49, 71, 85], [3, 41, 58, 66, 87, 96], [7, 29, 43, 69, 85, 106], [7, 15, 36, 45, 63, 101], [13, 34, 46, 68, 82, 108], [16, 39, 44, 53, 56, 80], [2, 49, 72, 77, 86, 89], [14, 28, 54, 70, 91, 103], [3, 21, 59, 76, 84, 105], [8, 55, 78, 83, 92, 95], [35, 52, 60, 81, 90, 108], [43, 66, 71, 80, 83, 107], [6, 15, 33, 71, 88, 96], [5, 52, 75, 80, 89, 92], [3, 19, 40, 52, 74, 88], [20, 37, 45, 66, 75, 93], [18, 27, 45, 83, 100, 108], [16, 28, 50, 64, 90, 106], [17, 64, 87, 92, 101, 104], [5, 19, 45, 61, 82, 94], [32, 49, 57, 78, 87, 105], [0, 18, 56, 73, 81, 102], [2, 26, 73, 96, 101, 110], [0, 21, 30, 48, 86, 103], [12, 21, 39, 77, 94, 102], [1, 23, 37, 63, 79, 100], [1, 22, 34, 56, 70, 96], [23, 40, 48, 69, 78, 96], [4, 16, 38, 52, 78, 94], [1, 9, 30, 39, 57, 95], [7, 28, 40, 62, 76, 102], [12, 28, 49, 61, 83, 97], [23, 70, 93, 98, 107, 110], [0, 5, 14, 17, 41, 88], [11, 58, 81, 86, 95, 98], [3, 8, 17, 20, 44, 91], [31, 54, 59, 68, 71, 95], [21, 26, 35, 38, 62, 109], [4, 25, 37, 59, 73, 99], [15, 20, 29, 32, 56, 103], [10, 36, 52, 73, 85, 107], [11, 28, 36, 57, 66, 84], [7, 30, 35, 44, 47, 71], [6, 44, 61, 69, 90, 99], [46, 69, 74, 83, 86, 110], [6, 27, 36, 54, 92, 109], [1, 13, 35, 49, 75, 91], [12, 50, 67, 75, 96, 105], [8, 11, 35, 82, 105, 110], [22, 45, 50, 59, 62, 86], [6, 22, 43, 55, 77, 91], [1, 24, 29, 38, 41, 65], [34, 57, 62, 71, 74, 98], [24, 40, 61, 73, 95, 109], [10, 31, 43, 65, 79, 105], [37, 60, 65, 74, 77, 101], [9, 18, 36, 74, 91, 99], [13, 39, 55, 76, 88, 110], [9, 14, 23, 26, 50, 97], [9, 25, 46, 58, 80, 94], [12, 17, 26, 29, 53, 100], [19, 42, 47, 56, 59, 83], [20, 34, 60, 76, 97, 109], [40, 63, 68, 77, 80, 104], [13, 36, 41, 50, 53, 77], [4, 30, 46, 67, 79, 101], [13, 25, 47, 61, 87, 103], [15, 24, 42, 80, 97, 105], [2, 11, 14, 38, 85, 108], [29, 46, 54, 75, 84, 102], [10, 18, 39, 48, 66, 104], [4, 27, 32, 41, 44, 68], [14, 31, 39, 60, 69, 87], [15, 53, 70, 78, 99, 108], [0, 9, 27, 65, 82, 90], [10, 32, 46, 72, 88, 109], [3, 12, 30, 68, 85, 93], [15, 31, 52, 64, 86, 100], [4, 12, 33, 42, 60, 98], [14, 61, 84, 89, 98, 101], [7, 33, 49, 70, 82, 104], [17, 31, 57, 73, 94, 106], [38, 44, 59, 87, 97, 110], [0, 40, 91, 97, 101, 106], [30, 42, 89, 96, 102, 109], [22, 28, 32, 37, 42, 82], [31, 82, 88, 92, 97, 102], [22, 73, 79, 83, 88, 93], [27, 39, 86, 93, 99, 106], [10, 16, 20, 25, 30, 70], [20, 26, 41, 69, 79, 92], [13, 64, 70, 74, 79, 84], [19, 70, 76, 80, 85, 90], [3, 10, 42, 54, 101, 108], [0, 7, 39, 51, 98, 105], [41, 48, 54, 61, 93, 105], [4, 36, 48, 95, 102, 108], [6, 18, 65, 72, 78, 85], [37, 88, 94, 98, 103, 108], [21, 31, 44, 83, 89, 104], [8, 15, 21, 28, 60, 72], [5, 20, 48, 58, 71, 110], [0, 47, 54, 60, 67, 99], [4, 8, 13, 18, 58, 109], [46, 52, 56, 61, 66, 106], [9, 21, 68, 75, 81, 88], [11, 50, 56, 71, 99, 109], [38, 45, 51, 58, 90, 102], [24, 36, 83, 90, 96, 103], [44, 51, 57, 64, 96, 108], [14, 42, 52, 65, 104, 110], [17, 23, 38, 66, 76, 89], [8, 47, 53, 68, 96, 106], [34, 40, 44, 49, 54, 94], [2, 41, 47, 62, 90, 100], [32, 38, 53, 81, 91, 104], [23, 29, 44, 72, 82, 95], [19, 25, 29, 34, 39, 79], [2, 8, 23, 51, 61, 74], [1, 6, 46, 97, 103, 107], [7, 13, 17, 22, 27, 67], [10, 61, 67, 71, 76, 81], [2, 9, 15, 22, 54, 66], [18, 30, 77, 84, 90, 97], [1, 5, 10, 15, 55, 106], [15, 27, 74, 81, 87, 94], [7, 58, 64, 68, 73, 78], [12, 24, 71, 78, 84, 91], [18, 28, 41, 80, 86, 101], [49, 55, 59, 64, 69, 109], [5, 44, 50, 65, 93, 103], [11, 39, 49, 62, 101, 107], [43, 49, 53, 58, 63, 103], [32, 39, 45, 52, 84, 96], [5, 11, 26, 54, 64, 77], [6, 53, 60, 66, 73, 105], [11, 17, 32, 60, 70, 83], [16, 67, 73, 77, 82, 87], [34, 85, 91, 95, 100, 105], [16, 22, 26, 31, 36, 76], [21, 33, 80, 87, 93, 100], [24, 34, 47, 86, 92, 107], [4, 55, 61, 65, 70, 75], [26, 33, 39, 46, 78, 90], [0, 10, 23, 62, 68, 83], [3, 13, 26, 65, 71, 86], [12, 22, 35, 74, 80, 95], [1, 7, 11, 16, 21, 61], [17, 24, 30, 37, 69, 81], [6, 12, 19, 51, 63, 110], [13, 19, 23, 28, 33, 73], [25, 76, 82, 86, 91, 96], [2, 17, 45, 55, 68, 107], [6, 16, 29, 68, 74, 89], [14, 20, 35, 63, 73, 86], [4, 17, 56, 62, 77, 105], [15, 25, 38, 77, 83, 98], [9, 19, 32, 71, 77, 92], [25, 31, 35, 40, 45, 85], [0, 6, 13, 45, 57, 104], [31, 37, 41, 46, 51, 91], [2, 30, 40, 53, 92, 98], [37, 43, 47, 52, 57, 97], [29, 35, 50, 78, 88, 101], [14, 21, 27, 34, 66, 78], [35, 41, 56, 84, 94, 107], [23, 30, 36, 43, 75, 87], [4, 10, 14, 19, 24, 64], [7, 20, 59, 65, 80, 108], [8, 36, 46, 59, 98, 104], [3, 9, 16, 48, 60, 107], [5, 12, 18, 25, 57, 69], [11, 18, 24, 31, 63, 75], [8, 14, 29, 57, 67, 80], [3, 50, 57, 63, 70, 102], [3, 43, 94, 100, 104, 109], [3, 15, 62, 69, 75, 82], [1, 33, 45, 92, 99, 105], [26, 32, 47, 75, 85, 98], [2, 7, 12, 52, 103, 109], [40, 46, 50, 55, 60, 100], [1, 14, 53, 59, 74, 102], [28, 34, 38, 43, 48, 88], [27, 37, 50, 89, 95, 110], [35, 42, 48, 55, 87, 99], [9, 56, 63, 69, 76, 108], [0, 12, 59, 66, 72, 79], [28, 79, 85, 89, 94, 99], [29, 36, 42, 49, 81, 93], [1, 52, 58, 62, 67, 72], [20, 27, 33, 40, 72, 84], [5, 33, 43, 56, 95, 101], [4, 9, 49, 100, 106, 110], [25, 41, 43, 78, 102, 110], [30, 54, 62, 88, 104, 106], [33, 57, 65, 91, 107, 109], [4, 20, 22, 57, 81, 89], [6, 14, 40, 56, 58, 93], [21, 29, 55, 71, 73, 108], [2, 4, 39, 63, 71, 97], [9, 33, 41, 67, 83, 85], [0, 8, 34, 50, 52, 87], [12, 36, 44, 70, 86, 88], [12, 20, 46, 62, 64, 99], [5, 7, 42, 66, 74, 100], [2, 28, 44, 46, 81, 105], [0, 24, 32, 58, 74, 76], [3, 27, 35, 61, 77, 79], [10, 26, 28, 63, 87, 95], [15, 39, 47, 73, 89, 91], [18, 42, 50, 76, 92, 94], [13, 29, 31, 66, 90, 98], [16, 32, 34, 69, 93, 101], [15, 23, 49, 65, 67, 102], [6, 30, 38, 64, 80, 82], [21, 45, 53, 79, 95, 97], [1, 17, 19, 54, 78, 86], [22, 38, 40, 75, 99, 107], [7, 23, 25, 60, 84, 92], [19, 35, 37, 72, 96, 104], [8, 10, 45, 69, 77, 103], [11, 13, 48, 72, 80, 106], [14, 16, 51, 75, 83, 109], [24, 48, 56, 82, 98, 100], [3, 11, 37, 53, 55, 90], [27, 51, 59, 85, 101, 103], [5, 31, 47, 49, 84, 108], [9, 17, 43, 59, 61, 96], [18, 26, 52, 68, 70, 105], [1, 36, 60, 68, 94, 110]]
\item 1 \{1=36408, 2=577422, 3=2175156, 4=2339214\} [[8, 39, 40, 41, 42, 70], [0, 1, 2, 3, 31, 80], [30, 31, 32, 33, 61, 110], [6, 7, 8, 9, 37, 86], [23, 54, 55, 56, 57, 85], [15, 16, 17, 18, 46, 95], [27, 28, 29, 30, 58, 107], [24, 25, 26, 27, 55, 104], [26, 57, 58, 59, 60, 88], [29, 60, 61, 62, 63, 91], [25, 74, 105, 106, 107, 108], [41, 72, 73, 74, 75, 103], [7, 56, 87, 88, 89, 90], [10, 59, 90, 91, 92, 93], [13, 62, 93, 94, 95, 96], [22, 71, 102, 103, 104, 105], [9, 10, 11, 12, 40, 89], [44, 75, 76, 77, 78, 106], [18, 19, 20, 21, 49, 98], [11, 42, 43, 44, 45, 73], [3, 4, 5, 6, 34, 83], [16, 65, 96, 97, 98, 99], [19, 68, 99, 100, 101, 102], [17, 48, 49, 50, 51, 79], [47, 78, 79, 80, 81, 109], [20, 51, 52, 53, 54, 82], [12, 13, 14, 15, 43, 92], [0, 28, 77, 108, 109, 110], [4, 53, 84, 85, 86, 87], [1, 50, 81, 82, 83, 84], [32, 63, 64, 65, 66, 94], [2, 33, 34, 35, 36, 64], [21, 22, 23, 24, 52, 101], [38, 69, 70, 71, 72, 100], [14, 45, 46, 47, 48, 76], [35, 66, 67, 68, 69, 97], [5, 36, 37, 38, 39, 67], [0, 52, 60, 95, 100, 106], [10, 27, 35, 41, 54, 99], [39, 61, 78, 86, 92, 105], [24, 28, 35, 50, 86, 94], [14, 22, 63, 67, 74, 89], [0, 4, 11, 26, 62, 70], [20, 28, 69, 73, 80, 95], [46, 54, 89, 94, 100, 105], [21, 25, 32, 47, 83, 91], [30, 52, 69, 77, 83, 96], [16, 24, 59, 64, 70, 75], [26, 34, 75, 79, 86, 101], [39, 43, 50, 65, 101, 109], [25, 33, 68, 73, 79, 84], [9, 54, 76, 93, 101, 107], [12, 16, 23, 38, 74, 82], [13, 21, 56, 61, 67, 72], [19, 27, 62, 67, 73, 78], [8, 44, 52, 93, 97, 104], [6, 14, 20, 33, 78, 100], [34, 42, 77, 82, 88, 93], [29, 34, 40, 45, 97, 105], [17, 25, 66, 70, 77, 92], [43, 51, 86, 91, 97, 102], [3, 7, 14, 29, 65, 73], [1, 9, 44, 49, 55, 60], [20, 25, 31, 36, 88, 96], [11, 47, 55, 96, 100, 107], [12, 20, 26, 39, 84, 106], [23, 28, 34, 39, 91, 99], [2, 10, 51, 55, 62, 77], [24, 46, 63, 71, 77, 90], [30, 34, 41, 56, 92, 100], [23, 31, 72, 76, 83, 98], [7, 48, 52, 59, 74, 110], [7, 15, 50, 55, 61, 66], [5, 20, 56, 64, 105, 109], [4, 21, 29, 35, 48, 93], [2, 15, 60, 82, 99, 107], [33, 37, 44, 59, 95, 103], [15, 37, 54, 62, 68, 81], [1, 18, 26, 32, 45, 90], [31, 39, 74, 79, 85, 90], [6, 28, 45, 53, 59, 72], [27, 31, 38, 53, 89, 97], [2, 17, 53, 61, 102, 106], [5, 13, 54, 58, 65, 80], [4, 45, 49, 56, 71, 107], [6, 41, 46, 52, 57, 109], [11, 19, 60, 64, 71, 86], [2, 7, 13, 18, 70, 78], [6, 51, 73, 90, 98, 104], [5, 18, 63, 85, 102, 110], [49, 57, 92, 97, 103, 108], [1, 7, 12, 64, 72, 107], [3, 11, 17, 30, 75, 97], [3, 38, 43, 49, 54, 106], [3, 48, 70, 87, 95, 101], [37, 45, 80, 85, 91, 96], [2, 38, 46, 87, 91, 98], [1, 8, 23, 59, 67, 108], [6, 10, 17, 32, 68, 76], [5, 10, 16, 21, 73, 81], [11, 16, 22, 27, 79, 87], [2, 8, 21, 66, 88, 105], [3, 25, 42, 50, 56, 69], [29, 37, 78, 82, 89, 104], [0, 35, 40, 46, 51, 103], [35, 43, 84, 88, 95, 110], [4, 10, 15, 67, 75, 110], [4, 12, 47, 52, 58, 63], [21, 43, 60, 68, 74, 87], [28, 36, 71, 76, 82, 87], [36, 40, 47, 62, 98, 106], [3, 55, 63, 98, 103, 109], [10, 18, 53, 58, 64, 69], [42, 64, 81, 89, 95, 108], [12, 57, 79, 96, 104, 110], [5, 11, 24, 69, 91, 108], [22, 30, 65, 70, 76, 81], [14, 19, 25, 30, 82, 90], [32, 40, 81, 85, 92, 107], [13, 30, 38, 44, 57, 102], [12, 34, 51, 59, 65, 78], [33, 55, 72, 80, 86, 99], [19, 36, 44, 50, 63, 108], [0, 22, 39, 47, 53, 66], [8, 13, 19, 24, 76, 84], [9, 13, 20, 35, 71, 79], [14, 50, 58, 99, 103, 110], [18, 40, 57, 65, 71, 84], [16, 33, 41, 47, 60, 105], [4, 9, 61, 69, 104, 109], [15, 23, 29, 42, 87, 109], [7, 24, 32, 38, 51, 96], [9, 31, 48, 56, 62, 75], [17, 22, 28, 33, 85, 93], [26, 31, 37, 42, 94, 102], [1, 6, 58, 66, 101, 106], [15, 19, 26, 41, 77, 85], [1, 42, 46, 53, 68, 104], [0, 8, 14, 27, 72, 94], [40, 48, 83, 88, 94, 99], [32, 37, 43, 48, 100, 108], [9, 17, 23, 36, 81, 103], [8, 16, 57, 61, 68, 83], [36, 58, 75, 83, 89, 102], [27, 49, 66, 74, 80, 93], [18, 22, 29, 44, 80, 88], [0, 45, 67, 84, 92, 98], [5, 41, 49, 90, 94, 101], [12, 17, 41, 67, 71, 83], [10, 28, 31, 63, 70, 104], [3, 12, 18, 36, 86, 100], [32, 46, 60, 69, 75, 93], [1, 19, 22, 54, 61, 95], [15, 22, 56, 73, 91, 94], [16, 20, 32, 72, 77, 101], [17, 34, 52, 55, 87, 94], [15, 65, 79, 93, 102, 108], [5, 31, 35, 47, 87, 92], [7, 10, 42, 49, 83, 100], [0, 50, 64, 78, 87, 93], [3, 53, 67, 81, 90, 96], [30, 37, 71, 88, 106, 109], [16, 34, 37, 69, 76, 110], [12, 21, 27, 45, 95, 109], [24, 29, 53, 79, 83, 95], [20, 34, 48, 57, 63, 81], [1, 15, 24, 30, 48, 98], [23, 37, 51, 60, 66, 84], [21, 28, 62, 79, 97, 100], [0, 9, 15, 33, 83, 97], [11, 51, 56, 80, 106, 110], [15, 20, 44, 70, 74, 86], [26, 43, 61, 64, 96, 103], [3, 21, 71, 85, 99, 108], [44, 58, 72, 81, 87, 105], [13, 16, 48, 55, 89, 106], [13, 31, 34, 66, 73, 107], [4, 8, 20, 60, 65, 89], [11, 37, 41, 53, 93, 98], [22, 26, 38, 78, 83, 107], [10, 14, 26, 66, 71, 95], [2, 28, 32, 44, 84, 89], [23, 40, 58, 61, 93, 100], [39, 44, 68, 94, 98, 110], [24, 31, 65, 82, 100, 103], [23, 49, 53, 65, 105, 110], [0, 7, 41, 58, 76, 79], [5, 22, 40, 43, 75, 82], [7, 25, 28, 60, 67, 101], [8, 22, 36, 45, 51, 69], [6, 15, 21, 39, 89, 103], [0, 18, 68, 82, 96, 105], [10, 13, 45, 52, 86, 103], [9, 16, 50, 67, 85, 88], [20, 46, 50, 62, 102, 107], [4, 22, 25, 57, 64, 98], [6, 11, 35, 61, 65, 77], [4, 18, 27, 33, 51, 101], [10, 24, 33, 39, 57, 107], [38, 52, 66, 75, 81, 99], [8, 34, 38, 50, 90, 95], [9, 59, 73, 87, 96, 102], [6, 12, 30, 80, 94, 108], [25, 29, 41, 81, 86, 110], [14, 31, 49, 52, 84, 91], [35, 49, 63, 72, 78, 96], [17, 43, 47, 59, 99, 104], [11, 28, 46, 49, 81, 88], [7, 21, 30, 36, 54, 104], [9, 14, 38, 64, 68, 80], [21, 26, 50, 76, 80, 92], [2, 16, 30, 39, 45, 63], [27, 32, 56, 82, 86, 98], [11, 25, 39, 48, 54, 72], [26, 40, 54, 63, 69, 87], [1, 5, 17, 57, 62, 86], [13, 17, 29, 69, 74, 98], [41, 55, 69, 78, 84, 102], [36, 41, 65, 91, 95, 107], [13, 27, 36, 42, 60, 110], [6, 13, 47, 64, 82, 85], [3, 9, 27, 77, 91, 105], [4, 7, 39, 46, 80, 97], [9, 18, 24, 42, 92, 106], [2, 26, 52, 56, 68, 108], [8, 48, 53, 77, 103, 107], [33, 38, 62, 88, 92, 104], [3, 10, 44, 61, 79, 82], [18, 25, 59, 76, 94, 97], [0, 5, 29, 55, 59, 71], [19, 23, 35, 75, 80, 104], [14, 28, 42, 51, 57, 75], [1, 33, 40, 74, 91, 109], [17, 31, 45, 54, 60, 78], [29, 46, 64, 67, 99, 106], [29, 43, 57, 66, 72, 90], [1, 35, 52, 70, 73, 105], [18, 23, 47, 73, 77, 89], [12, 62, 76, 90, 99, 105], [5, 45, 50, 74, 100, 104], [2, 42, 47, 71, 97, 101], [7, 11, 23, 63, 68, 92], [20, 37, 55, 58, 90, 97], [32, 49, 67, 70, 102, 109], [3, 8, 32, 58, 62, 74], [12, 19, 53, 70, 88, 91], [6, 56, 70, 84, 93, 99], [47, 61, 75, 84, 90, 108], [16, 19, 51, 58, 92, 109], [0, 6, 24, 74, 88, 102], [2, 19, 37, 40, 72, 79], [8, 25, 43, 46, 78, 85], [14, 40, 44, 56, 96, 101], [27, 34, 68, 85, 103, 106], [4, 38, 55, 73, 76, 108], [2, 14, 54, 59, 83, 109], [5, 19, 33, 42, 48, 66], [30, 35, 59, 85, 89, 101], [1, 4, 36, 43, 77, 94], [0, 13, 81, 91, 101, 104], [66, 76, 86, 89, 96, 109], [9, 19, 29, 32, 39, 52], [60, 70, 80, 83, 90, 103], [63, 73, 83, 86, 93, 106], [2, 9, 22, 90, 100, 110], [0, 10, 20, 23, 30, 43], [24, 34, 44, 47, 54, 67], [54, 64, 74, 77, 84, 97], [3, 16, 84, 94, 104, 107], [21, 31, 41, 44, 51, 64], [12, 22, 32, 35, 42, 55], [18, 28, 38, 41, 48, 61], [2, 5, 12, 25, 93, 103], [30, 40, 50, 53, 60, 73], [33, 43, 53, 56, 63, 76], [39, 49, 59, 62, 69, 82], [1, 69, 79, 89, 92, 99], [6, 16, 26, 29, 36, 49], [5, 8, 15, 28, 96, 106], [8, 11, 18, 31, 99, 109], [51, 61, 71, 74, 81, 94], [36, 46, 56, 59, 66, 79], [48, 58, 68, 71, 78, 91], [45, 55, 65, 68, 75, 88], [1, 11, 14, 21, 34, 102], [42, 52, 62, 65, 72, 85], [4, 14, 17, 24, 37, 105], [57, 67, 77, 80, 87, 100], [7, 75, 85, 95, 98, 105], [10, 78, 88, 98, 101, 108], [15, 25, 35, 38, 45, 58], [27, 37, 47, 50, 57, 70], [7, 17, 20, 27, 40, 108], [4, 72, 82, 92, 95, 102], [6, 19, 87, 97, 107, 110], [3, 13, 23, 26, 33, 46], [1, 27, 39, 71, 75, 96], [11, 33, 58, 67, 82, 94], [1, 13, 41, 63, 88, 97], [6, 31, 40, 55, 67, 95], [5, 51, 68, 70, 89, 107], [16, 25, 40, 52, 80, 102], [8, 10, 29, 47, 56, 102], [1, 29, 51, 76, 85, 100], [2, 48, 65, 67, 86, 104], [18, 35, 37, 56, 74, 83], [2, 6, 27, 43, 69, 81], [0, 21, 37, 63, 75, 107], [14, 23, 69, 86, 88, 107], [15, 32, 34, 53, 71, 80], [11, 29, 38, 84, 101, 103], [9, 25, 51, 63, 95, 99], [26, 30, 51, 67, 93, 105], [13, 25, 53, 75, 100, 109], [2, 11, 57, 74, 76, 95], [7, 35, 57, 82, 91, 106], [10, 19, 34, 46, 74, 96], [0, 12, 44, 48, 69, 85], [8, 12, 33, 49, 75, 87], [1, 20, 38, 47, 93, 110], [3, 15, 47, 51, 72, 88], [6, 18, 50, 54, 75, 91], [19, 28, 43, 55, 83, 105], [18, 34, 60, 72, 104, 108], [39, 56, 58, 77, 95, 104], [5, 23, 32, 78, 95, 97], [3, 19, 45, 57, 89, 93], [8, 30, 55, 64, 79, 91], [14, 18, 39, 55, 81, 93], [4, 19, 31, 59, 81, 106], [5, 27, 52, 61, 76, 88], [30, 47, 49, 68, 86, 95], [6, 22, 48, 60, 92, 96], [11, 13, 32, 50, 59, 105], [1, 16, 28, 56, 78, 103], [1, 10, 25, 37, 65, 87], [0, 25, 34, 49, 61, 89], [2, 24, 49, 58, 73, 85], [14, 32, 41, 87, 104, 106], [3, 28, 37, 52, 64, 92], [17, 39, 64, 73, 88, 100], [21, 46, 55, 70, 82, 110], [7, 22, 34, 62, 84, 109], [3, 20, 22, 41, 59, 68], [9, 34, 43, 58, 70, 98], [6, 23, 25, 44, 62, 71], [29, 33, 54, 70, 96, 108], [10, 36, 48, 80, 84, 105], [8, 54, 71, 73, 92, 110], [3, 35, 39, 60, 76, 102], [23, 45, 70, 79, 94, 106], [15, 40, 49, 64, 76, 104], [4, 32, 54, 79, 88, 103], [0, 17, 19, 38, 56, 65], [23, 27, 48, 64, 90, 102], [21, 38, 40, 59, 77, 86], [10, 38, 60, 85, 94, 109], [8, 17, 63, 80, 82, 101], [18, 43, 52, 67, 79, 107], [13, 22, 37, 49, 77, 99], [5, 14, 60, 77, 79, 98], [12, 29, 31, 50, 68, 77], [36, 53, 55, 74, 92, 101], [11, 15, 36, 52, 78, 90], [18, 30, 62, 66, 87, 103], [4, 30, 42, 74, 78, 99], [7, 33, 45, 77, 81, 102], [20, 42, 67, 76, 91, 103], [5, 9, 30, 46, 72, 84], [0, 32, 36, 57, 73, 99], [20, 24, 45, 61, 87, 99], [3, 24, 40, 66, 78, 110], [14, 16, 35, 53, 62, 108], [12, 24, 56, 60, 81, 97], [12, 28, 54, 66, 98, 102], [6, 38, 42, 63, 79, 105], [9, 21, 53, 57, 78, 94], [21, 33, 65, 69, 90, 106], [24, 36, 68, 72, 93, 109], [9, 41, 45, 66, 82, 108], [24, 41, 43, 62, 80, 89], [42, 59, 61, 80, 98, 107], [4, 13, 28, 40, 68, 90], [13, 39, 51, 83, 87, 108], [22, 31, 46, 58, 86, 108], [12, 37, 46, 61, 73, 101], [5, 7, 26, 44, 53, 99], [9, 26, 28, 47, 65, 74], [2, 20, 29, 75, 92, 94], [15, 27, 59, 63, 84, 100], [15, 31, 57, 69, 101, 105], [8, 26, 35, 81, 98, 100], [2, 4, 23, 41, 50, 96], [27, 44, 46, 65, 83, 92], [14, 36, 61, 70, 85, 97], [11, 20, 66, 83, 85, 104], [26, 48, 73, 82, 97, 109], [7, 16, 31, 43, 71, 93], [4, 16, 44, 66, 91, 100], [7, 19, 47, 69, 94, 103], [17, 21, 42, 58, 84, 96], [45, 62, 64, 83, 101, 110], [33, 50, 52, 71, 89, 98], [17, 35, 44, 90, 107, 109], [0, 16, 42, 54, 86, 90], [10, 22, 50, 72, 97, 106], [17, 26, 72, 89, 91, 110]]
\item 1 \{1=27528, 2=520146, 3=2186700, 4=2393826\} [[8, 39, 40, 41, 42, 70], [0, 1, 2, 3, 31, 80], [30, 31, 32, 33, 61, 110], [6, 7, 8, 9, 37, 86], [23, 54, 55, 56, 57, 85], [15, 16, 17, 18, 46, 95], [27, 28, 29, 30, 58, 107], [24, 25, 26, 27, 55, 104], [26, 57, 58, 59, 60, 88], [29, 60, 61, 62, 63, 91], [25, 74, 105, 106, 107, 108], [41, 72, 73, 74, 75, 103], [7, 56, 87, 88, 89, 90], [10, 59, 90, 91, 92, 93], [13, 62, 93, 94, 95, 96], [22, 71, 102, 103, 104, 105], [9, 10, 11, 12, 40, 89], [44, 75, 76, 77, 78, 106], [18, 19, 20, 21, 49, 98], [11, 42, 43, 44, 45, 73], [3, 4, 5, 6, 34, 83], [16, 65, 96, 97, 98, 99], [19, 68, 99, 100, 101, 102], [17, 48, 49, 50, 51, 79], [47, 78, 79, 80, 81, 109], [20, 51, 52, 53, 54, 82], [12, 13, 14, 15, 43, 92], [0, 28, 77, 108, 109, 110], [4, 53, 84, 85, 86, 87], [1, 50, 81, 82, 83, 84], [32, 63, 64, 65, 66, 94], [2, 33, 34, 35, 36, 64], [21, 22, 23, 24, 52, 101], [38, 69, 70, 71, 72, 100], [14, 45, 46, 47, 48, 76], [35, 66, 67, 68, 69, 97], [5, 36, 37, 38, 39, 67], [6, 55, 75, 79, 92, 110], [40, 60, 64, 77, 95, 102], [25, 45, 49, 62, 80, 87], [5, 12, 61, 81, 85, 98], [3, 7, 20, 38, 45, 94], [11, 29, 36, 85, 105, 109], [13, 33, 37, 50, 68, 75], [1, 14, 32, 39, 88, 108], [22, 42, 46, 59, 77, 84], [6, 10, 23, 41, 48, 97], [18, 22, 35, 53, 60, 109], [5, 23, 30, 79, 99, 103], [37, 57, 61, 74, 92, 99], [46, 66, 70, 83, 101, 108], [9, 13, 26, 44, 51, 100], [10, 30, 34, 47, 65, 72], [12, 16, 29, 47, 54, 103], [11, 18, 67, 87, 91, 104], [4, 24, 28, 41, 59, 66], [2, 9, 58, 78, 82, 95], [19, 39, 43, 56, 74, 81], [8, 26, 33, 82, 102, 106], [2, 20, 27, 76, 96, 100], [3, 52, 72, 76, 89, 107], [43, 63, 67, 80, 98, 105], [16, 36, 40, 53, 71, 78], [28, 48, 52, 65, 83, 90], [1, 21, 25, 38, 56, 63], [15, 19, 32, 50, 57, 106], [7, 27, 31, 44, 62, 69], [14, 21, 70, 90, 94, 107], [31, 51, 55, 68, 86, 93], [0, 4, 17, 35, 42, 91], [8, 15, 64, 84, 88, 101], [17, 24, 73, 93, 97, 110], [34, 54, 58, 71, 89, 96], [0, 49, 69, 73, 86, 104], [4, 31, 38, 50, 78, 99], [8, 13, 17, 34, 63, 99], [9, 45, 65, 70, 74, 91], [23, 46, 58, 81, 86, 90], [6, 42, 62, 67, 71, 88], [18, 38, 43, 47, 64, 93], [4, 33, 69, 89, 94, 98], [24, 44, 49, 53, 70, 99], [13, 20, 32, 60, 81, 97], [0, 20, 25, 29, 46, 75], [2, 19, 48, 84, 104, 109], [7, 34, 41, 53, 81, 102], [33, 53, 58, 62, 79, 108], [0, 5, 9, 53, 76, 88], [22, 29, 41, 69, 90, 106], [21, 26, 30, 74, 97, 109], [0, 44, 67, 79, 102, 107], [1, 8, 20, 48, 69, 85], [18, 34, 61, 68, 80, 108], [4, 16, 39, 44, 48, 92], [2, 14, 42, 63, 79, 106], [20, 43, 55, 78, 83, 87], [30, 50, 55, 59, 76, 105], [17, 22, 26, 43, 72, 108], [2, 25, 37, 60, 65, 69], [13, 42, 78, 98, 103, 107], [12, 32, 37, 41, 58, 87], [1, 30, 66, 86, 91, 95], [8, 31, 43, 66, 71, 75], [10, 22, 45, 50, 54, 98], [15, 51, 71, 76, 80, 97], [10, 39, 75, 95, 100, 104], [7, 19, 42, 47, 51, 95], [13, 40, 47, 59, 87, 108], [7, 36, 72, 92, 97, 101], [11, 34, 46, 69, 74, 78], [32, 55, 67, 90, 95, 99], [4, 27, 32, 36, 80, 103], [18, 39, 55, 82, 89, 101], [2, 6, 50, 73, 85, 108], [3, 19, 46, 53, 65, 93], [21, 57, 77, 82, 86, 103], [16, 23, 35, 63, 84, 100], [3, 23, 28, 32, 49, 78], [9, 14, 18, 62, 85, 97], [3, 39, 59, 64, 68, 85], [6, 22, 49, 56, 68, 96], [9, 29, 34, 38, 55, 84], [4, 8, 25, 54, 90, 110], [19, 31, 54, 59, 63, 107], [0, 16, 43, 50, 62, 90], [19, 26, 38, 66, 87, 103], [12, 28, 55, 62, 74, 102], [22, 34, 57, 62, 66, 110], [3, 24, 40, 67, 74, 86], [7, 14, 26, 54, 75, 91], [6, 26, 31, 35, 52, 81], [26, 49, 61, 84, 89, 93], [6, 11, 15, 59, 82, 94], [18, 54, 74, 79, 83, 100], [2, 30, 51, 67, 94, 101], [14, 19, 23, 40, 69, 105], [38, 61, 73, 96, 101, 105], [15, 35, 40, 44, 61, 90], [0, 21, 37, 64, 71, 83], [5, 10, 14, 31, 60, 96], [15, 20, 24, 68, 91, 103], [21, 41, 46, 50, 67, 96], [7, 30, 35, 39, 83, 106], [5, 28, 40, 63, 68, 72], [29, 52, 64, 87, 92, 96], [10, 37, 44, 56, 84, 105], [12, 33, 49, 76, 83, 95], [9, 25, 52, 59, 71, 99], [3, 8, 12, 56, 79, 91], [2, 7, 11, 28, 57, 93], [18, 23, 27, 71, 94, 106], [35, 58, 70, 93, 98, 102], [6, 27, 43, 70, 77, 89], [11, 39, 60, 76, 103, 110], [21, 42, 58, 85, 92, 104], [8, 36, 57, 73, 100, 107], [24, 45, 61, 88, 95, 107], [16, 45, 81, 101, 106, 110], [9, 30, 46, 73, 80, 92], [4, 11, 23, 51, 72, 88], [27, 63, 83, 88, 92, 109], [41, 64, 76, 99, 104, 108], [12, 17, 21, 65, 88, 100], [15, 31, 58, 65, 77, 105], [1, 24, 29, 33, 77, 100], [16, 28, 51, 56, 60, 104], [0, 36, 56, 61, 65, 82], [27, 48, 64, 91, 98, 110], [24, 60, 80, 85, 89, 106], [5, 33, 54, 70, 97, 104], [25, 32, 44, 72, 93, 109], [5, 17, 45, 66, 82, 109], [13, 25, 48, 53, 57, 101], [14, 37, 49, 72, 77, 81], [1, 5, 22, 51, 87, 107], [11, 16, 20, 37, 66, 102], [3, 47, 70, 82, 105, 110], [1, 13, 36, 41, 45, 89], [12, 48, 68, 73, 77, 94], [15, 36, 52, 79, 86, 98], [1, 28, 35, 47, 75, 96], [27, 47, 52, 56, 73, 102], [17, 40, 52, 75, 80, 84], [10, 17, 29, 57, 78, 94], [10, 33, 38, 42, 86, 109], [0, 6, 40, 54, 93, 101], [4, 20, 22, 58, 63, 73], [2, 17, 41, 54, 61, 86], [23, 38, 62, 75, 82, 107], [17, 44, 59, 83, 96, 103], [33, 41, 51, 57, 91, 105], [5, 7, 43, 48, 58, 100], [31, 45, 84, 92, 102, 108], [8, 32, 45, 52, 77, 104], [11, 35, 48, 55, 80, 107], [14, 27, 34, 59, 86, 101], [12, 51, 59, 69, 75, 109], [30, 38, 48, 54, 88, 102], [14, 29, 53, 66, 73, 98], [2, 12, 18, 52, 66, 105], [1, 37, 42, 52, 94, 110], [1, 15, 54, 62, 72, 78], [4, 29, 56, 71, 95, 108], [16, 21, 31, 73, 89, 91], [23, 36, 43, 68, 95, 110], [1, 43, 59, 61, 97, 102], [8, 10, 46, 51, 61, 103], [22, 38, 40, 76, 81, 91], [8, 35, 50, 74, 87, 94], [0, 7, 32, 59, 74, 98], [2, 26, 39, 46, 71, 98], [8, 23, 47, 60, 67, 92], [9, 17, 27, 33, 67, 81], [4, 9, 19, 61, 77, 79], [20, 35, 59, 72, 79, 104], [2, 29, 44, 68, 81, 88], [28, 44, 46, 82, 87, 97], [9, 15, 49, 63, 102, 110], [5, 20, 44, 57, 64, 89], [3, 11, 21, 27, 61, 75], [5, 32, 47, 71, 84, 91], [19, 35, 37, 73, 78, 88], [11, 38, 53, 77, 90, 97], [22, 27, 37, 79, 95, 97], [3, 37, 51, 90, 98, 108], [0, 39, 47, 57, 63, 97], [15, 23, 33, 39, 73, 87], [22, 36, 75, 83, 93, 99], [36, 44, 54, 60, 94, 108], [12, 19, 44, 71, 86, 110], [11, 26, 50, 63, 70, 95], [16, 30, 69, 77, 87, 93], [6, 12, 46, 60, 99, 107], [34, 39, 49, 91, 107, 109], [7, 12, 22, 64, 80, 82], [14, 38, 51, 58, 83, 110], [1, 17, 19, 55, 60, 70], [6, 14, 24, 30, 64, 78], [0, 8, 18, 24, 58, 72], [17, 30, 37, 62, 89, 104], [1, 6, 16, 58, 74, 76], [5, 18, 25, 50, 77, 92], [7, 21, 60, 68, 78, 84], [16, 32, 34, 70, 75, 85], [21, 29, 39, 45, 79, 93], [2, 4, 40, 45, 55, 97], [5, 15, 21, 55, 69, 108], [3, 42, 50, 60, 66, 100], [31, 36, 46, 88, 104, 106], [3, 9, 43, 57, 96, 104], [31, 47, 49, 85, 90, 100], [3, 13, 55, 71, 73, 109], [13, 18, 28, 70, 86, 88], [28, 33, 43, 85, 101, 103], [8, 21, 28, 53, 80, 95], [0, 10, 52, 68, 70, 106], [11, 13, 49, 54, 64, 106], [24, 32, 42, 48, 82, 96], [6, 45, 53, 63, 69, 103], [25, 39, 78, 86, 96, 102], [37, 53, 55, 91, 96, 106], [10, 24, 63, 71, 81, 87], [0, 34, 48, 87, 95, 105], [4, 46, 62, 64, 100, 105], [19, 24, 34, 76, 92, 94], [7, 23, 25, 61, 66, 76], [6, 13, 38, 65, 80, 104], [40, 56, 58, 94, 99, 109], [28, 42, 81, 89, 99, 105], [11, 24, 31, 56, 83, 98], [18, 26, 36, 42, 76, 90], [4, 18, 57, 65, 75, 81], [13, 29, 31, 67, 72, 82], [25, 41, 43, 79, 84, 94], [25, 30, 40, 82, 98, 100], [34, 50, 52, 88, 93, 103], [13, 27, 66, 74, 84, 90], [14, 41, 56, 80, 93, 100], [27, 35, 45, 51, 85, 99], [9, 48, 56, 66, 72, 106], [17, 32, 56, 69, 76, 101], [14, 16, 52, 57, 67, 109], [3, 10, 35, 62, 77, 101], [10, 26, 28, 64, 69, 79], [26, 41, 65, 78, 85, 110], [1, 26, 53, 68, 92, 105], [12, 20, 30, 36, 70, 84], [20, 33, 40, 65, 92, 107], [7, 49, 65, 67, 103, 108], [20, 47, 62, 86, 99, 106], [9, 16, 41, 68, 83, 107], [19, 33, 72, 80, 90, 96], [2, 15, 22, 47, 74, 89], [10, 15, 25, 67, 83, 85], [5, 29, 42, 49, 74, 101], [23, 50, 65, 89, 102, 109], [7, 33, 46, 52, 55, 63], [10, 36, 49, 55, 58, 66], [20, 31, 41, 95, 101, 109], [6, 61, 87, 100, 106, 109], [4, 10, 13, 21, 76, 102], [28, 54, 67, 73, 76, 84], [3, 54, 69, 81, 92, 95], [0, 51, 66, 78, 89, 92], [12, 27, 39, 50, 53, 72], [52, 78, 91, 97, 100, 108], [2, 21, 72, 87, 99, 110], [46, 72, 85, 91, 94, 102], [3, 16, 22, 25, 33, 88], [42, 57, 69, 80, 83, 102], [3, 58, 84, 97, 103, 106], [3, 14, 17, 36, 87, 102], [6, 17, 20, 39, 90, 105], [10, 20, 74, 80, 88, 110], [15, 30, 42, 53, 56, 75], [19, 45, 58, 64, 67, 75], [6, 21, 33, 44, 47, 66], [14, 20, 28, 50, 61, 71], [47, 53, 61, 83, 94, 104], [34, 60, 73, 79, 82, 90], [53, 59, 67, 89, 100, 110], [0, 13, 19, 22, 30, 85], [7, 29, 40, 50, 104, 110], [12, 25, 31, 34, 42, 97], [0, 11, 14, 33, 84, 99], [14, 25, 35, 89, 95, 103], [15, 28, 34, 37, 45, 100], [5, 11, 19, 41, 52, 62], [36, 51, 63, 74, 77, 96], [41, 47, 55, 77, 88, 98], [24, 39, 51, 62, 65, 84], [30, 45, 57, 68, 71, 90], [9, 24, 36, 47, 50, 69], [17, 23, 31, 53, 64, 74], [16, 42, 55, 61, 64, 72], [2, 5, 24, 75, 90, 102], [23, 29, 37, 59, 70, 80], [33, 48, 60, 71, 74, 93], [7, 13, 16, 24, 79, 105], [0, 12, 23, 26, 45, 96], [9, 60, 75, 87, 98, 101], [24, 37, 43, 46, 54, 109], [18, 69, 84, 96, 107, 110], [9, 21, 32, 35, 54, 105], [50, 56, 64, 86, 97, 107], [12, 24, 35, 38, 57, 108], [35, 41, 49, 71, 82, 92], [4, 7, 15, 70, 96, 109], [18, 33, 45, 56, 59, 78], [1, 9, 64, 90, 103, 109], [9, 22, 28, 31, 39, 94], [31, 57, 70, 76, 79, 87], [37, 63, 76, 82, 85, 93], [6, 19, 25, 28, 36, 91], [8, 14, 22, 44, 55, 65], [45, 60, 72, 83, 86, 105], [39, 54, 66, 77, 80, 99], [38, 44, 52, 74, 85, 95], [5, 13, 35, 46, 56, 110], [13, 39, 52, 58, 61, 69], [27, 42, 54, 65, 68, 87], [6, 18, 29, 32, 51, 102], [6, 57, 72, 84, 95, 98], [48, 63, 75, 86, 89, 108], [11, 17, 25, 47, 58, 68], [4, 30, 43, 49, 52, 60], [44, 50, 58, 80, 91, 101], [7, 17, 71, 77, 85, 107], [2, 8, 16, 38, 49, 59], [1, 11, 65, 71, 79, 101], [5, 8, 27, 78, 93, 105], [5, 59, 65, 73, 95, 106], [0, 15, 27, 38, 41, 60], [21, 36, 48, 59, 62, 81], [2, 13, 23, 77, 83, 91], [9, 20, 23, 42, 93, 108], [20, 26, 34, 56, 67, 77], [32, 38, 46, 68, 79, 89], [18, 31, 37, 40, 48, 103], [1, 27, 40, 46, 49, 57], [17, 28, 38, 92, 98, 106], [2, 10, 32, 43, 53, 107], [4, 26, 37, 47, 101, 107], [5, 16, 26, 80, 86, 94], [8, 11, 30, 81, 96, 108], [8, 19, 29, 83, 89, 97], [3, 18, 30, 41, 44, 63], [4, 14, 68, 74, 82, 104], [15, 66, 81, 93, 104, 107], [8, 62, 68, 76, 98, 109], [22, 48, 61, 67, 70, 78], [10, 16, 19, 27, 82, 108], [12, 63, 78, 90, 101, 104], [21, 34, 40, 43, 51, 106], [40, 66, 79, 85, 88, 96], [26, 32, 40, 62, 73, 83], [43, 69, 82, 88, 91, 99], [3, 15, 26, 29, 48, 99], [11, 22, 32, 86, 92, 100], [0, 55, 81, 94, 100, 103], [1, 4, 12, 67, 93, 106], [25, 51, 64, 70, 73, 81], [29, 35, 43, 65, 76, 86], [1, 7, 10, 18, 73, 99], [2, 56, 62, 70, 92, 103], [49, 75, 88, 94, 97, 105], [1, 23, 34, 44, 98, 104]]
\item 1 \{1=40848, 2=491508, 3=2207568, 4=2388276\} [[8, 39, 40, 41, 42, 70], [0, 1, 2, 3, 31, 80], [30, 31, 32, 33, 61, 110], [6, 7, 8, 9, 37, 86], [23, 54, 55, 56, 57, 85], [15, 16, 17, 18, 46, 95], [27, 28, 29, 30, 58, 107], [24, 25, 26, 27, 55, 104], [26, 57, 58, 59, 60, 88], [29, 60, 61, 62, 63, 91], [25, 74, 105, 106, 107, 108], [41, 72, 73, 74, 75, 103], [7, 56, 87, 88, 89, 90], [10, 59, 90, 91, 92, 93], [13, 62, 93, 94, 95, 96], [22, 71, 102, 103, 104, 105], [9, 10, 11, 12, 40, 89], [44, 75, 76, 77, 78, 106], [18, 19, 20, 21, 49, 98], [11, 42, 43, 44, 45, 73], [3, 4, 5, 6, 34, 83], [16, 65, 96, 97, 98, 99], [19, 68, 99, 100, 101, 102], [17, 48, 49, 50, 51, 79], [47, 78, 79, 80, 81, 109], [20, 51, 52, 53, 54, 82], [12, 13, 14, 15, 43, 92], [0, 28, 77, 108, 109, 110], [4, 53, 84, 85, 86, 87], [1, 50, 81, 82, 83, 84], [32, 63, 64, 65, 66, 94], [2, 33, 34, 35, 36, 64], [21, 22, 23, 24, 52, 101], [38, 69, 70, 71, 72, 100], [14, 45, 46, 47, 48, 76], [35, 66, 67, 68, 69, 97], [5, 36, 37, 38, 39, 67], [11, 23, 34, 49, 87, 92], [33, 37, 42, 53, 66, 92], [3, 8, 38, 50, 61, 76], [6, 10, 15, 26, 39, 65], [4, 7, 21, 38, 75, 97], [0, 17, 54, 76, 94, 97], [2, 32, 44, 55, 70, 108], [11, 48, 70, 88, 91, 105], [18, 23, 53, 65, 76, 91], [17, 69, 73, 78, 89, 102], [5, 17, 28, 43, 81, 86], [0, 22, 40, 43, 57, 74], [21, 43, 61, 64, 78, 95], [2, 14, 25, 40, 78, 83], [8, 45, 67, 85, 88, 102], [15, 37, 55, 58, 72, 89], [1, 15, 32, 69, 91, 109], [21, 26, 56, 68, 79, 94], [20, 72, 76, 81, 92, 105], [29, 41, 52, 67, 105, 110], [36, 41, 71, 83, 94, 109], [16, 34, 37, 51, 68, 105], [9, 35, 87, 91, 96, 107], [30, 35, 65, 77, 88, 103], [3, 20, 57, 79, 97, 100], [14, 51, 73, 91, 94, 108], [7, 45, 50, 80, 92, 103], [7, 25, 28, 42, 59, 96], [23, 35, 46, 61, 99, 104], [24, 29, 59, 71, 82, 97], [10, 13, 27, 44, 81, 103], [4, 19, 57, 62, 92, 104], [14, 66, 70, 75, 86, 99], [7, 10, 24, 41, 78, 100], [9, 13, 18, 29, 42, 68], [23, 75, 79, 84, 95, 108], [14, 26, 37, 52, 90, 95], [13, 31, 34, 48, 65, 102], [0, 11, 24, 50, 102, 106], [3, 14, 27, 53, 105, 109], [21, 25, 30, 41, 54, 80], [0, 26, 78, 82, 87, 98], [9, 26, 63, 85, 103, 106], [6, 11, 41, 53, 64, 79], [11, 22, 37, 75, 80, 110], [12, 16, 21, 32, 45, 71], [39, 43, 48, 59, 72, 98], [27, 49, 67, 70, 84, 101], [12, 17, 47, 59, 70, 85], [2, 15, 41, 93, 97, 102], [30, 52, 70, 73, 87, 104], [1, 19, 22, 36, 53, 90], [17, 29, 40, 55, 93, 98], [2, 54, 58, 63, 74, 87], [33, 55, 73, 76, 90, 107], [1, 6, 17, 30, 56, 108], [13, 16, 30, 47, 84, 106], [7, 22, 60, 65, 95, 107], [12, 38, 90, 94, 99, 110], [5, 42, 64, 82, 85, 99], [6, 23, 60, 82, 100, 103], [19, 37, 40, 54, 71, 108], [1, 4, 18, 35, 72, 94], [15, 19, 24, 35, 48, 74], [5, 57, 61, 66, 77, 90], [8, 21, 47, 99, 103, 108], [3, 7, 12, 23, 36, 62], [6, 32, 84, 88, 93, 104], [3, 25, 43, 46, 60, 77], [36, 58, 76, 79, 93, 110], [18, 22, 27, 38, 51, 77], [0, 5, 35, 47, 58, 73], [15, 20, 50, 62, 73, 88], [2, 39, 61, 79, 82, 96], [18, 40, 58, 61, 75, 92], [48, 52, 57, 68, 81, 107], [1, 16, 54, 59, 89, 101], [30, 34, 39, 50, 63, 89], [45, 49, 54, 65, 78, 104], [11, 63, 67, 72, 83, 96], [1, 39, 44, 74, 86, 97], [9, 14, 44, 56, 67, 82], [8, 19, 34, 72, 77, 107], [8, 20, 31, 46, 84, 89], [10, 48, 53, 83, 95, 106], [6, 28, 46, 49, 63, 80], [24, 46, 64, 67, 81, 98], [12, 29, 66, 88, 106, 109], [51, 55, 60, 71, 84, 110], [10, 28, 31, 45, 62, 99], [33, 38, 68, 80, 91, 106], [5, 18, 44, 96, 100, 105], [27, 32, 62, 74, 85, 100], [4, 42, 47, 77, 89, 100], [36, 40, 45, 56, 69, 95], [10, 25, 63, 68, 98, 110], [5, 16, 31, 69, 74, 104], [42, 46, 51, 62, 75, 101], [2, 13, 28, 66, 71, 101], [8, 60, 64, 69, 80, 93], [16, 19, 33, 50, 87, 109], [12, 34, 52, 55, 69, 86], [20, 32, 43, 58, 96, 101], [4, 22, 25, 39, 56, 93], [9, 31, 49, 52, 66, 83], [24, 28, 33, 44, 57, 83], [13, 51, 56, 86, 98, 109], [3, 29, 81, 85, 90, 101], [26, 38, 49, 64, 102, 107], [0, 4, 9, 20, 33, 59], [27, 31, 36, 47, 60, 86], [0, 10, 56, 61, 101, 105], [1, 41, 45, 51, 61, 107], [3, 9, 19, 65, 70, 110], [17, 21, 27, 37, 83, 88], [38, 43, 83, 87, 93, 103], [44, 49, 89, 93, 99, 109], [41, 46, 86, 90, 96, 106], [14, 18, 24, 34, 80, 85], [38, 42, 48, 58, 104, 109], [32, 36, 42, 52, 98, 103], [8, 13, 53, 57, 63, 73], [29, 34, 74, 78, 84, 94], [5, 9, 15, 25, 71, 76], [29, 33, 39, 49, 95, 100], [26, 30, 36, 46, 92, 97], [11, 15, 21, 31, 77, 82], [0, 6, 16, 62, 67, 107], [3, 13, 59, 64, 104, 108], [8, 12, 18, 28, 74, 79], [1, 47, 52, 92, 96, 102], [32, 37, 77, 81, 87, 97], [2, 7, 47, 51, 57, 67], [26, 31, 71, 75, 81, 91], [4, 50, 55, 95, 99, 105], [17, 22, 62, 66, 72, 82], [2, 6, 12, 22, 68, 73], [23, 28, 68, 72, 78, 88], [4, 44, 48, 54, 64, 110], [35, 40, 80, 84, 90, 100], [35, 39, 45, 55, 101, 106], [5, 10, 50, 54, 60, 70], [14, 19, 59, 63, 69, 79], [11, 16, 56, 60, 66, 76], [7, 53, 58, 98, 102, 108], [20, 24, 30, 40, 86, 91], [20, 25, 65, 69, 75, 85], [23, 27, 33, 43, 89, 94], [5, 48, 62, 80, 87, 108], [2, 4, 8, 17, 52, 91], [42, 56, 74, 81, 102, 110], [6, 14, 57, 71, 89, 96], [31, 70, 92, 94, 98, 107], [12, 39, 46, 54, 91, 103], [10, 32, 34, 38, 47, 82], [15, 29, 47, 54, 75, 83], [1, 40, 62, 64, 68, 77], [15, 22, 30, 67, 79, 99], [2, 11, 46, 85, 107, 109], [4, 16, 36, 63, 70, 78], [6, 13, 21, 58, 70, 90], [2, 45, 59, 77, 84, 105], [39, 53, 71, 78, 99, 107], [22, 61, 83, 85, 89, 98], [33, 47, 65, 72, 93, 101], [2, 37, 76, 98, 100, 104], [10, 22, 42, 69, 76, 84], [3, 40, 52, 72, 99, 106], [31, 53, 55, 59, 68, 103], [13, 52, 74, 76, 80, 89], [6, 20, 38, 45, 66, 74], [9, 36, 43, 51, 88, 100], [12, 26, 44, 51, 72, 80], [13, 35, 37, 41, 50, 85], [4, 43, 65, 67, 71, 80], [12, 19, 27, 64, 76, 96], [7, 29, 31, 35, 44, 79], [22, 34, 54, 81, 88, 96], [10, 49, 71, 73, 77, 86], [8, 26, 33, 54, 62, 105], [7, 19, 39, 66, 73, 81], [1, 23, 25, 29, 38, 73], [8, 15, 36, 44, 87, 101], [6, 43, 55, 75, 102, 109], [24, 31, 39, 76, 88, 108], [5, 40, 79, 101, 103, 107], [13, 25, 45, 72, 79, 87], [18, 26, 69, 83, 101, 108], [27, 41, 59, 66, 87, 95], [20, 22, 26, 35, 70, 109], [2, 20, 27, 48, 56, 99], [8, 10, 14, 23, 58, 97], [5, 7, 11, 20, 55, 94], [15, 23, 66, 80, 98, 105], [14, 16, 20, 29, 64, 103], [0, 14, 32, 39, 60, 68], [34, 46, 66, 93, 100, 108], [0, 27, 34, 42, 79, 91], [17, 19, 23, 32, 67, 106], [34, 73, 95, 97, 101, 110], [0, 37, 49, 69, 96, 103], [16, 55, 77, 79, 83, 92], [0, 7, 15, 52, 64, 84], [1, 21, 48, 55, 63, 100], [6, 27, 35, 78, 92, 110], [9, 23, 41, 48, 69, 77], [18, 25, 33, 70, 82, 102], [3, 17, 35, 42, 63, 71], [16, 38, 40, 44, 53, 88], [28, 40, 60, 87, 94, 102], [22, 44, 46, 50, 59, 94], [15, 42, 49, 57, 94, 106], [24, 38, 56, 63, 84, 92], [11, 13, 17, 26, 61, 100], [1, 5, 14, 49, 88, 110], [5, 23, 30, 51, 59, 102], [7, 46, 68, 70, 74, 83], [3, 30, 37, 45, 82, 94], [18, 32, 50, 57, 78, 86], [7, 27, 54, 61, 69, 106], [25, 64, 86, 88, 92, 101], [37, 59, 61, 65, 74, 109], [14, 21, 42, 50, 93, 107], [21, 28, 36, 73, 85, 105], [11, 29, 36, 57, 65, 108], [18, 45, 52, 60, 97, 109], [25, 37, 57, 84, 91, 99], [21, 35, 53, 60, 81, 89], [3, 10, 18, 55, 67, 87], [34, 56, 58, 62, 71, 106], [36, 50, 68, 75, 96, 104], [28, 50, 52, 56, 65, 100], [0, 21, 29, 72, 86, 104], [10, 30, 57, 64, 72, 109], [31, 43, 63, 90, 97, 105], [8, 43, 82, 104, 106, 110], [16, 28, 48, 75, 82, 90], [3, 11, 54, 68, 86, 93], [28, 67, 89, 91, 95, 104], [19, 41, 43, 47, 56, 91], [1, 13, 33, 60, 67, 75], [11, 18, 39, 47, 90, 104], [4, 24, 51, 58, 66, 103], [4, 26, 28, 32, 41, 76], [3, 24, 32, 75, 89, 107], [9, 16, 24, 61, 73, 93], [17, 24, 45, 53, 96, 110], [19, 31, 51, 78, 85, 93], [0, 8, 51, 65, 83, 90], [30, 44, 62, 69, 90, 98], [25, 47, 49, 53, 62, 97], [5, 12, 33, 41, 84, 98], [12, 20, 63, 77, 95, 102], [4, 12, 49, 61, 81, 108], [9, 17, 60, 74, 92, 99], [1, 9, 46, 58, 78, 105], [19, 58, 80, 82, 86, 95], [6, 33, 40, 48, 85, 97], [2, 9, 30, 38, 81, 95], [9, 53, 80, 94, 101, 104], [2, 29, 43, 50, 53, 69], [1, 7, 26, 34, 43, 99], [2, 24, 49, 60, 72, 90], [29, 56, 70, 77, 80, 96], [6, 18, 36, 59, 81, 106], [0, 12, 30, 53, 75, 100], [48, 61, 67, 86, 94, 103], [54, 67, 73, 92, 100, 109], [3, 21, 44, 66, 91, 102], [11, 25, 32, 35, 51, 95], [1, 10, 66, 79, 85, 104], [7, 13, 32, 40, 49, 105], [0, 13, 19, 38, 46, 55], [11, 33, 58, 69, 81, 99], [5, 32, 46, 53, 56, 72], [26, 53, 67, 74, 77, 93], [20, 47, 61, 68, 71, 87], [1, 20, 28, 37, 93, 106], [23, 37, 44, 47, 63, 107], [4, 23, 31, 40, 96, 109], [20, 42, 67, 78, 90, 108], [7, 14, 17, 33, 77, 104], [13, 20, 23, 39, 83, 110], [8, 22, 29, 32, 48, 92], [6, 31, 42, 54, 72, 95], [0, 44, 71, 85, 92, 95], [42, 55, 61, 80, 88, 97], [9, 34, 45, 57, 75, 98], [8, 30, 55, 66, 78, 96], [14, 41, 55, 62, 65, 81], [11, 19, 28, 84, 97, 103], [26, 40, 47, 50, 66, 110], [5, 27, 52, 63, 75, 93], [39, 52, 58, 77, 85, 94], [7, 18, 30, 48, 71, 93], [5, 19, 26, 29, 45, 89], [7, 63, 76, 82, 101, 109], [5, 13, 22, 78, 91, 97], [6, 19, 25, 44, 52, 61], [11, 38, 52, 59, 62, 78], [36, 49, 55, 74, 82, 91], [3, 15, 33, 56, 78, 103], [13, 24, 36, 54, 77, 99], [7, 16, 72, 85, 91, 110], [1, 8, 11, 27, 71, 98], [1, 12, 24, 42, 65, 87], [17, 25, 34, 90, 103, 109], [8, 35, 49, 56, 59, 75], [33, 46, 52, 71, 79, 88], [18, 43, 54, 66, 84, 107], [21, 46, 57, 69, 87, 110], [9, 22, 28, 47, 55, 64], [30, 43, 49, 68, 76, 85], [15, 28, 34, 53, 61, 70], [0, 18, 41, 63, 88, 99], [3, 26, 48, 73, 84, 96], [14, 22, 31, 87, 100, 106], [17, 44, 58, 65, 68, 84], [10, 17, 20, 36, 80, 107], [15, 59, 86, 100, 107, 110], [17, 31, 38, 41, 57, 101], [16, 27, 39, 57, 80, 102], [3, 28, 39, 51, 69, 92], [3, 16, 22, 41, 49, 58], [21, 34, 40, 59, 67, 76], [0, 23, 45, 70, 81, 93], [10, 16, 35, 43, 52, 108], [15, 38, 60, 85, 96, 108], [6, 24, 47, 69, 94, 105], [5, 8, 24, 68, 95, 109], [9, 21, 39, 62, 84, 109], [35, 62, 76, 83, 86, 102], [3, 47, 74, 88, 95, 98], [12, 25, 31, 50, 58, 67], [4, 60, 73, 79, 98, 106], [14, 28, 35, 38, 54, 98], [12, 35, 57, 82, 93, 105], [23, 50, 64, 71, 74, 90], [1, 57, 70, 76, 95, 103], [9, 27, 50, 72, 97, 108], [19, 30, 42, 60, 83, 105], [2, 18, 62, 89, 103, 110], [38, 65, 79, 86, 89, 105], [2, 10, 19, 75, 88, 94], [15, 40, 51, 63, 81, 104], [27, 40, 46, 65, 73, 82], [12, 56, 83, 97, 104, 107], [0, 25, 36, 48, 66, 89], [20, 34, 41, 44, 60, 104], [4, 15, 27, 45, 68, 90], [41, 68, 82, 89, 92, 108], [4, 13, 69, 82, 88, 107], [10, 21, 33, 51, 74, 96], [4, 11, 14, 30, 74, 101], [4, 10, 29, 37, 46, 102], [45, 58, 64, 83, 91, 100], [6, 50, 77, 91, 98, 101], [18, 31, 37, 56, 64, 73], [12, 37, 48, 60, 78, 101], [24, 37, 43, 62, 70, 79], [9, 32, 54, 79, 90, 102], [32, 59, 73, 80, 83, 99], [2, 16, 23, 26, 42, 86], [14, 36, 61, 72, 84, 102], [8, 16, 25, 81, 94, 100], [6, 29, 51, 76, 87, 99], [22, 33, 45, 63, 86, 108], [17, 39, 64, 75, 87, 105], [51, 64, 70, 89, 97, 106], [2, 5, 21, 65, 92, 106]]
\item 1 \{0=1110, 1=31968, 2=542790, 3=2167164, 4=2385168\} [[8, 39, 40, 41, 42, 70], [0, 1, 2, 3, 31, 80], [30, 31, 32, 33, 61, 110], [6, 7, 8, 9, 37, 86], [23, 54, 55, 56, 57, 85], [15, 16, 17, 18, 46, 95], [27, 28, 29, 30, 58, 107], [24, 25, 26, 27, 55, 104], [26, 57, 58, 59, 60, 88], [29, 60, 61, 62, 63, 91], [25, 74, 105, 106, 107, 108], [41, 72, 73, 74, 75, 103], [7, 56, 87, 88, 89, 90], [10, 59, 90, 91, 92, 93], [13, 62, 93, 94, 95, 96], [22, 71, 102, 103, 104, 105], [9, 10, 11, 12, 40, 89], [44, 75, 76, 77, 78, 106], [18, 19, 20, 21, 49, 98], [11, 42, 43, 44, 45, 73], [3, 4, 5, 6, 34, 83], [16, 65, 96, 97, 98, 99], [19, 68, 99, 100, 101, 102], [17, 48, 49, 50, 51, 79], [47, 78, 79, 80, 81, 109], [20, 51, 52, 53, 54, 82], [12, 13, 14, 15, 43, 92], [0, 28, 77, 108, 109, 110], [4, 53, 84, 85, 86, 87], [1, 50, 81, 82, 83, 84], [32, 63, 64, 65, 66, 94], [2, 33, 34, 35, 36, 64], [21, 22, 23, 24, 52, 101], [38, 69, 70, 71, 72, 100], [14, 45, 46, 47, 48, 76], [35, 66, 67, 68, 69, 97], [5, 36, 37, 38, 39, 67], [30, 53, 69, 73, 81, 96], [2, 38, 78, 88, 92, 98], [40, 49, 71, 81, 88, 94], [17, 27, 34, 40, 97, 106], [2, 18, 22, 30, 45, 90], [9, 54, 77, 93, 97, 105], [3, 26, 42, 46, 54, 69], [14, 54, 64, 68, 74, 89], [33, 56, 72, 76, 84, 99], [0, 23, 39, 43, 51, 66], [1, 5, 11, 26, 62, 102], [1, 9, 24, 69, 92, 108], [5, 21, 25, 33, 48, 93], [9, 16, 22, 79, 88, 110], [14, 24, 31, 37, 94, 103], [5, 45, 55, 59, 65, 80], [52, 61, 83, 93, 100, 106], [16, 25, 47, 57, 64, 70], [3, 7, 15, 30, 75, 98], [7, 29, 39, 46, 52, 109], [12, 22, 26, 32, 47, 83], [6, 51, 74, 90, 94, 102], [34, 43, 65, 75, 82, 88], [20, 30, 37, 43, 100, 109], [25, 34, 56, 66, 73, 79], [39, 62, 78, 82, 90, 105], [3, 18, 63, 86, 102, 106], [15, 19, 27, 42, 87, 110], [8, 24, 28, 36, 51, 96], [23, 63, 73, 77, 83, 98], [5, 15, 22, 28, 85, 94], [7, 11, 17, 32, 68, 108], [5, 20, 56, 96, 106, 110], [15, 38, 54, 58, 66, 81], [14, 50, 90, 100, 104, 110], [9, 13, 21, 36, 81, 104], [4, 13, 35, 45, 52, 58], [21, 31, 35, 41, 56, 92], [39, 49, 53, 59, 74, 110], [6, 21, 66, 89, 105, 109], [18, 41, 57, 61, 69, 84], [12, 16, 24, 39, 84, 107], [36, 59, 75, 79, 87, 102], [15, 25, 29, 35, 50, 86], [11, 47, 87, 97, 101, 107], [14, 30, 34, 42, 57, 102], [17, 57, 67, 71, 77, 92], [27, 50, 66, 70, 78, 93], [11, 51, 61, 65, 71, 86], [49, 58, 80, 90, 97, 103], [0, 15, 60, 83, 99, 103], [1, 7, 64, 73, 95, 105], [7, 16, 38, 48, 55, 61], [55, 64, 86, 96, 103, 109], [11, 27, 31, 39, 54, 99], [6, 16, 20, 26, 41, 77], [37, 46, 68, 78, 85, 91], [35, 75, 85, 89, 95, 110], [3, 10, 16, 73, 82, 104], [30, 40, 44, 50, 65, 101], [6, 13, 19, 76, 85, 107], [17, 33, 37, 45, 60, 105], [1, 10, 32, 42, 49, 55], [0, 10, 14, 20, 35, 71], [22, 31, 53, 63, 70, 76], [6, 10, 18, 33, 78, 101], [2, 17, 53, 93, 103, 107], [21, 44, 60, 64, 72, 87], [26, 66, 76, 80, 86, 101], [32, 72, 82, 86, 92, 107], [27, 37, 41, 47, 62, 98], [11, 21, 28, 34, 91, 100], [4, 61, 70, 92, 102, 109], [36, 46, 50, 56, 71, 107], [4, 8, 14, 29, 65, 105], [3, 13, 17, 23, 38, 74], [9, 19, 23, 29, 44, 80], [19, 28, 50, 60, 67, 73], [4, 10, 67, 76, 98, 108], [20, 36, 40, 48, 63, 108], [8, 18, 25, 31, 88, 97], [6, 29, 45, 49, 57, 72], [2, 12, 19, 25, 82, 91], [10, 19, 41, 51, 58, 64], [5, 41, 81, 91, 95, 101], [43, 52, 74, 84, 91, 97], [12, 35, 51, 55, 63, 78], [0, 7, 13, 70, 79, 101], [2, 8, 23, 59, 99, 109], [28, 37, 59, 69, 76, 82], [24, 47, 63, 67, 75, 90], [1, 23, 33, 40, 46, 103], [13, 22, 44, 54, 61, 67], [2, 42, 52, 56, 62, 77], [0, 4, 12, 27, 72, 95], [12, 57, 80, 96, 100, 108], [46, 55, 77, 87, 94, 100], [8, 44, 84, 94, 98, 104], [3, 48, 71, 87, 91, 99], [31, 40, 62, 72, 79, 85], [42, 65, 81, 85, 93, 108], [0, 45, 68, 84, 88, 96], [29, 69, 79, 83, 89, 104], [9, 32, 48, 52, 60, 75], [20, 60, 70, 74, 80, 95], [8, 48, 58, 62, 68, 83], [18, 28, 32, 38, 53, 89], [24, 34, 38, 44, 59, 95], [1, 58, 67, 89, 99, 106], [33, 43, 47, 53, 68, 104], [4, 26, 36, 43, 49, 106], [31, 55, 60, 89, 102, 107], [7, 12, 41, 54, 59, 94], [22, 46, 51, 80, 93, 98], [9, 14, 49, 73, 78, 107], [4, 9, 38, 51, 56, 91], [10, 34, 39, 68, 81, 86], [19, 24, 53, 66, 71, 106], [5, 18, 23, 58, 82, 87], [2, 37, 61, 66, 95, 108], [22, 27, 56, 69, 74, 109], [16, 21, 50, 63, 68, 103], [2, 15, 20, 55, 79, 84], [13, 37, 42, 71, 84, 89], [7, 31, 36, 65, 78, 83], [34, 58, 63, 92, 105, 110], [17, 30, 35, 70, 94, 99], [8, 21, 26, 61, 85, 90], [23, 36, 41, 76, 100, 105], [11, 24, 29, 64, 88, 93], [26, 39, 44, 79, 103, 108], [1, 6, 35, 48, 53, 88], [6, 11, 46, 70, 75, 104], [10, 15, 44, 57, 62, 97], [19, 43, 48, 77, 90, 95], [13, 18, 47, 60, 65, 100], [3, 8, 43, 67, 72, 101], [16, 40, 45, 74, 87, 92], [25, 49, 54, 83, 96, 101], [0, 5, 40, 64, 69, 98], [4, 28, 33, 62, 75, 80], [14, 27, 32, 67, 91, 96], [20, 33, 38, 73, 97, 102], [12, 17, 52, 76, 81, 110], [0, 29, 42, 47, 82, 106], [3, 32, 45, 50, 85, 109], [28, 52, 57, 86, 99, 104], [1, 25, 30, 59, 72, 77], [5, 16, 29, 31, 43, 71], [8, 15, 47, 56, 102, 108], [6, 27, 43, 61, 81, 103], [2, 9, 41, 50, 96, 102], [17, 28, 41, 43, 55, 83], [12, 28, 46, 66, 88, 102], [7, 20, 22, 34, 62, 107], [9, 15, 26, 33, 65, 74], [15, 21, 32, 39, 71, 80], [0, 22, 36, 57, 73, 91], [2, 4, 16, 44, 89, 100], [1, 21, 43, 57, 78, 94], [7, 27, 49, 63, 84, 100], [11, 13, 25, 53, 98, 109], [8, 53, 64, 77, 79, 91], [10, 23, 25, 37, 65, 110], [13, 27, 48, 64, 82, 102], [32, 43, 56, 58, 70, 98], [14, 23, 69, 75, 86, 93], [14, 59, 70, 83, 85, 97], [10, 28, 48, 70, 84, 105], [36, 42, 53, 60, 92, 101], [1, 29, 74, 85, 98, 100], [10, 30, 52, 66, 87, 103], [7, 35, 80, 91, 104, 106], [3, 9, 20, 27, 59, 68], [1, 14, 16, 28, 56, 101], [3, 25, 39, 60, 76, 94], [20, 29, 75, 81, 92, 99], [41, 52, 65, 67, 79, 107], [18, 24, 35, 42, 74, 83], [45, 51, 62, 69, 101, 110], [11, 20, 66, 72, 83, 90], [2, 11, 57, 63, 74, 81], [26, 35, 81, 87, 98, 105], [0, 32, 41, 87, 93, 104], [39, 45, 56, 63, 95, 104], [1, 13, 41, 86, 97, 110], [6, 22, 40, 60, 82, 96], [26, 71, 82, 95, 97, 109], [5, 51, 57, 68, 75, 107], [2, 13, 26, 28, 40, 68], [5, 12, 44, 53, 99, 105], [23, 68, 79, 92, 94, 106], [18, 34, 52, 72, 94, 108], [8, 17, 63, 69, 80, 87], [3, 19, 37, 57, 79, 93], [11, 22, 35, 37, 49, 77], [4, 32, 77, 88, 101, 103], [1, 19, 39, 61, 75, 96], [14, 25, 38, 40, 52, 80], [4, 22, 42, 64, 78, 99], [6, 12, 23, 30, 62, 71], [8, 54, 60, 71, 78, 110], [6, 28, 42, 63, 79, 97], [0, 16, 34, 54, 76, 90], [4, 17, 19, 31, 59, 104], [19, 33, 54, 70, 88, 108], [10, 38, 83, 94, 107, 109], [7, 25, 45, 67, 81, 102], [2, 47, 58, 71, 73, 85], [44, 55, 68, 70, 82, 110], [0, 6, 17, 24, 56, 65], [13, 33, 55, 69, 90, 106], [12, 34, 48, 69, 85, 103], [0, 21, 37, 55, 75, 97], [7, 21, 42, 58, 76, 96], [8, 19, 32, 34, 46, 74], [21, 27, 38, 45, 77, 86], [5, 7, 19, 47, 92, 103], [15, 37, 51, 72, 88, 106], [29, 40, 53, 55, 67, 95], [15, 31, 49, 69, 91, 105], [16, 36, 58, 72, 93, 109], [17, 62, 73, 86, 88, 100], [23, 32, 78, 84, 95, 102], [29, 38, 84, 90, 101, 108], [5, 14, 60, 66, 77, 84], [20, 65, 76, 89, 91, 103], [13, 31, 51, 73, 87, 108], [4, 24, 46, 60, 81, 97], [38, 49, 62, 64, 76, 104], [30, 36, 47, 54, 86, 95], [3, 14, 21, 53, 62, 108], [12, 18, 29, 36, 68, 77], [3, 35, 44, 90, 96, 107], [6, 38, 47, 93, 99, 110], [12, 33, 49, 67, 87, 109], [3, 24, 40, 58, 78, 100], [11, 56, 67, 80, 82, 94], [4, 18, 39, 55, 73, 93], [10, 24, 45, 61, 79, 99], [26, 37, 50, 52, 64, 92], [1, 15, 36, 52, 70, 90], [9, 30, 46, 64, 84, 106], [5, 50, 61, 74, 76, 88], [9, 31, 45, 66, 82, 100], [42, 48, 59, 66, 98, 107], [35, 46, 59, 61, 73, 101], [2, 48, 54, 65, 72, 104], [0, 11, 18, 50, 59, 105], [24, 30, 41, 48, 80, 89], [9, 25, 43, 63, 85, 99], [23, 34, 47, 49, 61, 89], [20, 31, 44, 46, 58, 86], [17, 26, 72, 78, 89, 96], [8, 10, 22, 50, 95, 106], [18, 40, 54, 75, 91, 109], [27, 33, 44, 51, 83, 92], [16, 30, 51, 67, 85, 105], [33, 39, 50, 57, 89, 98], [0, 9, 53, 58, 61, 94], [1, 18, 27, 71, 76, 79], [3, 51, 70, 77, 81, 89], [24, 72, 91, 98, 102, 110], [5, 32, 35, 54, 79, 100], [21, 69, 88, 95, 99, 107], [1, 8, 12, 20, 45, 93], [1, 4, 37, 54, 63, 107], [13, 20, 24, 32, 57, 105], [4, 11, 15, 23, 48, 96], [2, 5, 24, 49, 70, 86], [27, 46, 53, 57, 65, 90], [16, 23, 27, 35, 60, 108], [0, 8, 33, 81, 100, 107], [17, 20, 39, 64, 85, 101], [10, 26, 53, 56, 75, 100], [7, 14, 18, 26, 51, 99], [4, 7, 40, 57, 66, 110], [9, 28, 35, 39, 47, 72], [21, 40, 47, 51, 59, 84], [3, 11, 36, 84, 103, 110], [3, 22, 29, 33, 41, 66], [36, 55, 62, 66, 74, 99], [18, 37, 44, 48, 56, 81], [8, 11, 30, 55, 76, 92], [23, 28, 31, 64, 81, 90], [4, 25, 41, 68, 71, 90], [13, 29, 56, 59, 78, 103], [19, 36, 45, 89, 94, 97], [45, 64, 71, 75, 83, 108], [9, 57, 76, 83, 87, 95], [13, 30, 39, 83, 88, 91], [1, 17, 44, 47, 66, 91], [6, 54, 73, 80, 84, 92], [5, 8, 27, 52, 73, 89], [26, 31, 34, 67, 84, 93], [4, 21, 30, 74, 79, 82], [0, 19, 26, 30, 38, 63], [18, 66, 85, 92, 96, 104], [7, 24, 33, 77, 82, 85], [2, 27, 75, 94, 101, 105], [9, 34, 55, 71, 98, 101], [6, 31, 52, 68, 95, 98], [0, 44, 49, 52, 85, 102], [33, 52, 59, 63, 71, 96], [24, 43, 50, 54, 62, 87], [6, 50, 55, 58, 91, 108], [14, 19, 22, 55, 72, 81], [19, 40, 56, 83, 86, 105], [23, 26, 45, 70, 91, 107], [1, 34, 51, 60, 104, 109], [41, 46, 49, 82, 99, 108], [3, 28, 49, 65, 92, 95], [28, 45, 54, 98, 103, 106], [12, 60, 79, 86, 90, 98], [38, 43, 46, 79, 96, 105], [22, 39, 48, 92, 97, 100], [10, 27, 36, 80, 85, 88], [1, 22, 38, 65, 68, 87], [6, 15, 59, 64, 67, 100], [5, 9, 17, 42, 90, 109], [7, 23, 50, 53, 72, 97], [11, 38, 41, 60, 85, 106], [16, 33, 42, 86, 91, 94], [39, 58, 65, 69, 77, 102], [5, 30, 78, 97, 104, 108], [8, 13, 16, 49, 66, 75], [2, 29, 32, 51, 76, 97], [25, 42, 51, 95, 100, 103], [4, 20, 47, 50, 69, 94], [5, 10, 13, 46, 63, 72], [20, 25, 28, 61, 78, 87], [6, 25, 32, 36, 44, 69], [19, 35, 62, 65, 84, 109], [2, 21, 46, 67, 83, 110], [0, 48, 67, 74, 78, 86], [20, 23, 42, 67, 88, 104], [18, 43, 64, 80, 107, 110], [3, 12, 56, 61, 64, 97], [22, 43, 59, 86, 89, 108], [12, 37, 58, 74, 101, 104], [16, 32, 59, 62, 81, 106], [7, 28, 44, 71, 74, 93], [13, 34, 50, 77, 80, 99], [3, 47, 52, 55, 88, 105], [32, 37, 40, 73, 90, 99], [14, 41, 44, 63, 88, 109], [15, 34, 41, 45, 53, 78], [15, 40, 61, 77, 104, 107], [17, 22, 25, 58, 75, 84], [42, 61, 68, 72, 80, 105], [10, 31, 47, 74, 77, 96], [2, 7, 10, 43, 60, 69], [26, 29, 48, 73, 94, 110], [29, 34, 37, 70, 87, 96], [15, 24, 68, 73, 76, 109], [10, 17, 21, 29, 54, 102], [12, 21, 65, 70, 73, 106], [0, 25, 46, 62, 89, 92], [11, 16, 19, 52, 69, 78], [16, 37, 53, 80, 83, 102], [30, 49, 56, 60, 68, 93], [15, 63, 82, 89, 93, 101], [9, 18, 62, 67, 70, 103], [11, 14, 33, 58, 79, 95], [35, 40, 43, 76, 93, 102], [2, 6, 14, 39, 87, 106], [14, 17, 36, 61, 82, 98], [8, 35, 38, 57, 82, 103], [31, 48, 57, 101, 106, 109], [12, 31, 38, 42, 50, 75]]
\item 1 \{1=32856, 2=483516, 3=2213784, 4=2398044\} [[8, 39, 40, 41, 42, 70], [0, 1, 2, 3, 31, 80], [30, 31, 32, 33, 61, 110], [6, 7, 8, 9, 37, 86], [23, 54, 55, 56, 57, 85], [15, 16, 17, 18, 46, 95], [27, 28, 29, 30, 58, 107], [24, 25, 26, 27, 55, 104], [26, 57, 58, 59, 60, 88], [29, 60, 61, 62, 63, 91], [25, 74, 105, 106, 107, 108], [41, 72, 73, 74, 75, 103], [7, 56, 87, 88, 89, 90], [10, 59, 90, 91, 92, 93], [13, 62, 93, 94, 95, 96], [22, 71, 102, 103, 104, 105], [9, 10, 11, 12, 40, 89], [44, 75, 76, 77, 78, 106], [18, 19, 20, 21, 49, 98], [11, 42, 43, 44, 45, 73], [3, 4, 5, 6, 34, 83], [16, 65, 96, 97, 98, 99], [19, 68, 99, 100, 101, 102], [17, 48, 49, 50, 51, 79], [47, 78, 79, 80, 81, 109], [20, 51, 52, 53, 54, 82], [12, 13, 14, 15, 43, 92], [0, 28, 77, 108, 109, 110], [4, 53, 84, 85, 86, 87], [1, 50, 81, 82, 83, 84], [32, 63, 64, 65, 66, 94], [2, 33, 34, 35, 36, 64], [21, 22, 23, 24, 52, 101], [38, 69, 70, 71, 72, 100], [14, 45, 46, 47, 48, 76], [35, 66, 67, 68, 69, 97], [5, 36, 37, 38, 39, 67], [12, 16, 21, 54, 77, 88], [7, 38, 81, 92, 94, 97], [12, 23, 25, 28, 49, 80], [18, 29, 31, 34, 55, 86], [24, 39, 46, 56, 74, 86], [1, 36, 40, 45, 78, 101], [24, 47, 58, 93, 97, 102], [11, 22, 57, 61, 66, 99], [18, 33, 40, 50, 68, 80], [24, 28, 33, 66, 89, 100], [4, 39, 43, 48, 81, 104], [0, 7, 17, 35, 47, 96], [2, 14, 63, 78, 85, 95], [14, 57, 68, 70, 73, 94], [33, 44, 46, 49, 70, 101], [5, 7, 10, 31, 62, 105], [39, 54, 61, 71, 89, 101], [20, 63, 74, 76, 79, 100], [14, 25, 60, 64, 69, 102], [48, 63, 70, 80, 98, 110], [3, 14, 16, 19, 40, 71], [18, 41, 52, 87, 91, 96], [2, 20, 32, 81, 96, 103], [2, 13, 48, 52, 57, 90], [5, 23, 35, 84, 99, 106], [6, 21, 28, 38, 56, 68], [6, 10, 15, 48, 71, 82], [8, 51, 62, 64, 67, 88], [30, 45, 52, 62, 80, 92], [5, 16, 51, 55, 60, 93], [21, 36, 43, 53, 71, 83], [2, 45, 56, 58, 61, 82], [8, 10, 13, 34, 65, 108], [12, 27, 34, 44, 62, 74], [1, 6, 39, 62, 73, 108], [8, 26, 38, 87, 102, 109], [9, 24, 31, 41, 59, 71], [10, 45, 49, 54, 87, 110], [3, 7, 12, 45, 68, 79], [16, 47, 90, 101, 103, 106], [27, 31, 36, 69, 92, 103], [7, 42, 46, 51, 84, 107], [9, 13, 18, 51, 74, 85], [1, 32, 75, 86, 88, 91], [2, 51, 66, 73, 83, 101], [9, 32, 43, 78, 82, 87], [0, 33, 56, 67, 102, 106], [6, 17, 19, 22, 43, 74], [9, 20, 22, 25, 46, 77], [8, 19, 54, 58, 63, 96], [17, 60, 71, 73, 76, 97], [27, 38, 40, 43, 64, 95], [18, 22, 27, 60, 83, 94], [0, 15, 22, 32, 50, 62], [36, 47, 49, 52, 73, 104], [6, 29, 40, 75, 79, 84], [12, 35, 46, 81, 85, 90], [24, 35, 37, 40, 61, 92], [3, 26, 37, 72, 76, 81], [5, 48, 59, 61, 64, 85], [3, 18, 25, 35, 53, 65], [1, 22, 53, 96, 107, 109], [42, 53, 55, 58, 79, 110], [15, 30, 37, 47, 65, 77], [6, 13, 23, 41, 53, 102], [33, 48, 55, 65, 83, 95], [21, 44, 55, 90, 94, 99], [4, 35, 78, 89, 91, 94], [30, 53, 64, 99, 103, 108], [39, 50, 52, 55, 76, 107], [27, 50, 61, 96, 100, 105], [19, 50, 93, 104, 106, 109], [30, 34, 39, 72, 95, 106], [8, 20, 69, 84, 91, 101], [11, 23, 72, 87, 94, 104], [23, 66, 77, 79, 82, 103], [5, 54, 69, 76, 86, 104], [8, 57, 72, 79, 89, 107], [42, 57, 64, 74, 92, 104], [10, 41, 84, 95, 97, 100], [14, 26, 75, 90, 97, 107], [26, 69, 80, 82, 85, 106], [11, 60, 75, 82, 92, 110], [30, 41, 43, 46, 67, 98], [15, 26, 28, 31, 52, 83], [12, 19, 29, 47, 59, 108], [15, 19, 24, 57, 80, 91], [11, 54, 65, 67, 70, 91], [45, 60, 67, 77, 95, 107], [13, 44, 87, 98, 100, 103], [0, 4, 9, 42, 65, 76], [4, 14, 32, 44, 93, 108], [21, 32, 34, 37, 58, 89], [15, 38, 49, 84, 88, 93], [9, 16, 26, 44, 56, 105], [5, 17, 66, 81, 88, 98], [21, 25, 30, 63, 86, 97], [33, 37, 42, 75, 98, 109], [3, 36, 59, 70, 105, 109], [29, 72, 83, 85, 88, 109], [2, 4, 7, 28, 59, 102], [3, 10, 20, 38, 50, 99], [17, 29, 78, 93, 100, 110], [0, 11, 13, 16, 37, 68], [20, 31, 66, 70, 75, 108], [1, 4, 25, 56, 99, 110], [17, 28, 63, 67, 72, 105], [0, 23, 34, 69, 73, 78], [36, 51, 58, 68, 86, 98], [1, 11, 29, 41, 90, 105], [27, 42, 49, 59, 77, 89], [34, 42, 50, 63, 90, 102], [13, 66, 71, 80, 86, 107], [31, 43, 65, 79, 88, 102], [1, 10, 24, 64, 76, 98], [2, 19, 72, 77, 86, 92], [13, 21, 29, 42, 69, 81], [1, 15, 55, 67, 89, 103], [6, 46, 58, 80, 94, 103], [10, 22, 44, 58, 67, 81], [4, 26, 40, 49, 63, 103], [21, 33, 76, 84, 92, 105], [0, 43, 51, 59, 72, 99], [19, 31, 53, 67, 76, 90], [0, 40, 52, 74, 88, 97], [4, 12, 20, 33, 60, 72], [1, 9, 17, 30, 57, 69], [13, 25, 47, 61, 70, 84], [9, 52, 60, 68, 81, 108], [9, 21, 64, 72, 80, 93], [8, 22, 31, 45, 85, 97], [19, 27, 35, 48, 75, 87], [48, 53, 62, 68, 89, 106], [24, 36, 79, 87, 95, 108], [3, 15, 58, 66, 74, 87], [4, 57, 62, 71, 77, 98], [0, 5, 14, 20, 41, 58], [10, 18, 26, 39, 66, 78], [22, 30, 38, 51, 78, 90], [7, 19, 41, 55, 64, 78], [25, 33, 41, 54, 81, 93], [3, 46, 54, 62, 75, 102], [12, 17, 26, 32, 53, 70], [17, 31, 40, 54, 94, 106], [5, 11, 32, 49, 102, 107], [3, 11, 24, 51, 63, 106], [8, 14, 35, 52, 105, 110], [2, 8, 29, 46, 99, 104], [2, 11, 17, 38, 55, 108], [6, 14, 27, 54, 66, 109], [30, 35, 44, 50, 71, 88], [15, 20, 29, 35, 56, 73], [37, 45, 53, 66, 93, 105], [34, 46, 68, 82, 91, 105], [15, 27, 70, 78, 86, 99], [3, 43, 55, 77, 91, 100], [33, 38, 47, 53, 74, 91], [20, 37, 90, 95, 104, 110], [11, 25, 34, 48, 88, 100], [21, 26, 35, 41, 62, 79], [18, 23, 32, 38, 59, 76], [5, 19, 28, 42, 82, 94], [5, 22, 75, 80, 89, 95], [0, 27, 39, 82, 90, 98], [7, 60, 65, 74, 80, 101], [1, 23, 37, 46, 60, 100], [13, 22, 36, 76, 88, 110], [3, 8, 17, 23, 44, 61], [12, 39, 51, 94, 102, 110], [4, 16, 38, 52, 61, 75], [37, 49, 71, 85, 94, 108], [18, 30, 73, 81, 89, 102], [6, 11, 20, 26, 47, 64], [9, 49, 61, 83, 97, 106], [24, 29, 38, 44, 65, 82], [10, 32, 46, 55, 69, 109], [4, 18, 58, 70, 92, 106], [0, 8, 21, 48, 60, 103], [39, 44, 53, 59, 80, 97], [7, 21, 61, 73, 95, 109], [2, 16, 25, 39, 79, 91], [7, 16, 30, 70, 82, 104], [28, 36, 44, 57, 84, 96], [11, 28, 81, 86, 95, 101], [25, 37, 59, 73, 82, 96], [9, 36, 48, 91, 99, 107], [7, 29, 43, 52, 66, 106], [1, 13, 35, 49, 58, 72], [51, 56, 65, 71, 92, 109], [6, 49, 57, 65, 78, 105], [22, 34, 56, 70, 79, 93], [7, 15, 23, 36, 63, 75], [2, 23, 40, 93, 98, 107], [12, 24, 67, 75, 83, 96], [6, 33, 45, 88, 96, 104], [5, 26, 43, 96, 101, 110], [31, 39, 47, 60, 87, 99], [42, 47, 56, 62, 83, 100], [16, 69, 74, 83, 89, 110], [9, 14, 23, 29, 50, 67], [36, 41, 50, 56, 77, 94], [8, 25, 78, 83, 92, 98], [20, 34, 43, 57, 97, 109], [5, 18, 45, 57, 100, 108], [1, 54, 59, 68, 74, 95], [16, 24, 32, 45, 72, 84], [14, 31, 84, 89, 98, 104], [16, 28, 50, 64, 73, 87], [10, 63, 68, 77, 83, 104], [6, 18, 61, 69, 77, 90], [4, 13, 27, 67, 79, 101], [14, 28, 37, 51, 91, 103], [27, 32, 41, 47, 68, 85], [12, 52, 64, 86, 100, 109], [17, 34, 87, 92, 101, 107], [10, 19, 33, 73, 85, 107], [28, 40, 62, 76, 85, 99], [2, 15, 42, 54, 97, 105], [0, 12, 55, 63, 71, 84], [3, 30, 42, 85, 93, 101], [45, 50, 59, 65, 86, 103], [40, 48, 56, 69, 96, 108], [2, 5, 21, 50, 70, 74], [12, 22, 73, 91, 98, 106], [10, 14, 53, 56, 72, 101], [15, 61, 72, 98, 102, 108], [32, 35, 51, 80, 100, 104], [23, 43, 47, 86, 89, 105], [43, 54, 80, 84, 90, 108], [16, 34, 41, 49, 66, 76], [7, 18, 44, 48, 54, 72], [29, 32, 48, 77, 97, 101], [7, 24, 34, 85, 103, 110], [2, 10, 27, 37, 88, 106], [6, 35, 55, 59, 98, 101], [17, 21, 27, 45, 91, 102], [0, 26, 30, 36, 54, 100], [20, 23, 39, 68, 88, 92], [3, 21, 67, 78, 104, 108], [4, 21, 31, 82, 100, 107], [28, 46, 53, 61, 78, 88], [10, 21, 47, 51, 57, 75], [19, 37, 44, 52, 69, 79], [3, 32, 52, 56, 95, 98], [4, 55, 73, 80, 88, 105], [35, 38, 54, 83, 103, 107], [8, 12, 18, 36, 82, 93], [7, 11, 50, 53, 69, 98], [20, 40, 44, 83, 86, 102], [9, 35, 39, 45, 63, 109], [4, 8, 47, 50, 66, 95], [2, 41, 44, 60, 89, 109], [5, 8, 24, 53, 73, 77], [16, 23, 31, 48, 58, 109], [26, 46, 50, 89, 92, 108], [34, 52, 59, 67, 84, 94], [46, 64, 71, 79, 96, 106], [7, 58, 76, 83, 91, 108], [14, 18, 24, 42, 88, 99], [2, 18, 47, 67, 71, 110], [43, 61, 68, 76, 93, 103], [11, 15, 21, 39, 85, 96], [25, 43, 50, 58, 75, 85], [22, 33, 59, 63, 69, 87], [8, 28, 32, 71, 74, 90], [28, 39, 65, 69, 75, 93], [7, 25, 32, 40, 57, 67], [0, 46, 57, 83, 87, 93], [10, 17, 25, 42, 52, 103], [14, 17, 33, 62, 82, 86], [23, 27, 33, 51, 97, 108], [9, 19, 70, 88, 95, 103], [6, 52, 63, 89, 93, 99], [1, 19, 26, 34, 51, 61], [1, 18, 28, 79, 97, 104], [0, 18, 64, 75, 101, 105], [2, 6, 12, 30, 76, 87], [13, 17, 56, 59, 75, 104], [13, 24, 50, 54, 60, 78], [31, 49, 56, 64, 81, 91], [49, 67, 74, 82, 99, 109], [9, 55, 66, 92, 96, 102], [13, 31, 38, 46, 63, 73], [6, 32, 36, 42, 60, 106], [4, 22, 29, 37, 54, 64], [22, 40, 47, 55, 72, 82], [2, 22, 26, 65, 68, 84], [0, 10, 61, 79, 86, 94], [9, 38, 58, 62, 101, 104], [5, 13, 30, 40, 91, 109], [6, 16, 67, 85, 92, 100], [3, 9, 27, 73, 84, 110], [5, 25, 29, 68, 71, 87], [13, 20, 28, 45, 55, 106], [25, 36, 62, 66, 72, 90], [11, 31, 35, 74, 77, 93], [19, 23, 62, 65, 81, 110], [3, 13, 64, 82, 89, 97], [40, 51, 77, 81, 87, 105], [16, 27, 53, 57, 63, 81], [34, 45, 71, 75, 81, 99], [23, 26, 42, 71, 91, 95], [38, 41, 57, 86, 106, 110], [12, 41, 61, 65, 104, 107], [8, 11, 27, 56, 76, 80], [4, 11, 19, 36, 46, 97], [15, 44, 64, 68, 107, 110], [3, 49, 60, 86, 90, 96], [26, 29, 45, 74, 94, 98], [19, 30, 56, 60, 66, 84], [20, 24, 30, 48, 94, 105], [1, 5, 44, 47, 63, 92], [16, 20, 59, 62, 78, 107], [37, 48, 74, 78, 84, 102], [5, 9, 15, 33, 79, 90], [11, 14, 30, 59, 79, 83], [31, 42, 68, 72, 78, 96], [40, 58, 65, 73, 90, 100], [7, 14, 22, 39, 49, 100], [4, 15, 41, 45, 51, 69], [1, 8, 16, 33, 43, 94], [0, 6, 24, 70, 81, 107], [1, 12, 38, 42, 48, 66], [37, 55, 62, 70, 87, 97], [12, 58, 69, 95, 99, 105], [15, 25, 76, 94, 101, 109], [17, 37, 41, 80, 83, 99], [1, 52, 70, 77, 85, 102], [14, 34, 38, 77, 80, 96], [3, 29, 33, 39, 57, 103], [10, 28, 35, 43, 60, 70], [17, 20, 36, 65, 85, 89], [0, 29, 49, 53, 92, 95], [7, 13, 26, 33, 77, 99], [23, 45, 64, 70, 83, 90], [17, 39, 58, 64, 77, 84], [9, 28, 34, 47, 54, 98], [11, 33, 52, 58, 71, 78], [1, 7, 20, 27, 71, 93], [26, 48, 67, 73, 86, 93], [5, 12, 56, 78, 97, 103], [38, 60, 79, 85, 98, 105], [16, 22, 35, 42, 86, 108], [32, 54, 73, 79, 92, 99], [21, 40, 46, 59, 66, 110], [2, 9, 53, 75, 94, 100], [4, 10, 23, 30, 74, 96], [14, 36, 55, 61, 74, 81], [20, 42, 61, 67, 80, 87], [12, 31, 37, 50, 57, 101], [8, 15, 59, 81, 100, 106], [0, 44, 66, 85, 91, 104], [11, 18, 62, 84, 103, 109], [10, 16, 29, 36, 80, 102], [15, 34, 40, 53, 60, 104], [5, 27, 46, 52, 65, 72], [3, 47, 69, 88, 94, 107], [0, 19, 25, 38, 45, 89], [18, 37, 43, 56, 63, 107], [2, 24, 43, 49, 62, 69], [13, 19, 32, 39, 83, 105], [41, 63, 82, 88, 101, 108], [3, 22, 28, 41, 48, 92], [1, 14, 21, 65, 87, 106], [8, 30, 49, 55, 68, 75], [6, 25, 31, 44, 51, 95], [29, 51, 70, 76, 89, 96], [35, 57, 76, 82, 95, 102], [6, 50, 72, 91, 97, 110], [4, 17, 24, 68, 90, 109]]
\item 1 \{1=35520, 2=550116, 3=2187144, 4=2355420\} [[8, 39, 40, 41, 42, 70], [0, 1, 2, 3, 31, 80], [30, 31, 32, 33, 61, 110], [6, 7, 8, 9, 37, 86], [23, 54, 55, 56, 57, 85], [15, 16, 17, 18, 46, 95], [27, 28, 29, 30, 58, 107], [24, 25, 26, 27, 55, 104], [26, 57, 58, 59, 60, 88], [29, 60, 61, 62, 63, 91], [25, 74, 105, 106, 107, 108], [41, 72, 73, 74, 75, 103], [7, 56, 87, 88, 89, 90], [10, 59, 90, 91, 92, 93], [13, 62, 93, 94, 95, 96], [22, 71, 102, 103, 104, 105], [9, 10, 11, 12, 40, 89], [44, 75, 76, 77, 78, 106], [18, 19, 20, 21, 49, 98], [11, 42, 43, 44, 45, 73], [3, 4, 5, 6, 34, 83], [16, 65, 96, 97, 98, 99], [19, 68, 99, 100, 101, 102], [17, 48, 49, 50, 51, 79], [47, 78, 79, 80, 81, 109], [20, 51, 52, 53, 54, 82], [12, 13, 14, 15, 43, 92], [0, 28, 77, 108, 109, 110], [4, 53, 84, 85, 86, 87], [1, 50, 81, 82, 83, 84], [32, 63, 64, 65, 66, 94], [2, 33, 34, 35, 36, 64], [21, 22, 23, 24, 52, 101], [38, 69, 70, 71, 72, 100], [14, 45, 46, 47, 48, 76], [35, 66, 67, 68, 69, 97], [5, 36, 37, 38, 39, 67], [7, 43, 49, 64, 80, 110], [22, 28, 43, 59, 89, 97], [11, 35, 50, 54, 62, 87], [1, 12, 51, 60, 105, 109], [14, 22, 58, 64, 79, 95], [8, 38, 46, 82, 88, 103], [24, 59, 83, 98, 102, 110], [0, 45, 49, 52, 63, 102], [12, 47, 71, 86, 90, 98], [2, 27, 62, 86, 101, 105], [9, 48, 57, 102, 106, 109], [33, 37, 40, 51, 90, 99], [0, 9, 54, 58, 61, 72], [1, 4, 15, 54, 63, 108], [21, 25, 28, 39, 78, 87], [24, 33, 78, 82, 85, 96], [21, 30, 75, 79, 82, 93], [0, 4, 7, 18, 57, 66], [21, 56, 80, 95, 99, 107], [3, 11, 36, 71, 95, 110], [14, 29, 33, 41, 66, 101], [27, 31, 34, 45, 84, 93], [0, 8, 33, 68, 92, 107], [4, 10, 25, 41, 71, 79], [29, 53, 68, 72, 80, 105], [36, 40, 43, 54, 93, 102], [5, 20, 24, 32, 57, 92], [3, 48, 52, 55, 66, 105], [20, 28, 64, 70, 85, 101], [15, 19, 22, 33, 72, 81], [8, 12, 20, 45, 80, 104], [14, 44, 52, 88, 94, 109], [8, 23, 27, 35, 60, 95], [1, 16, 32, 62, 70, 106], [23, 47, 62, 66, 74, 99], [13, 19, 34, 50, 80, 88], [31, 37, 52, 68, 98, 106], [2, 10, 46, 52, 67, 83], [7, 13, 28, 44, 74, 82], [23, 31, 67, 73, 88, 104], [39, 43, 46, 57, 96, 105], [13, 29, 59, 67, 103, 109], [4, 19, 35, 65, 73, 109], [2, 32, 40, 76, 82, 97], [9, 44, 68, 83, 87, 95], [34, 40, 55, 71, 101, 109], [26, 50, 65, 69, 77, 102], [14, 38, 53, 57, 65, 90], [16, 22, 37, 53, 83, 91], [32, 56, 71, 75, 83, 108], [15, 50, 74, 89, 93, 101], [0, 39, 48, 93, 97, 100], [15, 24, 69, 73, 76, 87], [11, 15, 23, 48, 83, 107], [9, 18, 63, 67, 70, 81], [3, 38, 62, 77, 81, 89], [26, 34, 70, 76, 91, 107], [14, 18, 26, 51, 86, 110], [23, 38, 42, 50, 75, 110], [33, 42, 87, 91, 94, 105], [30, 34, 37, 48, 87, 96], [28, 34, 49, 65, 95, 103], [5, 29, 44, 48, 56, 81], [11, 26, 30, 38, 63, 98], [6, 10, 13, 24, 63, 72], [6, 15, 60, 64, 67, 78], [2, 6, 14, 39, 74, 98], [5, 35, 43, 79, 85, 100], [10, 16, 31, 47, 77, 85], [1, 7, 22, 38, 68, 76], [42, 46, 49, 60, 99, 108], [17, 25, 61, 67, 82, 98], [30, 39, 84, 88, 91, 102], [6, 41, 65, 80, 84, 92], [0, 35, 59, 74, 78, 86], [27, 36, 81, 85, 88, 99], [18, 53, 77, 92, 96, 104], [11, 41, 49, 85, 91, 106], [29, 37, 73, 79, 94, 110], [5, 30, 65, 89, 104, 108], [7, 23, 53, 61, 97, 103], [4, 20, 50, 58, 94, 100], [36, 45, 90, 94, 97, 108], [3, 42, 51, 96, 100, 103], [20, 44, 59, 63, 71, 96], [11, 19, 55, 61, 76, 92], [3, 7, 10, 21, 60, 69], [4, 40, 46, 61, 77, 107], [3, 12, 57, 61, 64, 75], [18, 27, 72, 76, 79, 90], [24, 28, 31, 42, 81, 90], [20, 35, 39, 47, 72, 107], [6, 45, 54, 99, 103, 106], [1, 17, 47, 55, 91, 97], [8, 16, 52, 58, 73, 89], [17, 32, 36, 44, 69, 104], [10, 26, 56, 64, 100, 106], [25, 31, 46, 62, 92, 100], [2, 26, 41, 45, 53, 78], [6, 51, 55, 58, 69, 108], [1, 37, 43, 58, 74, 104], [5, 13, 49, 55, 70, 86], [2, 17, 21, 29, 54, 89], [9, 13, 16, 27, 66, 75], [17, 41, 56, 60, 68, 93], [12, 21, 66, 70, 73, 84], [8, 32, 47, 51, 59, 84], [12, 16, 19, 30, 69, 78], [19, 25, 40, 56, 86, 94], [18, 22, 25, 36, 75, 84], [5, 9, 17, 42, 77, 101], [21, 27, 32, 43, 77, 103], [25, 54, 60, 65, 76, 110], [17, 43, 72, 78, 83, 94], [7, 41, 67, 96, 102, 107], [10, 44, 70, 99, 105, 110], [3, 9, 14, 25, 59, 85], [1, 35, 61, 90, 96, 101], [23, 49, 78, 84, 89, 100], [5, 31, 60, 66, 71, 82], [7, 36, 42, 47, 58, 92], [11, 37, 66, 72, 77, 88], [18, 24, 29, 40, 74, 100], [12, 18, 23, 34, 68, 94], [24, 30, 35, 46, 80, 106], [26, 52, 81, 87, 92, 103], [9, 15, 20, 31, 65, 91], [6, 12, 17, 28, 62, 88], [10, 39, 45, 50, 61, 95], [13, 42, 48, 53, 64, 98], [15, 21, 26, 37, 71, 97], [2, 28, 57, 63, 68, 79], [29, 55, 84, 90, 95, 106], [19, 48, 54, 59, 70, 104], [0, 5, 16, 50, 76, 105], [1, 30, 36, 41, 52, 86], [3, 8, 19, 53, 79, 108], [32, 58, 87, 93, 98, 109], [0, 6, 11, 22, 56, 82], [14, 40, 69, 75, 80, 91], [4, 33, 39, 44, 55, 89], [2, 13, 47, 73, 102, 108], [16, 45, 51, 56, 67, 101], [22, 51, 57, 62, 73, 107], [27, 33, 38, 49, 83, 109], [20, 46, 75, 81, 86, 97], [4, 38, 64, 93, 99, 104], [8, 34, 63, 69, 74, 85], [7, 14, 24, 50, 70, 108], [3, 29, 49, 87, 97, 104], [0, 26, 46, 84, 94, 101], [6, 16, 23, 33, 59, 79], [5, 25, 63, 73, 80, 90], [4, 42, 52, 59, 69, 95], [12, 22, 29, 39, 65, 85], [16, 54, 64, 71, 81, 107], [11, 31, 69, 79, 86, 96], [5, 15, 41, 61, 99, 109], [13, 51, 61, 68, 78, 104], [20, 40, 78, 88, 95, 105], [9, 35, 55, 93, 103, 110], [1, 39, 49, 56, 66, 92], [2, 22, 60, 70, 77, 87], [1, 8, 18, 44, 64, 102], [17, 37, 75, 85, 92, 102], [6, 32, 52, 90, 100, 107], [19, 57, 67, 74, 84, 110], [33, 43, 50, 60, 86, 106], [7, 45, 55, 62, 72, 98], [18, 28, 35, 45, 71, 91], [23, 43, 81, 91, 98, 108], [24, 34, 41, 51, 77, 97], [4, 11, 21, 47, 67, 105], [3, 13, 20, 30, 56, 76], [9, 19, 26, 36, 62, 82], [8, 28, 66, 76, 83, 93], [14, 34, 72, 82, 89, 99], [30, 40, 47, 57, 83, 103], [36, 46, 53, 63, 89, 109], [27, 37, 44, 54, 80, 100], [10, 48, 58, 65, 75, 101], [0, 10, 17, 27, 53, 73], [21, 31, 38, 48, 74, 94], [2, 12, 38, 58, 96, 106], [15, 25, 32, 42, 68, 88], [7, 31, 39, 54, 75, 109], [70, 74, 92, 95, 97, 109], [43, 47, 65, 68, 70, 82], [12, 26, 42, 54, 74, 83], [12, 46, 55, 79, 87, 102], [16, 20, 38, 41, 43, 55], [22, 30, 45, 66, 100, 109], [64, 68, 86, 89, 91, 103], [1, 73, 77, 95, 98, 100], [12, 33, 67, 76, 100, 108], [31, 35, 53, 56, 58, 70], [5, 21, 33, 53, 62, 102], [12, 24, 44, 53, 93, 107], [5, 45, 59, 75, 87, 107], [13, 21, 36, 57, 91, 100], [19, 28, 52, 60, 75, 96], [2, 42, 56, 72, 84, 104], [52, 56, 74, 77, 79, 91], [2, 4, 16, 88, 92, 110], [28, 32, 50, 53, 55, 67], [4, 28, 36, 51, 72, 106], [19, 27, 42, 63, 97, 106], [10, 18, 33, 54, 88, 97], [19, 23, 41, 44, 46, 58], [31, 40, 64, 72, 87, 108], [10, 19, 43, 51, 66, 87], [55, 59, 77, 80, 82, 94], [2, 18, 30, 50, 59, 99], [14, 23, 63, 77, 93, 105], [67, 71, 89, 92, 94, 106], [13, 17, 35, 38, 40, 52], [18, 52, 61, 85, 93, 108], [11, 27, 39, 59, 68, 108], [7, 16, 40, 48, 63, 84], [4, 76, 80, 98, 101, 103], [21, 35, 51, 63, 83, 92], [1, 5, 23, 26, 28, 40], [49, 53, 71, 74, 76, 88], [2, 5, 7, 19, 91, 95], [5, 8, 10, 22, 94, 98], [0, 34, 43, 67, 75, 90], [5, 14, 54, 68, 84, 96], [7, 15, 30, 51, 85, 94], [0, 14, 30, 42, 62, 71], [7, 79, 83, 101, 104, 106], [2, 20, 23, 25, 37, 109], [6, 26, 35, 75, 89, 105], [1, 10, 34, 42, 57, 78], [14, 17, 19, 31, 103, 107], [18, 32, 48, 60, 80, 89], [1, 25, 33, 48, 69, 103], [0, 21, 55, 64, 88, 96], [0, 15, 36, 70, 79, 103], [6, 27, 61, 70, 94, 102], [6, 21, 42, 76, 85, 109], [9, 29, 38, 78, 92, 108], [25, 34, 58, 66, 81, 102], [24, 38, 54, 66, 86, 95], [10, 82, 86, 104, 107, 109], [22, 31, 55, 63, 78, 99], [6, 20, 36, 48, 68, 77], [15, 49, 58, 82, 90, 105], [15, 29, 45, 57, 77, 86], [4, 13, 37, 45, 60, 81], [4, 12, 27, 48, 82, 91], [61, 65, 83, 86, 88, 100], [16, 25, 49, 57, 72, 93], [9, 21, 41, 50, 90, 104], [37, 41, 59, 62, 64, 76], [39, 53, 69, 81, 101, 110], [4, 8, 26, 29, 31, 43], [8, 24, 36, 56, 65, 105], [3, 17, 33, 45, 65, 74], [9, 23, 39, 51, 71, 80], [25, 29, 47, 50, 52, 64], [3, 24, 58, 67, 91, 99], [3, 18, 39, 73, 82, 106], [15, 27, 47, 56, 96, 110], [6, 40, 49, 73, 81, 96], [9, 30, 64, 73, 97, 105], [0, 12, 32, 41, 81, 95], [2, 11, 51, 65, 81, 93], [9, 43, 52, 76, 84, 99], [40, 44, 62, 65, 67, 79], [36, 50, 66, 78, 98, 107], [22, 26, 44, 47, 49, 61], [28, 37, 61, 69, 84, 105], [3, 37, 46, 70, 78, 93], [17, 20, 22, 34, 106, 110], [7, 11, 29, 32, 34, 46], [3, 23, 32, 72, 86, 102], [1, 13, 85, 89, 107, 110], [8, 48, 62, 78, 90, 110], [16, 24, 39, 60, 94, 103], [30, 44, 60, 72, 92, 101], [3, 15, 35, 44, 84, 98], [27, 41, 57, 69, 89, 98], [8, 11, 13, 25, 97, 101], [6, 18, 38, 47, 87, 101], [13, 22, 46, 54, 69, 90], [58, 62, 80, 83, 85, 97], [33, 47, 63, 75, 95, 104], [11, 20, 60, 74, 90, 102], [0, 20, 29, 69, 83, 99], [46, 50, 68, 71, 73, 85], [11, 14, 16, 28, 100, 104], [8, 17, 57, 71, 87, 99], [17, 26, 66, 80, 96, 108], [1, 9, 24, 45, 79, 88], [34, 38, 56, 59, 61, 73], [10, 14, 32, 35, 37, 49], [3, 22, 27, 40, 50, 92], [9, 28, 33, 46, 56, 98], [7, 17, 59, 81, 100, 105], [4, 9, 22, 32, 74, 96], [38, 60, 79, 84, 97, 107], [2, 24, 43, 48, 61, 71], [21, 40, 45, 58, 68, 110], [20, 42, 61, 66, 79, 89], [8, 30, 49, 54, 67, 77], [35, 57, 76, 81, 94, 104], [8, 50, 72, 91, 96, 109], [14, 36, 55, 60, 73, 83], [17, 39, 58, 63, 76, 86], [32, 54, 73, 78, 91, 101], [10, 20, 62, 84, 103, 108], [15, 34, 39, 52, 62, 104], [7, 12, 25, 35, 77, 99], [1, 6, 19, 29, 71, 93], [0, 19, 24, 37, 47, 89], [13, 18, 31, 41, 83, 105], [5, 47, 69, 88, 93, 106], [11, 33, 52, 57, 70, 80], [10, 15, 28, 38, 80, 102], [41, 63, 82, 87, 100, 110], [3, 16, 26, 68, 90, 109], [5, 27, 46, 51, 64, 74], [18, 37, 42, 55, 65, 107], [2, 44, 66, 85, 90, 103], [6, 25, 30, 43, 53, 95], [0, 13, 23, 65, 87, 106], [23, 45, 64, 69, 82, 92], [1, 11, 53, 75, 94, 99], [12, 31, 36, 49, 59, 101], [26, 48, 67, 72, 85, 95], [4, 14, 56, 78, 97, 102], [16, 21, 34, 44, 86, 108], [29, 51, 70, 75, 88, 98], [34, 54, 79, 92, 98, 105], [25, 45, 70, 83, 89, 96], [11, 17, 24, 64, 84, 109], [6, 46, 66, 91, 104, 110], [7, 20, 26, 33, 73, 93], [16, 36, 61, 74, 80, 87], [28, 48, 73, 86, 92, 99], [18, 43, 56, 62, 69, 109], [9, 34, 47, 53, 60, 100], [8, 14, 21, 61, 81, 106], [15, 40, 53, 59, 66, 106], [10, 23, 29, 36, 76, 96], [7, 27, 52, 65, 71, 78], [19, 32, 38, 45, 85, 105], [2, 8, 15, 55, 75, 100], [19, 39, 64, 77, 83, 90], [0, 40, 60, 85, 98, 104], [37, 57, 82, 95, 101, 108], [22, 42, 67, 80, 86, 93], [4, 17, 23, 30, 70, 90], [5, 12, 52, 72, 97, 110], [10, 30, 55, 68, 74, 81], [3, 43, 63, 88, 101, 107], [0, 25, 38, 44, 51, 91], [16, 29, 35, 42, 82, 102], [5, 11, 18, 58, 78, 103], [12, 37, 50, 56, 63, 103], [2, 9, 49, 69, 94, 107], [1, 21, 46, 59, 65, 72], [22, 35, 41, 48, 88, 108], [13, 33, 58, 71, 77, 84], [13, 26, 32, 39, 79, 99], [1, 14, 20, 27, 67, 87], [3, 28, 41, 47, 54, 94], [31, 51, 76, 89, 95, 102], [6, 31, 44, 50, 57, 97], [4, 24, 49, 62, 68, 75]]
\item 1 \{1=47952, 2=516150, 3=2201796, 4=2362302\} [[8, 39, 40, 41, 42, 70], [0, 1, 2, 3, 31, 80], [30, 31, 32, 33, 61, 110], [6, 7, 8, 9, 37, 86], [23, 54, 55, 56, 57, 85], [15, 16, 17, 18, 46, 95], [27, 28, 29, 30, 58, 107], [24, 25, 26, 27, 55, 104], [26, 57, 58, 59, 60, 88], [29, 60, 61, 62, 63, 91], [25, 74, 105, 106, 107, 108], [41, 72, 73, 74, 75, 103], [7, 56, 87, 88, 89, 90], [10, 59, 90, 91, 92, 93], [13, 62, 93, 94, 95, 96], [22, 71, 102, 103, 104, 105], [9, 10, 11, 12, 40, 89], [44, 75, 76, 77, 78, 106], [18, 19, 20, 21, 49, 98], [11, 42, 43, 44, 45, 73], [3, 4, 5, 6, 34, 83], [16, 65, 96, 97, 98, 99], [19, 68, 99, 100, 101, 102], [17, 48, 49, 50, 51, 79], [47, 78, 79, 80, 81, 109], [20, 51, 52, 53, 54, 82], [12, 13, 14, 15, 43, 92], [0, 28, 77, 108, 109, 110], [4, 53, 84, 85, 86, 87], [1, 50, 81, 82, 83, 84], [32, 63, 64, 65, 66, 94], [2, 33, 34, 35, 36, 64], [21, 22, 23, 24, 52, 101], [38, 69, 70, 71, 72, 100], [14, 45, 46, 47, 48, 76], [35, 66, 67, 68, 69, 97], [5, 36, 37, 38, 39, 67], [8, 14, 19, 33, 58, 89], [27, 38, 60, 76, 87, 95], [6, 31, 62, 92, 98, 103], [1, 45, 60, 64, 67, 104], [15, 30, 34, 37, 74, 82], [0, 11, 33, 49, 60, 68], [7, 38, 68, 74, 79, 93], [9, 20, 42, 58, 69, 77], [29, 59, 65, 70, 84, 109], [2, 32, 38, 43, 57, 82], [23, 53, 59, 64, 78, 103], [33, 44, 66, 82, 93, 101], [10, 21, 29, 72, 83, 105], [4, 35, 65, 71, 76, 90], [14, 22, 66, 81, 85, 88], [0, 15, 19, 22, 59, 67], [29, 37, 81, 96, 100, 103], [4, 15, 23, 66, 77, 99], [39, 50, 72, 88, 99, 107], [16, 47, 77, 83, 88, 102], [9, 25, 36, 44, 87, 98], [26, 32, 37, 51, 76, 107], [2, 45, 56, 78, 94, 105], [6, 10, 13, 50, 58, 102], [18, 33, 37, 40, 77, 85], [23, 31, 75, 90, 94, 97], [10, 41, 71, 77, 82, 96], [27, 42, 46, 49, 86, 94], [12, 23, 45, 61, 72, 80], [36, 51, 55, 58, 95, 103], [15, 31, 42, 50, 93, 104], [11, 17, 22, 36, 61, 92], [26, 34, 78, 93, 97, 100], [24, 39, 43, 46, 83, 91], [21, 36, 40, 43, 80, 88], [21, 37, 48, 56, 99, 110], [2, 7, 21, 46, 77, 107], [9, 34, 65, 95, 101, 106], [5, 27, 43, 54, 62, 105], [18, 34, 45, 53, 96, 107], [6, 21, 25, 28, 65, 73], [22, 53, 83, 89, 94, 108], [6, 22, 33, 41, 84, 95], [12, 28, 39, 47, 90, 101], [17, 23, 28, 42, 67, 98], [19, 50, 80, 86, 91, 105], [32, 40, 84, 99, 103, 106], [35, 43, 87, 102, 106, 109], [29, 35, 40, 54, 79, 110], [3, 11, 54, 65, 87, 103], [7, 51, 66, 70, 73, 110], [18, 29, 51, 67, 78, 86], [5, 11, 16, 30, 55, 86], [9, 17, 60, 71, 93, 109], [33, 48, 52, 55, 92, 100], [5, 48, 59, 81, 97, 108], [21, 32, 54, 70, 81, 89], [3, 19, 30, 38, 81, 92], [39, 54, 58, 61, 98, 106], [2, 24, 40, 51, 59, 102], [1, 15, 40, 71, 101, 107], [6, 14, 57, 68, 90, 106], [2, 10, 54, 69, 73, 76], [5, 10, 24, 49, 80, 110], [0, 25, 56, 86, 92, 97], [23, 29, 34, 48, 73, 104], [0, 4, 7, 44, 52, 96], [13, 24, 32, 75, 86, 108], [5, 13, 57, 72, 76, 79], [24, 35, 57, 73, 84, 92], [1, 32, 62, 68, 73, 87], [17, 25, 69, 84, 88, 91], [5, 35, 41, 46, 60, 85], [20, 28, 72, 87, 91, 94], [14, 20, 25, 39, 64, 95], [12, 37, 68, 98, 104, 109], [1, 38, 46, 90, 105, 109], [6, 17, 39, 55, 66, 74], [20, 50, 56, 61, 75, 100], [12, 27, 31, 34, 71, 79], [2, 8, 13, 27, 52, 83], [8, 38, 44, 49, 63, 88], [9, 13, 16, 53, 61, 105], [8, 30, 46, 57, 65, 108], [3, 18, 22, 25, 62, 70], [11, 41, 47, 52, 66, 91], [30, 41, 63, 79, 90, 98], [0, 16, 27, 35, 78, 89], [3, 28, 59, 89, 95, 100], [26, 56, 62, 67, 81, 106], [4, 18, 43, 74, 104, 110], [3, 7, 10, 47, 55, 99], [9, 24, 28, 31, 68, 76], [1, 4, 41, 49, 93, 108], [20, 26, 31, 45, 70, 101], [12, 16, 19, 56, 64, 108], [7, 18, 26, 69, 80, 102], [11, 19, 63, 78, 82, 85], [3, 14, 36, 52, 63, 71], [42, 53, 75, 91, 102, 110], [42, 57, 61, 64, 101, 109], [8, 16, 60, 75, 79, 82], [14, 44, 50, 55, 69, 94], [4, 48, 63, 67, 70, 107], [17, 47, 53, 58, 72, 97], [13, 44, 74, 80, 85, 99], [1, 12, 20, 63, 74, 96], [30, 45, 49, 52, 89, 97], [15, 26, 48, 64, 75, 83], [0, 8, 51, 62, 84, 100], [36, 47, 69, 85, 96, 104], [25, 30, 54, 60, 83, 99], [20, 36, 73, 78, 102, 108], [13, 20, 23, 38, 84, 89], [23, 36, 49, 82, 100, 106], [39, 44, 79, 86, 89, 104], [31, 49, 55, 83, 96, 109], [2, 37, 44, 47, 62, 108], [2, 15, 28, 61, 79, 85], [14, 30, 67, 72, 96, 102], [4, 22, 28, 56, 69, 82], [19, 24, 48, 54, 77, 93], [21, 27, 50, 66, 103, 108], [8, 54, 59, 94, 101, 104], [5, 18, 31, 64, 82, 88], [0, 6, 29, 45, 82, 87], [6, 19, 52, 70, 76, 104], [19, 37, 43, 71, 84, 97], [22, 29, 32, 47, 93, 98], [7, 14, 17, 32, 78, 83], [0, 13, 46, 64, 70, 98], [2, 18, 55, 60, 84, 90], [7, 12, 36, 42, 65, 81], [19, 26, 29, 44, 90, 95], [28, 46, 52, 80, 93, 106], [16, 21, 45, 51, 74, 90], [5, 51, 56, 91, 98, 101], [10, 28, 34, 62, 75, 88], [34, 41, 44, 59, 105, 110], [15, 52, 57, 81, 87, 110], [5, 8, 23, 69, 74, 109], [24, 29, 64, 71, 74, 89], [7, 40, 58, 64, 92, 105], [15, 20, 55, 62, 65, 80], [8, 21, 34, 67, 85, 91], [0, 5, 40, 47, 50, 65], [9, 15, 38, 54, 91, 96], [2, 5, 20, 66, 71, 106], [4, 37, 55, 61, 89, 102], [18, 24, 47, 63, 100, 105], [13, 18, 42, 48, 71, 87], [45, 50, 85, 92, 95, 110], [7, 25, 31, 59, 72, 85], [28, 35, 38, 53, 99, 104], [1, 7, 35, 48, 61, 94], [16, 34, 40, 68, 81, 94], [4, 9, 33, 39, 62, 78], [21, 26, 61, 68, 71, 86], [10, 15, 39, 45, 68, 84], [20, 33, 46, 79, 97, 103], [17, 33, 70, 75, 99, 105], [22, 27, 51, 57, 80, 96], [25, 43, 49, 77, 90, 103], [28, 33, 57, 63, 86, 102], [1, 29, 42, 55, 88, 106], [6, 12, 35, 51, 88, 93], [6, 43, 48, 72, 78, 101], [9, 22, 55, 73, 79, 107], [25, 32, 35, 50, 96, 101], [9, 14, 49, 56, 59, 74], [4, 10, 38, 51, 64, 97], [3, 8, 43, 50, 53, 68], [22, 40, 46, 74, 87, 100], [36, 41, 76, 83, 86, 101], [16, 22, 50, 63, 76, 109], [7, 13, 41, 54, 67, 100], [10, 16, 44, 57, 70, 103], [13, 31, 37, 65, 78, 91], [3, 26, 42, 79, 84, 108], [1, 6, 30, 36, 59, 75], [1, 34, 52, 58, 86, 99], [31, 36, 60, 66, 89, 105], [12, 18, 41, 57, 94, 99], [4, 11, 14, 29, 75, 80], [3, 16, 49, 67, 73, 101], [2, 48, 53, 88, 95, 98], [5, 21, 58, 63, 87, 93], [31, 38, 41, 56, 102, 107], [11, 57, 62, 97, 104, 107], [9, 46, 51, 75, 81, 104], [14, 27, 40, 73, 91, 97], [0, 37, 42, 66, 72, 95], [8, 24, 61, 66, 90, 96], [12, 25, 58, 76, 82, 110], [10, 43, 61, 67, 95, 108], [10, 17, 20, 35, 81, 86], [33, 38, 73, 80, 83, 98], [13, 19, 47, 60, 73, 106], [26, 39, 52, 85, 103, 109], [12, 49, 54, 78, 84, 107], [4, 32, 45, 58, 91, 109], [42, 47, 82, 89, 92, 107], [18, 23, 58, 65, 68, 83], [3, 40, 45, 69, 75, 98], [1, 19, 25, 53, 66, 79], [14, 60, 65, 100, 107, 110], [34, 39, 63, 69, 92, 108], [16, 23, 26, 41, 87, 92], [0, 24, 30, 53, 69, 106], [11, 24, 37, 70, 88, 94], [0, 23, 39, 76, 81, 105], [3, 9, 32, 48, 85, 90], [27, 32, 67, 74, 77, 92], [6, 11, 46, 53, 56, 71], [11, 27, 64, 69, 93, 99], [17, 30, 43, 76, 94, 100], [12, 17, 52, 59, 62, 77], [15, 21, 44, 60, 97, 102], [2, 17, 63, 68, 103, 110], [3, 27, 33, 56, 72, 109], [30, 35, 70, 77, 80, 95], [1, 8, 11, 26, 72, 77], [12, 48, 60, 69, 86, 103], [2, 49, 65, 72, 92, 104], [6, 23, 40, 60, 96, 108], [1, 18, 76, 89, 91, 103], [43, 56, 58, 70, 79, 96], [8, 25, 45, 81, 93, 102], [7, 20, 22, 34, 43, 60], [13, 33, 69, 81, 90, 107], [15, 27, 36, 53, 70, 90], [8, 10, 22, 31, 48, 106], [37, 53, 60, 80, 92, 101], [2, 4, 16, 25, 42, 100], [19, 32, 34, 46, 55, 72], [10, 23, 25, 37, 46, 63], [12, 70, 83, 85, 97, 106], [5, 12, 32, 44, 53, 100], [11, 20, 67, 83, 90, 110], [0, 20, 32, 41, 88, 104], [1, 10, 27, 85, 98, 100], [8, 20, 29, 76, 92, 99], [8, 15, 35, 47, 56, 103], [1, 14, 16, 28, 37, 54], [52, 65, 67, 79, 88, 105], [9, 21, 30, 47, 64, 84], [18, 30, 39, 56, 73, 93], [7, 23, 30, 50, 62, 71], [34, 50, 57, 77, 89, 98], [1, 17, 24, 44, 56, 65], [4, 20, 27, 47, 59, 68], [3, 23, 35, 44, 91, 107], [28, 44, 51, 71, 83, 92], [0, 58, 71, 73, 85, 94], [28, 41, 43, 55, 64, 81], [5, 52, 68, 75, 95, 107], [6, 42, 54, 63, 80, 97], [2, 9, 29, 41, 50, 97], [1, 13, 22, 39, 97, 110], [10, 30, 66, 78, 87, 104], [7, 24, 82, 95, 97, 109], [10, 19, 36, 94, 107, 109], [9, 67, 80, 82, 94, 103], [4, 21, 79, 92, 94, 106], [4, 13, 30, 88, 101, 103], [5, 7, 19, 28, 45, 103], [49, 62, 64, 76, 85, 102], [2, 11, 58, 74, 81, 101], [7, 16, 33, 91, 104, 106], [0, 17, 34, 54, 90, 102], [22, 35, 37, 49, 58, 75], [5, 22, 42, 78, 90, 99], [3, 20, 37, 57, 93, 105], [25, 38, 40, 52, 61, 78], [7, 27, 63, 75, 84, 101], [8, 17, 64, 80, 87, 107], [11, 18, 38, 50, 59, 106], [43, 59, 66, 86, 98, 107], [4, 17, 19, 31, 40, 57], [13, 26, 28, 40, 49, 66], [22, 38, 45, 65, 77, 86], [6, 15, 32, 49, 69, 105], [6, 64, 77, 79, 91, 100], [55, 68, 70, 82, 91, 108], [5, 17, 26, 73, 89, 96], [16, 36, 72, 84, 93, 110], [15, 73, 86, 88, 100, 109], [27, 39, 48, 65, 82, 102], [11, 23, 32, 79, 95, 102], [3, 15, 24, 41, 58, 78], [40, 53, 55, 67, 76, 93], [34, 47, 49, 61, 70, 87], [14, 26, 35, 82, 98, 105], [14, 31, 51, 87, 99, 108], [5, 14, 61, 77, 84, 104], [4, 24, 60, 72, 81, 98], [19, 35, 42, 62, 74, 83], [2, 14, 23, 70, 86, 93], [3, 61, 74, 76, 88, 97], [33, 45, 54, 71, 88, 108], [24, 36, 45, 62, 79, 99], [16, 32, 39, 59, 71, 80], [11, 13, 25, 34, 51, 109], [3, 39, 51, 60, 77, 94], [0, 12, 21, 38, 55, 75], [46, 59, 61, 73, 82, 99], [13, 29, 36, 56, 68, 77], [8, 55, 71, 78, 98, 110], [15, 51, 63, 72, 89, 106], [2, 19, 39, 75, 87, 96], [14, 21, 41, 53, 62, 109], [30, 42, 51, 68, 85, 105], [40, 56, 63, 83, 95, 104], [46, 62, 69, 89, 101, 110], [6, 26, 38, 47, 94, 110], [37, 50, 52, 64, 73, 90], [31, 47, 54, 74, 86, 95], [0, 36, 48, 57, 74, 91], [31, 44, 46, 58, 67, 84], [9, 18, 35, 52, 72, 108], [9, 45, 57, 66, 83, 100], [21, 33, 42, 59, 76, 96], [10, 26, 33, 53, 65, 74], [0, 9, 26, 43, 63, 99], [11, 28, 48, 84, 96, 105], [16, 29, 31, 43, 52, 69], [25, 41, 48, 68, 80, 89], [3, 12, 29, 46, 66, 102], [18, 54, 66, 75, 92, 109], [12, 24, 33, 50, 67, 87], [6, 18, 27, 44, 61, 81], [17, 29, 38, 85, 101, 108], [1, 21, 57, 69, 78, 95], [4, 26, 36, 46, 50, 54], [0, 10, 14, 18, 79, 101], [16, 38, 48, 58, 62, 66], [40, 62, 72, 82, 86, 90], [52, 74, 84, 94, 98, 102], [49, 71, 81, 91, 95, 99], [10, 32, 42, 52, 56, 60], [3, 13, 17, 21, 82, 104], [31, 53, 63, 73, 77, 81], [11, 21, 31, 35, 39, 100], [6, 16, 20, 24, 85, 107], [8, 18, 28, 32, 36, 97], [4, 8, 12, 73, 95, 105], [9, 19, 23, 27, 88, 110], [13, 35, 45, 55, 59, 63], [25, 47, 57, 67, 71, 75], [1, 23, 33, 43, 47, 51], [37, 59, 69, 79, 83, 87], [5, 15, 25, 29, 33, 94], [20, 30, 40, 44, 48, 109], [34, 56, 66, 76, 80, 84], [58, 80, 90, 100, 104, 108], [22, 44, 54, 64, 68, 72], [7, 11, 15, 76, 98, 108], [46, 68, 78, 88, 92, 96], [1, 5, 9, 70, 92, 102], [7, 29, 39, 49, 53, 57], [43, 65, 75, 85, 89, 93], [17, 27, 37, 41, 45, 106], [19, 41, 51, 61, 65, 69], [55, 77, 87, 97, 101, 105], [0, 61, 83, 93, 103, 107], [14, 24, 34, 38, 42, 103], [3, 64, 86, 96, 106, 110], [2, 12, 22, 26, 30, 91], [28, 50, 60, 70, 74, 78], [2, 6, 67, 89, 99, 109]]
\item 1 \{1=39960, 2=535464, 3=2202240, 4=2350536\} [[8, 39, 40, 41, 42, 70], [0, 1, 2, 3, 31, 80], [30, 31, 32, 33, 61, 110], [6, 7, 8, 9, 37, 86], [23, 54, 55, 56, 57, 85], [15, 16, 17, 18, 46, 95], [27, 28, 29, 30, 58, 107], [24, 25, 26, 27, 55, 104], [26, 57, 58, 59, 60, 88], [29, 60, 61, 62, 63, 91], [25, 74, 105, 106, 107, 108], [41, 72, 73, 74, 75, 103], [7, 56, 87, 88, 89, 90], [10, 59, 90, 91, 92, 93], [13, 62, 93, 94, 95, 96], [22, 71, 102, 103, 104, 105], [9, 10, 11, 12, 40, 89], [44, 75, 76, 77, 78, 106], [18, 19, 20, 21, 49, 98], [11, 42, 43, 44, 45, 73], [3, 4, 5, 6, 34, 83], [16, 65, 96, 97, 98, 99], [19, 68, 99, 100, 101, 102], [17, 48, 49, 50, 51, 79], [47, 78, 79, 80, 81, 109], [20, 51, 52, 53, 54, 82], [12, 13, 14, 15, 43, 92], [0, 28, 77, 108, 109, 110], [4, 53, 84, 85, 86, 87], [1, 50, 81, 82, 83, 84], [32, 63, 64, 65, 66, 94], [2, 33, 34, 35, 36, 64], [21, 22, 23, 24, 52, 101], [38, 69, 70, 71, 72, 100], [14, 45, 46, 47, 48, 76], [35, 66, 67, 68, 69, 97], [5, 36, 37, 38, 39, 67], [7, 17, 29, 44, 93, 104], [9, 54, 72, 76, 81, 95], [5, 54, 65, 79, 89, 101], [45, 56, 70, 80, 92, 107], [1, 32, 37, 43, 91, 99], [3, 17, 42, 87, 105, 109], [33, 44, 58, 68, 80, 95], [8, 33, 78, 96, 100, 105], [21, 32, 46, 56, 68, 83], [22, 30, 43, 74, 79, 85], [18, 36, 40, 45, 59, 84], [8, 22, 32, 44, 59, 108], [36, 47, 61, 71, 83, 98], [10, 41, 46, 52, 100, 108], [42, 53, 67, 77, 89, 104], [43, 51, 64, 95, 100, 106], [12, 23, 37, 47, 59, 74], [5, 20, 69, 80, 94, 104], [36, 54, 58, 63, 77, 102], [21, 66, 84, 88, 93, 107], [8, 13, 19, 67, 75, 88], [3, 21, 25, 30, 44, 69], [23, 28, 34, 82, 90, 103], [11, 36, 81, 99, 103, 108], [7, 38, 43, 49, 97, 105], [29, 34, 40, 88, 96, 109], [4, 52, 60, 73, 104, 109], [30, 41, 55, 65, 77, 92], [39, 57, 61, 66, 80, 105], [2, 27, 72, 90, 94, 99], [24, 35, 49, 59, 71, 86], [15, 26, 40, 50, 62, 77], [8, 23, 72, 83, 97, 107], [1, 6, 20, 45, 90, 108], [30, 48, 52, 57, 71, 96], [0, 4, 9, 23, 48, 93], [9, 20, 34, 44, 56, 71], [6, 10, 15, 29, 54, 99], [5, 19, 29, 41, 56, 105], [6, 51, 69, 73, 78, 92], [40, 48, 61, 92, 97, 103], [12, 30, 34, 39, 53, 78], [18, 63, 81, 85, 90, 104], [2, 17, 66, 77, 91, 101], [11, 60, 71, 85, 95, 107], [1, 49, 57, 70, 101, 106], [24, 69, 87, 91, 96, 110], [27, 38, 52, 62, 74, 89], [15, 60, 78, 82, 87, 101], [6, 19, 50, 55, 61, 109], [2, 51, 62, 76, 86, 98], [33, 51, 55, 60, 74, 99], [5, 17, 32, 81, 92, 106], [26, 31, 37, 85, 93, 106], [15, 19, 24, 38, 63, 108], [46, 54, 67, 98, 103, 109], [7, 15, 28, 59, 64, 70], [6, 17, 31, 41, 53, 68], [11, 26, 75, 86, 100, 110], [13, 21, 34, 65, 70, 76], [5, 10, 16, 64, 72, 85], [18, 29, 43, 53, 65, 80], [9, 13, 18, 32, 57, 102], [1, 9, 22, 53, 58, 64], [20, 25, 31, 79, 87, 100], [0, 13, 44, 49, 55, 103], [3, 16, 47, 52, 58, 106], [0, 45, 63, 67, 72, 86], [16, 24, 37, 68, 73, 79], [6, 24, 28, 33, 47, 72], [37, 45, 58, 89, 94, 100], [4, 10, 58, 66, 79, 110], [2, 14, 29, 78, 89, 103], [10, 18, 31, 62, 67, 73], [0, 14, 39, 84, 102, 106], [4, 12, 25, 56, 61, 67], [8, 20, 35, 84, 95, 109], [12, 16, 21, 35, 60, 105], [31, 39, 52, 83, 88, 94], [19, 27, 40, 71, 76, 82], [3, 7, 12, 26, 51, 96], [14, 19, 25, 73, 81, 94], [13, 23, 35, 50, 99, 110], [21, 39, 43, 48, 62, 87], [1, 11, 23, 38, 87, 98], [2, 16, 26, 38, 53, 102], [24, 42, 46, 51, 65, 90], [0, 11, 25, 35, 47, 62], [9, 27, 31, 36, 50, 75], [15, 33, 37, 42, 56, 81], [17, 22, 28, 76, 84, 97], [14, 63, 74, 88, 98, 110], [27, 45, 49, 54, 68, 93], [8, 57, 68, 82, 92, 104], [12, 57, 75, 79, 84, 98], [0, 18, 22, 27, 41, 66], [4, 14, 26, 41, 90, 101], [2, 7, 13, 61, 69, 82], [39, 50, 64, 74, 86, 101], [10, 20, 32, 47, 96, 107], [11, 16, 22, 70, 78, 91], [3, 48, 66, 70, 75, 89], [5, 30, 75, 93, 97, 102], [48, 59, 73, 83, 95, 110], [42, 60, 64, 69, 83, 108], [4, 35, 40, 46, 94, 102], [28, 36, 49, 80, 85, 91], [25, 33, 46, 77, 82, 88], [3, 14, 28, 38, 50, 65], [34, 42, 55, 86, 91, 97], [1, 7, 55, 63, 76, 107], [1, 26, 39, 65, 95, 103], [8, 61, 76, 79, 90, 102], [6, 23, 76, 91, 94, 105], [0, 16, 33, 71, 87, 92], [7, 32, 45, 71, 101, 109], [3, 20, 73, 88, 91, 102], [8, 38, 46, 55, 80, 93], [19, 34, 37, 48, 60, 77], [25, 40, 43, 54, 66, 83], [15, 53, 69, 74, 93, 109], [34, 49, 52, 63, 75, 92], [22, 37, 40, 51, 63, 80], [11, 19, 28, 53, 66, 92], [16, 31, 34, 45, 57, 74], [32, 48, 53, 72, 88, 105], [5, 18, 44, 74, 82, 91], [14, 22, 31, 56, 69, 95], [5, 58, 73, 76, 87, 99], [20, 28, 37, 62, 75, 101], [12, 50, 66, 71, 90, 106], [29, 37, 46, 71, 84, 110], [14, 27, 53, 83, 91, 100], [15, 20, 39, 55, 72, 110], [17, 30, 56, 86, 94, 103], [2, 10, 19, 44, 57, 83], [23, 36, 62, 92, 100, 109], [20, 36, 41, 60, 76, 93], [35, 51, 56, 75, 91, 108], [23, 31, 40, 65, 78, 104], [9, 14, 33, 49, 66, 104], [5, 21, 26, 45, 61, 78], [46, 61, 64, 75, 87, 104], [3, 15, 32, 85, 100, 103], [18, 34, 51, 89, 105, 110], [26, 34, 43, 68, 81, 107], [3, 19, 36, 74, 90, 95], [4, 19, 22, 33, 45, 62], [10, 13, 24, 36, 53, 106], [4, 21, 59, 75, 80, 99], [8, 16, 25, 50, 63, 89], [0, 38, 54, 59, 78, 94], [9, 21, 38, 91, 106, 109], [23, 53, 61, 70, 95, 108], [13, 16, 27, 39, 56, 109], [40, 55, 58, 69, 81, 98], [5, 35, 43, 52, 77, 90], [15, 31, 48, 86, 102, 107], [0, 17, 70, 85, 88, 99], [1, 18, 56, 72, 77, 96], [17, 25, 34, 59, 72, 98], [17, 33, 38, 57, 73, 90], [10, 27, 65, 81, 86, 105], [6, 44, 60, 65, 84, 100], [31, 46, 49, 60, 72, 89], [13, 28, 31, 42, 54, 71], [11, 41, 49, 58, 83, 96], [6, 18, 35, 88, 103, 106], [7, 24, 62, 78, 83, 102], [4, 13, 38, 51, 77, 107], [7, 16, 41, 54, 80, 110], [7, 10, 21, 33, 50, 103], [20, 50, 58, 67, 92, 105], [6, 11, 30, 46, 63, 101], [10, 25, 28, 39, 51, 68], [28, 43, 46, 57, 69, 86], [2, 15, 41, 71, 79, 88], [5, 13, 22, 47, 60, 86], [12, 38, 68, 76, 85, 110], [3, 41, 57, 62, 81, 97], [3, 8, 27, 43, 60, 98], [2, 21, 37, 54, 92, 108], [11, 24, 50, 80, 88, 97], [14, 30, 35, 54, 70, 87], [29, 45, 50, 69, 85, 102], [9, 35, 65, 73, 82, 107], [0, 5, 24, 40, 57, 95], [23, 39, 44, 63, 79, 96], [2, 18, 23, 42, 58, 75], [13, 30, 68, 84, 89, 108], [9, 26, 79, 94, 97, 108], [12, 28, 45, 83, 99, 104], [4, 29, 42, 68, 98, 106], [6, 22, 39, 77, 93, 98], [11, 64, 79, 82, 93, 105], [4, 7, 18, 30, 47, 100], [8, 21, 47, 77, 85, 94], [17, 47, 55, 64, 89, 102], [37, 52, 55, 66, 78, 95], [26, 42, 47, 66, 82, 99], [9, 47, 63, 68, 87, 103], [11, 27, 32, 51, 67, 84], [1, 4, 15, 27, 44, 97], [6, 32, 62, 70, 79, 104], [14, 67, 82, 85, 96, 108], [3, 29, 59, 67, 76, 101], [2, 32, 40, 49, 74, 87], [8, 24, 29, 48, 64, 81], [9, 25, 42, 80, 96, 101], [0, 26, 56, 64, 73, 98], [1, 16, 19, 30, 42, 59], [1, 10, 35, 48, 74, 104], [7, 22, 25, 36, 48, 65], [2, 55, 70, 73, 84, 96], [14, 44, 52, 61, 86, 99], [20, 33, 59, 89, 97, 106], [43, 58, 61, 72, 84, 101], [49, 64, 67, 78, 90, 107], [0, 12, 29, 82, 97, 100], [52, 67, 70, 81, 93, 110], [1, 12, 24, 41, 94, 109], [12, 17, 36, 52, 69, 107], [6, 13, 59, 66, 81, 87], [22, 55, 67, 83, 87, 106], [0, 20, 61, 65, 68, 74], [17, 24, 39, 45, 75, 82], [38, 45, 60, 66, 96, 103], [10, 14, 17, 23, 60, 80], [0, 15, 21, 51, 58, 104], [6, 26, 67, 71, 74, 80], [9, 29, 70, 74, 77, 83], [1, 17, 21, 40, 67, 100], [16, 20, 23, 29, 66, 86], [19, 31, 47, 51, 70, 97], [2, 5, 11, 48, 68, 109], [26, 33, 48, 54, 84, 91], [15, 22, 68, 75, 90, 96], [5, 12, 27, 33, 63, 70], [5, 42, 62, 103, 107, 110], [14, 18, 37, 64, 97, 109], [11, 15, 34, 61, 94, 106], [40, 44, 47, 53, 90, 110], [5, 9, 28, 55, 88, 100], [22, 34, 50, 54, 73, 100], [1, 13, 29, 33, 52, 79], [18, 25, 71, 78, 93, 99], [2, 8, 45, 65, 106, 110], [10, 22, 38, 42, 61, 88], [25, 58, 70, 86, 90, 109], [8, 12, 31, 58, 91, 103], [4, 16, 32, 36, 55, 82], [2, 9, 24, 30, 60, 67], [18, 38, 79, 83, 86, 92], [3, 9, 39, 46, 92, 99], [4, 8, 11, 17, 54, 74], [14, 55, 59, 62, 68, 105], [1, 34, 46, 62, 66, 85], [14, 21, 36, 42, 72, 79], [37, 41, 44, 50, 87, 107], [6, 21, 27, 57, 64, 110], [19, 23, 26, 32, 69, 89], [22, 26, 29, 35, 72, 92], [17, 58, 62, 65, 71, 108], [6, 12, 42, 49, 95, 102], [4, 20, 24, 43, 70, 103], [0, 7, 53, 60, 75, 81], [3, 18, 24, 54, 61, 107], [4, 37, 49, 65, 69, 88], [16, 49, 61, 77, 81, 100], [0, 6, 36, 43, 89, 96], [31, 35, 38, 44, 81, 101], [1, 5, 8, 14, 51, 71], [25, 29, 32, 38, 75, 95], [10, 37, 70, 82, 98, 102], [15, 35, 76, 80, 83, 89], [7, 19, 35, 39, 58, 85], [7, 40, 52, 68, 72, 91], [32, 39, 54, 60, 90, 97], [12, 32, 73, 77, 80, 86], [3, 10, 56, 63, 78, 84], [2, 43, 47, 50, 56, 93], [4, 31, 64, 76, 92, 96], [3, 33, 40, 86, 93, 108], [30, 50, 91, 95, 98, 104], [7, 11, 14, 20, 57, 77], [8, 15, 30, 36, 66, 73], [12, 19, 65, 72, 87, 93], [3, 22, 49, 82, 94, 110], [19, 52, 64, 80, 84, 103], [23, 30, 45, 51, 81, 88], [7, 34, 67, 79, 95, 99], [10, 43, 55, 71, 75, 94], [5, 46, 50, 53, 59, 96], [25, 37, 53, 57, 76, 103], [28, 40, 56, 60, 79, 106], [13, 25, 41, 45, 64, 91], [33, 53, 94, 98, 101, 107], [31, 43, 59, 63, 82, 109], [13, 40, 73, 85, 101, 105], [9, 15, 45, 52, 98, 105], [7, 23, 27, 46, 73, 106], [16, 28, 44, 48, 67, 94], [13, 46, 58, 74, 78, 97], [0, 30, 37, 83, 90, 105], [29, 36, 51, 57, 87, 94], [35, 42, 57, 63, 93, 100], [3, 23, 64, 68, 71, 77], [2, 39, 59, 100, 104, 107], [4, 50, 57, 72, 78, 108], [9, 16, 62, 69, 84, 90], [1, 47, 54, 69, 75, 105], [36, 56, 97, 101, 104, 110], [8, 49, 53, 56, 62, 99], [24, 31, 77, 84, 99, 105], [28, 32, 35, 41, 78, 98], [16, 43, 76, 88, 104, 108], [0, 19, 46, 79, 91, 107], [21, 28, 74, 81, 96, 102], [11, 18, 33, 39, 69, 76], [34, 38, 41, 47, 84, 104], [24, 44, 85, 89, 92, 98], [21, 41, 82, 86, 89, 95], [20, 27, 42, 48, 78, 85], [27, 47, 88, 92, 95, 101], [2, 6, 25, 52, 85, 97], [11, 52, 56, 59, 65, 102], [41, 48, 63, 69, 99, 106], [13, 17, 20, 26, 63, 83], [12, 18, 48, 55, 101, 108], [27, 34, 80, 87, 102, 108], [1, 28, 61, 73, 89, 93], [10, 26, 30, 49, 76, 109], [44, 51, 66, 72, 102, 109], [12, 44, 54, 62, 64, 88], [2, 12, 20, 22, 46, 81], [22, 57, 89, 99, 107, 109], [29, 39, 47, 49, 73, 108], [18, 50, 60, 68, 70, 94], [3, 11, 13, 37, 72, 104], [5, 7, 31, 66, 98, 108], [11, 21, 29, 31, 55, 90], [17, 27, 35, 37, 61, 96], [3, 35, 45, 53, 55, 79], [6, 14, 16, 40, 75, 107], [33, 65, 75, 83, 85, 109], [2, 4, 28, 63, 95, 105], [8, 18, 26, 28, 52, 87], [0, 32, 42, 50, 52, 76], [4, 39, 71, 81, 89, 91], [19, 54, 86, 96, 104, 106], [1, 25, 60, 92, 102, 110], [1, 36, 68, 78, 86, 88], [21, 53, 63, 71, 73, 97], [30, 62, 72, 80, 82, 106], [23, 33, 41, 43, 67, 102], [10, 45, 77, 87, 95, 97], [15, 47, 57, 65, 67, 91], [7, 42, 74, 84, 92, 94], [24, 56, 66, 74, 76, 100], [9, 41, 51, 59, 61, 85], [14, 24, 32, 34, 58, 93], [0, 8, 10, 34, 69, 101], [27, 59, 69, 77, 79, 103], [6, 38, 48, 56, 58, 82], [9, 17, 19, 43, 78, 110], [13, 48, 80, 90, 98, 100], [26, 36, 44, 46, 70, 105], [20, 30, 38, 40, 64, 99], [5, 15, 23, 25, 49, 84], [16, 51, 83, 93, 101, 103]]
\item 1 \{1=17760, 2=512154, 3=2168940, 4=2429346\} [[8, 39, 40, 41, 42, 70], [0, 1, 2, 3, 31, 80], [30, 31, 32, 33, 61, 110], [6, 7, 8, 9, 37, 86], [23, 54, 55, 56, 57, 85], [15, 16, 17, 18, 46, 95], [27, 28, 29, 30, 58, 107], [24, 25, 26, 27, 55, 104], [26, 57, 58, 59, 60, 88], [29, 60, 61, 62, 63, 91], [25, 74, 105, 106, 107, 108], [41, 72, 73, 74, 75, 103], [7, 56, 87, 88, 89, 90], [10, 59, 90, 91, 92, 93], [13, 62, 93, 94, 95, 96], [22, 71, 102, 103, 104, 105], [9, 10, 11, 12, 40, 89], [44, 75, 76, 77, 78, 106], [18, 19, 20, 21, 49, 98], [11, 42, 43, 44, 45, 73], [3, 4, 5, 6, 34, 83], [16, 65, 96, 97, 98, 99], [19, 68, 99, 100, 101, 102], [17, 48, 49, 50, 51, 79], [47, 78, 79, 80, 81, 109], [20, 51, 52, 53, 54, 82], [12, 13, 14, 15, 43, 92], [0, 28, 77, 108, 109, 110], [4, 53, 84, 85, 86, 87], [1, 50, 81, 82, 83, 84], [32, 63, 64, 65, 66, 94], [2, 33, 34, 35, 36, 64], [21, 22, 23, 24, 52, 101], [38, 69, 70, 71, 72, 100], [14, 45, 46, 47, 48, 76], [35, 66, 67, 68, 69, 97], [5, 36, 37, 38, 39, 67], [11, 18, 50, 54, 71, 87], [8, 24, 59, 66, 98, 102], [23, 73, 79, 86, 88, 104], [14, 18, 35, 51, 86, 93], [15, 19, 22, 60, 72, 106], [40, 46, 53, 55, 71, 101], [9, 43, 63, 67, 70, 108], [14, 64, 70, 77, 79, 95], [0, 17, 33, 68, 75, 107], [26, 33, 65, 69, 86, 102], [8, 38, 88, 94, 101, 103], [3, 7, 10, 48, 60, 94], [19, 25, 32, 34, 50, 80], [0, 12, 46, 66, 70, 73], [16, 22, 29, 31, 47, 77], [1, 4, 42, 54, 88, 108], [25, 45, 49, 52, 90, 102], [22, 28, 35, 37, 53, 83], [5, 55, 61, 68, 70, 86], [0, 32, 36, 53, 69, 104], [20, 27, 59, 63, 80, 96], [9, 13, 16, 54, 66, 100], [46, 52, 59, 61, 77, 107], [2, 9, 41, 45, 62, 78], [20, 70, 76, 83, 85, 101], [20, 24, 41, 57, 92, 99], [13, 19, 26, 28, 44, 74], [28, 34, 41, 43, 59, 89], [3, 35, 39, 56, 72, 107], [8, 58, 64, 71, 73, 89], [9, 44, 51, 83, 87, 104], [23, 30, 62, 66, 83, 99], [8, 12, 29, 45, 80, 87], [5, 21, 56, 63, 95, 99], [15, 50, 57, 89, 93, 110], [28, 48, 52, 55, 93, 105], [18, 30, 64, 84, 88, 91], [22, 42, 46, 49, 87, 99], [2, 4, 20, 50, 100, 106], [19, 39, 43, 46, 84, 96], [30, 42, 76, 96, 100, 103], [0, 4, 7, 45, 57, 91], [27, 39, 73, 93, 97, 100], [2, 32, 82, 88, 95, 97], [17, 21, 38, 54, 89, 96], [29, 33, 50, 66, 101, 108], [2, 18, 53, 60, 92, 96], [21, 33, 67, 87, 91, 94], [1, 39, 51, 85, 105, 109], [34, 40, 47, 49, 65, 95], [12, 16, 19, 57, 69, 103], [3, 38, 45, 77, 81, 98], [49, 55, 62, 64, 80, 110], [6, 10, 13, 51, 63, 97], [23, 27, 44, 60, 95, 102], [5, 12, 44, 48, 65, 81], [33, 45, 79, 99, 103, 106], [6, 40, 60, 64, 67, 105], [5, 35, 85, 91, 98, 100], [17, 67, 73, 80, 82, 98], [6, 41, 48, 80, 84, 101], [16, 36, 40, 43, 81, 93], [3, 15, 49, 69, 73, 76], [4, 10, 17, 19, 35, 65], [10, 30, 34, 37, 75, 87], [7, 13, 20, 22, 38, 68], [0, 34, 54, 58, 61, 99], [9, 21, 55, 75, 79, 82], [31, 37, 44, 46, 62, 92], [7, 27, 31, 34, 72, 84], [5, 9, 26, 42, 77, 84], [4, 11, 13, 29, 59, 109], [1, 7, 14, 16, 32, 62], [31, 51, 55, 58, 96, 108], [2, 52, 58, 65, 67, 83], [25, 31, 38, 40, 56, 86], [29, 36, 68, 72, 89, 105], [11, 27, 62, 69, 101, 105], [5, 7, 23, 53, 103, 109], [6, 38, 42, 59, 75, 110], [11, 41, 91, 97, 104, 106], [14, 30, 65, 72, 104, 108], [26, 76, 82, 89, 91, 107], [14, 21, 53, 57, 74, 90], [11, 61, 67, 74, 76, 92], [26, 30, 47, 63, 98, 105], [1, 17, 47, 97, 103, 110], [12, 47, 54, 86, 90, 107], [36, 48, 82, 102, 106, 109], [4, 24, 28, 31, 69, 81], [11, 15, 32, 48, 83, 90], [0, 35, 42, 74, 78, 95], [14, 44, 94, 100, 107, 109], [8, 15, 47, 51, 68, 84], [37, 43, 50, 52, 68, 98], [3, 20, 36, 71, 78, 110], [3, 37, 57, 61, 64, 102], [13, 33, 37, 40, 78, 90], [10, 16, 23, 25, 41, 71], [32, 39, 71, 75, 92, 108], [29, 79, 85, 92, 94, 110], [1, 21, 25, 28, 66, 78], [12, 24, 58, 78, 82, 85], [6, 18, 52, 72, 76, 79], [15, 27, 61, 81, 85, 88], [24, 36, 70, 90, 94, 97], [18, 22, 25, 63, 75, 109], [43, 49, 56, 58, 74, 104], [17, 24, 56, 60, 77, 93], [2, 6, 23, 39, 74, 81], [1, 8, 10, 26, 56, 106], [15, 20, 33, 44, 55, 97], [28, 57, 62, 75, 86, 97], [7, 49, 78, 83, 96, 107], [25, 54, 59, 72, 83, 94], [2, 15, 26, 37, 79, 108], [21, 26, 39, 50, 61, 103], [1, 43, 72, 77, 90, 101], [24, 29, 42, 53, 64, 106], [0, 11, 22, 64, 93, 98], [1, 30, 35, 48, 59, 70], [19, 48, 53, 66, 77, 88], [37, 66, 71, 84, 95, 106], [4, 46, 75, 80, 93, 104], [8, 19, 61, 90, 95, 108], [5, 16, 58, 87, 92, 105], [12, 23, 34, 76, 105, 110], [0, 5, 18, 29, 40, 82], [9, 14, 27, 38, 49, 91], [18, 23, 36, 47, 58, 100], [34, 63, 68, 81, 92, 103], [10, 39, 44, 57, 68, 79], [31, 60, 65, 78, 89, 100], [13, 42, 47, 60, 71, 82], [4, 33, 38, 51, 62, 73], [22, 51, 56, 69, 80, 91], [3, 8, 21, 32, 43, 85], [40, 69, 74, 87, 98, 109], [3, 14, 25, 67, 96, 101], [27, 32, 45, 56, 67, 109], [2, 13, 55, 84, 89, 102], [6, 11, 24, 35, 46, 88], [12, 17, 30, 41, 52, 94], [9, 20, 31, 73, 102, 107], [7, 36, 41, 54, 65, 76], [6, 17, 28, 70, 99, 104], [10, 52, 81, 86, 99, 110], [16, 45, 50, 63, 74, 85], [0, 10, 72, 85, 96, 102], [7, 18, 24, 33, 43, 105], [12, 35, 38, 50, 60, 104], [6, 19, 30, 36, 45, 55], [1, 34, 44, 52, 71, 91], [1, 20, 40, 61, 94, 104], [32, 51, 74, 77, 89, 99], [0, 9, 19, 81, 94, 105], [29, 48, 71, 74, 86, 96], [4, 15, 21, 30, 40, 102], [4, 37, 47, 55, 74, 94], [30, 43, 54, 60, 69, 79], [0, 13, 24, 30, 39, 49], [7, 17, 25, 44, 64, 85], [9, 53, 72, 95, 98, 110], [8, 16, 35, 55, 76, 109], [51, 64, 75, 81, 90, 100], [5, 24, 47, 50, 62, 72], [11, 30, 53, 56, 68, 78], [22, 32, 40, 59, 79, 100], [9, 15, 24, 34, 96, 109], [11, 21, 65, 84, 107, 110], [5, 13, 32, 52, 73, 106], [20, 39, 62, 65, 77, 87], [10, 43, 53, 61, 80, 100], [4, 23, 43, 64, 97, 107], [10, 20, 28, 47, 67, 88], [35, 54, 77, 80, 92, 102], [2, 22, 43, 76, 86, 94], [25, 35, 43, 62, 82, 103], [27, 40, 51, 57, 66, 76], [0, 6, 15, 25, 87, 100], [7, 40, 50, 58, 77, 97], [0, 23, 26, 38, 48, 92], [0, 44, 63, 86, 89, 101], [1, 22, 55, 65, 73, 92], [3, 26, 29, 41, 51, 95], [3, 9, 18, 28, 90, 103], [1, 11, 19, 38, 58, 79], [42, 55, 66, 72, 81, 91], [12, 25, 36, 42, 51, 61], [7, 26, 46, 67, 100, 110], [13, 46, 56, 64, 83, 103], [5, 8, 20, 30, 74, 93], [11, 31, 52, 85, 95, 103], [23, 42, 65, 68, 80, 90], [18, 31, 42, 48, 57, 67], [16, 49, 59, 67, 86, 106], [10, 21, 27, 36, 46, 108], [26, 45, 68, 71, 83, 93], [19, 29, 37, 56, 76, 97], [48, 61, 72, 78, 87, 97], [6, 12, 21, 31, 93, 106], [2, 12, 56, 75, 98, 101], [13, 34, 67, 77, 85, 104], [6, 29, 32, 44, 54, 98], [3, 12, 22, 84, 97, 108], [7, 69, 82, 93, 99, 108], [15, 38, 41, 53, 63, 107], [14, 33, 56, 59, 71, 81], [4, 14, 22, 41, 61, 82], [24, 37, 48, 54, 63, 73], [57, 70, 81, 87, 96, 106], [8, 11, 23, 33, 77, 96], [38, 57, 80, 83, 95, 105], [41, 60, 83, 86, 98, 108], [17, 20, 32, 42, 86, 105], [15, 28, 39, 45, 54, 64], [39, 52, 63, 69, 78, 88], [8, 28, 49, 82, 92, 100], [54, 67, 78, 84, 93, 103], [1, 12, 18, 27, 37, 99], [8, 27, 50, 53, 65, 75], [9, 22, 33, 39, 48, 58], [2, 5, 17, 27, 71, 90], [28, 38, 46, 65, 85, 106], [21, 34, 45, 51, 60, 70], [5, 25, 46, 79, 89, 97], [2, 21, 44, 47, 59, 69], [11, 14, 26, 36, 80, 99], [3, 47, 66, 89, 92, 104], [4, 66, 79, 90, 96, 105], [18, 41, 44, 56, 66, 110], [8, 18, 62, 81, 104, 107], [1, 63, 76, 87, 93, 102], [4, 25, 58, 68, 76, 95], [16, 37, 70, 80, 88, 107], [19, 40, 73, 83, 91, 110], [16, 26, 34, 53, 73, 94], [3, 13, 75, 88, 99, 105], [2, 10, 29, 49, 70, 103], [2, 14, 24, 68, 87, 110], [7, 28, 61, 71, 79, 98], [14, 34, 55, 88, 98, 106], [20, 23, 35, 45, 89, 108], [14, 17, 29, 39, 83, 102], [5, 15, 59, 78, 101, 104], [60, 73, 84, 90, 99, 109], [6, 50, 69, 92, 95, 107], [6, 16, 78, 91, 102, 108], [10, 31, 64, 74, 82, 101], [19, 52, 62, 70, 89, 109], [36, 49, 60, 66, 75, 85], [31, 41, 49, 68, 88, 109], [45, 58, 69, 75, 84, 94], [33, 46, 57, 63, 72, 82], [9, 32, 35, 47, 57, 101], [3, 16, 27, 33, 42, 52], [13, 23, 31, 50, 70, 91], [17, 36, 59, 62, 74, 84], [17, 37, 58, 91, 101, 109], [6, 43, 57, 65, 71, 109], [1, 9, 46, 60, 68, 74], [1, 15, 23, 29, 67, 75], [3, 40, 54, 62, 68, 106], [22, 30, 67, 81, 89, 95], [22, 36, 44, 50, 88, 96], [0, 8, 14, 52, 60, 97], [5, 43, 51, 88, 102, 110], [5, 11, 49, 57, 94, 108], [31, 45, 53, 59, 97, 105], [16, 30, 38, 44, 82, 90], [25, 39, 47, 53, 91, 99], [37, 45, 82, 96, 104, 110], [12, 20, 26, 64, 72, 109], [10, 24, 32, 38, 76, 84], [0, 37, 51, 59, 65, 103], [13, 27, 35, 41, 79, 87], [7, 15, 52, 66, 74, 80], [31, 39, 76, 90, 98, 104], [28, 36, 73, 87, 95, 101], [34, 48, 56, 62, 100, 108], [25, 33, 70, 84, 92, 98], [19, 33, 41, 47, 85, 93], [7, 21, 29, 35, 73, 81], [34, 42, 79, 93, 101, 107], [9, 17, 23, 61, 69, 106], [10, 18, 55, 69, 77, 83], [13, 21, 58, 72, 80, 86], [16, 24, 61, 75, 83, 89], [4, 18, 26, 32, 70, 78], [4, 12, 49, 63, 71, 77], [6, 14, 20, 58, 66, 103], [3, 11, 17, 55, 63, 100], [2, 8, 46, 54, 91, 105], [19, 27, 64, 78, 86, 92], [2, 40, 48, 85, 99, 107], [28, 42, 50, 56, 94, 102], [2, 38, 61, 66, 93, 109], [22, 27, 54, 70, 74, 110], [17, 40, 45, 72, 88, 92], [11, 34, 39, 66, 82, 86], [1, 6, 33, 49, 53, 89], [13, 18, 45, 61, 65, 101], [4, 9, 36, 52, 56, 92], [16, 21, 48, 64, 68, 104], [4, 8, 44, 67, 72, 99], [8, 31, 36, 63, 79, 83], [35, 58, 63, 90, 106, 110], [2, 25, 30, 57, 73, 77], [7, 12, 39, 55, 59, 95], [24, 40, 44, 80, 103, 108], [9, 25, 29, 65, 88, 93], [18, 34, 38, 74, 97, 102], [3, 19, 23, 59, 82, 87], [10, 14, 50, 73, 78, 105], [23, 46, 51, 78, 94, 98], [19, 24, 51, 67, 71, 107], [10, 15, 42, 58, 62, 98], [15, 31, 35, 71, 94, 99], [7, 11, 47, 70, 75, 102], [20, 43, 48, 75, 91, 95], [21, 37, 41, 77, 100, 105], [29, 52, 57, 84, 100, 104], [5, 28, 33, 60, 76, 80], [32, 55, 60, 87, 103, 107], [6, 22, 26, 62, 85, 90], [13, 17, 53, 76, 81, 108], [26, 49, 54, 81, 97, 101], [14, 37, 42, 69, 85, 89], [12, 28, 32, 68, 91, 96], [0, 27, 43, 47, 83, 106], [3, 30, 46, 50, 86, 109], [0, 16, 20, 56, 79, 84], [1, 5, 41, 64, 69, 96], [0, 21, 62, 71, 76, 88], [12, 33, 74, 83, 88, 100], [29, 38, 43, 55, 78, 99], [6, 47, 56, 61, 73, 96], [15, 56, 65, 70, 82, 105], [0, 41, 50, 55, 67, 90], [6, 27, 68, 77, 82, 94], [3, 44, 53, 58, 70, 93], [7, 30, 51, 92, 101, 106], [5, 14, 19, 31, 54, 75], [4, 27, 48, 89, 98, 103], [2, 11, 16, 28, 51, 72], [20, 29, 34, 46, 69, 90], [32, 41, 46, 58, 81, 102], [21, 42, 83, 92, 97, 109], [9, 50, 59, 64, 76, 99], [3, 24, 65, 74, 79, 91], [18, 39, 80, 89, 94, 106], [10, 33, 54, 95, 104, 109], [8, 13, 25, 48, 69, 110], [8, 17, 22, 34, 57, 78], [2, 7, 19, 42, 63, 104], [26, 35, 40, 52, 75, 96], [1, 24, 45, 86, 95, 100], [18, 59, 68, 73, 85, 108], [4, 16, 39, 60, 101, 110], [11, 20, 25, 37, 60, 81], [14, 23, 28, 40, 63, 84], [5, 10, 22, 45, 66, 107], [9, 30, 71, 80, 85, 97], [12, 53, 62, 67, 79, 102], [1, 13, 36, 57, 98, 107], [35, 44, 49, 61, 84, 105], [17, 26, 31, 43, 66, 87], [38, 47, 52, 64, 87, 108], [15, 36, 77, 86, 91, 103], [23, 32, 37, 49, 72, 93]]
\item 1 \{0=1110, 1=38184, 2=572094, 3=2184036, 4=2332776\} [[8, 39, 40, 41, 42, 70], [0, 1, 2, 3, 31, 80], [30, 31, 32, 33, 61, 110], [6, 7, 8, 9, 37, 86], [23, 54, 55, 56, 57, 85], [15, 16, 17, 18, 46, 95], [27, 28, 29, 30, 58, 107], [24, 25, 26, 27, 55, 104], [26, 57, 58, 59, 60, 88], [29, 60, 61, 62, 63, 91], [25, 74, 105, 106, 107, 108], [41, 72, 73, 74, 75, 103], [7, 56, 87, 88, 89, 90], [10, 59, 90, 91, 92, 93], [13, 62, 93, 94, 95, 96], [22, 71, 102, 103, 104, 105], [9, 10, 11, 12, 40, 89], [44, 75, 76, 77, 78, 106], [18, 19, 20, 21, 49, 98], [11, 42, 43, 44, 45, 73], [3, 4, 5, 6, 34, 83], [16, 65, 96, 97, 98, 99], [19, 68, 99, 100, 101, 102], [17, 48, 49, 50, 51, 79], [47, 78, 79, 80, 81, 109], [20, 51, 52, 53, 54, 82], [12, 13, 14, 15, 43, 92], [0, 28, 77, 108, 109, 110], [4, 53, 84, 85, 86, 87], [1, 50, 81, 82, 83, 84], [32, 63, 64, 65, 66, 94], [2, 33, 34, 35, 36, 64], [21, 22, 23, 24, 52, 101], [38, 69, 70, 71, 72, 100], [14, 45, 46, 47, 48, 76], [35, 66, 67, 68, 69, 97], [5, 36, 37, 38, 39, 67], [26, 54, 73, 80, 93, 106], [37, 54, 69, 94, 104, 107], [9, 13, 17, 32, 39, 85], [13, 20, 33, 46, 77, 105], [31, 66, 70, 74, 89, 96], [13, 30, 45, 70, 80, 83], [19, 54, 58, 62, 77, 84], [28, 63, 67, 71, 86, 93], [15, 19, 23, 38, 45, 91], [17, 45, 64, 71, 84, 97], [19, 29, 32, 73, 90, 105], [22, 32, 35, 76, 93, 108], [21, 25, 29, 44, 51, 97], [30, 34, 38, 53, 60, 106], [2, 30, 49, 56, 69, 82], [7, 38, 66, 85, 92, 105], [9, 34, 44, 47, 88, 105], [12, 37, 47, 50, 91, 108], [1, 18, 33, 58, 68, 71], [27, 46, 53, 66, 79, 110], [21, 40, 47, 60, 73, 104], [1, 5, 20, 27, 73, 108], [8, 36, 55, 62, 75, 88], [12, 27, 52, 62, 65, 106], [11, 18, 64, 99, 103, 107], [3, 22, 29, 42, 55, 86], [10, 17, 30, 43, 74, 102], [24, 28, 32, 47, 54, 100], [4, 39, 43, 47, 62, 69], [22, 39, 54, 79, 89, 92], [11, 39, 58, 65, 78, 91], [2, 17, 24, 70, 105, 109], [27, 31, 35, 50, 57, 103], [2, 43, 60, 75, 100, 110], [0, 4, 8, 23, 30, 76], [13, 23, 26, 67, 84, 99], [13, 48, 52, 56, 71, 78], [12, 25, 56, 84, 103, 110], [20, 48, 67, 74, 87, 100], [16, 51, 55, 59, 74, 81], [15, 30, 55, 65, 68, 109], [16, 26, 29, 70, 87, 102], [34, 69, 73, 77, 92, 99], [1, 32, 60, 79, 86, 99], [7, 17, 20, 61, 78, 93], [16, 23, 36, 49, 80, 108], [9, 22, 53, 81, 100, 107], [3, 7, 11, 26, 33, 79], [6, 21, 46, 56, 59, 100], [12, 31, 38, 51, 64, 95], [18, 37, 44, 57, 70, 101], [6, 10, 14, 29, 36, 82], [9, 24, 49, 59, 62, 103], [22, 57, 61, 65, 80, 87], [3, 28, 38, 41, 82, 99], [15, 34, 41, 54, 67, 98], [0, 15, 40, 50, 53, 94], [0, 19, 26, 39, 52, 83], [2, 5, 46, 63, 78, 103], [0, 25, 35, 38, 79, 96], [34, 51, 66, 91, 101, 104], [18, 22, 26, 41, 48, 94], [1, 8, 21, 34, 65, 93], [4, 21, 36, 61, 71, 74], [4, 11, 24, 37, 68, 96], [6, 25, 32, 45, 58, 89], [43, 78, 82, 86, 101, 108], [9, 28, 35, 48, 61, 92], [12, 16, 20, 35, 42, 88], [2, 15, 28, 59, 87, 106], [10, 20, 23, 64, 81, 96], [16, 33, 48, 73, 83, 86], [24, 43, 50, 63, 76, 107], [2, 9, 55, 90, 94, 98], [7, 24, 39, 64, 74, 77], [5, 12, 58, 93, 97, 101], [40, 57, 72, 97, 107, 110], [25, 42, 57, 82, 92, 95], [7, 42, 46, 50, 65, 72], [10, 27, 42, 67, 77, 80], [0, 46, 81, 85, 89, 104], [29, 57, 76, 83, 96, 109], [40, 75, 79, 83, 98, 105], [7, 14, 27, 40, 71, 99], [5, 33, 52, 59, 72, 85], [25, 60, 64, 68, 83, 90], [14, 42, 61, 68, 81, 94], [5, 18, 31, 62, 90, 109], [8, 15, 61, 96, 100, 104], [33, 37, 41, 56, 63, 109], [8, 11, 52, 69, 84, 109], [37, 72, 76, 80, 95, 102], [10, 41, 69, 88, 95, 108], [0, 13, 44, 72, 91, 98], [14, 21, 67, 102, 106, 110], [6, 31, 41, 44, 85, 102], [5, 8, 49, 66, 81, 106], [28, 45, 60, 85, 95, 98], [3, 18, 43, 53, 56, 97], [4, 35, 63, 82, 89, 102], [1, 36, 40, 44, 59, 66], [19, 36, 51, 76, 86, 89], [6, 52, 87, 91, 95, 110], [10, 45, 49, 53, 68, 75], [4, 14, 17, 58, 75, 90], [1, 11, 14, 55, 72, 87], [3, 16, 47, 75, 94, 101], [6, 19, 50, 78, 97, 104], [3, 49, 84, 88, 92, 107], [23, 51, 70, 77, 90, 103], [31, 48, 63, 88, 98, 101], [36, 42, 52, 60, 102, 107], [8, 28, 73, 79, 97, 102], [1, 19, 24, 41, 61, 106], [0, 14, 56, 73, 101, 107], [38, 55, 83, 89, 93, 107], [2, 6, 20, 62, 79, 107], [4, 49, 55, 73, 78, 95], [7, 12, 29, 49, 94, 100], [30, 36, 46, 54, 96, 101], [17, 23, 27, 41, 83, 100], [0, 5, 45, 51, 61, 69], [21, 26, 66, 72, 82, 90], [8, 25, 53, 59, 63, 77], [3, 9, 19, 27, 69, 74], [31, 37, 55, 60, 77, 97], [16, 22, 40, 45, 62, 82], [5, 9, 23, 65, 82, 110], [3, 13, 21, 63, 68, 108], [37, 43, 61, 66, 83, 103], [5, 11, 15, 29, 71, 88], [36, 41, 81, 87, 97, 105], [17, 34, 62, 68, 72, 86], [13, 18, 35, 55, 100, 106], [12, 18, 28, 36, 78, 83], [0, 6, 16, 24, 66, 71], [1, 6, 23, 43, 88, 94], [26, 43, 71, 77, 81, 95], [39, 44, 84, 90, 100, 108], [6, 12, 22, 30, 72, 77], [2, 44, 61, 89, 95, 99], [13, 58, 64, 82, 87, 104], [34, 40, 58, 63, 80, 100], [23, 40, 68, 74, 78, 92], [4, 10, 28, 33, 50, 70], [2, 19, 47, 53, 57, 71], [13, 19, 37, 42, 59, 79], [7, 35, 41, 45, 59, 101], [2, 42, 48, 58, 66, 108], [10, 15, 32, 52, 97, 103], [8, 14, 18, 32, 74, 91], [3, 8, 48, 54, 64, 72], [22, 28, 46, 51, 68, 88], [2, 8, 12, 26, 68, 85], [6, 11, 51, 57, 67, 75], [20, 26, 30, 44, 86, 103], [23, 29, 33, 47, 89, 106], [15, 21, 31, 39, 81, 86], [4, 22, 27, 44, 64, 109], [10, 55, 61, 79, 84, 101], [6, 48, 53, 93, 99, 109], [4, 12, 54, 59, 99, 105], [14, 20, 24, 38, 80, 97], [12, 17, 57, 63, 73, 81], [0, 17, 37, 82, 88, 106], [41, 58, 86, 92, 96, 110], [2, 22, 67, 73, 91, 96], [3, 17, 59, 76, 104, 110], [1, 7, 25, 30, 47, 67], [35, 52, 80, 86, 90, 104], [16, 44, 50, 54, 68, 110], [7, 15, 57, 62, 102, 108], [4, 32, 38, 42, 56, 98], [14, 34, 79, 85, 103, 108], [8, 50, 67, 95, 101, 105], [1, 29, 35, 39, 53, 95], [18, 23, 63, 69, 79, 87], [5, 22, 50, 56, 60, 74], [0, 10, 18, 60, 65, 105], [10, 38, 44, 48, 62, 104], [10, 16, 34, 39, 56, 76], [40, 46, 64, 69, 86, 106], [11, 17, 21, 35, 77, 94], [3, 20, 40, 85, 91, 109], [26, 32, 36, 50, 92, 109], [18, 24, 34, 42, 84, 89], [7, 52, 58, 76, 81, 98], [1, 9, 51, 56, 96, 102], [33, 39, 49, 57, 99, 104], [4, 9, 26, 46, 91, 97], [25, 31, 49, 54, 71, 91], [24, 29, 69, 75, 85, 93], [21, 27, 37, 45, 87, 92], [30, 35, 75, 81, 91, 99], [11, 28, 56, 62, 66, 80], [3, 45, 50, 90, 96, 106], [39, 45, 55, 63, 105, 110], [32, 49, 77, 83, 87, 101], [1, 46, 52, 70, 75, 92], [13, 41, 47, 51, 65, 107], [9, 14, 54, 60, 70, 78], [15, 20, 60, 66, 76, 84], [28, 34, 52, 57, 74, 94], [7, 13, 31, 36, 53, 73], [20, 37, 65, 71, 75, 89], [11, 31, 76, 82, 100, 105], [29, 46, 74, 80, 84, 98], [24, 30, 40, 48, 90, 95], [27, 32, 72, 78, 88, 96], [11, 53, 70, 98, 104, 108], [5, 25, 70, 76, 94, 99], [5, 47, 64, 92, 98, 102], [14, 31, 59, 65, 69, 83], [33, 38, 78, 84, 94, 102], [9, 15, 25, 33, 75, 80], [0, 42, 47, 87, 93, 103], [19, 25, 43, 48, 65, 85], [43, 49, 67, 72, 89, 109], [16, 61, 67, 85, 90, 107], [19, 64, 70, 88, 93, 110], [27, 33, 43, 51, 93, 98], [16, 21, 38, 58, 103, 109], [1, 62, 64, 76, 85, 101], [12, 23, 66, 86, 98, 107], [1, 16, 57, 69, 78, 105], [56, 58, 70, 79, 95, 106], [12, 24, 33, 60, 67, 82], [8, 10, 22, 31, 47, 58], [2, 11, 27, 38, 81, 101], [5, 21, 32, 75, 95, 107], [29, 31, 43, 52, 68, 79], [0, 20, 32, 41, 57, 68], [2, 4, 16, 25, 41, 52], [9, 21, 30, 57, 64, 79], [23, 25, 37, 46, 62, 73], [53, 55, 67, 76, 92, 103], [11, 22, 83, 85, 97, 106], [10, 51, 63, 72, 99, 106], [21, 28, 43, 84, 96, 105], [9, 20, 63, 83, 95, 104], [18, 30, 39, 66, 73, 88], [39, 51, 60, 87, 94, 109], [2, 45, 65, 77, 86, 102], [2, 13, 74, 76, 88, 97], [3, 10, 25, 66, 78, 87], [8, 20, 29, 45, 56, 99], [7, 23, 34, 95, 97, 109], [8, 17, 33, 44, 87, 107], [6, 33, 40, 55, 96, 108], [7, 16, 32, 43, 104, 106], [41, 43, 55, 64, 80, 91], [8, 51, 71, 83, 92, 108], [1, 17, 28, 89, 91, 103], [2, 18, 29, 72, 92, 104], [4, 20, 31, 92, 94, 106], [14, 26, 35, 51, 62, 105], [12, 19, 34, 75, 87, 96], [4, 19, 60, 72, 81, 108], [24, 31, 46, 87, 99, 108], [4, 13, 29, 40, 101, 103], [32, 34, 46, 55, 71, 82], [5, 7, 19, 28, 44, 55], [27, 39, 48, 75, 82, 97], [10, 71, 73, 85, 94, 110], [15, 22, 37, 78, 90, 99], [3, 30, 37, 52, 93, 105], [3, 12, 39, 46, 61, 102], [38, 40, 52, 61, 77, 88], [17, 29, 38, 54, 65, 108], [18, 25, 40, 81, 93, 102], [47, 49, 61, 70, 86, 97], [5, 48, 68, 80, 89, 105], [7, 68, 70, 82, 91, 107], [14, 25, 86, 88, 100, 109], [6, 17, 60, 80, 92, 101], [36, 48, 57, 84, 91, 106], [0, 12, 21, 48, 55, 70], [0, 27, 34, 49, 90, 102], [6, 13, 28, 69, 81, 90], [14, 16, 28, 37, 53, 64], [5, 16, 77, 79, 91, 100], [0, 7, 22, 63, 75, 84], [11, 20, 36, 47, 90, 110], [15, 35, 47, 56, 72, 83], [24, 36, 45, 72, 79, 94], [5, 14, 30, 41, 84, 104], [9, 18, 45, 52, 67, 108], [10, 19, 35, 46, 107, 109], [9, 29, 41, 50, 66, 77], [30, 50, 62, 71, 87, 98], [3, 14, 57, 77, 89, 98], [33, 45, 54, 81, 88, 103], [4, 65, 67, 79, 88, 104], [8, 19, 80, 82, 94, 103], [5, 17, 26, 42, 53, 96], [8, 24, 35, 78, 98, 110], [1, 13, 22, 38, 49, 110], [6, 18, 27, 54, 61, 76], [11, 13, 25, 34, 50, 61], [2, 14, 23, 39, 50, 93], [3, 15, 24, 51, 58, 73], [50, 52, 64, 73, 89, 100], [11, 23, 32, 48, 59, 102], [3, 23, 35, 44, 60, 71], [24, 44, 56, 65, 81, 92], [18, 38, 50, 59, 75, 86], [21, 41, 53, 62, 78, 89], [1, 42, 54, 63, 90, 97], [33, 53, 65, 74, 90, 101], [13, 54, 66, 75, 102, 109], [36, 56, 68, 77, 93, 104], [59, 61, 73, 82, 98, 109], [42, 62, 74, 83, 99, 110], [39, 59, 71, 80, 96, 107], [4, 45, 57, 66, 93, 100], [35, 37, 49, 58, 74, 85], [30, 42, 51, 78, 85, 100], [6, 15, 42, 49, 64, 105], [20, 22, 34, 43, 59, 70], [15, 26, 69, 89, 101, 110], [6, 26, 38, 47, 63, 74], [27, 47, 59, 68, 84, 95], [44, 46, 58, 67, 83, 94], [0, 9, 36, 43, 58, 99], [15, 27, 36, 63, 70, 85], [0, 11, 54, 74, 86, 95], [12, 32, 44, 53, 69, 80], [9, 16, 31, 72, 84, 93], [1, 10, 26, 37, 98, 100], [26, 28, 40, 49, 65, 76], [7, 48, 60, 69, 96, 103], [17, 19, 31, 40, 56, 67], [21, 33, 42, 69, 76, 91], [27, 56, 86, 91, 94, 105], [11, 41, 46, 49, 60, 93], [8, 38, 43, 46, 57, 90], [9, 38, 68, 73, 76, 87], [20, 25, 28, 39, 72, 101], [2, 7, 10, 21, 54, 83], [6, 39, 68, 98, 103, 106], [23, 53, 58, 61, 72, 105], [3, 32, 62, 67, 70, 81], [6, 35, 65, 70, 73, 84], [17, 22, 25, 36, 69, 98], [0, 33, 62, 92, 97, 100], [2, 32, 37, 40, 51, 84], [21, 50, 80, 85, 88, 99], [26, 31, 34, 45, 78, 107], [3, 36, 65, 95, 100, 103], [14, 19, 22, 33, 66, 95], [29, 34, 37, 48, 81, 110], [1, 4, 15, 48, 77, 107], [23, 28, 31, 42, 75, 104], [26, 56, 61, 64, 75, 108], [12, 41, 71, 76, 79, 90], [30, 59, 89, 94, 97, 108], [4, 7, 18, 51, 80, 110], [20, 50, 55, 58, 69, 102], [15, 44, 74, 79, 82, 93], [11, 16, 19, 30, 63, 92], [5, 35, 40, 43, 54, 87], [9, 42, 71, 101, 106, 109], [8, 13, 16, 27, 60, 89], [5, 10, 13, 24, 57, 86], [0, 29, 59, 64, 67, 78], [24, 53, 83, 88, 91, 102], [18, 47, 77, 82, 85, 96], [17, 47, 52, 55, 66, 99], [14, 44, 49, 52, 63, 96], [1, 12, 45, 74, 104, 109]]
\item 1 \{0=1110, 1=24864, 2=544122, 3=2208900, 4=2349204\} [[8, 39, 40, 41, 42, 70], [0, 1, 2, 3, 31, 80], [30, 31, 32, 33, 61, 110], [6, 7, 8, 9, 37, 86], [23, 54, 55, 56, 57, 85], [15, 16, 17, 18, 46, 95], [27, 28, 29, 30, 58, 107], [24, 25, 26, 27, 55, 104], [26, 57, 58, 59, 60, 88], [29, 60, 61, 62, 63, 91], [25, 74, 105, 106, 107, 108], [41, 72, 73, 74, 75, 103], [7, 56, 87, 88, 89, 90], [10, 59, 90, 91, 92, 93], [13, 62, 93, 94, 95, 96], [22, 71, 102, 103, 104, 105], [9, 10, 11, 12, 40, 89], [44, 75, 76, 77, 78, 106], [18, 19, 20, 21, 49, 98], [11, 42, 43, 44, 45, 73], [3, 4, 5, 6, 34, 83], [16, 65, 96, 97, 98, 99], [19, 68, 99, 100, 101, 102], [17, 48, 49, 50, 51, 79], [47, 78, 79, 80, 81, 109], [20, 51, 52, 53, 54, 82], [12, 13, 14, 15, 43, 92], [0, 28, 77, 108, 109, 110], [4, 53, 84, 85, 86, 87], [1, 50, 81, 82, 83, 84], [32, 63, 64, 65, 66, 94], [2, 33, 34, 35, 36, 64], [21, 22, 23, 24, 52, 101], [38, 69, 70, 71, 72, 100], [14, 45, 46, 47, 48, 76], [35, 66, 67, 68, 69, 97], [5, 36, 37, 38, 39, 67], [35, 51, 55, 61, 86, 96], [29, 45, 49, 55, 80, 90], [44, 60, 64, 70, 95, 105], [6, 56, 72, 76, 82, 107], [41, 57, 61, 67, 92, 102], [23, 39, 43, 49, 74, 84], [23, 33, 83, 99, 103, 109], [5, 15, 65, 81, 85, 91], [1, 7, 32, 42, 92, 108], [11, 27, 31, 37, 62, 72], [2, 18, 22, 28, 53, 63], [11, 21, 71, 87, 91, 97], [20, 30, 80, 96, 100, 106], [12, 16, 22, 47, 57, 107], [4, 29, 39, 89, 105, 109], [6, 10, 16, 41, 51, 101], [14, 30, 34, 40, 65, 75], [0, 50, 66, 70, 76, 101], [26, 42, 46, 52, 77, 87], [3, 53, 69, 73, 79, 104], [2, 12, 62, 78, 82, 88], [15, 19, 25, 50, 60, 110], [47, 63, 67, 73, 98, 108], [5, 21, 25, 31, 56, 66], [20, 36, 40, 46, 71, 81], [38, 54, 58, 64, 89, 99], [0, 4, 10, 35, 45, 95], [17, 33, 37, 43, 68, 78], [14, 24, 74, 90, 94, 100], [8, 18, 68, 84, 88, 94], [8, 24, 28, 34, 59, 69], [3, 7, 13, 38, 48, 98], [9, 59, 75, 79, 85, 110], [1, 26, 36, 86, 102, 106], [32, 48, 52, 58, 83, 93], [17, 27, 77, 93, 97, 103], [9, 13, 19, 44, 54, 104], [14, 22, 44, 51, 62, 89], [7, 10, 18, 23, 34, 70], [32, 40, 62, 69, 80, 107], [22, 70, 73, 81, 86, 97], [21, 43, 63, 69, 77, 96], [26, 62, 70, 92, 99, 110], [19, 67, 70, 78, 83, 94], [12, 34, 54, 60, 68, 87], [2, 13, 49, 97, 100, 108], [18, 24, 32, 51, 87, 109], [5, 41, 49, 71, 78, 89], [1, 37, 85, 88, 96, 101], [12, 48, 70, 90, 96, 104], [23, 31, 53, 60, 71, 98], [8, 15, 26, 53, 89, 97], [31, 34, 42, 47, 58, 94], [15, 37, 57, 63, 71, 90], [1, 49, 52, 60, 65, 76], [8, 44, 52, 74, 81, 92], [2, 9, 20, 47, 83, 91], [7, 29, 36, 47, 74, 110], [18, 54, 76, 96, 102, 110], [0, 22, 42, 48, 56, 75], [6, 11, 22, 58, 106, 109], [3, 8, 19, 55, 103, 106], [9, 15, 23, 42, 78, 100], [25, 28, 36, 41, 52, 88], [7, 55, 58, 66, 71, 82], [10, 13, 21, 26, 37, 73], [3, 25, 45, 51, 59, 78], [2, 38, 46, 68, 75, 86], [16, 36, 42, 50, 69, 105], [23, 59, 67, 89, 96, 107], [37, 40, 48, 53, 64, 100], [20, 27, 38, 65, 101, 109], [1, 21, 27, 35, 54, 90], [8, 35, 71, 79, 101, 108], [26, 34, 56, 63, 74, 101], [0, 5, 16, 52, 100, 103], [34, 37, 45, 50, 61, 97], [6, 28, 48, 54, 62, 81], [11, 47, 55, 77, 84, 95], [34, 82, 85, 93, 98, 109], [27, 49, 69, 75, 83, 102], [9, 45, 67, 87, 93, 101], [4, 26, 33, 44, 71, 107], [4, 24, 30, 38, 57, 93], [19, 22, 30, 35, 46, 82], [5, 13, 35, 42, 53, 80], [33, 55, 75, 81, 89, 108], [2, 29, 65, 73, 95, 102], [7, 27, 33, 41, 60, 96], [7, 43, 91, 94, 102, 107], [17, 24, 35, 62, 98, 106], [46, 49, 57, 62, 73, 109], [14, 21, 32, 59, 95, 103], [1, 23, 30, 41, 68, 104], [3, 14, 41, 77, 85, 107], [6, 42, 64, 84, 90, 98], [3, 9, 17, 36, 72, 94], [11, 19, 41, 48, 59, 86], [5, 12, 23, 50, 86, 94], [19, 39, 45, 53, 72, 108], [2, 10, 32, 39, 50, 77], [22, 25, 33, 38, 49, 85], [2, 21, 57, 79, 99, 105], [29, 37, 59, 66, 77, 104], [0, 11, 38, 74, 82, 104], [8, 16, 38, 45, 56, 83], [11, 18, 29, 56, 92, 100], [40, 43, 51, 56, 67, 103], [4, 52, 55, 63, 68, 79], [35, 43, 65, 72, 83, 110], [18, 40, 60, 66, 74, 93], [4, 40, 88, 91, 99, 104], [10, 30, 36, 44, 63, 99], [3, 39, 61, 81, 87, 95], [3, 11, 30, 66, 88, 108], [20, 28, 50, 57, 68, 95], [0, 8, 27, 63, 85, 105], [6, 17, 44, 80, 88, 110], [0, 36, 58, 78, 84, 92], [14, 50, 58, 80, 87, 98], [10, 58, 61, 69, 74, 85], [25, 73, 76, 84, 89, 100], [28, 76, 79, 87, 92, 103], [12, 18, 26, 45, 81, 103], [28, 31, 39, 44, 55, 91], [24, 46, 66, 72, 80, 99], [31, 79, 82, 90, 95, 106], [4, 7, 15, 20, 31, 67], [13, 61, 64, 72, 77, 88], [17, 53, 61, 83, 90, 101], [43, 46, 54, 59, 70, 106], [9, 31, 51, 57, 65, 84], [10, 46, 94, 97, 105, 110], [17, 25, 47, 54, 65, 92], [15, 51, 73, 93, 99, 107], [1, 4, 12, 17, 28, 64], [13, 33, 39, 47, 66, 102], [5, 32, 68, 76, 98, 105], [15, 21, 29, 48, 84, 106], [13, 16, 24, 29, 40, 76], [5, 24, 60, 82, 102, 108], [16, 64, 67, 75, 80, 91], [6, 12, 20, 39, 75, 97], [30, 52, 72, 78, 86, 105], [20, 56, 64, 86, 93, 104], [1, 9, 14, 25, 61, 109], [16, 19, 27, 32, 43, 79], [0, 6, 14, 33, 69, 91], [6, 29, 35, 78, 87, 99], [14, 20, 63, 72, 84, 102], [33, 46, 79, 88, 93, 100], [32, 41, 46, 53, 56, 91], [10, 62, 71, 76, 83, 86], [19, 28, 33, 40, 84, 97], [1, 8, 11, 46, 98, 107], [1, 45, 58, 91, 100, 105], [22, 74, 83, 88, 95, 98], [34, 86, 95, 100, 107, 110], [5, 10, 17, 20, 55, 107], [6, 18, 36, 59, 65, 108], [3, 10, 54, 67, 100, 109], [20, 26, 69, 78, 90, 108], [15, 28, 61, 70, 75, 82], [26, 35, 40, 47, 50, 85], [15, 38, 44, 87, 96, 108], [8, 17, 22, 29, 32, 67], [1, 34, 43, 48, 55, 99], [12, 21, 33, 51, 74, 80], [21, 30, 42, 60, 83, 89], [4, 37, 46, 51, 58, 102], [4, 48, 61, 94, 103, 108], [41, 50, 55, 62, 65, 100], [10, 43, 52, 57, 64, 108], [31, 40, 45, 52, 96, 109], [10, 19, 24, 31, 75, 88], [3, 12, 24, 42, 65, 71], [25, 34, 39, 46, 90, 103], [42, 51, 63, 81, 104, 110], [2, 5, 40, 92, 101, 106], [50, 59, 64, 71, 74, 109], [27, 40, 73, 82, 87, 94], [4, 13, 18, 25, 69, 82], [23, 32, 37, 44, 47, 82], [7, 59, 68, 73, 80, 83], [19, 71, 80, 85, 92, 95], [39, 52, 85, 94, 99, 106], [5, 14, 19, 26, 29, 64], [5, 48, 57, 69, 87, 110], [7, 16, 21, 28, 72, 85], [18, 31, 64, 73, 78, 85], [47, 56, 61, 68, 71, 106], [5, 11, 54, 63, 75, 93], [0, 23, 29, 72, 81, 93], [0, 18, 41, 47, 90, 99], [0, 13, 46, 55, 60, 67], [15, 24, 36, 54, 77, 83], [6, 24, 47, 53, 96, 105], [6, 19, 52, 61, 66, 73], [36, 49, 82, 91, 96, 103], [1, 6, 13, 57, 70, 103], [38, 47, 52, 59, 62, 97], [39, 48, 60, 78, 101, 107], [16, 25, 30, 37, 81, 94], [22, 31, 36, 43, 87, 100], [6, 15, 27, 45, 68, 74], [42, 55, 88, 97, 102, 109], [35, 44, 49, 56, 59, 94], [13, 65, 74, 79, 86, 89], [2, 7, 14, 17, 52, 104], [0, 9, 21, 39, 62, 68], [2, 45, 54, 66, 84, 107], [4, 11, 14, 49, 101, 110], [0, 7, 51, 64, 97, 106], [30, 39, 51, 69, 92, 98], [28, 37, 42, 49, 93, 106], [5, 8, 43, 95, 104, 109], [2, 37, 89, 98, 103, 110], [2, 11, 16, 23, 26, 61], [9, 22, 55, 64, 69, 76], [3, 21, 44, 50, 93, 102], [3, 16, 49, 58, 63, 70], [25, 77, 86, 91, 98, 101], [21, 34, 67, 76, 81, 88], [11, 17, 60, 69, 81, 99], [27, 36, 48, 66, 89, 95], [12, 35, 41, 84, 93, 105], [4, 9, 16, 60, 73, 106], [44, 53, 58, 65, 68, 103], [9, 27, 50, 56, 99, 108], [8, 13, 20, 23, 58, 110], [36, 45, 57, 75, 98, 104], [7, 12, 19, 63, 76, 109], [24, 33, 45, 63, 86, 92], [7, 40, 49, 54, 61, 105], [17, 26, 31, 38, 41, 76], [0, 12, 30, 53, 59, 102], [2, 8, 51, 60, 72, 90], [3, 26, 32, 75, 84, 96], [31, 83, 92, 97, 104, 107], [12, 25, 58, 67, 72, 79], [8, 14, 57, 66, 78, 96], [20, 29, 34, 41, 44, 79], [1, 10, 15, 22, 66, 79], [9, 18, 30, 48, 71, 77], [29, 38, 43, 50, 53, 88], [33, 42, 54, 72, 95, 101], [28, 80, 89, 94, 101, 104], [9, 32, 38, 81, 90, 102], [30, 43, 76, 85, 90, 97], [1, 53, 62, 67, 74, 77], [11, 20, 25, 32, 35, 70], [24, 37, 70, 79, 84, 91], [4, 56, 65, 70, 77, 80], [16, 68, 77, 82, 89, 92], [3, 15, 33, 56, 62, 105], [13, 22, 27, 34, 78, 91], [17, 23, 66, 75, 87, 105], [18, 27, 39, 57, 80, 86], [14, 23, 28, 35, 38, 73], [8, 25, 62, 64, 87, 102], [18, 33, 50, 67, 104, 106], [10, 47, 49, 72, 87, 104], [14, 31, 68, 70, 93, 108], [35, 37, 60, 75, 92, 109], [6, 21, 38, 55, 92, 94], [8, 10, 33, 48, 65, 82], [5, 7, 30, 45, 62, 79], [14, 16, 39, 54, 71, 88], [2, 4, 27, 42, 59, 76], [3, 18, 35, 52, 89, 91], [13, 50, 52, 75, 90, 107], [4, 41, 43, 66, 81, 98], [16, 53, 55, 78, 93, 110], [1, 38, 40, 63, 78, 95], [6, 23, 40, 77, 79, 102], [9, 24, 41, 58, 95, 97], [2, 19, 56, 58, 81, 96], [9, 26, 43, 80, 82, 105], [11, 28, 65, 67, 90, 105], [0, 15, 32, 49, 86, 88], [3, 20, 37, 74, 76, 99], [12, 27, 44, 61, 98, 100], [26, 28, 51, 66, 83, 100], [7, 44, 46, 69, 84, 101], [32, 34, 57, 72, 89, 106], [11, 13, 36, 51, 68, 85], [29, 31, 54, 69, 86, 103], [15, 30, 47, 64, 101, 103], [23, 25, 48, 63, 80, 97], [0, 17, 34, 71, 73, 96], [12, 29, 46, 83, 85, 108], [21, 36, 53, 70, 107, 109], [20, 22, 45, 60, 77, 94], [17, 19, 42, 57, 74, 91], [5, 22, 59, 61, 84, 99], [1, 24, 39, 56, 73, 110], [2, 24, 44, 48, 67, 85], [10, 28, 56, 78, 98, 102], [1, 29, 51, 71, 75, 94], [3, 22, 40, 68, 90, 110], [8, 12, 31, 49, 77, 99], [4, 22, 50, 72, 92, 96], [10, 38, 60, 80, 84, 103], [0, 20, 24, 43, 61, 89], [4, 32, 54, 74, 78, 97], [16, 34, 62, 84, 104, 108], [17, 39, 59, 63, 82, 100], [1, 19, 47, 69, 89, 93], [11, 33, 53, 57, 76, 94], [17, 21, 40, 58, 86, 108], [16, 44, 66, 86, 90, 109], [7, 25, 53, 75, 95, 99], [5, 9, 28, 46, 74, 96], [14, 18, 37, 55, 83, 105], [21, 41, 45, 64, 82, 110], [9, 29, 33, 52, 70, 98], [15, 35, 39, 58, 76, 104], [6, 26, 30, 49, 67, 95], [0, 19, 37, 65, 87, 107], [14, 36, 56, 60, 79, 97], [3, 23, 27, 46, 64, 92], [7, 35, 57, 77, 81, 100], [13, 41, 63, 83, 87, 106], [26, 48, 68, 72, 91, 109], [20, 42, 62, 66, 85, 103], [23, 45, 65, 69, 88, 106], [13, 31, 59, 81, 101, 105], [5, 27, 47, 51, 70, 88], [12, 32, 36, 55, 73, 101], [2, 6, 25, 43, 71, 93], [11, 15, 34, 52, 80, 102], [18, 38, 42, 61, 79, 107], [8, 30, 50, 54, 73, 91], [6, 31, 46, 50, 63, 89], [7, 22, 26, 39, 65, 93], [12, 38, 66, 91, 106, 110], [14, 42, 67, 82, 86, 99], [6, 32, 60, 85, 100, 104], [24, 49, 64, 68, 81, 107], [13, 17, 30, 56, 84, 109], [1, 5, 18, 44, 72, 97], [18, 43, 58, 62, 75, 101], [22, 37, 41, 54, 80, 108], [17, 45, 70, 85, 89, 102], [0, 26, 54, 79, 94, 98], [10, 14, 27, 53, 81, 106], [21, 46, 61, 65, 78, 104], [4, 8, 21, 47, 75, 100], [12, 37, 52, 56, 69, 95], [8, 36, 61, 76, 80, 93], [4, 19, 23, 36, 62, 90], [23, 51, 76, 91, 95, 108], [9, 35, 63, 88, 103, 107], [7, 11, 24, 50, 78, 103], [11, 39, 64, 79, 83, 96], [27, 52, 67, 71, 84, 110], [15, 40, 55, 59, 72, 98], [1, 16, 20, 33, 59, 87], [20, 48, 73, 88, 92, 105], [2, 15, 41, 69, 94, 109], [10, 25, 29, 42, 68, 96], [9, 34, 49, 53, 66, 92], [5, 33, 58, 73, 77, 90], [3, 29, 57, 82, 97, 101], [13, 28, 32, 45, 71, 99], [3, 28, 43, 47, 60, 86], [0, 25, 40, 44, 57, 83], [16, 31, 35, 48, 74, 102], [19, 34, 38, 51, 77, 105], [2, 30, 55, 70, 74, 87]]
\item 1 \{0=1110, 1=26640, 2=534132, 3=2198688, 4=2367630\} [[8, 39, 40, 41, 42, 70], [0, 1, 2, 3, 31, 80], [30, 31, 32, 33, 61, 110], [6, 7, 8, 9, 37, 86], [23, 54, 55, 56, 57, 85], [15, 16, 17, 18, 46, 95], [27, 28, 29, 30, 58, 107], [24, 25, 26, 27, 55, 104], [26, 57, 58, 59, 60, 88], [29, 60, 61, 62, 63, 91], [25, 74, 105, 106, 107, 108], [41, 72, 73, 74, 75, 103], [7, 56, 87, 88, 89, 90], [10, 59, 90, 91, 92, 93], [13, 62, 93, 94, 95, 96], [22, 71, 102, 103, 104, 105], [9, 10, 11, 12, 40, 89], [44, 75, 76, 77, 78, 106], [18, 19, 20, 21, 49, 98], [11, 42, 43, 44, 45, 73], [3, 4, 5, 6, 34, 83], [16, 65, 96, 97, 98, 99], [19, 68, 99, 100, 101, 102], [17, 48, 49, 50, 51, 79], [47, 78, 79, 80, 81, 109], [20, 51, 52, 53, 54, 82], [12, 13, 14, 15, 43, 92], [0, 28, 77, 108, 109, 110], [4, 53, 84, 85, 86, 87], [1, 50, 81, 82, 83, 84], [32, 63, 64, 65, 66, 94], [2, 33, 34, 35, 36, 64], [21, 22, 23, 24, 52, 101], [38, 69, 70, 71, 72, 100], [14, 45, 46, 47, 48, 76], [35, 66, 67, 68, 69, 97], [5, 36, 37, 38, 39, 67], [20, 32, 79, 96, 100, 105], [37, 54, 58, 63, 89, 101], [11, 23, 70, 87, 91, 96], [3, 7, 12, 38, 50, 97], [1, 18, 22, 27, 53, 65], [2, 49, 66, 70, 75, 101], [4, 21, 25, 30, 56, 68], [2, 14, 61, 78, 82, 87], [10, 27, 31, 36, 62, 74], [5, 52, 69, 73, 78, 104], [19, 36, 40, 45, 71, 83], [3, 29, 41, 88, 105, 109], [13, 30, 34, 39, 65, 77], [17, 29, 76, 93, 97, 102], [0, 26, 38, 85, 102, 106], [25, 42, 46, 51, 77, 89], [5, 17, 64, 81, 85, 90], [15, 19, 24, 50, 62, 109], [9, 13, 18, 44, 56, 103], [28, 45, 49, 54, 80, 92], [8, 55, 72, 76, 81, 107], [23, 35, 82, 99, 103, 108], [40, 57, 61, 66, 92, 104], [46, 63, 67, 72, 98, 110], [0, 4, 9, 35, 47, 94], [8, 20, 67, 84, 88, 93], [14, 26, 73, 90, 94, 99], [16, 33, 37, 42, 68, 80], [11, 58, 75, 79, 84, 110], [12, 16, 21, 47, 59, 106], [31, 48, 52, 57, 83, 95], [1, 6, 32, 44, 91, 108], [6, 10, 15, 41, 53, 100], [22, 39, 43, 48, 74, 86], [43, 60, 64, 69, 95, 107], [34, 51, 55, 60, 86, 98], [7, 24, 28, 33, 59, 71], [15, 51, 56, 70, 97, 104], [24, 29, 43, 70, 77, 99], [4, 11, 33, 69, 74, 88], [33, 38, 52, 79, 86, 108], [6, 42, 47, 61, 88, 95], [4, 31, 38, 60, 96, 101], [3, 39, 44, 58, 85, 92], [0, 5, 19, 46, 53, 75], [25, 32, 54, 90, 95, 109], [2, 16, 43, 50, 72, 108], [21, 57, 62, 76, 103, 110], [1, 8, 30, 66, 71, 85], [27, 32, 46, 73, 80, 102], [21, 26, 40, 67, 74, 96], [15, 20, 34, 61, 68, 90], [5, 27, 63, 68, 82, 109], [16, 23, 45, 81, 86, 100], [7, 14, 36, 72, 77, 91], [2, 24, 60, 65, 79, 106], [13, 20, 42, 78, 83, 97], [7, 34, 41, 63, 99, 104], [9, 45, 50, 64, 91, 98], [3, 8, 22, 49, 56, 78], [0, 36, 41, 55, 82, 89], [9, 14, 28, 55, 62, 84], [18, 23, 37, 64, 71, 93], [19, 26, 48, 84, 89, 103], [10, 37, 44, 66, 102, 107], [12, 48, 53, 67, 94, 101], [18, 54, 59, 73, 100, 107], [6, 11, 25, 52, 59, 81], [13, 40, 47, 69, 105, 110], [22, 29, 51, 87, 92, 106], [1, 28, 35, 57, 93, 98], [12, 17, 31, 58, 65, 87], [10, 17, 39, 75, 80, 94], [30, 35, 49, 76, 83, 105], [16, 27, 34, 79, 85, 88], [23, 38, 44, 80, 88, 104], [4, 10, 13, 52, 63, 70], [24, 32, 36, 51, 69, 75], [37, 43, 46, 85, 96, 103], [4, 49, 55, 58, 97, 108], [54, 62, 66, 81, 99, 105], [36, 44, 48, 63, 81, 87], [9, 16, 61, 67, 70, 109], [0, 8, 12, 27, 45, 51], [9, 27, 33, 93, 101, 105], [22, 33, 40, 85, 91, 94], [26, 34, 50, 80, 95, 101], [16, 22, 25, 64, 75, 82], [45, 53, 57, 72, 90, 96], [10, 21, 28, 73, 79, 82], [25, 31, 34, 73, 84, 91], [1, 4, 43, 54, 61, 106], [43, 49, 52, 91, 102, 109], [37, 48, 55, 100, 106, 109], [15, 21, 81, 89, 93, 108], [20, 28, 44, 74, 89, 95], [7, 23, 53, 68, 74, 110], [14, 29, 35, 71, 79, 95], [57, 65, 69, 84, 102, 108], [4, 7, 46, 57, 64, 109], [3, 63, 71, 75, 90, 108], [14, 20, 56, 64, 80, 110], [17, 25, 41, 71, 86, 92], [32, 40, 56, 86, 101, 107], [1, 17, 47, 62, 68, 104], [8, 14, 50, 58, 74, 104], [8, 23, 29, 65, 73, 89], [34, 40, 43, 82, 93, 100], [23, 31, 47, 77, 92, 98], [15, 23, 27, 42, 60, 66], [34, 45, 52, 97, 103, 106], [29, 37, 53, 83, 98, 104], [39, 47, 51, 66, 84, 90], [28, 39, 46, 91, 97, 100], [3, 11, 15, 30, 48, 54], [42, 50, 54, 69, 87, 93], [19, 30, 37, 82, 88, 91], [14, 22, 38, 68, 83, 89], [4, 20, 50, 65, 71, 107], [0, 18, 24, 84, 92, 96], [9, 17, 21, 36, 54, 60], [31, 42, 49, 94, 100, 103], [0, 60, 68, 72, 87, 105], [4, 15, 22, 67, 73, 76], [13, 19, 22, 61, 72, 79], [12, 18, 78, 86, 90, 105], [29, 44, 50, 86, 94, 110], [0, 7, 52, 58, 61, 100], [11, 19, 35, 65, 80, 86], [30, 38, 42, 57, 75, 81], [2, 17, 23, 59, 67, 83], [1, 40, 51, 58, 103, 109], [1, 46, 52, 55, 94, 105], [6, 13, 58, 64, 67, 106], [40, 46, 49, 88, 99, 106], [3, 18, 36, 42, 102, 110], [9, 15, 75, 83, 87, 102], [6, 12, 72, 80, 84, 99], [28, 34, 37, 76, 87, 94], [5, 41, 49, 65, 95, 110], [2, 6, 21, 39, 45, 105], [18, 26, 30, 45, 63, 69], [3, 21, 27, 87, 95, 99], [31, 37, 40, 79, 90, 97], [2, 8, 44, 52, 68, 98], [0, 6, 66, 74, 78, 93], [5, 13, 29, 59, 74, 80], [0, 15, 33, 39, 99, 107], [12, 30, 36, 96, 104, 108], [2, 32, 47, 53, 89, 97], [12, 20, 24, 39, 57, 63], [25, 36, 43, 88, 94, 97], [35, 43, 59, 89, 104, 110], [1, 7, 10, 49, 60, 67], [1, 12, 19, 64, 70, 73], [3, 10, 55, 61, 64, 103], [2, 10, 26, 56, 71, 77], [5, 9, 24, 42, 48, 108], [22, 28, 31, 70, 81, 88], [11, 41, 56, 62, 98, 106], [6, 14, 18, 33, 51, 57], [6, 24, 30, 90, 98, 102], [8, 16, 32, 62, 77, 83], [20, 35, 41, 77, 85, 101], [19, 25, 28, 67, 78, 85], [5, 11, 47, 55, 71, 101], [5, 20, 26, 62, 70, 86], [5, 35, 50, 56, 92, 100], [48, 56, 60, 75, 93, 99], [26, 41, 47, 83, 91, 107], [3, 9, 69, 77, 81, 96], [11, 17, 53, 61, 77, 107], [8, 38, 53, 59, 95, 103], [13, 24, 31, 76, 82, 85], [17, 32, 38, 74, 82, 98], [7, 13, 16, 55, 66, 73], [14, 44, 59, 65, 101, 109], [33, 41, 45, 60, 78, 84], [11, 26, 32, 68, 76, 92], [10, 16, 19, 58, 69, 76], [27, 35, 39, 54, 72, 78], [2, 38, 46, 62, 92, 107], [51, 59, 63, 78, 96, 102], [21, 29, 33, 48, 66, 72], [7, 18, 25, 70, 76, 79], [0, 10, 25, 29, 57, 101], [10, 14, 42, 86, 96, 106], [32, 42, 52, 67, 71, 99], [3, 47, 57, 67, 82, 86], [14, 24, 34, 49, 53, 81], [38, 48, 58, 73, 77, 105], [17, 27, 37, 52, 56, 84], [9, 53, 63, 73, 88, 92], [1, 5, 33, 77, 87, 97], [11, 21, 31, 46, 50, 78], [0, 44, 54, 64, 79, 83], [2, 30, 74, 84, 94, 109], [5, 15, 25, 40, 44, 72], [4, 19, 23, 51, 95, 105], [20, 30, 40, 55, 59, 87], [26, 36, 46, 61, 65, 93], [41, 51, 61, 76, 80, 108], [35, 45, 55, 70, 74, 102], [27, 71, 81, 91, 106, 110], [7, 22, 26, 54, 98, 108], [12, 56, 66, 76, 91, 95], [7, 11, 39, 83, 93, 103], [8, 18, 28, 43, 47, 75], [1, 16, 20, 48, 92, 102], [4, 8, 36, 80, 90, 100], [24, 68, 78, 88, 103, 107], [21, 65, 75, 85, 100, 104], [3, 13, 28, 32, 60, 104], [6, 50, 60, 70, 85, 89], [18, 62, 72, 82, 97, 101], [6, 16, 31, 35, 63, 107], [29, 39, 49, 64, 68, 96], [15, 59, 69, 79, 94, 98], [13, 17, 45, 89, 99, 109], [9, 19, 34, 38, 66, 110], [2, 12, 22, 37, 41, 69], [23, 33, 43, 58, 62, 90], [22, 63, 77, 84, 95, 100], [6, 23, 26, 75, 97, 109], [8, 10, 24, 54, 97, 110], [11, 20, 22, 36, 66, 109], [8, 11, 60, 82, 94, 102], [0, 14, 21, 32, 37, 70], [5, 14, 16, 30, 60, 103], [27, 70, 83, 92, 94, 108], [3, 33, 76, 89, 98, 100], [1, 15, 45, 88, 101, 110], [33, 47, 54, 65, 70, 103], [1, 9, 26, 29, 78, 100], [6, 36, 79, 92, 101, 103], [1, 14, 23, 25, 39, 69], [3, 46, 59, 68, 70, 84], [30, 44, 51, 62, 67, 100], [39, 53, 60, 71, 76, 109], [7, 48, 62, 69, 80, 85], [40, 53, 62, 64, 78, 108], [48, 70, 82, 90, 107, 110], [21, 64, 77, 86, 88, 102], [2, 51, 73, 85, 93, 110], [6, 20, 27, 38, 43, 76], [30, 52, 64, 72, 89, 92], [9, 52, 65, 74, 76, 90], [11, 18, 29, 34, 67, 108], [19, 60, 74, 81, 92, 97], [9, 39, 82, 95, 104, 106], [7, 15, 32, 35, 84, 106], [2, 11, 13, 27, 57, 100], [9, 23, 30, 41, 46, 79], [7, 19, 27, 44, 47, 96], [45, 67, 79, 87, 104, 107], [16, 28, 36, 53, 56, 105], [12, 55, 68, 77, 79, 93], [4, 12, 29, 32, 81, 103], [3, 14, 19, 52, 93, 107], [15, 58, 71, 80, 82, 96], [1, 13, 21, 38, 41, 90], [16, 57, 71, 78, 89, 94], [5, 12, 23, 28, 61, 102], [13, 25, 33, 50, 53, 102], [13, 26, 35, 37, 51, 81], [39, 61, 73, 81, 98, 101], [8, 13, 46, 87, 101, 108], [18, 61, 74, 83, 85, 99], [16, 29, 38, 40, 54, 84], [22, 35, 44, 46, 60, 90], [34, 47, 56, 58, 72, 102], [12, 34, 46, 54, 71, 74], [1, 34, 75, 89, 96, 107], [31, 72, 86, 93, 104, 109], [10, 22, 30, 47, 50, 99], [24, 38, 45, 56, 61, 94], [10, 18, 35, 38, 87, 109], [33, 55, 67, 75, 92, 95], [31, 44, 53, 55, 69, 99], [24, 46, 58, 66, 83, 86], [6, 28, 40, 48, 65, 68], [4, 17, 26, 28, 42, 72], [4, 16, 24, 41, 44, 93], [14, 17, 66, 88, 100, 108], [3, 25, 37, 45, 62, 65], [19, 32, 41, 43, 57, 87], [4, 37, 78, 92, 99, 110], [0, 22, 34, 42, 59, 62], [27, 49, 61, 69, 86, 89], [15, 29, 36, 47, 52, 85], [0, 30, 73, 86, 95, 97], [9, 31, 43, 51, 68, 71], [21, 43, 55, 63, 80, 83], [5, 8, 57, 79, 91, 99], [28, 41, 50, 52, 66, 96], [18, 40, 52, 60, 77, 80], [0, 17, 20, 69, 91, 103], [7, 20, 29, 31, 45, 75], [5, 10, 43, 84, 98, 105], [2, 9, 20, 25, 58, 99], [2, 4, 18, 48, 91, 104], [19, 31, 39, 56, 59, 108], [21, 35, 42, 53, 58, 91], [0, 11, 16, 49, 90, 104], [12, 26, 33, 44, 49, 82], [36, 50, 57, 68, 73, 106], [4, 45, 59, 66, 77, 82], [6, 49, 62, 71, 73, 87], [25, 66, 80, 87, 98, 103], [10, 23, 32, 34, 48, 78], [10, 51, 65, 72, 83, 88], [24, 67, 80, 89, 91, 105], [8, 17, 19, 33, 63, 106], [28, 69, 83, 90, 101, 106], [3, 20, 23, 72, 94, 106], [25, 38, 47, 49, 63, 93], [8, 15, 26, 31, 64, 105], [13, 54, 68, 75, 86, 91], [27, 41, 48, 59, 64, 97], [37, 50, 59, 61, 75, 105], [2, 5, 54, 76, 88, 96], [18, 32, 39, 50, 55, 88], [12, 42, 85, 98, 107, 109], [0, 43, 56, 65, 67, 81], [11, 14, 63, 85, 97, 105], [15, 37, 49, 57, 74, 77], [42, 64, 76, 84, 101, 104], [36, 58, 70, 78, 95, 98], [5, 7, 21, 51, 94, 107], [6, 17, 22, 55, 96, 110], [1, 42, 56, 63, 74, 79], [3, 17, 24, 35, 40, 73], [2, 7, 40, 81, 95, 102], [9, 22, 57, 80, 97, 107], [5, 22, 32, 45, 58, 93], [10, 20, 33, 46, 81, 104], [33, 56, 73, 83, 96, 109], [15, 38, 55, 65, 78, 91], [0, 23, 40, 50, 63, 76], [8, 25, 35, 48, 61, 96], [6, 29, 46, 56, 69, 82], [18, 41, 58, 68, 81, 94], [4, 14, 27, 40, 75, 98], [21, 44, 61, 71, 84, 97], [27, 50, 67, 77, 90, 103], [1, 36, 59, 76, 86, 99], [6, 19, 54, 77, 94, 104], [8, 21, 34, 69, 92, 109], [9, 32, 49, 59, 72, 85], [14, 31, 41, 54, 67, 102], [17, 34, 44, 57, 70, 105], [2, 19, 29, 42, 55, 90], [0, 13, 48, 71, 88, 98], [7, 42, 65, 82, 92, 105], [12, 35, 52, 62, 75, 88], [3, 16, 51, 74, 91, 101], [16, 26, 39, 52, 87, 110], [5, 18, 31, 66, 89, 106], [12, 25, 60, 83, 100, 110], [3, 26, 43, 53, 66, 79], [24, 47, 64, 74, 87, 100], [1, 11, 24, 37, 72, 95], [7, 17, 30, 43, 78, 101], [20, 37, 47, 60, 73, 108], [13, 23, 36, 49, 84, 107], [30, 53, 70, 80, 93, 106], [4, 39, 62, 79, 89, 102], [2, 15, 28, 63, 86, 103], [11, 28, 38, 51, 64, 99], [10, 45, 68, 85, 95, 108]]
\item 1 \{1=31968, 2=588078, 3=2177820, 4=2330334\} [[8, 39, 40, 41, 42, 70], [0, 1, 2, 3, 31, 80], [30, 31, 32, 33, 61, 110], [6, 7, 8, 9, 37, 86], [23, 54, 55, 56, 57, 85], [15, 16, 17, 18, 46, 95], [27, 28, 29, 30, 58, 107], [24, 25, 26, 27, 55, 104], [26, 57, 58, 59, 60, 88], [29, 60, 61, 62, 63, 91], [25, 74, 105, 106, 107, 108], [41, 72, 73, 74, 75, 103], [7, 56, 87, 88, 89, 90], [10, 59, 90, 91, 92, 93], [13, 62, 93, 94, 95, 96], [22, 71, 102, 103, 104, 105], [9, 10, 11, 12, 40, 89], [44, 75, 76, 77, 78, 106], [18, 19, 20, 21, 49, 98], [11, 42, 43, 44, 45, 73], [3, 4, 5, 6, 34, 83], [16, 65, 96, 97, 98, 99], [19, 68, 99, 100, 101, 102], [17, 48, 49, 50, 51, 79], [47, 78, 79, 80, 81, 109], [20, 51, 52, 53, 54, 82], [12, 13, 14, 15, 43, 92], [0, 28, 77, 108, 109, 110], [4, 53, 84, 85, 86, 87], [1, 50, 81, 82, 83, 84], [32, 63, 64, 65, 66, 94], [2, 33, 34, 35, 36, 64], [21, 22, 23, 24, 52, 101], [38, 69, 70, 71, 72, 100], [14, 45, 46, 47, 48, 76], [35, 66, 67, 68, 69, 97], [5, 36, 37, 38, 39, 67], [5, 65, 71, 85, 88, 92], [1, 10, 33, 66, 70, 76], [13, 22, 45, 78, 82, 88], [15, 19, 25, 61, 70, 93], [8, 18, 38, 66, 78, 84], [14, 42, 54, 60, 95, 105], [0, 20, 48, 60, 66, 101], [17, 27, 47, 75, 87, 93], [7, 10, 14, 38, 98, 104], [1, 4, 8, 32, 92, 98], [12, 45, 49, 55, 91, 100], [2, 30, 42, 48, 83, 93], [15, 27, 33, 68, 78, 98], [30, 34, 40, 76, 85, 108], [32, 42, 62, 90, 102, 108], [53, 59, 73, 76, 80, 104], [21, 25, 31, 67, 76, 99], [9, 21, 27, 62, 72, 92], [28, 37, 60, 93, 97, 103], [12, 16, 22, 58, 67, 90], [18, 30, 36, 71, 81, 101], [24, 36, 42, 77, 87, 107], [18, 51, 55, 61, 97, 106], [44, 50, 64, 67, 71, 95], [1, 37, 46, 69, 102, 106], [9, 42, 46, 52, 88, 97], [3, 9, 44, 54, 74, 102], [24, 28, 34, 70, 79, 102], [9, 13, 19, 55, 64, 87], [20, 26, 40, 43, 47, 71], [22, 31, 54, 87, 91, 97], [27, 39, 45, 80, 90, 110], [11, 71, 77, 91, 94, 98], [6, 26, 54, 66, 72, 107], [4, 7, 11, 35, 95, 101], [7, 30, 63, 67, 73, 109], [9, 15, 50, 60, 80, 108], [25, 34, 57, 90, 94, 100], [0, 4, 10, 46, 55, 78], [50, 56, 70, 73, 77, 101], [2, 12, 32, 60, 72, 78], [6, 39, 43, 49, 85, 94], [15, 48, 52, 58, 94, 103], [20, 80, 86, 100, 103, 107], [23, 33, 53, 81, 93, 99], [8, 36, 48, 54, 89, 99], [2, 62, 68, 82, 85, 89], [17, 23, 37, 40, 44, 68], [21, 54, 58, 64, 100, 109], [11, 39, 51, 57, 92, 102], [35, 41, 55, 58, 62, 86], [3, 7, 13, 49, 58, 81], [2, 26, 86, 92, 106, 109], [11, 21, 41, 69, 81, 87], [38, 44, 58, 61, 65, 89], [17, 45, 57, 63, 98, 108], [0, 12, 18, 53, 63, 83], [1, 24, 57, 61, 67, 103], [59, 65, 79, 82, 86, 110], [5, 15, 35, 63, 75, 81], [1, 7, 43, 52, 75, 108], [21, 33, 39, 74, 84, 104], [4, 40, 49, 72, 105, 109], [6, 12, 47, 57, 77, 105], [3, 15, 21, 56, 66, 86], [5, 19, 22, 26, 50, 110], [3, 38, 48, 68, 96, 108], [17, 77, 83, 97, 100, 104], [29, 35, 49, 52, 56, 80], [2, 16, 19, 23, 47, 107], [2, 8, 22, 25, 29, 53], [6, 10, 16, 52, 61, 84], [3, 23, 51, 63, 69, 104], [3, 36, 40, 46, 82, 91], [10, 13, 17, 41, 101, 107], [27, 31, 37, 73, 82, 105], [20, 30, 50, 78, 90, 96], [8, 14, 28, 31, 35, 59], [5, 11, 25, 28, 32, 56], [11, 17, 31, 34, 38, 62], [12, 24, 30, 65, 75, 95], [23, 29, 43, 46, 50, 74], [26, 32, 46, 49, 53, 77], [34, 43, 66, 99, 103, 109], [29, 39, 59, 87, 99, 105], [56, 62, 76, 79, 83, 107], [23, 83, 89, 103, 106, 110], [0, 6, 41, 51, 71, 99], [19, 28, 51, 84, 88, 94], [8, 68, 74, 88, 91, 95], [31, 40, 63, 96, 100, 106], [13, 16, 20, 44, 104, 110], [14, 24, 44, 72, 84, 90], [32, 38, 52, 55, 59, 83], [0, 35, 45, 65, 93, 105], [10, 19, 42, 75, 79, 85], [5, 33, 45, 51, 86, 96], [6, 18, 24, 59, 69, 89], [41, 47, 61, 64, 68, 92], [1, 5, 29, 89, 95, 109], [9, 29, 57, 69, 75, 110], [26, 36, 56, 84, 96, 102], [47, 53, 67, 70, 74, 98], [7, 16, 39, 72, 76, 82], [16, 25, 48, 81, 85, 91], [4, 13, 36, 69, 73, 79], [14, 20, 34, 37, 41, 65], [14, 74, 80, 94, 97, 101], [4, 27, 60, 64, 70, 106], [18, 22, 28, 64, 73, 96], [0, 33, 37, 43, 79, 88], [12, 61, 73, 98, 102, 107], [10, 22, 47, 51, 56, 72], [9, 58, 70, 95, 99, 104], [11, 15, 20, 36, 85, 97], [7, 19, 44, 48, 53, 69], [3, 8, 24, 73, 85, 110], [22, 34, 59, 63, 68, 84], [0, 49, 61, 86, 90, 95], [25, 37, 62, 66, 71, 87], [5, 9, 14, 30, 79, 91], [4, 16, 41, 45, 50, 66], [6, 55, 67, 92, 96, 101], [1, 13, 38, 42, 47, 63], [20, 24, 29, 45, 94, 106], [40, 52, 77, 81, 86, 102], [34, 46, 71, 75, 80, 96], [2, 18, 67, 79, 104, 108], [19, 31, 56, 60, 65, 81], [10, 35, 39, 44, 60, 109], [8, 12, 17, 33, 82, 94], [16, 28, 53, 57, 62, 78], [46, 58, 83, 87, 92, 108], [3, 52, 64, 89, 93, 98], [7, 32, 36, 41, 57, 106], [31, 43, 68, 72, 77, 93], [43, 55, 80, 84, 89, 105], [14, 18, 23, 39, 88, 100], [17, 21, 26, 42, 91, 103], [28, 40, 65, 69, 74, 90], [0, 5, 21, 70, 82, 107], [15, 64, 76, 101, 105, 110], [4, 29, 33, 38, 54, 103], [37, 49, 74, 78, 83, 99], [2, 6, 11, 27, 76, 88], [23, 27, 32, 48, 97, 109], [13, 25, 50, 54, 59, 75], [1, 26, 30, 35, 51, 100], [11, 22, 46, 61, 79, 99], [39, 54, 61, 71, 78, 108], [21, 63, 78, 85, 95, 102], [11, 54, 65, 67, 80, 83], [12, 23, 25, 38, 41, 80], [4, 22, 42, 65, 76, 100], [13, 33, 56, 67, 91, 106], [18, 29, 31, 44, 47, 86], [33, 44, 46, 59, 62, 101], [30, 45, 52, 62, 69, 99], [1, 25, 40, 58, 78, 101], [6, 21, 28, 38, 45, 75], [11, 14, 53, 96, 107, 109], [2, 41, 84, 95, 97, 110], [3, 10, 20, 27, 57, 99], [20, 63, 74, 76, 89, 92], [42, 53, 55, 68, 71, 110], [21, 32, 34, 47, 50, 89], [36, 51, 58, 68, 75, 105], [8, 51, 62, 64, 77, 80], [8, 19, 43, 58, 76, 96], [16, 31, 49, 69, 92, 103], [2, 9, 39, 81, 96, 103], [5, 16, 40, 55, 73, 93], [17, 60, 71, 73, 86, 89], [39, 50, 52, 65, 68, 107], [32, 75, 86, 88, 101, 104], [1, 14, 17, 56, 99, 110], [6, 13, 23, 30, 60, 102], [2, 4, 17, 20, 59, 102], [8, 10, 23, 26, 65, 108], [23, 66, 77, 79, 92, 95], [24, 39, 46, 56, 63, 93], [18, 60, 75, 82, 92, 99], [14, 25, 49, 64, 82, 102], [5, 8, 47, 90, 101, 103], [12, 19, 29, 36, 66, 108], [12, 35, 46, 70, 85, 103], [5, 7, 20, 23, 62, 105], [7, 31, 46, 64, 84, 107], [14, 57, 68, 70, 83, 86], [0, 30, 72, 87, 94, 104], [6, 48, 63, 70, 80, 87], [27, 38, 40, 53, 56, 95], [10, 28, 48, 71, 82, 106], [1, 11, 18, 48, 90, 105], [4, 19, 37, 57, 80, 91], [4, 14, 21, 51, 93, 108], [15, 30, 37, 47, 54, 84], [20, 31, 55, 70, 88, 108], [1, 19, 39, 62, 73, 97], [3, 14, 16, 29, 32, 71], [18, 33, 40, 50, 57, 87], [6, 36, 78, 93, 100, 110], [10, 25, 43, 63, 86, 97], [38, 81, 92, 94, 107, 110], [6, 29, 40, 64, 79, 97], [17, 28, 52, 67, 85, 105], [7, 22, 40, 60, 83, 94], [8, 15, 45, 87, 102, 109], [3, 26, 37, 61, 76, 94], [1, 21, 44, 55, 79, 94], [0, 7, 17, 24, 54, 96], [4, 24, 47, 58, 82, 97], [10, 30, 53, 64, 88, 103], [1, 16, 34, 54, 77, 88], [18, 41, 52, 76, 91, 109], [0, 15, 22, 32, 39, 69], [2, 13, 37, 52, 70, 90], [6, 17, 19, 32, 35, 74], [2, 5, 44, 87, 98, 100], [21, 36, 43, 53, 60, 90], [4, 28, 43, 61, 81, 104], [35, 78, 89, 91, 104, 107], [5, 12, 42, 84, 99, 106], [8, 11, 50, 93, 104, 106], [3, 18, 25, 35, 42, 72], [2, 45, 56, 58, 71, 74], [5, 48, 59, 61, 74, 77], [24, 66, 81, 88, 98, 105], [30, 41, 43, 56, 59, 98], [9, 51, 66, 73, 83, 90], [24, 35, 37, 50, 53, 92], [12, 27, 34, 44, 51, 81], [9, 20, 22, 35, 38, 77], [15, 38, 49, 73, 88, 106], [29, 72, 83, 85, 98, 101], [7, 27, 50, 61, 85, 100], [27, 69, 84, 91, 101, 108], [15, 57, 72, 79, 89, 96], [7, 25, 45, 68, 79, 103], [19, 34, 52, 72, 95, 106], [3, 45, 60, 67, 77, 84], [15, 26, 28, 41, 44, 83], [0, 23, 34, 58, 73, 91], [13, 31, 51, 74, 85, 109], [9, 24, 31, 41, 48, 78], [13, 28, 46, 66, 89, 100], [33, 48, 55, 65, 72, 102], [9, 16, 26, 33, 63, 105], [26, 69, 80, 82, 95, 98], [0, 42, 57, 64, 74, 81], [12, 54, 69, 76, 86, 93], [22, 37, 55, 75, 98, 109], [16, 36, 59, 70, 94, 109], [0, 11, 13, 26, 29, 68], [27, 42, 49, 59, 66, 96], [3, 33, 75, 90, 97, 107], [10, 34, 49, 67, 87, 110], [9, 32, 43, 67, 82, 100], [36, 47, 49, 62, 65, 104], [22, 33, 41, 49, 89, 108], [12, 37, 48, 56, 64, 104], [26, 45, 70, 81, 89, 97], [23, 42, 67, 78, 86, 94], [7, 47, 66, 91, 102, 110], [9, 17, 25, 65, 84, 109], [17, 36, 61, 72, 80, 88], [5, 24, 49, 60, 68, 76], [32, 51, 76, 87, 95, 103], [1, 41, 60, 85, 96, 104], [2, 10, 50, 69, 94, 105], [0, 25, 36, 44, 52, 92], [38, 57, 82, 93, 101, 109], [8, 27, 52, 63, 71, 79], [3, 28, 39, 47, 55, 95], [4, 44, 63, 88, 99, 107], [14, 33, 58, 69, 77, 85], [9, 34, 45, 53, 61, 101], [20, 39, 64, 75, 83, 91], [18, 43, 54, 62, 70, 110], [4, 15, 23, 31, 71, 90], [19, 30, 38, 46, 86, 105], [3, 11, 19, 59, 78, 103], [5, 13, 53, 72, 97, 108], [11, 30, 55, 66, 74, 82], [7, 18, 26, 34, 74, 93], [10, 21, 29, 37, 77, 96], [6, 14, 22, 62, 81, 106], [15, 40, 51, 59, 67, 107], [2, 21, 46, 57, 65, 73], [16, 27, 35, 43, 83, 102], [13, 24, 32, 40, 80, 99], [0, 8, 16, 56, 75, 100], [6, 31, 42, 50, 58, 98], [35, 54, 79, 90, 98, 106], [1, 12, 20, 28, 68, 87], [29, 48, 73, 84, 92, 100], [17, 29, 55, 76, 81, 90], [6, 44, 56, 82, 103, 108], [29, 41, 67, 88, 93, 102], [0, 38, 50, 76, 97, 102], [32, 44, 70, 91, 96, 105], [11, 37, 58, 63, 72, 110], [7, 12, 21, 59, 71, 97], [22, 43, 48, 57, 95, 107], [13, 34, 39, 48, 86, 98], [2, 28, 49, 54, 63, 101], [3, 41, 53, 79, 100, 105], [26, 38, 64, 85, 90, 99], [1, 6, 15, 53, 65, 91], [8, 20, 46, 67, 72, 81], [4, 25, 30, 39, 77, 89], [14, 26, 52, 73, 78, 87], [4, 9, 18, 56, 68, 94], [10, 31, 36, 45, 83, 95], [23, 35, 61, 82, 87, 96], [10, 15, 24, 62, 74, 100], [13, 18, 27, 65, 77, 103], [25, 46, 51, 60, 98, 110], [2, 14, 40, 61, 66, 75], [3, 12, 50, 62, 88, 109], [11, 23, 49, 70, 75, 84], [35, 47, 73, 94, 99, 108], [0, 9, 47, 59, 85, 106], [16, 21, 30, 68, 80, 106], [5, 31, 52, 57, 66, 104], [16, 37, 42, 51, 89, 101], [8, 34, 55, 60, 69, 107], [20, 32, 58, 79, 84, 93], [1, 22, 27, 36, 74, 86], [5, 17, 43, 64, 69, 78], [19, 40, 45, 54, 92, 104], [7, 28, 33, 42, 80, 92], [19, 24, 33, 71, 83, 109], [8, 30, 44, 49, 57, 97], [1, 23, 45, 59, 64, 72], [2, 24, 38, 43, 51, 91], [3, 17, 22, 30, 70, 92], [11, 33, 47, 52, 60, 100], [20, 42, 56, 61, 69, 109], [6, 46, 68, 90, 104, 109], [0, 14, 19, 27, 67, 89], [5, 10, 18, 58, 80, 102], [16, 38, 60, 74, 79, 87], [11, 16, 24, 64, 86, 108], [4, 26, 48, 62, 67, 75], [22, 44, 66, 80, 85, 93], [15, 29, 34, 42, 82, 104], [13, 35, 57, 71, 76, 84], [8, 13, 21, 61, 83, 105], [14, 36, 50, 55, 63, 103], [10, 32, 54, 68, 73, 81], [37, 59, 81, 95, 100, 108], [9, 23, 28, 36, 76, 98], [19, 41, 63, 77, 82, 90], [34, 56, 78, 92, 97, 105], [4, 12, 52, 74, 96, 110], [1, 9, 49, 71, 93, 107], [21, 35, 40, 48, 88, 110], [2, 7, 15, 55, 77, 99], [6, 20, 25, 33, 73, 95], [3, 43, 65, 87, 101, 106], [17, 39, 53, 58, 66, 106], [18, 32, 37, 45, 85, 107], [31, 53, 75, 89, 94, 102], [12, 26, 31, 39, 79, 101], [5, 27, 41, 46, 54, 94], [0, 40, 62, 84, 98, 103], [25, 47, 69, 83, 88, 96], [28, 50, 72, 86, 91, 99], [7, 29, 51, 65, 70, 78]]
\item 1 \{0=1110, 1=26640, 2=524808, 3=2201352, 4=2374290\} [[8, 39, 40, 41, 42, 70], [0, 1, 2, 3, 31, 80], [30, 31, 32, 33, 61, 110], [6, 7, 8, 9, 37, 86], [23, 54, 55, 56, 57, 85], [15, 16, 17, 18, 46, 95], [27, 28, 29, 30, 58, 107], [24, 25, 26, 27, 55, 104], [26, 57, 58, 59, 60, 88], [29, 60, 61, 62, 63, 91], [25, 74, 105, 106, 107, 108], [41, 72, 73, 74, 75, 103], [7, 56, 87, 88, 89, 90], [10, 59, 90, 91, 92, 93], [13, 62, 93, 94, 95, 96], [22, 71, 102, 103, 104, 105], [9, 10, 11, 12, 40, 89], [44, 75, 76, 77, 78, 106], [18, 19, 20, 21, 49, 98], [11, 42, 43, 44, 45, 73], [3, 4, 5, 6, 34, 83], [16, 65, 96, 97, 98, 99], [19, 68, 99, 100, 101, 102], [17, 48, 49, 50, 51, 79], [47, 78, 79, 80, 81, 109], [20, 51, 52, 53, 54, 82], [12, 13, 14, 15, 43, 92], [0, 28, 77, 108, 109, 110], [4, 53, 84, 85, 86, 87], [1, 50, 81, 82, 83, 84], [32, 63, 64, 65, 66, 94], [2, 33, 34, 35, 36, 64], [21, 22, 23, 24, 52, 101], [38, 69, 70, 71, 72, 100], [14, 45, 46, 47, 48, 76], [35, 66, 67, 68, 69, 97], [5, 36, 37, 38, 39, 67], [30, 41, 47, 54, 67, 98], [32, 43, 78, 82, 88, 101], [42, 57, 64, 82, 92, 102], [7, 38, 81, 92, 98, 105], [9, 20, 26, 33, 46, 77], [1, 36, 40, 46, 59, 101], [11, 22, 57, 61, 67, 80], [14, 25, 60, 64, 70, 83], [17, 60, 71, 77, 84, 97], [4, 39, 43, 49, 62, 104], [4, 22, 32, 42, 93, 108], [10, 45, 49, 55, 68, 110], [3, 54, 69, 76, 94, 104], [48, 63, 70, 88, 98, 108], [6, 10, 16, 29, 71, 82], [33, 44, 50, 57, 70, 101], [42, 53, 59, 66, 79, 110], [8, 51, 62, 68, 75, 88], [15, 30, 37, 55, 65, 75], [5, 47, 58, 93, 97, 103], [5, 16, 51, 55, 61, 74], [6, 13, 31, 41, 51, 102], [39, 50, 56, 63, 76, 107], [6, 17, 23, 30, 43, 74], [1, 14, 56, 67, 102, 106], [24, 28, 34, 47, 89, 100], [21, 32, 38, 45, 58, 89], [3, 7, 13, 26, 68, 79], [0, 4, 10, 23, 65, 76], [24, 39, 46, 64, 74, 84], [12, 16, 22, 35, 77, 88], [33, 48, 55, 73, 83, 93], [0, 7, 25, 35, 45, 96], [2, 45, 56, 62, 69, 82], [27, 31, 37, 50, 92, 103], [36, 47, 53, 60, 73, 104], [0, 13, 44, 87, 98, 104], [6, 21, 28, 46, 56, 66], [1, 32, 75, 86, 92, 99], [15, 26, 32, 39, 52, 83], [8, 18, 69, 84, 91, 109], [16, 26, 36, 87, 102, 109], [7, 17, 27, 78, 93, 100], [12, 23, 29, 36, 49, 80], [2, 8, 15, 28, 59, 102], [18, 33, 40, 58, 68, 78], [27, 42, 49, 67, 77, 87], [8, 19, 54, 58, 64, 77], [27, 38, 44, 51, 64, 95], [11, 54, 65, 71, 78, 91], [6, 57, 72, 79, 97, 107], [18, 22, 28, 41, 83, 94], [18, 29, 35, 42, 55, 86], [1, 19, 29, 39, 90, 105], [15, 19, 25, 38, 80, 91], [39, 54, 61, 79, 89, 99], [36, 51, 58, 76, 86, 96], [1, 7, 20, 62, 73, 108], [35, 46, 81, 85, 91, 104], [13, 23, 33, 84, 99, 106], [5, 48, 59, 65, 72, 85], [2, 9, 22, 53, 96, 107], [3, 10, 28, 38, 48, 99], [2, 12, 63, 78, 85, 103], [1, 11, 21, 72, 87, 94], [20, 63, 74, 80, 87, 100], [26, 69, 80, 86, 93, 106], [26, 37, 72, 76, 82, 95], [24, 35, 41, 48, 61, 92], [12, 19, 37, 47, 57, 108], [5, 12, 25, 56, 99, 110], [0, 15, 22, 40, 50, 60], [11, 53, 64, 99, 103, 109], [8, 50, 61, 96, 100, 106], [2, 13, 48, 52, 58, 71], [10, 20, 30, 81, 96, 103], [41, 52, 87, 91, 97, 110], [5, 15, 66, 81, 88, 106], [4, 14, 24, 75, 90, 97], [17, 28, 63, 67, 73, 86], [7, 42, 46, 52, 65, 107], [29, 72, 83, 89, 96, 109], [12, 27, 34, 52, 62, 72], [2, 44, 55, 90, 94, 100], [23, 66, 77, 83, 90, 103], [3, 14, 20, 27, 40, 71], [38, 49, 84, 88, 94, 107], [4, 17, 59, 70, 105, 109], [21, 36, 43, 61, 71, 81], [9, 60, 75, 82, 100, 110], [8, 14, 21, 34, 65, 108], [14, 57, 68, 74, 81, 94], [5, 11, 18, 31, 62, 105], [9, 16, 34, 44, 54, 105], [30, 34, 40, 53, 95, 106], [23, 34, 69, 73, 79, 92], [9, 24, 31, 49, 59, 69], [21, 25, 31, 44, 86, 97], [33, 37, 43, 56, 98, 109], [3, 18, 25, 43, 53, 63], [29, 40, 75, 79, 85, 98], [4, 35, 78, 89, 95, 102], [45, 60, 67, 85, 95, 105], [0, 11, 17, 24, 37, 68], [10, 41, 84, 95, 101, 108], [6, 19, 50, 93, 104, 110], [20, 31, 66, 70, 76, 89], [30, 45, 52, 70, 80, 90], [3, 16, 47, 90, 101, 107], [9, 13, 19, 32, 74, 85], [0, 51, 66, 73, 91, 101], [40, 49, 54, 73, 81, 102], [19, 28, 33, 52, 60, 81], [7, 15, 36, 85, 94, 99], [2, 5, 43, 79, 91, 95], [0, 47, 69, 74, 83, 99], [32, 54, 59, 68, 84, 96], [6, 55, 64, 69, 88, 96], [11, 14, 52, 88, 100, 104], [10, 22, 26, 44, 47, 85], [14, 36, 41, 50, 66, 78], [19, 55, 67, 71, 89, 92], [9, 21, 68, 90, 95, 104], [13, 21, 42, 91, 100, 105], [4, 8, 26, 29, 67, 103], [29, 51, 56, 65, 81, 93], [2, 18, 30, 77, 99, 104], [34, 70, 82, 86, 104, 107], [3, 24, 73, 82, 87, 106], [2, 40, 76, 88, 92, 110], [26, 48, 53, 62, 78, 90], [7, 19, 23, 41, 44, 82], [22, 34, 38, 56, 59, 97], [25, 34, 39, 58, 66, 87], [28, 40, 44, 62, 65, 103], [1, 37, 49, 53, 71, 74], [4, 40, 52, 56, 74, 77], [17, 39, 44, 53, 69, 81], [12, 24, 71, 93, 98, 107], [44, 66, 71, 80, 96, 108], [3, 50, 72, 77, 86, 102], [13, 49, 61, 65, 83, 86], [13, 25, 29, 47, 50, 88], [34, 46, 50, 68, 71, 109], [2, 20, 23, 61, 97, 109], [0, 21, 70, 79, 84, 103], [6, 18, 65, 87, 92, 101], [6, 53, 75, 80, 89, 105], [6, 11, 20, 36, 48, 95], [21, 26, 35, 51, 63, 110], [19, 31, 35, 53, 56, 94], [16, 24, 45, 94, 103, 108], [34, 43, 48, 67, 75, 96], [1, 5, 23, 26, 64, 100], [3, 22, 30, 51, 100, 109], [1, 10, 15, 34, 42, 63], [4, 16, 20, 38, 41, 79], [7, 16, 21, 40, 48, 69], [37, 73, 85, 89, 107, 110], [14, 17, 55, 91, 103, 107], [16, 28, 32, 50, 53, 91], [16, 25, 30, 49, 57, 78], [16, 52, 64, 68, 86, 89], [7, 11, 29, 32, 70, 106], [6, 27, 76, 85, 90, 109], [0, 49, 58, 63, 82, 90], [3, 15, 62, 84, 89, 98], [0, 12, 59, 81, 86, 95], [18, 67, 76, 81, 100, 108], [3, 52, 61, 66, 85, 93], [5, 8, 46, 82, 94, 98], [11, 33, 38, 47, 63, 75], [25, 37, 41, 59, 62, 100], [3, 8, 17, 33, 45, 92], [8, 24, 36, 83, 105, 110], [8, 30, 35, 44, 60, 72], [8, 11, 49, 85, 97, 101], [23, 45, 50, 59, 75, 87], [10, 19, 24, 43, 51, 72], [43, 52, 57, 76, 84, 105], [12, 61, 70, 75, 94, 102], [2, 24, 29, 38, 54, 66], [35, 57, 62, 71, 87, 99], [4, 9, 28, 36, 57, 106], [1, 6, 25, 33, 54, 103], [28, 37, 42, 61, 69, 90], [9, 58, 67, 72, 91, 99], [22, 58, 70, 74, 92, 95], [15, 20, 29, 45, 57, 104], [7, 12, 31, 39, 60, 109], [38, 60, 65, 74, 90, 102], [31, 43, 47, 65, 68, 106], [7, 43, 55, 59, 77, 80], [9, 14, 23, 39, 51, 98], [25, 61, 73, 77, 95, 98], [31, 40, 45, 64, 72, 93], [0, 5, 14, 30, 42, 89], [9, 56, 78, 83, 92, 108], [46, 55, 60, 79, 87, 108], [13, 22, 27, 46, 54, 75], [15, 64, 73, 78, 97, 105], [10, 14, 32, 35, 73, 109], [10, 46, 58, 62, 80, 83], [20, 42, 47, 56, 72, 84], [41, 63, 68, 77, 93, 105], [5, 27, 32, 41, 57, 69], [12, 17, 26, 42, 54, 101], [22, 31, 36, 55, 63, 84], [31, 67, 79, 83, 101, 104], [28, 64, 76, 80, 98, 101], [4, 13, 18, 37, 45, 66], [37, 46, 51, 70, 78, 99], [15, 27, 74, 96, 101, 110], [17, 20, 58, 94, 106, 110], [4, 12, 33, 82, 91, 96], [0, 19, 27, 48, 97, 106], [1, 13, 17, 35, 38, 76], [18, 23, 32, 48, 60, 107], [2, 11, 27, 39, 86, 108], [5, 21, 33, 80, 102, 107], [1, 9, 30, 79, 88, 93], [10, 18, 39, 88, 97, 102], [43, 58, 69, 85, 102, 108], [10, 21, 37, 54, 60, 106], [1, 41, 58, 65, 104, 109], [23, 38, 42, 78, 96, 104], [5, 10, 13, 53, 70, 77], [33, 51, 59, 89, 104, 108], [23, 40, 47, 86, 91, 94], [0, 36, 54, 62, 92, 107], [2, 6, 42, 60, 68, 98], [37, 52, 63, 79, 96, 102], [12, 30, 38, 68, 83, 87], [2, 17, 21, 57, 75, 83], [0, 18, 26, 56, 71, 75], [1, 16, 27, 43, 60, 66], [10, 25, 36, 52, 69, 75], [38, 55, 62, 101, 106, 109], [4, 19, 30, 46, 63, 69], [4, 15, 31, 48, 54, 100], [15, 33, 41, 71, 86, 90], [29, 44, 48, 84, 102, 110], [19, 34, 45, 61, 78, 84], [29, 46, 53, 92, 97, 100], [13, 20, 59, 64, 67, 107], [15, 21, 67, 82, 93, 109], [3, 39, 57, 65, 95, 110], [29, 34, 37, 77, 94, 101], [0, 8, 38, 53, 57, 93], [26, 31, 34, 74, 91, 98], [28, 43, 54, 70, 87, 93], [8, 23, 27, 63, 81, 89], [5, 9, 45, 63, 71, 101], [9, 17, 47, 62, 66, 102], [38, 43, 46, 86, 103, 110], [7, 18, 34, 51, 57, 103], [4, 11, 50, 55, 58, 98], [40, 55, 66, 82, 99, 105], [32, 49, 56, 95, 100, 103], [0, 6, 52, 67, 78, 94], [20, 37, 44, 83, 88, 91], [4, 21, 27, 73, 88, 99], [3, 21, 29, 59, 74, 78], [8, 25, 32, 71, 76, 79], [17, 34, 41, 80, 85, 88], [5, 22, 29, 68, 73, 76], [3, 9, 55, 70, 81, 97], [2, 41, 46, 49, 89, 106], [34, 49, 60, 76, 93, 99], [7, 22, 33, 49, 66, 72], [6, 12, 58, 73, 84, 100], [6, 14, 44, 59, 63, 99], [1, 8, 47, 52, 55, 95], [7, 24, 30, 76, 91, 102], [14, 18, 54, 72, 80, 110], [3, 49, 64, 75, 91, 108], [7, 14, 53, 58, 61, 101], [2, 19, 26, 65, 70, 73], [12, 20, 50, 65, 69, 105], [0, 16, 33, 39, 85, 100], [5, 44, 49, 52, 92, 109], [5, 20, 24, 60, 78, 86], [16, 23, 62, 67, 70, 110], [14, 29, 33, 69, 87, 95], [6, 24, 32, 62, 77, 81], [5, 35, 50, 54, 90, 108], [20, 35, 39, 75, 93, 101], [9, 27, 35, 65, 80, 84], [14, 31, 38, 77, 82, 85], [21, 39, 47, 77, 92, 96], [13, 30, 36, 82, 97, 108], [24, 42, 50, 80, 95, 99], [2, 7, 10, 50, 67, 74], [3, 11, 41, 56, 60, 96], [12, 18, 64, 79, 90, 106], [11, 28, 35, 74, 79, 82], [4, 7, 47, 64, 71, 110], [15, 23, 53, 68, 72, 108], [30, 48, 56, 86, 101, 105], [26, 43, 50, 89, 94, 97], [35, 40, 43, 83, 100, 107], [14, 19, 22, 62, 79, 86], [17, 22, 25, 65, 82, 89], [20, 25, 28, 68, 85, 92], [23, 28, 31, 71, 88, 95], [25, 40, 51, 67, 84, 90], [6, 22, 39, 45, 91, 106], [32, 37, 40, 80, 97, 104], [22, 37, 48, 64, 81, 87], [13, 24, 40, 57, 63, 109], [13, 28, 39, 55, 72, 78], [31, 46, 57, 73, 90, 96], [18, 36, 44, 74, 89, 93], [35, 52, 59, 98, 103, 106], [9, 15, 61, 76, 87, 103], [26, 41, 45, 81, 99, 107], [3, 19, 36, 42, 88, 103], [11, 26, 30, 66, 84, 92], [2, 32, 47, 51, 87, 105], [1, 18, 24, 70, 85, 96], [9, 25, 42, 48, 94, 109], [16, 31, 42, 58, 75, 81], [11, 16, 19, 59, 76, 83], [8, 12, 48, 66, 74, 104], [8, 13, 16, 56, 73, 80], [0, 46, 61, 72, 88, 105], [11, 15, 51, 69, 77, 107], [17, 32, 36, 72, 90, 98], [27, 45, 53, 83, 98, 102], [1, 12, 28, 45, 51, 97], [10, 17, 56, 61, 64, 104], [1, 4, 44, 61, 68, 107], [10, 27, 33, 79, 94, 105], [27, 36, 56, 68, 70, 91], [33, 42, 62, 74, 76, 97], [6, 15, 35, 47, 49, 70], [14, 26, 28, 49, 96, 105], [11, 23, 25, 46, 93, 102], [7, 54, 63, 83, 95, 97], [0, 20, 32, 34, 55, 102], [10, 57, 66, 86, 98, 100], [5, 17, 19, 40, 87, 96], [3, 23, 35, 37, 58, 105], [2, 14, 16, 37, 84, 93], [5, 7, 28, 75, 84, 104], [11, 13, 34, 81, 90, 110], [19, 66, 75, 95, 107, 109], [21, 30, 50, 62, 64, 85], [6, 26, 38, 40, 61, 108], [30, 39, 59, 71, 73, 94], [8, 10, 31, 78, 87, 107], [13, 60, 69, 89, 101, 103], [12, 21, 41, 53, 55, 76], [45, 54, 74, 86, 88, 109], [36, 45, 65, 77, 79, 100], [4, 51, 60, 80, 92, 94], [18, 27, 47, 59, 61, 82], [2, 4, 25, 72, 81, 101], [15, 24, 44, 56, 58, 79], [16, 63, 72, 92, 104, 106], [1, 48, 57, 77, 89, 91], [0, 9, 29, 41, 43, 64], [42, 51, 71, 83, 85, 106], [1, 22, 69, 78, 98, 110], [3, 12, 32, 44, 46, 67], [17, 29, 31, 52, 99, 108], [24, 33, 53, 65, 67, 88], [9, 18, 38, 50, 52, 73], [39, 48, 68, 80, 82, 103], [8, 20, 22, 43, 90, 99]]
\item 1 \{1=46176, 2=538128, 3=2211120, 4=2332776\} [[8, 39, 40, 41, 42, 70], [0, 1, 2, 3, 31, 80], [30, 31, 32, 33, 61, 110], [6, 7, 8, 9, 37, 86], [23, 54, 55, 56, 57, 85], [15, 16, 17, 18, 46, 95], [27, 28, 29, 30, 58, 107], [24, 25, 26, 27, 55, 104], [26, 57, 58, 59, 60, 88], [29, 60, 61, 62, 63, 91], [25, 74, 105, 106, 107, 108], [41, 72, 73, 74, 75, 103], [7, 56, 87, 88, 89, 90], [10, 59, 90, 91, 92, 93], [13, 62, 93, 94, 95, 96], [22, 71, 102, 103, 104, 105], [9, 10, 11, 12, 40, 89], [44, 75, 76, 77, 78, 106], [18, 19, 20, 21, 49, 98], [11, 42, 43, 44, 45, 73], [3, 4, 5, 6, 34, 83], [16, 65, 96, 97, 98, 99], [19, 68, 99, 100, 101, 102], [17, 48, 49, 50, 51, 79], [47, 78, 79, 80, 81, 109], [20, 51, 52, 53, 54, 82], [12, 13, 14, 15, 43, 92], [0, 28, 77, 108, 109, 110], [4, 53, 84, 85, 86, 87], [1, 50, 81, 82, 83, 84], [32, 63, 64, 65, 66, 94], [2, 33, 34, 35, 36, 64], [21, 22, 23, 24, 52, 101], [38, 69, 70, 71, 72, 100], [14, 45, 46, 47, 48, 76], [35, 66, 67, 68, 69, 97], [5, 36, 37, 38, 39, 67], [11, 54, 58, 61, 87, 92], [24, 29, 59, 102, 106, 109], [33, 37, 40, 66, 71, 101], [42, 46, 49, 75, 80, 110], [2, 45, 49, 52, 78, 83], [6, 10, 13, 39, 44, 74], [12, 17, 47, 90, 94, 97], [23, 66, 70, 73, 99, 104], [18, 23, 53, 96, 100, 103], [15, 19, 22, 48, 53, 83], [6, 11, 41, 84, 88, 91], [0, 5, 35, 78, 82, 85], [1, 27, 32, 62, 105, 109], [26, 69, 73, 76, 102, 107], [15, 20, 50, 93, 97, 100], [3, 7, 10, 36, 41, 71], [21, 26, 56, 99, 103, 106], [30, 34, 37, 63, 68, 98], [18, 22, 25, 51, 56, 86], [27, 31, 34, 60, 65, 95], [24, 28, 31, 57, 62, 92], [29, 72, 76, 79, 105, 110], [14, 57, 61, 64, 90, 95], [9, 14, 44, 87, 91, 94], [8, 51, 55, 58, 84, 89], [39, 43, 46, 72, 77, 107], [36, 40, 43, 69, 74, 104], [12, 16, 19, 45, 50, 80], [2, 32, 75, 79, 82, 108], [20, 63, 67, 70, 96, 101], [3, 8, 38, 81, 85, 88], [0, 4, 7, 33, 38, 68], [21, 25, 28, 54, 59, 89], [1, 4, 30, 35, 65, 108], [5, 48, 52, 55, 81, 86], [17, 60, 64, 67, 93, 98], [9, 13, 16, 42, 47, 77], [53, 56, 76, 80, 91, 101], [9, 22, 31, 63, 73, 106], [24, 44, 54, 60, 72, 83], [3, 23, 33, 39, 51, 62], [5, 15, 21, 33, 44, 96], [1, 10, 42, 52, 85, 99], [29, 32, 52, 56, 67, 77], [1, 11, 74, 77, 97, 101], [5, 8, 28, 32, 43, 53], [0, 10, 43, 57, 70, 79], [5, 57, 77, 87, 93, 105], [3, 9, 21, 32, 84, 104], [14, 24, 30, 42, 53, 105], [9, 19, 52, 66, 79, 88], [5, 68, 71, 91, 95, 106], [4, 18, 31, 40, 72, 82], [25, 39, 52, 61, 93, 103], [8, 18, 24, 36, 47, 99], [21, 31, 64, 78, 91, 100], [51, 71, 81, 87, 99, 110], [21, 41, 51, 57, 69, 80], [22, 36, 49, 58, 90, 100], [1, 33, 43, 76, 90, 103], [12, 32, 42, 48, 60, 71], [7, 16, 48, 58, 91, 105], [28, 42, 55, 64, 96, 106], [23, 26, 46, 50, 61, 71], [2, 12, 18, 30, 41, 93], [12, 25, 34, 66, 76, 109], [8, 71, 74, 94, 98, 109], [19, 23, 34, 44, 107, 110], [4, 8, 19, 29, 92, 95], [20, 23, 43, 47, 58, 68], [30, 40, 73, 87, 100, 109], [9, 29, 39, 45, 57, 68], [44, 47, 67, 71, 82, 92], [2, 5, 25, 29, 40, 50], [0, 12, 23, 75, 95, 105], [1, 34, 48, 61, 70, 102], [4, 37, 51, 64, 73, 105], [35, 38, 58, 62, 73, 83], [4, 13, 45, 55, 88, 102], [17, 20, 40, 44, 55, 65], [7, 39, 49, 82, 96, 109], [13, 27, 40, 49, 81, 91], [7, 17, 80, 83, 103, 107], [59, 62, 82, 86, 97, 107], [7, 40, 54, 67, 76, 108], [27, 47, 57, 63, 75, 86], [9, 15, 27, 38, 90, 110], [14, 17, 37, 41, 52, 62], [6, 17, 69, 89, 99, 105], [32, 35, 55, 59, 70, 80], [50, 53, 73, 77, 88, 98], [12, 22, 55, 69, 82, 91], [36, 56, 66, 72, 84, 95], [7, 11, 22, 32, 95, 98], [10, 24, 37, 46, 78, 88], [47, 50, 70, 74, 85, 95], [2, 54, 74, 84, 90, 102], [8, 11, 31, 35, 46, 56], [15, 35, 45, 51, 63, 74], [3, 15, 26, 78, 98, 108], [6, 26, 36, 42, 54, 65], [3, 13, 46, 60, 73, 82], [30, 50, 60, 66, 78, 89], [6, 19, 28, 60, 70, 103], [3, 14, 66, 86, 96, 102], [11, 14, 34, 38, 49, 59], [7, 21, 34, 43, 75, 85], [0, 6, 18, 29, 81, 101], [17, 27, 33, 45, 56, 108], [18, 28, 61, 75, 88, 97], [1, 5, 16, 26, 89, 92], [10, 19, 51, 61, 94, 108], [38, 41, 61, 65, 76, 86], [8, 60, 80, 90, 96, 108], [31, 45, 58, 67, 99, 109], [45, 65, 75, 81, 93, 104], [24, 34, 67, 81, 94, 103], [4, 36, 46, 79, 93, 106], [56, 59, 79, 83, 94, 104], [4, 14, 77, 80, 100, 104], [3, 16, 25, 57, 67, 100], [16, 20, 31, 41, 104, 107], [2, 65, 68, 88, 92, 103], [41, 44, 64, 68, 79, 89], [2, 22, 26, 37, 47, 110], [10, 14, 25, 35, 98, 101], [11, 21, 27, 39, 50, 102], [26, 29, 49, 53, 64, 74], [18, 38, 48, 54, 66, 77], [16, 30, 43, 52, 84, 94], [0, 20, 30, 36, 48, 59], [0, 13, 22, 54, 64, 97], [0, 11, 63, 83, 93, 99], [13, 17, 28, 38, 101, 104], [42, 62, 72, 78, 90, 101], [39, 59, 69, 75, 87, 98], [15, 25, 58, 72, 85, 94], [1, 15, 28, 37, 69, 79], [9, 20, 72, 92, 102, 108], [27, 37, 70, 84, 97, 106], [62, 65, 85, 89, 100, 110], [10, 20, 83, 86, 106, 110], [2, 13, 23, 86, 89, 109], [6, 12, 24, 35, 87, 107], [6, 16, 49, 63, 76, 85], [33, 53, 63, 69, 81, 92], [48, 68, 78, 84, 96, 107], [19, 33, 46, 55, 87, 97], [4, 26, 52, 70, 75, 91], [0, 9, 56, 61, 69, 96], [13, 37, 59, 85, 103, 108], [15, 30, 39, 86, 91, 99], [41, 46, 54, 81, 96, 105], [1, 6, 22, 46, 68, 94], [2, 6, 14, 20, 56, 73], [4, 9, 25, 49, 71, 97], [1, 25, 47, 73, 91, 96], [1, 9, 36, 51, 60, 107], [29, 33, 41, 47, 83, 100], [3, 30, 45, 54, 101, 106], [1, 41, 45, 53, 59, 95], [5, 31, 49, 54, 70, 94], [5, 9, 17, 23, 59, 76], [14, 19, 27, 54, 69, 78], [22, 40, 45, 61, 85, 107], [17, 21, 29, 35, 71, 88], [11, 37, 55, 60, 76, 100], [8, 12, 20, 26, 62, 79], [12, 21, 68, 73, 81, 108], [38, 43, 51, 78, 93, 102], [1, 19, 24, 40, 64, 86], [23, 27, 35, 41, 77, 94], [32, 49, 89, 93, 101, 107], [5, 10, 18, 45, 60, 69], [32, 37, 45, 72, 87, 96], [3, 19, 43, 65, 91, 109], [16, 34, 39, 55, 79, 101], [17, 34, 74, 78, 86, 92], [18, 33, 42, 89, 94, 102], [3, 50, 55, 63, 90, 105], [2, 19, 59, 63, 71, 77], [3, 18, 27, 74, 79, 87], [13, 31, 36, 52, 76, 98], [26, 31, 39, 66, 81, 90], [14, 40, 58, 63, 79, 103], [10, 34, 56, 82, 100, 105], [35, 52, 92, 96, 104, 110], [29, 34, 42, 69, 84, 93], [14, 31, 71, 75, 83, 89], [24, 39, 48, 95, 100, 108], [20, 24, 32, 38, 74, 91], [21, 36, 45, 92, 97, 105], [14, 18, 26, 32, 68, 85], [20, 37, 77, 81, 89, 95], [23, 40, 80, 84, 92, 98], [11, 15, 23, 29, 65, 82], [19, 37, 42, 58, 82, 104], [0, 47, 52, 60, 87, 102], [2, 8, 44, 61, 101, 105], [7, 31, 53, 79, 97, 102], [7, 47, 51, 59, 65, 101], [5, 22, 62, 66, 74, 80], [6, 53, 58, 66, 93, 108], [6, 21, 30, 77, 82, 90], [26, 43, 83, 87, 95, 101], [0, 8, 14, 50, 67, 107], [4, 22, 27, 43, 67, 89], [9, 18, 65, 70, 78, 105], [0, 15, 24, 71, 76, 84], [26, 30, 38, 44, 80, 97], [29, 46, 86, 90, 98, 104], [2, 38, 55, 95, 99, 107], [7, 25, 30, 46, 70, 92], [17, 43, 61, 66, 82, 106], [11, 28, 68, 72, 80, 86], [4, 44, 48, 56, 62, 98], [44, 49, 57, 84, 99, 108], [0, 27, 42, 51, 98, 103], [13, 53, 57, 65, 71, 107], [6, 33, 48, 57, 104, 109], [16, 38, 64, 82, 87, 103], [7, 12, 28, 52, 74, 100], [16, 21, 37, 61, 83, 109], [23, 28, 36, 63, 78, 87], [5, 11, 47, 64, 104, 108], [1, 23, 49, 67, 72, 88], [2, 28, 46, 51, 67, 91], [12, 27, 36, 83, 88, 96], [13, 18, 34, 58, 80, 106], [17, 22, 30, 57, 72, 81], [7, 29, 55, 73, 78, 94], [10, 28, 33, 49, 73, 95], [8, 13, 21, 48, 63, 72], [22, 44, 70, 88, 93, 109], [3, 12, 59, 64, 72, 99], [4, 28, 50, 76, 94, 99], [6, 15, 62, 67, 75, 102], [16, 56, 60, 68, 74, 110], [8, 34, 52, 57, 73, 97], [2, 7, 15, 42, 57, 66], [19, 41, 67, 85, 90, 106], [10, 15, 31, 55, 77, 103], [35, 39, 47, 53, 89, 106], [11, 16, 24, 51, 66, 75], [13, 35, 61, 79, 84, 100], [35, 40, 48, 75, 90, 99], [9, 24, 33, 80, 85, 93], [10, 32, 58, 76, 81, 97], [25, 43, 48, 64, 88, 110], [5, 41, 58, 98, 102, 110], [3, 11, 17, 53, 70, 110], [8, 25, 65, 69, 77, 83], [32, 36, 44, 50, 86, 103], [10, 50, 54, 62, 68, 104], [4, 12, 39, 54, 63, 110], [20, 46, 64, 69, 85, 109], [38, 42, 50, 56, 92, 109], [20, 25, 33, 60, 75, 84], [0, 16, 40, 62, 88, 106], [21, 38, 40, 53, 60, 94], [31, 69, 86, 88, 101, 108], [8, 10, 23, 30, 64, 102], [11, 13, 26, 33, 67, 105], [6, 40, 78, 95, 97, 110], [25, 63, 80, 82, 95, 102], [22, 60, 77, 79, 92, 99], [8, 15, 49, 87, 104, 106], [2, 4, 17, 24, 58, 96], [13, 51, 68, 70, 83, 90], [1, 14, 21, 55, 93, 110], [18, 35, 37, 50, 57, 91], [3, 37, 75, 92, 94, 107], [15, 32, 34, 47, 54, 88], [36, 53, 55, 68, 75, 109], [30, 47, 49, 62, 69, 103], [11, 18, 52, 90, 107, 109], [4, 42, 59, 61, 74, 81], [19, 57, 74, 76, 89, 96], [24, 41, 43, 56, 63, 97], [6, 23, 25, 38, 45, 79], [2, 9, 43, 81, 98, 100], [7, 45, 62, 64, 77, 84], [14, 16, 29, 36, 70, 108], [5, 12, 46, 84, 101, 103], [5, 7, 20, 27, 61, 99], [12, 29, 31, 44, 51, 85], [27, 44, 46, 59, 66, 100], [0, 34, 72, 89, 91, 104], [16, 54, 71, 73, 86, 93], [9, 26, 28, 41, 48, 82], [28, 66, 83, 85, 98, 105], [10, 48, 65, 67, 80, 87], [33, 50, 52, 65, 72, 106], [3, 20, 22, 35, 42, 76], [0, 17, 19, 32, 39, 73], [1, 39, 56, 58, 71, 78], [22, 34, 41, 50, 87, 108], [24, 45, 70, 82, 89, 98], [0, 21, 46, 58, 65, 74], [10, 17, 26, 63, 84, 109], [19, 31, 38, 47, 84, 105], [2, 39, 60, 85, 97, 104], [15, 36, 61, 73, 80, 89], [0, 25, 37, 44, 53, 90], [9, 34, 46, 53, 62, 99], [15, 40, 52, 59, 68, 105], [5, 42, 63, 88, 100, 107], [3, 28, 40, 47, 56, 93], [30, 51, 76, 88, 95, 104], [18, 43, 55, 62, 71, 108], [4, 11, 20, 57, 78, 103], [8, 45, 66, 91, 103, 110], [9, 30, 55, 67, 74, 83], [7, 19, 26, 35, 72, 93], [10, 22, 29, 38, 75, 96], [33, 54, 79, 91, 98, 107], [18, 39, 64, 76, 83, 92], [7, 14, 23, 60, 81, 106], [36, 57, 82, 94, 101, 110], [16, 28, 35, 44, 81, 102], [5, 14, 51, 72, 97, 109], [12, 37, 49, 56, 65, 102], [27, 48, 73, 85, 92, 101], [21, 42, 67, 79, 86, 95], [13, 25, 32, 41, 78, 99], [4, 16, 23, 32, 69, 90], [1, 13, 20, 29, 66, 87], [1, 8, 17, 54, 75, 100], [6, 27, 52, 64, 71, 80], [3, 24, 49, 61, 68, 77], [2, 11, 48, 69, 94, 106], [12, 33, 58, 70, 77, 86], [6, 31, 43, 50, 59, 96], [1, 7, 18, 44, 63, 104], [6, 47, 55, 61, 72, 98], [14, 33, 74, 82, 88, 99], [15, 56, 64, 70, 81, 107], [14, 22, 28, 39, 65, 84], [2, 21, 62, 70, 76, 87], [9, 50, 58, 64, 75, 101], [8, 16, 22, 33, 59, 78], [5, 24, 65, 73, 79, 90], [3, 29, 48, 89, 97, 103], [3, 44, 52, 58, 69, 95], [12, 53, 61, 67, 78, 104], [17, 36, 77, 85, 91, 102], [32, 40, 46, 57, 83, 102], [20, 39, 80, 88, 94, 105], [0, 26, 45, 86, 94, 100], [23, 31, 37, 48, 74, 93], [2, 10, 16, 27, 53, 72], [1, 12, 38, 57, 98, 106], [4, 10, 21, 47, 66, 107], [18, 59, 67, 73, 84, 110], [26, 34, 40, 51, 77, 96], [35, 43, 49, 60, 86, 105], [6, 32, 51, 92, 100, 106], [11, 19, 25, 36, 62, 81], [23, 42, 83, 91, 97, 108], [20, 28, 34, 45, 71, 90], [8, 27, 68, 76, 82, 93], [5, 13, 19, 30, 56, 75], [0, 41, 49, 55, 66, 92], [4, 15, 41, 60, 101, 109], [7, 13, 24, 50, 69, 110], [11, 30, 71, 79, 85, 96], [29, 37, 43, 54, 80, 99], [17, 25, 31, 42, 68, 87], [38, 46, 52, 63, 89, 108], [9, 35, 54, 95, 103, 109]]
\item 1 \{1=26640, 2=516816, 3=2203128, 4=2381616\} [[8, 39, 40, 41, 42, 70], [0, 1, 2, 3, 31, 80], [30, 31, 32, 33, 61, 110], [6, 7, 8, 9, 37, 86], [23, 54, 55, 56, 57, 85], [15, 16, 17, 18, 46, 95], [27, 28, 29, 30, 58, 107], [24, 25, 26, 27, 55, 104], [26, 57, 58, 59, 60, 88], [29, 60, 61, 62, 63, 91], [25, 74, 105, 106, 107, 108], [41, 72, 73, 74, 75, 103], [7, 56, 87, 88, 89, 90], [10, 59, 90, 91, 92, 93], [13, 62, 93, 94, 95, 96], [22, 71, 102, 103, 104, 105], [9, 10, 11, 12, 40, 89], [44, 75, 76, 77, 78, 106], [18, 19, 20, 21, 49, 98], [11, 42, 43, 44, 45, 73], [3, 4, 5, 6, 34, 83], [16, 65, 96, 97, 98, 99], [19, 68, 99, 100, 101, 102], [17, 48, 49, 50, 51, 79], [47, 78, 79, 80, 81, 109], [20, 51, 52, 53, 54, 82], [12, 13, 14, 15, 43, 92], [0, 28, 77, 108, 109, 110], [4, 53, 84, 85, 86, 87], [1, 50, 81, 82, 83, 84], [32, 63, 64, 65, 66, 94], [2, 33, 34, 35, 36, 64], [21, 22, 23, 24, 52, 101], [38, 69, 70, 71, 72, 100], [14, 45, 46, 47, 48, 76], [35, 66, 67, 68, 69, 97], [5, 36, 37, 38, 39, 67], [10, 25, 35, 41, 54, 99], [17, 22, 28, 33, 56, 93], [0, 5, 13, 54, 58, 65], [39, 61, 76, 86, 92, 105], [9, 14, 22, 63, 67, 74], [30, 65, 70, 76, 81, 104], [24, 59, 64, 70, 75, 98], [4, 14, 20, 33, 78, 100], [21, 25, 32, 78, 83, 91], [9, 44, 49, 55, 60, 83], [0, 4, 11, 57, 62, 70], [11, 16, 22, 27, 50, 87], [21, 26, 34, 75, 79, 86], [14, 19, 25, 30, 53, 90], [39, 43, 50, 96, 101, 109], [33, 68, 73, 79, 84, 107], [9, 54, 76, 91, 101, 107], [2, 8, 21, 66, 88, 103], [3, 48, 70, 85, 95, 101], [0, 35, 40, 46, 51, 74], [21, 56, 61, 67, 72, 95], [2, 15, 60, 82, 97, 107], [5, 18, 63, 85, 100, 110], [17, 54, 89, 94, 100, 105], [30, 52, 67, 77, 83, 96], [1, 16, 26, 32, 45, 90], [3, 7, 14, 60, 65, 73], [14, 51, 86, 91, 97, 102], [2, 7, 13, 18, 41, 78], [3, 26, 63, 98, 103, 109], [1, 11, 17, 30, 75, 97], [10, 20, 26, 39, 84, 106], [20, 25, 31, 36, 59, 96], [42, 47, 55, 96, 100, 107], [33, 38, 46, 87, 91, 98], [12, 16, 23, 69, 74, 82], [23, 28, 34, 39, 62, 99], [2, 10, 51, 55, 62, 108], [12, 17, 25, 66, 70, 77], [24, 46, 61, 71, 77, 90], [29, 34, 40, 45, 68, 105], [19, 34, 44, 50, 63, 108], [7, 48, 52, 59, 105, 110], [18, 40, 55, 65, 71, 84], [5, 51, 56, 64, 105, 109], [33, 37, 44, 90, 95, 103], [6, 28, 43, 53, 59, 72], [32, 37, 43, 48, 71, 108], [4, 10, 15, 38, 75, 110], [0, 22, 37, 47, 53, 66], [0, 23, 60, 95, 100, 106], [7, 17, 23, 36, 81, 103], [1, 7, 12, 35, 72, 107], [6, 41, 46, 52, 57, 80], [6, 51, 73, 88, 98, 104], [12, 34, 49, 59, 65, 78], [12, 47, 52, 58, 63, 86], [2, 48, 53, 61, 102, 106], [3, 38, 43, 49, 54, 77], [33, 55, 70, 80, 86, 99], [27, 31, 38, 84, 89, 97], [7, 22, 32, 38, 51, 96], [2, 39, 74, 79, 85, 90], [5, 10, 16, 21, 44, 81], [3, 25, 40, 50, 56, 69], [5, 42, 77, 82, 88, 93], [24, 29, 37, 78, 82, 89], [30, 35, 43, 84, 88, 95], [6, 10, 17, 63, 68, 76], [15, 37, 52, 62, 68, 81], [4, 45, 49, 56, 102, 107], [9, 13, 20, 66, 71, 79], [21, 43, 58, 68, 74, 87], [39, 44, 52, 93, 97, 104], [20, 57, 92, 97, 103, 108], [36, 71, 76, 82, 87, 110], [36, 40, 47, 93, 98, 106], [4, 19, 29, 35, 48, 93], [18, 53, 58, 64, 69, 92], [5, 11, 24, 69, 91, 106], [27, 32, 40, 81, 85, 92], [42, 64, 79, 89, 95, 108], [12, 57, 79, 94, 104, 110], [13, 28, 38, 44, 57, 102], [18, 23, 31, 72, 76, 83], [13, 23, 29, 42, 87, 109], [1, 8, 54, 59, 67, 108], [9, 31, 46, 56, 62, 75], [8, 14, 27, 72, 94, 109], [45, 50, 58, 99, 103, 110], [30, 34, 41, 87, 92, 100], [15, 50, 55, 61, 66, 89], [3, 8, 16, 57, 61, 68], [8, 13, 19, 24, 47, 84], [15, 20, 28, 69, 73, 80], [26, 31, 37, 42, 65, 102], [1, 6, 29, 66, 101, 106], [15, 19, 26, 72, 77, 85], [8, 45, 80, 85, 91, 96], [27, 62, 67, 73, 78, 101], [0, 45, 67, 82, 92, 98], [11, 48, 83, 88, 94, 99], [16, 31, 41, 47, 60, 105], [1, 42, 46, 53, 99, 104], [36, 41, 49, 90, 94, 101], [4, 9, 32, 69, 104, 109], [36, 58, 73, 83, 89, 102], [24, 28, 35, 81, 86, 94], [6, 11, 19, 60, 64, 71], [18, 22, 29, 75, 80, 88], [27, 49, 64, 74, 80, 93], [19, 37, 46, 54, 69, 79], [22, 41, 44, 61, 65, 107], [32, 58, 77, 80, 97, 101], [16, 35, 38, 55, 59, 101], [15, 27, 41, 98, 102, 108], [12, 22, 73, 91, 100, 108], [3, 21, 33, 47, 104, 108], [9, 21, 35, 92, 96, 102], [41, 45, 51, 69, 81, 95], [4, 23, 26, 43, 47, 89], [53, 57, 63, 81, 93, 107], [5, 8, 25, 29, 71, 97], [4, 13, 21, 36, 46, 97], [8, 12, 18, 36, 48, 62], [13, 32, 35, 52, 56, 98], [11, 14, 31, 35, 77, 103], [4, 22, 31, 39, 54, 64], [44, 48, 54, 72, 84, 98], [4, 55, 73, 82, 90, 105], [35, 39, 45, 63, 75, 89], [17, 43, 62, 65, 82, 86], [29, 55, 74, 77, 94, 98], [2, 19, 23, 65, 91, 110], [28, 46, 55, 63, 78, 88], [38, 64, 83, 86, 103, 107], [1, 20, 23, 40, 44, 86], [2, 28, 47, 50, 67, 71], [7, 15, 30, 40, 91, 109], [9, 23, 80, 84, 90, 108], [16, 34, 43, 51, 66, 76], [34, 52, 61, 69, 84, 94], [17, 21, 27, 45, 57, 71], [7, 58, 76, 85, 93, 108], [0, 14, 71, 75, 81, 99], [47, 51, 57, 75, 87, 101], [0, 12, 26, 83, 87, 93], [0, 15, 25, 76, 94, 103], [26, 30, 36, 54, 66, 80], [6, 21, 31, 82, 100, 109], [19, 38, 41, 58, 62, 104], [0, 18, 30, 44, 101, 105], [2, 44, 70, 89, 92, 109], [43, 61, 70, 78, 93, 103], [13, 22, 30, 45, 55, 106], [7, 11, 53, 79, 98, 101], [31, 49, 58, 66, 81, 91], [23, 27, 33, 51, 63, 77], [12, 24, 38, 95, 99, 105], [41, 67, 86, 89, 106, 110], [7, 16, 24, 39, 49, 100], [49, 67, 76, 84, 99, 109], [1, 52, 70, 79, 87, 102], [32, 36, 42, 60, 72, 86], [4, 12, 27, 37, 88, 106], [14, 18, 24, 42, 54, 68], [29, 33, 39, 57, 69, 83], [17, 20, 37, 41, 83, 109], [38, 42, 48, 66, 78, 92], [13, 17, 59, 85, 104, 107], [11, 37, 56, 59, 76, 80], [37, 55, 64, 72, 87, 97], [10, 19, 27, 42, 52, 103], [11, 15, 21, 39, 51, 65], [26, 52, 71, 74, 91, 95], [3, 17, 74, 78, 84, 102], [14, 17, 34, 38, 80, 106], [1, 10, 18, 33, 43, 94], [5, 31, 50, 53, 70, 74], [35, 61, 80, 83, 100, 104], [7, 26, 29, 46, 50, 92], [22, 40, 49, 57, 72, 82], [25, 44, 47, 64, 68, 110], [20, 46, 65, 68, 85, 89], [7, 25, 34, 42, 57, 67], [3, 13, 64, 82, 91, 99], [3, 18, 28, 79, 97, 106], [46, 64, 73, 81, 96, 106], [4, 8, 50, 76, 95, 98], [5, 9, 15, 33, 45, 59], [20, 24, 30, 48, 60, 74], [16, 20, 62, 88, 107, 110], [25, 43, 52, 60, 75, 85], [1, 5, 47, 73, 92, 95], [14, 40, 59, 62, 79, 83], [11, 68, 72, 78, 96, 108], [6, 18, 32, 89, 93, 99], [56, 60, 66, 84, 96, 110], [3, 9, 27, 39, 53, 110], [40, 58, 67, 75, 90, 100], [1, 9, 24, 34, 85, 103], [8, 34, 53, 56, 73, 77], [0, 10, 61, 79, 88, 96], [13, 31, 40, 48, 63, 73], [8, 65, 69, 75, 93, 105], [10, 28, 37, 45, 60, 70], [2, 59, 63, 69, 87, 99], [5, 62, 66, 72, 90, 102], [0, 6, 24, 36, 50, 107], [10, 14, 56, 82, 101, 104], [2, 5, 22, 26, 68, 94], [1, 19, 28, 36, 51, 61], [23, 49, 68, 71, 88, 92], [8, 11, 28, 32, 74, 100], [6, 20, 77, 81, 87, 105], [3, 15, 29, 86, 90, 96], [6, 16, 67, 85, 94, 102], [9, 19, 70, 88, 97, 105], [50, 54, 60, 78, 90, 104], [10, 29, 32, 49, 53, 95], [2, 6, 12, 30, 42, 56], [16, 25, 33, 48, 58, 109], [2, 16, 37, 40, 77, 104], [28, 40, 54, 87, 96, 103], [12, 45, 54, 61, 97, 109], [18, 39, 47, 59, 77, 102], [2, 29, 38, 52, 73, 76], [8, 22, 43, 46, 83, 110], [0, 21, 29, 41, 59, 84], [10, 13, 50, 77, 86, 100], [35, 62, 71, 85, 106, 109], [1, 13, 27, 60, 69, 76], [8, 33, 60, 81, 89, 101], [6, 15, 22, 58, 70, 84], [6, 39, 48, 55, 91, 103], [21, 48, 69, 77, 89, 107], [12, 21, 28, 64, 76, 90], [7, 21, 54, 63, 70, 106], [14, 23, 37, 58, 61, 98], [5, 17, 35, 60, 87, 108], [24, 45, 53, 65, 83, 108], [8, 26, 51, 78, 99, 107], [10, 22, 36, 69, 78, 85], [26, 53, 62, 76, 97, 100], [24, 51, 72, 80, 92, 110], [13, 34, 37, 74, 101, 110], [9, 42, 51, 58, 94, 106], [19, 22, 59, 86, 95, 109], [20, 29, 43, 64, 67, 104], [3, 11, 23, 41, 66, 93], [5, 23, 48, 75, 96, 104], [9, 18, 25, 61, 73, 87], [1, 15, 48, 57, 64, 100], [24, 33, 40, 76, 88, 102], [23, 50, 59, 73, 94, 97], [5, 32, 41, 55, 76, 79], [0, 7, 43, 55, 69, 102], [13, 16, 53, 80, 89, 103], [13, 25, 39, 72, 81, 88], [30, 39, 46, 82, 94, 108], [3, 36, 45, 52, 88, 100], [2, 14, 32, 57, 84, 105], [12, 20, 32, 50, 75, 102], [11, 20, 34, 55, 58, 95], [4, 18, 51, 60, 67, 103], [4, 40, 52, 66, 99, 108], [0, 8, 20, 38, 63, 90], [7, 28, 31, 68, 95, 104], [29, 56, 65, 79, 100, 103], [12, 39, 60, 68, 80, 98], [3, 24, 32, 44, 62, 87], [5, 19, 40, 43, 80, 107], [7, 10, 47, 74, 83, 97], [12, 33, 41, 53, 71, 96], [14, 41, 50, 64, 85, 88], [8, 17, 31, 52, 55, 92], [6, 27, 35, 47, 65, 90], [32, 59, 68, 82, 103, 106], [15, 42, 63, 71, 83, 101], [18, 45, 66, 74, 86, 104], [34, 46, 60, 93, 102, 109], [0, 33, 42, 49, 85, 97], [6, 33, 54, 62, 74, 92], [7, 19, 33, 66, 75, 82], [9, 17, 29, 47, 72, 99], [5, 14, 28, 49, 52, 89], [21, 42, 50, 62, 80, 105], [11, 36, 63, 84, 92, 104], [15, 36, 44, 56, 74, 99], [21, 30, 37, 73, 85, 99], [6, 14, 26, 44, 69, 96], [2, 20, 45, 72, 93, 101], [4, 16, 30, 63, 72, 79], [23, 32, 46, 67, 70, 107], [0, 27, 48, 56, 68, 86], [27, 36, 43, 79, 91, 105], [3, 10, 46, 58, 72, 105], [4, 7, 44, 71, 80, 94], [9, 36, 57, 65, 77, 95], [11, 38, 47, 61, 82, 85], [25, 37, 51, 84, 93, 100], [3, 12, 19, 55, 67, 81], [2, 11, 25, 46, 49, 86], [19, 31, 45, 78, 87, 94], [8, 35, 44, 58, 79, 82], [10, 31, 34, 71, 98, 107], [3, 30, 51, 59, 71, 89], [1, 38, 65, 74, 88, 109], [11, 29, 54, 81, 102, 110], [14, 39, 66, 87, 95, 107], [1, 22, 25, 62, 89, 98], [6, 13, 49, 61, 75, 108], [10, 24, 57, 66, 73, 109], [26, 35, 49, 70, 73, 110], [31, 43, 57, 90, 99, 106], [4, 25, 28, 65, 92, 101], [15, 23, 35, 53, 78, 105], [17, 44, 53, 67, 88, 91], [15, 24, 31, 67, 79, 93], [17, 42, 69, 90, 98, 110], [17, 26, 40, 61, 64, 101], [22, 34, 48, 81, 90, 97], [1, 37, 49, 63, 96, 105], [16, 19, 56, 83, 92, 106], [5, 30, 57, 78, 86, 98], [9, 30, 38, 50, 68, 93], [0, 9, 16, 52, 64, 78], [16, 28, 42, 75, 84, 91], [2, 27, 54, 75, 83, 95], [18, 27, 34, 70, 82, 96], [20, 47, 56, 70, 91, 94], [18, 26, 38, 56, 81, 108], [1, 4, 41, 68, 77, 91], [21, 38, 40, 53, 60, 94], [31, 69, 86, 88, 101, 108], [8, 10, 23, 30, 64, 102], [11, 13, 26, 33, 67, 105], [6, 40, 78, 95, 97, 110], [25, 63, 80, 82, 95, 102], [22, 60, 77, 79, 92, 99], [8, 15, 49, 87, 104, 106], [2, 4, 17, 24, 58, 96], [13, 51, 68, 70, 83, 90], [1, 14, 21, 55, 93, 110], [18, 35, 37, 50, 57, 91], [3, 37, 75, 92, 94, 107], [15, 32, 34, 47, 54, 88], [36, 53, 55, 68, 75, 109], [30, 47, 49, 62, 69, 103], [11, 18, 52, 90, 107, 109], [4, 42, 59, 61, 74, 81], [19, 57, 74, 76, 89, 96], [24, 41, 43, 56, 63, 97], [6, 23, 25, 38, 45, 79], [2, 9, 43, 81, 98, 100], [7, 45, 62, 64, 77, 84], [14, 16, 29, 36, 70, 108], [5, 12, 46, 84, 101, 103], [5, 7, 20, 27, 61, 99], [12, 29, 31, 44, 51, 85], [27, 44, 46, 59, 66, 100], [0, 34, 72, 89, 91, 104], [16, 54, 71, 73, 86, 93], [9, 26, 28, 41, 48, 82], [28, 66, 83, 85, 98, 105], [10, 48, 65, 67, 80, 87], [33, 50, 52, 65, 72, 106], [3, 20, 22, 35, 42, 76], [0, 17, 19, 32, 39, 73], [1, 39, 56, 58, 71, 78]]
\item 1 \{1=28416, 2=510156, 3=2203128, 4=2386500\} [[8, 39, 40, 41, 42, 70], [0, 1, 2, 3, 31, 80], [30, 31, 32, 33, 61, 110], [6, 7, 8, 9, 37, 86], [23, 54, 55, 56, 57, 85], [15, 16, 17, 18, 46, 95], [27, 28, 29, 30, 58, 107], [24, 25, 26, 27, 55, 104], [26, 57, 58, 59, 60, 88], [29, 60, 61, 62, 63, 91], [25, 74, 105, 106, 107, 108], [41, 72, 73, 74, 75, 103], [7, 56, 87, 88, 89, 90], [10, 59, 90, 91, 92, 93], [13, 62, 93, 94, 95, 96], [22, 71, 102, 103, 104, 105], [9, 10, 11, 12, 40, 89], [44, 75, 76, 77, 78, 106], [18, 19, 20, 21, 49, 98], [11, 42, 43, 44, 45, 73], [3, 4, 5, 6, 34, 83], [16, 65, 96, 97, 98, 99], [19, 68, 99, 100, 101, 102], [17, 48, 49, 50, 51, 79], [47, 78, 79, 80, 81, 109], [20, 51, 52, 53, 54, 82], [12, 13, 14, 15, 43, 92], [0, 28, 77, 108, 109, 110], [4, 53, 84, 85, 86, 87], [1, 50, 81, 82, 83, 84], [32, 63, 64, 65, 66, 94], [2, 33, 34, 35, 36, 64], [21, 22, 23, 24, 52, 101], [38, 69, 70, 71, 72, 100], [14, 45, 46, 47, 48, 76], [35, 66, 67, 68, 69, 97], [5, 36, 37, 38, 39, 67], [6, 41, 79, 93, 97, 110], [6, 10, 23, 30, 65, 103], [26, 64, 78, 82, 95, 102], [8, 15, 50, 88, 102, 106], [9, 13, 26, 33, 68, 106], [12, 16, 29, 36, 71, 109], [22, 36, 40, 53, 60, 95], [4, 18, 22, 35, 42, 77], [5, 12, 47, 85, 99, 103], [1, 14, 21, 56, 94, 108], [1, 15, 19, 32, 39, 74], [11, 18, 53, 91, 105, 109], [17, 55, 69, 73, 86, 93], [31, 45, 49, 62, 69, 104], [5, 43, 57, 61, 74, 81], [19, 33, 37, 50, 57, 92], [34, 48, 52, 65, 72, 107], [20, 58, 72, 76, 89, 96], [3, 7, 20, 27, 62, 100], [25, 39, 43, 56, 63, 98], [2, 9, 44, 82, 96, 100], [23, 61, 75, 79, 92, 99], [37, 51, 55, 68, 75, 110], [32, 70, 84, 88, 101, 108], [16, 30, 34, 47, 54, 89], [0, 35, 73, 87, 91, 104], [8, 46, 60, 64, 77, 84], [7, 21, 25, 38, 45, 80], [10, 24, 28, 41, 48, 83], [0, 4, 17, 24, 59, 97], [28, 42, 46, 59, 66, 101], [29, 67, 81, 85, 98, 105], [11, 49, 63, 67, 80, 87], [3, 38, 76, 90, 94, 107], [14, 52, 66, 70, 83, 90], [13, 27, 31, 44, 51, 86], [2, 40, 54, 58, 71, 78], [10, 62, 67, 70, 76, 99], [26, 45, 86, 92, 103, 107], [22, 27, 33, 38, 66, 78], [38, 44, 55, 59, 89, 108], [14, 19, 22, 28, 51, 73], [12, 34, 86, 91, 94, 100], [0, 5, 33, 45, 100, 105], [6, 18, 73, 78, 84, 89], [2, 13, 17, 47, 66, 107], [9, 50, 56, 67, 71, 101], [8, 13, 16, 22, 45, 67], [41, 46, 49, 55, 78, 100], [3, 58, 63, 69, 74, 102], [26, 31, 34, 40, 63, 85], [16, 21, 27, 32, 60, 72], [29, 34, 37, 43, 66, 88], [31, 36, 42, 47, 75, 87], [2, 7, 10, 16, 39, 61], [52, 57, 63, 68, 96, 108], [0, 22, 74, 79, 82, 88], [23, 28, 31, 37, 60, 82], [15, 37, 89, 94, 97, 103], [1, 4, 10, 33, 55, 107], [17, 23, 34, 38, 68, 87], [3, 8, 36, 48, 103, 108], [1, 5, 35, 54, 95, 101], [34, 39, 45, 50, 78, 90], [49, 54, 60, 65, 93, 105], [11, 30, 71, 77, 88, 92], [19, 24, 30, 35, 63, 75], [23, 29, 40, 44, 74, 93], [4, 7, 13, 36, 58, 110], [29, 35, 46, 50, 80, 99], [7, 59, 64, 67, 73, 96], [27, 39, 94, 99, 105, 110], [5, 16, 20, 50, 69, 110], [47, 52, 55, 61, 84, 106], [20, 39, 80, 86, 97, 101], [15, 27, 82, 87, 93, 98], [5, 10, 13, 19, 42, 64], [23, 42, 83, 89, 100, 104], [19, 71, 76, 79, 85, 108], [7, 11, 41, 60, 101, 107], [5, 11, 22, 26, 56, 75], [2, 8, 19, 23, 53, 72], [5, 24, 65, 71, 82, 86], [10, 15, 21, 26, 54, 66], [28, 33, 39, 44, 72, 84], [44, 49, 52, 58, 81, 103], [3, 9, 14, 42, 54, 109], [4, 8, 38, 57, 98, 104], [6, 47, 53, 64, 68, 98], [38, 43, 46, 52, 75, 97], [18, 30, 85, 90, 96, 101], [2, 30, 42, 97, 102, 108], [1, 24, 46, 98, 103, 106], [6, 61, 66, 72, 77, 105], [35, 40, 43, 49, 72, 94], [35, 41, 52, 56, 86, 105], [1, 7, 30, 52, 104, 109], [43, 48, 54, 59, 87, 99], [2, 32, 51, 92, 98, 109], [3, 25, 77, 82, 85, 91], [13, 18, 24, 29, 57, 69], [37, 42, 48, 53, 81, 93], [11, 17, 28, 32, 62, 81], [6, 28, 80, 85, 88, 94], [9, 21, 76, 81, 87, 92], [2, 21, 62, 68, 79, 83], [3, 15, 70, 75, 81, 86], [4, 9, 15, 20, 48, 60], [18, 40, 92, 97, 100, 106], [21, 43, 95, 100, 103, 109], [7, 12, 18, 23, 51, 63], [14, 20, 31, 35, 65, 84], [50, 55, 58, 64, 87, 109], [18, 59, 65, 76, 80, 110], [32, 38, 49, 53, 83, 102], [3, 44, 50, 61, 65, 95], [13, 65, 70, 73, 79, 102], [8, 14, 25, 29, 59, 78], [10, 14, 44, 63, 104, 110], [9, 31, 83, 88, 91, 97], [20, 25, 28, 34, 57, 79], [40, 45, 51, 56, 84, 96], [15, 56, 62, 73, 77, 107], [20, 26, 37, 41, 71, 90], [25, 30, 36, 41, 69, 81], [12, 24, 79, 84, 90, 95], [0, 6, 11, 39, 51, 106], [29, 48, 89, 95, 106, 110], [12, 53, 59, 70, 74, 104], [0, 12, 67, 72, 78, 83], [9, 64, 69, 75, 80, 108], [32, 37, 40, 46, 69, 91], [1, 53, 58, 61, 67, 90], [21, 33, 88, 93, 99, 104], [0, 55, 60, 66, 71, 99], [0, 41, 47, 58, 62, 92], [4, 27, 49, 101, 106, 109], [8, 27, 68, 74, 85, 89], [26, 32, 43, 47, 77, 96], [4, 56, 61, 64, 70, 93], [16, 68, 73, 76, 82, 105], [24, 36, 91, 96, 102, 107], [17, 36, 77, 83, 94, 98], [14, 33, 74, 80, 91, 95], [11, 16, 19, 25, 48, 70], [1, 6, 12, 17, 45, 57], [46, 51, 57, 62, 90, 102], [17, 22, 25, 31, 54, 76], [13, 28, 49, 61, 71, 89], [2, 5, 48, 77, 90, 104], [33, 60, 73, 81, 90, 97], [7, 22, 43, 55, 65, 83], [8, 11, 54, 83, 96, 110], [2, 11, 14, 57, 86, 99], [11, 24, 38, 47, 50, 93], [2, 37, 52, 73, 85, 95], [4, 12, 21, 28, 75, 102], [1, 11, 29, 64, 79, 100], [6, 13, 60, 87, 100, 108], [12, 41, 54, 68, 77, 80], [18, 47, 60, 74, 83, 86], [27, 54, 67, 75, 84, 91], [1, 48, 75, 88, 96, 105], [5, 23, 58, 73, 94, 106], [33, 62, 75, 89, 98, 101], [22, 37, 58, 70, 80, 98], [24, 53, 66, 80, 89, 92], [5, 40, 55, 76, 88, 98], [45, 72, 85, 93, 102, 109], [7, 17, 35, 70, 85, 106], [0, 13, 21, 30, 37, 84], [13, 34, 46, 56, 74, 109], [6, 33, 46, 54, 63, 70], [24, 51, 64, 72, 81, 88], [3, 30, 43, 51, 60, 67], [19, 31, 41, 59, 94, 109], [19, 34, 55, 67, 77, 95], [5, 18, 32, 41, 44, 87], [34, 49, 70, 82, 92, 110], [12, 25, 33, 42, 49, 96], [4, 51, 78, 91, 99, 108], [10, 18, 27, 34, 81, 108], [3, 12, 19, 66, 93, 106], [18, 45, 58, 66, 75, 82], [39, 66, 79, 87, 96, 103], [8, 21, 35, 44, 47, 90], [7, 19, 29, 47, 82, 97], [23, 36, 50, 59, 62, 105], [7, 15, 24, 31, 78, 105], [6, 20, 29, 32, 75, 104], [15, 42, 55, 63, 72, 79], [4, 14, 32, 67, 82, 103], [39, 68, 81, 95, 104, 107], [42, 71, 84, 98, 107, 110], [5, 8, 51, 80, 93, 107], [4, 19, 40, 52, 62, 80], [9, 23, 32, 35, 78, 107], [15, 44, 57, 71, 80, 83], [21, 50, 63, 77, 86, 89], [0, 7, 54, 81, 94, 102], [24, 37, 45, 54, 61, 108], [3, 16, 24, 33, 40, 87], [3, 32, 45, 59, 68, 71], [6, 15, 22, 69, 96, 109], [10, 20, 38, 73, 88, 109], [1, 9, 18, 25, 72, 99], [9, 38, 51, 65, 74, 77], [1, 13, 23, 41, 76, 91], [12, 26, 35, 38, 81, 110], [0, 9, 16, 63, 90, 103], [30, 57, 70, 78, 87, 94], [1, 16, 37, 49, 59, 77], [0, 14, 23, 26, 69, 98], [28, 43, 64, 76, 86, 104], [30, 59, 72, 86, 95, 98], [6, 19, 27, 36, 43, 90], [9, 22, 30, 39, 46, 93], [9, 36, 49, 57, 66, 73], [2, 20, 55, 70, 91, 103], [4, 16, 26, 44, 79, 94], [12, 39, 52, 60, 69, 76], [10, 25, 46, 58, 68, 86], [18, 31, 39, 48, 55, 102], [16, 28, 38, 56, 91, 106], [10, 31, 43, 53, 71, 106], [14, 49, 64, 85, 97, 107], [26, 39, 53, 62, 65, 108], [5, 14, 17, 60, 89, 102], [14, 27, 41, 50, 53, 96], [3, 10, 57, 84, 97, 105], [4, 25, 37, 47, 65, 100], [10, 22, 32, 50, 85, 100], [7, 28, 40, 50, 68, 103], [11, 20, 23, 66, 95, 108], [15, 28, 36, 45, 52, 99], [1, 22, 34, 44, 62, 97], [20, 33, 47, 56, 59, 102], [31, 46, 67, 79, 89, 107], [8, 26, 61, 76, 97, 109], [16, 31, 52, 64, 74, 92], [11, 46, 61, 82, 94, 104], [36, 63, 76, 84, 93, 100], [36, 65, 78, 92, 101, 104], [17, 30, 44, 53, 56, 99], [2, 45, 74, 87, 101, 110], [13, 25, 35, 53, 88, 103], [8, 43, 58, 79, 91, 101], [2, 15, 29, 38, 41, 84], [21, 34, 42, 51, 58, 105], [6, 35, 48, 62, 71, 74], [3, 17, 26, 29, 72, 101], [25, 40, 61, 73, 83, 101], [27, 56, 69, 83, 92, 95], [42, 69, 82, 90, 99, 106], [8, 17, 20, 63, 92, 105], [0, 29, 42, 56, 65, 68], [17, 52, 67, 88, 100, 110], [0, 27, 40, 48, 57, 64], [21, 48, 61, 69, 78, 85], [5, 21, 29, 31, 70, 96], [3, 22, 41, 64, 89, 99], [1, 40, 66, 86, 102, 110], [7, 26, 49, 74, 84, 99], [6, 14, 16, 55, 81, 101], [22, 47, 57, 72, 91, 110], [9, 28, 47, 70, 95, 105], [8, 18, 33, 52, 71, 94], [24, 44, 60, 68, 70, 109], [28, 54, 74, 90, 98, 100], [0, 18, 43, 50, 70, 107], [14, 37, 62, 72, 87, 106], [5, 7, 46, 72, 92, 108], [21, 41, 57, 65, 67, 106], [16, 53, 57, 75, 100, 107], [0, 25, 32, 52, 89, 93], [11, 15, 33, 58, 65, 85], [2, 12, 27, 46, 65, 88], [2, 18, 26, 28, 67, 93], [20, 30, 45, 64, 83, 106], [20, 24, 42, 67, 74, 94], [7, 14, 34, 71, 75, 93], [25, 51, 71, 87, 95, 97], [35, 39, 57, 82, 89, 109], [18, 38, 54, 62, 64, 103], [19, 45, 65, 81, 89, 91], [22, 48, 68, 84, 92, 94], [13, 20, 40, 77, 81, 99], [34, 60, 80, 96, 104, 106], [2, 25, 50, 60, 75, 94], [3, 28, 35, 55, 92, 96], [14, 18, 36, 61, 68, 88], [5, 28, 53, 63, 78, 97], [19, 26, 46, 83, 87, 105], [8, 12, 30, 55, 62, 82], [3, 11, 13, 52, 78, 98], [2, 4, 43, 69, 89, 105], [15, 23, 25, 64, 90, 110], [5, 9, 27, 52, 59, 79], [14, 24, 39, 58, 77, 100], [19, 44, 54, 69, 88, 107], [7, 44, 48, 66, 91, 98], [7, 32, 42, 57, 76, 95], [31, 57, 77, 93, 101, 103], [8, 24, 32, 34, 73, 99], [13, 32, 55, 80, 90, 105], [6, 21, 40, 59, 82, 107], [13, 39, 59, 75, 83, 85], [1, 27, 47, 63, 71, 73], [10, 29, 52, 77, 87, 102], [17, 27, 42, 61, 80, 103], [7, 33, 53, 69, 77, 79], [29, 33, 51, 76, 83, 103], [9, 34, 41, 61, 98, 102], [9, 17, 19, 58, 84, 104], [13, 38, 48, 63, 82, 101], [23, 27, 45, 70, 77, 97], [14, 30, 38, 40, 79, 105], [0, 15, 34, 53, 76, 101], [0, 19, 38, 61, 86, 96], [10, 17, 37, 74, 78, 96], [16, 42, 62, 78, 86, 88], [16, 35, 58, 83, 93, 108], [1, 38, 42, 60, 85, 92], [1, 8, 28, 65, 69, 87], [37, 63, 83, 99, 107, 109], [9, 29, 45, 53, 55, 94], [16, 23, 43, 80, 84, 102], [3, 21, 46, 53, 73, 110], [10, 36, 56, 72, 80, 82], [4, 41, 45, 63, 88, 95], [5, 25, 62, 66, 84, 109], [19, 56, 60, 78, 103, 110], [10, 35, 45, 60, 79, 98], [4, 23, 46, 71, 81, 96], [12, 31, 50, 73, 98, 108], [17, 33, 41, 43, 82, 108], [4, 11, 31, 68, 72, 90], [12, 32, 48, 56, 58, 97], [12, 37, 44, 64, 101, 105], [2, 6, 24, 49, 56, 76], [10, 47, 51, 69, 94, 101], [6, 25, 44, 67, 92, 102], [0, 20, 36, 44, 46, 85], [6, 31, 38, 58, 95, 99], [4, 30, 50, 66, 74, 76], [15, 40, 47, 67, 104, 108], [11, 34, 59, 69, 84, 103], [15, 35, 51, 59, 61, 100], [32, 36, 54, 79, 86, 106], [11, 21, 36, 55, 74, 97], [1, 26, 36, 51, 70, 89], [8, 31, 56, 66, 81, 100], [23, 33, 48, 67, 86, 109], [11, 27, 35, 37, 76, 102], [17, 21, 39, 64, 71, 91], [0, 8, 10, 49, 75, 95], [9, 24, 43, 62, 85, 110], [3, 18, 37, 56, 79, 104], [6, 26, 42, 50, 52, 91], [12, 20, 22, 61, 87, 107], [4, 29, 39, 54, 73, 92], [26, 30, 48, 73, 80, 100], [16, 41, 51, 66, 85, 104], [3, 23, 39, 47, 49, 88], [5, 15, 30, 49, 68, 91], [13, 50, 54, 72, 97, 104], [17, 40, 65, 75, 90, 109], [1, 20, 43, 68, 78, 93], [2, 22, 59, 63, 81, 106], [22, 29, 49, 86, 90, 108]]
\item 1 \{1=32856, 2=535464, 3=2217336, 4=2342544\} [[8, 39, 40, 41, 42, 70], [0, 1, 2, 3, 31, 80], [30, 31, 32, 33, 61, 110], [6, 7, 8, 9, 37, 86], [23, 54, 55, 56, 57, 85], [15, 16, 17, 18, 46, 95], [27, 28, 29, 30, 58, 107], [24, 25, 26, 27, 55, 104], [26, 57, 58, 59, 60, 88], [29, 60, 61, 62, 63, 91], [25, 74, 105, 106, 107, 108], [41, 72, 73, 74, 75, 103], [7, 56, 87, 88, 89, 90], [10, 59, 90, 91, 92, 93], [13, 62, 93, 94, 95, 96], [22, 71, 102, 103, 104, 105], [9, 10, 11, 12, 40, 89], [44, 75, 76, 77, 78, 106], [18, 19, 20, 21, 49, 98], [11, 42, 43, 44, 45, 73], [3, 4, 5, 6, 34, 83], [16, 65, 96, 97, 98, 99], [19, 68, 99, 100, 101, 102], [17, 48, 49, 50, 51, 79], [47, 78, 79, 80, 81, 109], [20, 51, 52, 53, 54, 82], [12, 13, 14, 15, 43, 92], [0, 28, 77, 108, 109, 110], [4, 53, 84, 85, 86, 87], [1, 50, 81, 82, 83, 84], [32, 63, 64, 65, 66, 94], [2, 33, 34, 35, 36, 64], [21, 22, 23, 24, 52, 101], [38, 69, 70, 71, 72, 100], [14, 45, 46, 47, 48, 76], [35, 66, 67, 68, 69, 97], [5, 36, 37, 38, 39, 67], [26, 54, 70, 102, 106, 110], [6, 22, 54, 58, 62, 89], [9, 13, 17, 44, 72, 88], [6, 10, 14, 41, 69, 85], [7, 39, 43, 47, 74, 102], [0, 16, 48, 52, 56, 83], [21, 25, 29, 56, 84, 100], [21, 37, 69, 73, 77, 104], [17, 45, 61, 93, 97, 101], [8, 36, 52, 84, 88, 92], [1, 33, 37, 41, 68, 96], [23, 51, 67, 99, 103, 107], [27, 31, 35, 62, 90, 106], [2, 29, 57, 73, 105, 109], [5, 33, 49, 81, 85, 89], [11, 39, 55, 87, 91, 95], [24, 40, 72, 76, 80, 107], [15, 19, 23, 50, 78, 94], [10, 42, 46, 50, 77, 105], [15, 31, 63, 67, 71, 98], [9, 25, 57, 61, 65, 92], [0, 4, 8, 35, 63, 79], [27, 43, 75, 79, 83, 110], [4, 36, 40, 44, 71, 99], [14, 42, 58, 90, 94, 98], [12, 28, 60, 64, 68, 95], [18, 22, 26, 53, 81, 97], [18, 34, 66, 70, 74, 101], [20, 48, 64, 96, 100, 104], [13, 45, 49, 53, 80, 108], [24, 28, 32, 59, 87, 103], [12, 16, 20, 47, 75, 91], [2, 30, 46, 78, 82, 86], [3, 7, 11, 38, 66, 82], [3, 19, 51, 55, 59, 86], [1, 5, 32, 60, 76, 108], [30, 34, 38, 65, 93, 109], [10, 25, 31, 53, 60, 101], [23, 43, 58, 64, 86, 93], [9, 46, 54, 59, 63, 73], [4, 24, 38, 50, 60, 89], [17, 37, 52, 58, 80, 87], [22, 42, 56, 68, 78, 107], [2, 28, 48, 62, 74, 84], [2, 6, 16, 63, 100, 108], [3, 8, 12, 22, 69, 106], [12, 26, 38, 48, 77, 103], [0, 10, 57, 94, 102, 107], [1, 16, 22, 44, 51, 92], [23, 49, 69, 83, 95, 105], [11, 37, 57, 71, 83, 93], [4, 51, 88, 96, 101, 105], [19, 27, 32, 36, 46, 93], [13, 33, 47, 59, 69, 98], [29, 49, 64, 70, 92, 99], [34, 42, 47, 51, 61, 108], [25, 33, 38, 42, 52, 99], [6, 43, 51, 56, 60, 70], [11, 31, 46, 52, 74, 81], [13, 28, 34, 56, 63, 104], [3, 40, 48, 53, 57, 67], [5, 12, 53, 73, 88, 94], [42, 79, 87, 92, 96, 106], [21, 58, 66, 71, 75, 85], [25, 45, 59, 71, 81, 110], [16, 36, 50, 62, 72, 101], [8, 18, 47, 73, 93, 107], [9, 23, 35, 45, 74, 100], [31, 39, 44, 48, 58, 105], [4, 10, 32, 39, 80, 100], [1, 7, 29, 36, 77, 97], [28, 36, 41, 45, 55, 102], [4, 19, 25, 47, 54, 95], [3, 17, 29, 39, 68, 94], [7, 22, 28, 50, 57, 98], [1, 9, 14, 18, 28, 75], [6, 35, 61, 81, 95, 107], [1, 48, 85, 93, 98, 102], [35, 55, 70, 76, 98, 105], [11, 23, 33, 62, 88, 108], [8, 28, 43, 49, 71, 78], [0, 29, 55, 75, 89, 101], [36, 73, 81, 86, 90, 100], [12, 49, 57, 62, 66, 76], [14, 40, 60, 74, 86, 96], [7, 54, 91, 99, 104, 108], [8, 34, 54, 68, 80, 90], [17, 24, 65, 85, 100, 106], [6, 11, 15, 25, 72, 109], [4, 12, 17, 21, 31, 78], [9, 38, 64, 84, 98, 110], [3, 44, 64, 79, 85, 107], [6, 47, 67, 82, 88, 110], [19, 39, 53, 65, 75, 104], [11, 21, 50, 76, 96, 110], [7, 27, 41, 53, 63, 92], [38, 58, 73, 79, 101, 108], [20, 27, 68, 88, 103, 109], [5, 17, 27, 56, 82, 102], [13, 21, 26, 30, 40, 87], [19, 34, 40, 62, 69, 110], [18, 32, 44, 54, 83, 109], [14, 21, 62, 82, 97, 103], [11, 18, 59, 79, 94, 100], [2, 14, 24, 53, 79, 99], [10, 18, 23, 27, 37, 84], [7, 13, 35, 42, 83, 103], [5, 25, 40, 46, 68, 75], [20, 40, 55, 61, 83, 90], [16, 31, 37, 59, 66, 107], [0, 5, 9, 19, 66, 103], [7, 15, 20, 24, 34, 81], [3, 13, 60, 97, 105, 110], [27, 64, 72, 77, 81, 91], [45, 82, 90, 95, 99, 109], [22, 30, 35, 39, 49, 96], [10, 30, 44, 56, 66, 95], [15, 29, 41, 51, 80, 106], [33, 70, 78, 83, 87, 97], [8, 15, 56, 76, 91, 97], [26, 46, 61, 67, 89, 96], [8, 20, 30, 59, 85, 105], [18, 55, 63, 68, 72, 82], [39, 76, 84, 89, 93, 103], [6, 20, 32, 42, 71, 97], [0, 14, 26, 36, 65, 91], [1, 21, 35, 47, 57, 86], [32, 52, 67, 73, 95, 102], [30, 67, 75, 80, 84, 94], [3, 32, 58, 78, 92, 104], [13, 19, 41, 48, 89, 109], [5, 31, 51, 65, 77, 87], [2, 22, 37, 43, 65, 72], [16, 24, 29, 33, 43, 90], [24, 61, 69, 74, 78, 88], [2, 9, 50, 70, 85, 91], [10, 16, 38, 45, 86, 106], [15, 52, 60, 65, 69, 79], [0, 41, 61, 76, 82, 104], [2, 12, 41, 67, 87, 101], [14, 34, 49, 55, 77, 84], [4, 26, 33, 74, 94, 109], [5, 15, 44, 70, 90, 104], [26, 52, 72, 86, 98, 108], [17, 43, 63, 77, 89, 99], [20, 46, 66, 80, 92, 102], [0, 37, 45, 50, 54, 64], [1, 23, 30, 71, 91, 106], [32, 43, 82, 98, 101, 106], [14, 17, 22, 59, 70, 109], [3, 9, 20, 33, 77, 95], [2, 5, 10, 47, 58, 97], [25, 36, 63, 78, 85, 103], [15, 21, 32, 45, 89, 107], [1, 19, 52, 63, 90, 105], [35, 53, 72, 78, 89, 102], [8, 27, 33, 44, 57, 101], [12, 19, 37, 70, 81, 108], [1, 40, 56, 59, 64, 101], [6, 33, 48, 55, 73, 106], [18, 24, 35, 48, 92, 110], [11, 29, 48, 54, 65, 78], [11, 22, 61, 77, 80, 85], [2, 21, 27, 38, 51, 95], [15, 30, 37, 55, 88, 99], [2, 20, 39, 45, 56, 69], [12, 27, 34, 52, 85, 96], [5, 8, 13, 50, 61, 100], [20, 31, 70, 86, 89, 94], [17, 35, 54, 60, 71, 84], [19, 30, 57, 72, 79, 97], [8, 26, 45, 51, 62, 75], [8, 19, 58, 74, 77, 82], [13, 24, 51, 66, 73, 91], [16, 27, 54, 69, 76, 94], [5, 23, 42, 48, 59, 72], [2, 13, 52, 68, 71, 76], [3, 14, 27, 71, 89, 108], [0, 15, 22, 40, 73, 84], [2, 7, 44, 55, 94, 110], [4, 41, 52, 91, 107, 110], [8, 11, 16, 53, 64, 103], [25, 41, 44, 49, 86, 97], [18, 33, 40, 58, 91, 102], [24, 39, 46, 64, 97, 108], [3, 18, 25, 43, 76, 87], [37, 53, 56, 61, 98, 109], [23, 41, 60, 66, 77, 90], [29, 40, 79, 95, 98, 103], [28, 39, 66, 81, 88, 106], [9, 36, 51, 58, 76, 109], [31, 47, 50, 55, 92, 103], [17, 28, 67, 83, 86, 91], [22, 38, 41, 46, 83, 94], [21, 36, 43, 61, 94, 105], [9, 16, 34, 67, 78, 105], [34, 50, 53, 58, 95, 106], [5, 24, 30, 41, 54, 98], [0, 44, 62, 81, 87, 98], [6, 21, 28, 46, 79, 90], [6, 50, 68, 87, 93, 104], [2, 15, 59, 77, 96, 102], [9, 15, 26, 39, 83, 101], [12, 18, 29, 42, 86, 104], [10, 43, 54, 81, 96, 103], [0, 7, 25, 58, 69, 96], [1, 38, 49, 88, 104, 107], [19, 35, 38, 43, 80, 91], [38, 56, 75, 81, 92, 105], [9, 24, 31, 49, 82, 93], [7, 46, 62, 65, 70, 107], [3, 30, 45, 52, 70, 103], [4, 37, 48, 75, 90, 97], [0, 6, 17, 30, 74, 92], [11, 30, 36, 47, 60, 104], [9, 53, 71, 90, 96, 107], [4, 15, 42, 57, 64, 82], [7, 40, 51, 78, 93, 100], [17, 36, 42, 53, 66, 110], [0, 11, 24, 68, 86, 105], [26, 37, 76, 92, 95, 100], [32, 50, 69, 75, 86, 99], [10, 26, 29, 34, 71, 82], [4, 43, 59, 62, 67, 104], [1, 12, 39, 54, 61, 79], [5, 16, 55, 71, 74, 79], [23, 34, 73, 89, 92, 97], [3, 10, 28, 61, 72, 99], [6, 13, 31, 64, 75, 102], [4, 20, 23, 28, 65, 76], [3, 47, 65, 84, 90, 101], [0, 27, 42, 49, 67, 100], [16, 32, 35, 40, 77, 88], [11, 14, 19, 56, 67, 106], [22, 33, 60, 75, 82, 100], [12, 56, 74, 93, 99, 110], [28, 44, 47, 52, 89, 100], [1, 17, 20, 25, 62, 73], [7, 23, 26, 31, 68, 79], [29, 47, 66, 72, 83, 96], [41, 59, 78, 84, 95, 108], [13, 46, 57, 84, 99, 106], [14, 25, 64, 80, 83, 88], [5, 18, 62, 80, 99, 105], [7, 18, 45, 60, 67, 85], [16, 49, 60, 87, 102, 109], [13, 29, 32, 37, 74, 85], [26, 44, 63, 69, 80, 93], [8, 21, 65, 83, 102, 108], [10, 49, 65, 68, 73, 110], [31, 42, 69, 84, 91, 109], [14, 32, 51, 57, 68, 81], [20, 38, 57, 63, 74, 87], [4, 22, 55, 66, 93, 108], [6, 12, 23, 36, 80, 98], [1, 34, 45, 72, 87, 94], [14, 33, 39, 50, 63, 107], [10, 21, 48, 63, 70, 88], [35, 46, 85, 101, 104, 109], [1, 42, 65, 74, 89, 95], [0, 21, 34, 39, 59, 99], [10, 13, 20, 22, 36, 79], [1, 8, 10, 24, 67, 109], [28, 70, 73, 80, 82, 96], [25, 67, 70, 77, 79, 93], [22, 64, 67, 74, 76, 90], [11, 51, 63, 84, 97, 102], [0, 20, 60, 72, 93, 106], [15, 28, 33, 53, 93, 105], [4, 9, 29, 69, 81, 102], [4, 7, 14, 16, 30, 73], [33, 45, 66, 79, 84, 104], [8, 17, 32, 38, 55, 96], [6, 49, 91, 94, 101, 103], [0, 13, 18, 38, 78, 90], [6, 29, 38, 53, 59, 76], [10, 15, 35, 75, 87, 108], [1, 4, 11, 13, 27, 70], [11, 20, 35, 41, 58, 99], [9, 52, 94, 97, 104, 106], [14, 20, 37, 78, 101, 110], [30, 53, 62, 77, 83, 100], [3, 23, 63, 75, 96, 109], [19, 61, 64, 71, 73, 87], [2, 19, 60, 83, 92, 107], [2, 8, 25, 66, 89, 98], [34, 76, 79, 86, 88, 102], [16, 58, 61, 68, 70, 84], [28, 31, 38, 40, 54, 97], [18, 41, 50, 65, 71, 88], [27, 39, 60, 73, 78, 98], [14, 23, 38, 44, 61, 102], [7, 10, 17, 19, 33, 76], [31, 34, 41, 43, 57, 100], [13, 55, 58, 65, 67, 81], [24, 47, 56, 71, 77, 94], [12, 55, 97, 100, 107, 109], [39, 51, 72, 85, 90, 110], [17, 57, 69, 90, 103, 108], [25, 28, 35, 37, 51, 94], [2, 11, 26, 32, 49, 90], [6, 18, 39, 52, 57, 77], [5, 14, 29, 35, 52, 93], [0, 23, 32, 47, 53, 70], [7, 49, 52, 59, 61, 75], [36, 48, 69, 82, 87, 107], [7, 12, 32, 72, 84, 105], [1, 43, 46, 53, 55, 69], [19, 22, 29, 31, 45, 88], [12, 25, 30, 50, 90, 102], [11, 17, 34, 75, 98, 107], [1, 15, 58, 100, 103, 110], [39, 62, 71, 86, 92, 109], [22, 25, 32, 34, 48, 91], [9, 22, 27, 47, 87, 99], [6, 27, 40, 45, 65, 105], [6, 19, 24, 44, 84, 96], [17, 26, 41, 47, 64, 105], [14, 54, 66, 87, 100, 105], [34, 37, 44, 46, 60, 103], [18, 30, 51, 64, 69, 89], [0, 12, 33, 46, 51, 71], [9, 21, 42, 55, 60, 80], [5, 20, 26, 43, 84, 107], [5, 22, 63, 86, 95, 110], [5, 11, 28, 69, 92, 101], [3, 15, 36, 49, 54, 74], [37, 79, 82, 89, 91, 105], [5, 45, 57, 78, 91, 96], [21, 33, 54, 67, 72, 92], [20, 29, 44, 50, 67, 108], [40, 43, 50, 52, 66, 109], [15, 27, 48, 61, 66, 86], [15, 38, 47, 62, 68, 85], [10, 52, 55, 62, 64, 78], [3, 16, 21, 41, 81, 93], [3, 26, 35, 50, 56, 73], [27, 50, 59, 74, 80, 97], [12, 35, 44, 59, 65, 82], [33, 56, 65, 80, 86, 103], [4, 45, 68, 77, 92, 98], [7, 48, 71, 80, 95, 101], [0, 43, 85, 88, 95, 97], [5, 7, 21, 64, 106, 109], [10, 51, 74, 83, 98, 104], [8, 14, 31, 72, 95, 104], [8, 23, 29, 46, 87, 110], [36, 59, 68, 83, 89, 106], [3, 46, 88, 91, 98, 100], [1, 6, 26, 66, 78, 99], [12, 24, 45, 58, 63, 83], [18, 31, 36, 56, 96, 108], [21, 44, 53, 68, 74, 91], [2, 4, 18, 61, 103, 106], [3, 24, 37, 42, 62, 102], [9, 32, 41, 56, 62, 79], [31, 73, 76, 83, 85, 99], [16, 57, 80, 89, 104, 110], [30, 42, 63, 76, 81, 101], [40, 82, 85, 92, 94, 108], [13, 54, 77, 86, 101, 107], [13, 16, 23, 25, 39, 82], [8, 48, 60, 81, 94, 99], [37, 40, 47, 49, 63, 106], [24, 36, 57, 70, 75, 95], [9, 30, 43, 48, 68, 108], [4, 46, 49, 56, 58, 72], [16, 19, 26, 28, 42, 85], [2, 17, 23, 40, 81, 104], [2, 42, 54, 75, 88, 93]]
\item 1 \{1=39960, 2=552780, 3=2185368, 4=2350092\} [[8, 39, 40, 41, 42, 70], [0, 1, 2, 3, 31, 80], [30, 31, 32, 33, 61, 110], [6, 7, 8, 9, 37, 86], [23, 54, 55, 56, 57, 85], [15, 16, 17, 18, 46, 95], [27, 28, 29, 30, 58, 107], [24, 25, 26, 27, 55, 104], [26, 57, 58, 59, 60, 88], [29, 60, 61, 62, 63, 91], [25, 74, 105, 106, 107, 108], [41, 72, 73, 74, 75, 103], [7, 56, 87, 88, 89, 90], [10, 59, 90, 91, 92, 93], [13, 62, 93, 94, 95, 96], [22, 71, 102, 103, 104, 105], [9, 10, 11, 12, 40, 89], [44, 75, 76, 77, 78, 106], [18, 19, 20, 21, 49, 98], [11, 42, 43, 44, 45, 73], [3, 4, 5, 6, 34, 83], [16, 65, 96, 97, 98, 99], [19, 68, 99, 100, 101, 102], [17, 48, 49, 50, 51, 79], [47, 78, 79, 80, 81, 109], [20, 51, 52, 53, 54, 82], [12, 13, 14, 15, 43, 92], [0, 28, 77, 108, 109, 110], [4, 53, 84, 85, 86, 87], [1, 50, 81, 82, 83, 84], [32, 63, 64, 65, 66, 94], [2, 33, 34, 35, 36, 64], [21, 22, 23, 24, 52, 101], [38, 69, 70, 71, 72, 100], [14, 45, 46, 47, 48, 76], [35, 66, 67, 68, 69, 97], [5, 36, 37, 38, 39, 67], [28, 43, 53, 61, 74, 101], [23, 61, 76, 86, 94, 107], [9, 32, 45, 101, 104, 108], [9, 66, 70, 79, 90, 106], [1, 10, 21, 37, 51, 108], [54, 58, 67, 78, 94, 108], [26, 29, 33, 45, 68, 81], [48, 52, 61, 72, 88, 102], [0, 4, 13, 24, 40, 54], [3, 7, 16, 27, 43, 57], [37, 52, 62, 70, 83, 110], [0, 57, 61, 70, 81, 97], [32, 35, 39, 51, 74, 87], [10, 20, 28, 41, 68, 106], [53, 56, 60, 72, 95, 108], [7, 18, 34, 48, 105, 109], [1, 15, 72, 76, 85, 96], [20, 23, 27, 39, 62, 75], [8, 21, 77, 80, 84, 96], [27, 31, 40, 51, 67, 81], [25, 40, 50, 58, 71, 98], [20, 58, 73, 83, 91, 104], [14, 17, 21, 33, 56, 69], [11, 24, 80, 83, 87, 99], [2, 10, 23, 50, 88, 103], [5, 13, 26, 53, 91, 106], [23, 26, 30, 42, 65, 78], [8, 46, 61, 71, 79, 92], [42, 46, 55, 66, 82, 96], [8, 11, 15, 27, 50, 63], [14, 27, 83, 86, 90, 102], [47, 50, 54, 66, 89, 102], [7, 20, 47, 85, 100, 110], [2, 29, 67, 82, 92, 100], [9, 25, 39, 96, 100, 109], [1, 16, 26, 34, 47, 74], [0, 23, 36, 92, 95, 99], [33, 37, 46, 57, 73, 87], [7, 21, 78, 82, 91, 102], [0, 16, 30, 87, 91, 100], [34, 49, 59, 67, 80, 107], [3, 26, 39, 95, 98, 102], [4, 14, 22, 35, 62, 100], [3, 15, 38, 51, 107, 110], [29, 32, 36, 48, 71, 84], [11, 49, 64, 74, 82, 95], [44, 47, 51, 63, 86, 99], [6, 10, 19, 30, 46, 60], [38, 41, 45, 57, 80, 93], [1, 12, 28, 42, 99, 103], [45, 49, 58, 69, 85, 99], [6, 63, 67, 76, 87, 103], [36, 40, 49, 60, 76, 90], [3, 19, 33, 90, 94, 103], [16, 31, 41, 49, 62, 89], [17, 55, 70, 80, 88, 101], [17, 20, 24, 36, 59, 72], [5, 32, 70, 85, 95, 103], [11, 38, 76, 91, 101, 109], [11, 14, 18, 30, 53, 66], [6, 62, 65, 69, 81, 104], [12, 16, 25, 36, 52, 66], [6, 29, 42, 98, 101, 105], [30, 34, 43, 54, 70, 84], [10, 25, 35, 43, 56, 83], [39, 43, 52, 63, 79, 93], [5, 8, 12, 24, 47, 60], [5, 43, 58, 68, 76, 89], [2, 15, 71, 74, 78, 90], [6, 22, 36, 93, 97, 106], [35, 38, 42, 54, 77, 90], [51, 55, 64, 75, 91, 105], [7, 17, 25, 38, 65, 103], [24, 28, 37, 48, 64, 78], [8, 35, 73, 88, 98, 106], [50, 53, 57, 69, 92, 105], [17, 30, 86, 89, 93, 105], [7, 22, 32, 40, 53, 80], [2, 6, 18, 41, 54, 110], [9, 65, 68, 72, 84, 107], [13, 23, 31, 44, 71, 109], [26, 64, 79, 89, 97, 110], [31, 46, 56, 64, 77, 104], [4, 18, 75, 79, 88, 99], [13, 28, 38, 46, 59, 86], [13, 27, 84, 88, 97, 108], [12, 68, 71, 75, 87, 110], [19, 34, 44, 52, 65, 92], [1, 11, 19, 32, 59, 97], [3, 60, 64, 73, 84, 100], [14, 52, 67, 77, 85, 98], [3, 59, 62, 66, 78, 101], [21, 25, 34, 45, 61, 75], [8, 16, 29, 56, 94, 109], [5, 18, 74, 77, 81, 93], [4, 15, 31, 45, 102, 106], [10, 24, 81, 85, 94, 105], [0, 56, 59, 63, 75, 98], [4, 19, 29, 37, 50, 77], [15, 19, 28, 39, 55, 69], [12, 69, 73, 82, 93, 109], [41, 44, 48, 60, 83, 96], [20, 33, 89, 92, 96, 108], [4, 17, 44, 82, 97, 107], [1, 14, 41, 79, 94, 104], [0, 12, 35, 48, 104, 107], [9, 13, 22, 33, 49, 63], [18, 22, 31, 42, 58, 72], [22, 37, 47, 55, 68, 95], [2, 5, 9, 21, 44, 57], [2, 40, 55, 65, 73, 86], [2, 14, 70, 91, 99, 108], [28, 49, 57, 66, 71, 83], [43, 64, 72, 81, 86, 98], [1, 9, 18, 23, 35, 91], [10, 18, 27, 32, 44, 100], [6, 11, 23, 79, 100, 108], [4, 12, 21, 26, 38, 94], [25, 46, 54, 63, 68, 80], [8, 64, 85, 93, 102, 107], [16, 24, 33, 38, 50, 106], [37, 58, 66, 75, 80, 92], [5, 61, 82, 90, 99, 104], [13, 34, 42, 51, 56, 68], [52, 73, 81, 90, 95, 107], [10, 31, 39, 48, 53, 65], [1, 22, 30, 39, 44, 56], [34, 55, 63, 72, 77, 89], [0, 5, 17, 73, 94, 102], [49, 70, 78, 87, 92, 104], [16, 37, 45, 54, 59, 71], [13, 21, 30, 35, 47, 103], [3, 8, 20, 76, 97, 105], [11, 67, 88, 96, 105, 110], [7, 15, 24, 29, 41, 97], [7, 28, 36, 45, 50, 62], [31, 52, 60, 69, 74, 86], [55, 76, 84, 93, 98, 110], [19, 40, 48, 57, 62, 74], [2, 58, 79, 87, 96, 101], [6, 15, 20, 32, 88, 109], [3, 12, 17, 29, 85, 106], [19, 27, 36, 41, 53, 109], [40, 61, 69, 78, 83, 95], [0, 9, 14, 26, 82, 103], [22, 43, 51, 60, 65, 77], [46, 67, 75, 84, 89, 101], [4, 25, 33, 42, 47, 59], [22, 29, 46, 74, 99, 110], [6, 12, 39, 59, 61, 77], [8, 33, 44, 67, 74, 91], [27, 33, 60, 80, 82, 98], [24, 44, 46, 62, 102, 108], [8, 31, 38, 55, 83, 108], [34, 39, 71, 82, 85, 108], [0, 11, 34, 41, 58, 86], [1, 8, 25, 53, 78, 89], [2, 42, 48, 75, 95, 97], [12, 32, 34, 50, 90, 96], [21, 58, 63, 95, 106, 109], [5, 30, 41, 64, 71, 88], [11, 13, 29, 69, 75, 102], [21, 41, 43, 59, 99, 105], [2, 27, 38, 61, 68, 85], [13, 20, 37, 65, 90, 101], [2, 19, 47, 72, 83, 106], [18, 24, 51, 71, 73, 89], [4, 7, 30, 67, 72, 104], [17, 28, 31, 54, 91, 96], [0, 6, 33, 53, 55, 71], [11, 22, 25, 48, 85, 90], [14, 25, 28, 51, 88, 93], [13, 18, 50, 61, 64, 87], [28, 33, 65, 76, 79, 102], [12, 18, 45, 65, 67, 83], [16, 44, 69, 80, 103, 110], [29, 40, 43, 66, 103, 108], [39, 45, 72, 92, 94, 110], [5, 45, 51, 78, 98, 100], [16, 23, 40, 68, 93, 104], [9, 29, 31, 47, 87, 93], [15, 35, 37, 53, 93, 99], [7, 14, 31, 59, 84, 95], [15, 21, 48, 68, 70, 86], [0, 20, 22, 38, 78, 84], [21, 27, 54, 74, 76, 92], [23, 48, 59, 82, 89, 106], [5, 28, 35, 52, 80, 105], [12, 23, 46, 53, 70, 98], [20, 45, 56, 79, 86, 103], [8, 48, 54, 81, 101, 103], [21, 32, 55, 62, 79, 107], [8, 19, 22, 45, 82, 87], [17, 19, 35, 75, 81, 108], [10, 38, 63, 74, 97, 104], [3, 35, 46, 49, 72, 109], [14, 54, 60, 87, 107, 109], [5, 16, 19, 42, 79, 84], [3, 14, 37, 44, 61, 89], [6, 26, 28, 44, 84, 90], [2, 13, 16, 39, 76, 81], [6, 43, 48, 80, 91, 94], [2, 25, 32, 49, 77, 102], [9, 20, 43, 50, 67, 95], [3, 30, 50, 52, 68, 108], [3, 9, 36, 56, 58, 74], [15, 26, 49, 56, 73, 101], [18, 29, 52, 59, 76, 104], [15, 52, 57, 89, 100, 103], [0, 27, 47, 49, 65, 105], [26, 37, 40, 63, 100, 105], [4, 11, 28, 56, 81, 92], [10, 13, 36, 73, 78, 110], [20, 31, 34, 57, 94, 99], [19, 26, 43, 71, 96, 107], [3, 23, 25, 41, 81, 87], [3, 40, 45, 77, 88, 91], [18, 55, 60, 92, 103, 106], [18, 38, 40, 56, 96, 102], [7, 10, 33, 70, 75, 107], [24, 35, 58, 65, 82, 110], [24, 30, 57, 77, 79, 95], [4, 9, 41, 52, 55, 78], [1, 17, 57, 63, 90, 110], [2, 4, 20, 60, 66, 93], [1, 29, 54, 65, 88, 95], [5, 7, 23, 63, 69, 96], [22, 27, 59, 70, 73, 96], [33, 39, 66, 86, 88, 104], [7, 12, 44, 55, 58, 81], [31, 36, 68, 79, 82, 105], [14, 16, 32, 72, 78, 105], [14, 39, 50, 73, 80, 97], [11, 51, 57, 84, 104, 106], [25, 30, 62, 73, 76, 99], [9, 46, 51, 83, 94, 97], [8, 10, 26, 66, 72, 99], [6, 17, 40, 47, 64, 92], [36, 42, 69, 89, 91, 107], [1, 6, 38, 49, 52, 75], [10, 15, 47, 58, 61, 84], [13, 41, 66, 77, 100, 107], [9, 15, 42, 62, 64, 80], [16, 21, 53, 64, 67, 90], [26, 51, 62, 85, 92, 109], [7, 35, 60, 71, 94, 101], [5, 22, 50, 75, 86, 109], [23, 34, 37, 60, 97, 102], [0, 37, 42, 74, 85, 88], [11, 36, 47, 70, 77, 94], [4, 32, 57, 68, 91, 98], [0, 32, 43, 46, 69, 106], [12, 49, 54, 86, 97, 100], [1, 4, 27, 64, 69, 101], [19, 24, 56, 67, 70, 93], [17, 42, 53, 76, 83, 100], [30, 36, 63, 83, 85, 101], [1, 24, 61, 66, 98, 109], [10, 17, 34, 62, 87, 98], [12, 22, 57, 64, 76, 108], [38, 44, 66, 87, 95, 105], [33, 40, 52, 84, 99, 109], [10, 42, 57, 67, 102, 109], [7, 39, 54, 64, 99, 106], [18, 26, 36, 80, 86, 108], [4, 8, 23, 43, 49, 110], [17, 23, 45, 66, 74, 84], [0, 10, 45, 52, 64, 96], [3, 13, 48, 55, 67, 99], [29, 34, 38, 53, 73, 79], [5, 11, 33, 54, 62, 72], [2, 17, 37, 43, 104, 109], [11, 31, 37, 98, 103, 107], [16, 22, 83, 88, 92, 107], [26, 32, 54, 75, 83, 93], [1, 5, 20, 40, 46, 107], [21, 36, 46, 81, 88, 100], [2, 12, 56, 62, 84, 105], [5, 15, 59, 65, 87, 108], [50, 55, 59, 74, 94, 100], [1, 36, 43, 55, 87, 102], [4, 10, 71, 76, 80, 95], [4, 39, 46, 58, 90, 105], [11, 16, 20, 35, 55, 61], [17, 22, 26, 41, 61, 67], [12, 19, 31, 63, 78, 88], [1, 13, 45, 60, 70, 105], [2, 7, 11, 26, 46, 52], [30, 45, 55, 90, 97, 109], [7, 13, 74, 79, 83, 98], [18, 25, 37, 69, 84, 94], [8, 28, 34, 95, 100, 104], [14, 34, 40, 101, 106, 110], [3, 47, 53, 75, 96, 104], [59, 64, 68, 83, 103, 109], [9, 24, 34, 69, 76, 88], [0, 7, 19, 51, 66, 76], [15, 36, 44, 54, 98, 104], [5, 10, 14, 29, 49, 55], [23, 28, 32, 47, 67, 73], [41, 46, 50, 65, 85, 91], [15, 23, 33, 77, 83, 105], [11, 17, 39, 60, 68, 78], [20, 26, 48, 69, 77, 87], [9, 17, 27, 71, 77, 99], [9, 19, 54, 61, 73, 105], [3, 18, 28, 63, 70, 82], [6, 27, 35, 45, 89, 95], [0, 44, 50, 72, 93, 101], [7, 42, 49, 61, 93, 108], [2, 22, 28, 89, 94, 98], [10, 16, 77, 82, 86, 101], [6, 21, 31, 66, 73, 85], [26, 31, 35, 50, 70, 76], [6, 13, 25, 57, 72, 82], [0, 21, 29, 39, 83, 89], [0, 8, 18, 62, 68, 90], [1, 7, 68, 73, 77, 92], [24, 39, 49, 84, 91, 103], [15, 22, 34, 66, 81, 91], [6, 14, 24, 68, 74, 96], [38, 43, 47, 62, 82, 88], [6, 16, 51, 58, 70, 102], [3, 24, 32, 42, 86, 92], [9, 30, 38, 48, 92, 98], [35, 40, 44, 59, 79, 85], [18, 33, 43, 78, 85, 97], [21, 42, 50, 60, 104, 110], [21, 28, 40, 72, 87, 97], [12, 27, 37, 72, 79, 91], [53, 58, 62, 77, 97, 103], [2, 8, 30, 51, 59, 69], [32, 38, 60, 81, 89, 99], [27, 42, 52, 87, 94, 106], [4, 16, 48, 63, 73, 108], [8, 14, 36, 57, 65, 75], [27, 34, 46, 78, 93, 103], [23, 29, 51, 72, 80, 90], [5, 27, 48, 56, 66, 110], [4, 36, 51, 61, 96, 103], [14, 20, 42, 63, 71, 81], [47, 52, 56, 71, 91, 97], [5, 25, 31, 92, 97, 101], [32, 37, 41, 56, 76, 82], [1, 62, 67, 71, 86, 106], [1, 33, 48, 58, 93, 100], [44, 49, 53, 68, 88, 94], [24, 31, 43, 75, 90, 100], [8, 13, 17, 32, 52, 58], [15, 30, 40, 75, 82, 94], [56, 61, 65, 80, 100, 106], [3, 10, 22, 54, 69, 79], [3, 11, 21, 65, 71, 93], [13, 19, 80, 85, 89, 104], [2, 24, 45, 53, 63, 107], [30, 37, 49, 81, 96, 106], [9, 53, 59, 81, 102, 110], [12, 33, 41, 51, 95, 101], [35, 41, 63, 84, 92, 102], [41, 47, 69, 90, 98, 108], [18, 39, 47, 57, 101, 107], [4, 65, 70, 74, 89, 109], [20, 25, 29, 44, 64, 70], [12, 20, 30, 74, 80, 102], [6, 50, 56, 78, 99, 107], [0, 15, 25, 60, 67, 79], [19, 25, 86, 91, 95, 110], [29, 35, 57, 78, 86, 96], [14, 19, 23, 38, 58, 64], [9, 16, 28, 60, 75, 85]]
\item 1 \{1=34632, 2=525474, 3=2205348, 4=2362746\} [[8, 39, 40, 41, 42, 70], [0, 1, 2, 3, 31, 80], [30, 31, 32, 33, 61, 110], [6, 7, 8, 9, 37, 86], [23, 54, 55, 56, 57, 85], [15, 16, 17, 18, 46, 95], [27, 28, 29, 30, 58, 107], [24, 25, 26, 27, 55, 104], [26, 57, 58, 59, 60, 88], [29, 60, 61, 62, 63, 91], [25, 74, 105, 106, 107, 108], [41, 72, 73, 74, 75, 103], [7, 56, 87, 88, 89, 90], [10, 59, 90, 91, 92, 93], [13, 62, 93, 94, 95, 96], [22, 71, 102, 103, 104, 105], [9, 10, 11, 12, 40, 89], [44, 75, 76, 77, 78, 106], [18, 19, 20, 21, 49, 98], [11, 42, 43, 44, 45, 73], [3, 4, 5, 6, 34, 83], [16, 65, 96, 97, 98, 99], [19, 68, 99, 100, 101, 102], [17, 48, 49, 50, 51, 79], [47, 78, 79, 80, 81, 109], [20, 51, 52, 53, 54, 82], [12, 13, 14, 15, 43, 92], [0, 28, 77, 108, 109, 110], [4, 53, 84, 85, 86, 87], [1, 50, 81, 82, 83, 84], [32, 63, 64, 65, 66, 94], [2, 33, 34, 35, 36, 64], [21, 22, 23, 24, 52, 101], [38, 69, 70, 71, 72, 100], [14, 45, 46, 47, 48, 76], [35, 66, 67, 68, 69, 97], [5, 36, 37, 38, 39, 67], [8, 22, 57, 61, 75, 89], [12, 26, 56, 70, 105, 109], [15, 30, 65, 67, 70, 89], [15, 19, 33, 47, 77, 91], [3, 28, 36, 44, 70, 101], [7, 38, 51, 76, 84, 92], [6, 31, 39, 47, 73, 104], [2, 15, 40, 48, 56, 82], [1, 20, 57, 72, 107, 109], [8, 10, 13, 32, 69, 84], [33, 37, 51, 65, 95, 109], [6, 14, 40, 71, 84, 109], [1, 32, 45, 70, 78, 86], [11, 13, 16, 35, 72, 87], [0, 15, 50, 52, 55, 74], [17, 54, 69, 104, 106, 109], [25, 56, 69, 94, 102, 110], [23, 37, 72, 76, 90, 104], [30, 34, 48, 62, 92, 106], [14, 27, 52, 60, 68, 94], [5, 7, 10, 29, 66, 81], [0, 25, 33, 41, 67, 98], [18, 33, 68, 70, 73, 92], [5, 18, 43, 51, 59, 85], [17, 19, 22, 41, 78, 93], [32, 34, 37, 56, 93, 108], [2, 28, 59, 72, 97, 105], [16, 24, 32, 58, 89, 102], [1, 9, 17, 43, 74, 87], [11, 41, 55, 90, 94, 108], [19, 27, 35, 61, 92, 105], [9, 34, 42, 50, 76, 107], [27, 31, 45, 59, 89, 103], [22, 30, 38, 64, 95, 108], [26, 39, 64, 72, 80, 106], [24, 39, 74, 76, 79, 98], [12, 37, 45, 53, 79, 110], [4, 12, 20, 46, 77, 90], [7, 42, 46, 60, 74, 104], [11, 48, 63, 98, 100, 103], [21, 25, 39, 53, 83, 97], [20, 34, 69, 73, 87, 101], [2, 32, 46, 81, 85, 99], [17, 31, 66, 70, 84, 98], [23, 25, 28, 47, 84, 99], [5, 19, 54, 58, 72, 86], [6, 21, 56, 58, 61, 80], [9, 13, 27, 41, 71, 85], [3, 18, 53, 55, 58, 77], [8, 21, 46, 54, 62, 88], [6, 20, 50, 64, 99, 103], [9, 23, 53, 67, 102, 106], [13, 21, 29, 55, 86, 99], [22, 53, 66, 91, 99, 107], [1, 4, 23, 60, 75, 110], [18, 22, 36, 50, 80, 94], [20, 22, 25, 44, 81, 96], [2, 39, 54, 89, 91, 94], [12, 47, 49, 52, 71, 108], [10, 41, 54, 79, 87, 95], [3, 38, 40, 43, 62, 99], [10, 45, 49, 63, 77, 107], [3, 17, 47, 61, 96, 100], [4, 35, 48, 73, 81, 89], [36, 51, 86, 88, 91, 110], [12, 27, 62, 64, 67, 86], [7, 15, 23, 49, 80, 93], [24, 28, 42, 56, 86, 100], [9, 44, 46, 49, 68, 105], [5, 35, 49, 84, 88, 102], [14, 51, 66, 101, 103, 106], [14, 16, 19, 38, 75, 90], [3, 7, 21, 35, 65, 79], [5, 42, 57, 92, 94, 97], [3, 11, 37, 68, 81, 106], [2, 4, 7, 26, 63, 78], [16, 47, 60, 85, 93, 101], [5, 31, 62, 75, 100, 108], [19, 50, 63, 88, 96, 104], [10, 18, 26, 52, 83, 96], [33, 48, 83, 85, 88, 107], [0, 35, 37, 40, 59, 96], [14, 28, 63, 67, 81, 95], [26, 40, 75, 79, 93, 107], [29, 42, 67, 75, 83, 109], [2, 16, 51, 55, 69, 83], [1, 36, 40, 54, 68, 98], [11, 25, 60, 64, 78, 92], [29, 43, 78, 82, 96, 110], [6, 41, 43, 46, 65, 102], [9, 24, 59, 61, 64, 83], [13, 44, 57, 82, 90, 98], [0, 4, 18, 32, 62, 76], [4, 39, 43, 57, 71, 101], [13, 48, 52, 66, 80, 110], [21, 36, 71, 73, 76, 95], [23, 36, 61, 69, 77, 103], [6, 10, 24, 38, 68, 82], [8, 45, 60, 95, 97, 100], [29, 31, 34, 53, 90, 105], [12, 16, 30, 44, 74, 88], [0, 8, 34, 65, 78, 103], [20, 33, 58, 66, 74, 100], [27, 42, 77, 79, 82, 101], [30, 45, 80, 82, 85, 104], [0, 14, 44, 58, 93, 97], [1, 15, 29, 59, 73, 108], [8, 38, 52, 87, 91, 105], [26, 28, 31, 50, 87, 102], [11, 24, 49, 57, 65, 91], [17, 30, 55, 63, 71, 97], [22, 40, 63, 69, 76, 85], [42, 54, 65, 81, 105, 110], [4, 22, 45, 51, 58, 67], [7, 55, 73, 96, 102, 109], [21, 26, 69, 81, 92, 108], [30, 42, 53, 69, 93, 98], [1, 35, 41, 52, 62, 77], [1, 10, 58, 76, 99, 105], [4, 52, 70, 93, 99, 106], [10, 20, 35, 70, 104, 110], [13, 31, 54, 60, 67, 76], [22, 56, 62, 73, 83, 98], [21, 33, 44, 60, 84, 89], [0, 11, 27, 51, 56, 99], [9, 15, 22, 31, 79, 97], [5, 11, 22, 32, 47, 82], [8, 14, 25, 35, 50, 85], [9, 21, 32, 48, 72, 77], [1, 11, 26, 61, 95, 101], [37, 55, 78, 84, 91, 100], [0, 6, 13, 22, 70, 88], [33, 45, 56, 72, 96, 101], [3, 15, 26, 42, 66, 71], [15, 21, 28, 37, 85, 103], [23, 29, 40, 50, 65, 100], [13, 47, 53, 64, 74, 89], [14, 20, 31, 41, 56, 91], [18, 23, 66, 78, 89, 105], [27, 39, 50, 66, 90, 95], [10, 28, 51, 57, 64, 73], [8, 24, 48, 53, 96, 108], [15, 20, 63, 75, 86, 102], [29, 35, 46, 56, 71, 106], [3, 14, 30, 54, 59, 102], [24, 36, 47, 63, 87, 92], [3, 8, 51, 63, 74, 90], [0, 7, 16, 64, 82, 105], [9, 20, 36, 60, 65, 108], [12, 24, 35, 51, 75, 80], [6, 11, 54, 66, 77, 93], [8, 23, 58, 92, 98, 109], [7, 41, 47, 58, 68, 83], [31, 65, 71, 82, 92, 107], [16, 39, 45, 52, 61, 109], [16, 50, 56, 67, 77, 92], [46, 64, 87, 93, 100, 109], [15, 27, 38, 54, 78, 83], [19, 37, 60, 66, 73, 82], [20, 26, 37, 47, 62, 97], [3, 27, 32, 75, 87, 98], [7, 25, 48, 54, 61, 70], [2, 45, 57, 68, 84, 108], [16, 34, 57, 63, 70, 79], [34, 68, 74, 85, 95, 110], [10, 33, 39, 46, 55, 103], [39, 51, 62, 78, 102, 107], [6, 12, 19, 28, 76, 94], [2, 13, 23, 38, 73, 107], [18, 30, 41, 57, 81, 86], [6, 17, 33, 57, 62, 105], [17, 23, 34, 44, 59, 94], [11, 46, 80, 86, 97, 107], [12, 17, 60, 72, 83, 99], [26, 32, 43, 53, 68, 103], [1, 19, 42, 48, 55, 64], [6, 30, 35, 78, 90, 101], [11, 17, 28, 38, 53, 88], [32, 38, 49, 59, 74, 109], [31, 49, 72, 78, 85, 94], [14, 49, 83, 89, 100, 110], [12, 36, 41, 84, 96, 107], [2, 37, 71, 77, 88, 98], [25, 59, 65, 76, 86, 101], [3, 9, 16, 25, 73, 91], [36, 48, 59, 75, 99, 104], [7, 30, 36, 43, 52, 100], [4, 14, 29, 64, 98, 104], [1, 49, 67, 90, 96, 103], [7, 17, 32, 67, 101, 107], [9, 33, 38, 81, 93, 104], [18, 24, 31, 40, 88, 106], [5, 16, 26, 41, 76, 110], [5, 40, 74, 80, 91, 101], [40, 58, 81, 87, 94, 103], [25, 43, 66, 72, 79, 88], [2, 8, 19, 29, 44, 79], [8, 43, 77, 83, 94, 104], [12, 18, 25, 34, 82, 100], [34, 52, 75, 81, 88, 97], [4, 27, 33, 40, 49, 97], [28, 62, 68, 79, 89, 104], [5, 20, 55, 89, 95, 106], [43, 61, 84, 90, 97, 106], [2, 17, 52, 86, 92, 103], [1, 24, 30, 37, 46, 94], [15, 39, 44, 87, 99, 110], [19, 53, 59, 70, 80, 95], [3, 10, 19, 67, 85, 108], [21, 27, 34, 43, 91, 109], [6, 18, 29, 45, 69, 74], [0, 5, 48, 60, 71, 87], [28, 46, 69, 75, 82, 91], [9, 14, 57, 69, 80, 96], [0, 24, 29, 72, 84, 95], [4, 38, 44, 55, 65, 80], [2, 18, 42, 47, 90, 102], [0, 12, 23, 39, 63, 68], [13, 36, 42, 49, 58, 106], [4, 13, 61, 79, 102, 108], [5, 21, 45, 50, 93, 105], [10, 44, 50, 61, 71, 86], [24, 60, 70, 81, 90, 107], [6, 55, 59, 98, 107, 110], [24, 34, 45, 54, 71, 99], [5, 33, 69, 79, 90, 99], [19, 25, 32, 40, 52, 57], [32, 41, 44, 51, 100, 104], [38, 47, 50, 57, 106, 110], [21, 31, 42, 51, 68, 96], [40, 44, 83, 92, 95, 102], [0, 10, 21, 30, 47, 75], [22, 26, 65, 74, 77, 84], [7, 11, 50, 59, 62, 69], [0, 49, 53, 92, 101, 104], [2, 5, 12, 61, 65, 104], [5, 13, 25, 30, 103, 109], [3, 39, 49, 60, 69, 86], [40, 46, 53, 61, 73, 78], [10, 14, 53, 62, 65, 72], [0, 9, 26, 54, 90, 100], [31, 37, 44, 52, 64, 69], [14, 23, 26, 33, 82, 86], [8, 11, 18, 67, 71, 110], [4, 9, 82, 88, 95, 103], [13, 19, 26, 34, 46, 51], [3, 12, 29, 57, 93, 103], [4, 11, 19, 31, 36, 109], [30, 40, 51, 60, 77, 105], [7, 19, 24, 97, 103, 110], [19, 23, 62, 71, 74, 81], [33, 43, 54, 63, 80, 108], [4, 16, 21, 94, 100, 107], [1, 5, 44, 53, 56, 63], [7, 18, 27, 44, 72, 108], [10, 16, 23, 31, 43, 48], [9, 45, 55, 66, 75, 92], [58, 64, 71, 79, 91, 96], [49, 55, 62, 70, 82, 87], [9, 19, 30, 39, 56, 84], [3, 52, 56, 95, 104, 107], [25, 29, 68, 77, 80, 87], [55, 61, 68, 76, 88, 93], [23, 32, 35, 42, 91, 95], [2, 9, 58, 62, 101, 110], [16, 20, 59, 68, 71, 78], [2, 41, 50, 53, 60, 109], [11, 39, 75, 85, 96, 105], [0, 17, 45, 81, 91, 102], [46, 50, 89, 98, 101, 108], [7, 12, 85, 91, 98, 106], [28, 32, 71, 80, 83, 90], [6, 16, 27, 36, 53, 81], [6, 42, 52, 63, 72, 89], [3, 13, 24, 33, 50, 78], [6, 23, 51, 87, 97, 108], [1, 13, 18, 91, 97, 104], [26, 35, 38, 45, 94, 98], [28, 34, 41, 49, 61, 66], [43, 49, 56, 64, 76, 81], [25, 31, 38, 46, 58, 63], [70, 76, 83, 91, 103, 108], [2, 10, 22, 27, 100, 106], [20, 29, 32, 39, 88, 92], [13, 17, 56, 65, 68, 75], [2, 11, 14, 21, 70, 74], [43, 47, 86, 95, 98, 105], [37, 43, 50, 58, 70, 75], [0, 73, 79, 86, 94, 106], [15, 25, 36, 45, 62, 90], [3, 20, 48, 84, 94, 105], [37, 41, 80, 89, 92, 99], [21, 57, 67, 78, 87, 104], [6, 15, 32, 60, 96, 106], [18, 54, 64, 75, 84, 101], [4, 10, 17, 25, 37, 42], [7, 13, 20, 28, 40, 45], [11, 20, 23, 30, 79, 83], [4, 8, 47, 56, 59, 66], [1, 7, 14, 22, 34, 39], [29, 38, 41, 48, 97, 101], [9, 18, 35, 63, 99, 109], [12, 22, 33, 42, 59, 87], [5, 8, 15, 64, 68, 107], [34, 38, 77, 86, 89, 96], [3, 76, 82, 89, 97, 109], [12, 48, 58, 69, 78, 95], [46, 52, 59, 67, 79, 84], [14, 42, 78, 88, 99, 108], [2, 30, 66, 76, 87, 96], [61, 67, 74, 82, 94, 99], [34, 40, 47, 55, 67, 72], [15, 51, 61, 72, 81, 98], [27, 63, 73, 84, 93, 110], [52, 58, 65, 73, 85, 90], [18, 28, 39, 48, 65, 93], [67, 73, 80, 88, 100, 105], [0, 36, 46, 57, 66, 83], [4, 15, 24, 41, 69, 105], [1, 8, 16, 28, 33, 106], [1, 12, 21, 38, 66, 102], [5, 14, 17, 24, 73, 77], [64, 70, 77, 85, 97, 102], [35, 44, 47, 54, 103, 107], [8, 17, 20, 27, 76, 80], [17, 26, 29, 36, 85, 89], [1, 6, 79, 85, 92, 100], [31, 35, 74, 83, 86, 93], [10, 15, 88, 94, 101, 109], [8, 36, 72, 82, 93, 102], [22, 28, 35, 43, 55, 60], [16, 22, 29, 37, 49, 54], [27, 37, 48, 57, 74, 102], [10, 34, 60, 80, 98, 102], [9, 29, 47, 51, 70, 94], [1, 25, 51, 71, 89, 93], [16, 40, 66, 86, 104, 108], [11, 15, 34, 58, 84, 104], [7, 33, 53, 71, 75, 94], [14, 18, 37, 61, 87, 107], [21, 41, 59, 63, 82, 106], [15, 35, 53, 57, 76, 100], [6, 26, 44, 48, 67, 91], [1, 27, 47, 65, 69, 88], [0, 20, 38, 42, 61, 85], [10, 36, 56, 74, 78, 97], [2, 20, 24, 43, 67, 93], [4, 30, 50, 68, 72, 91], [3, 23, 41, 45, 64, 88], [19, 45, 65, 83, 87, 106], [2, 6, 25, 49, 75, 95], [7, 31, 57, 77, 95, 99], [13, 37, 63, 83, 101, 105], [11, 29, 33, 52, 76, 102], [22, 48, 68, 86, 90, 109], [5, 9, 28, 52, 78, 98], [13, 39, 59, 77, 81, 100], [14, 32, 36, 55, 79, 105], [17, 21, 40, 64, 90, 110], [12, 32, 50, 54, 73, 97], [5, 23, 27, 46, 70, 96], [8, 12, 31, 55, 81, 101], [24, 44, 62, 66, 85, 109], [3, 22, 46, 72, 92, 110], [16, 42, 62, 80, 84, 103], [4, 28, 54, 74, 92, 96], [8, 26, 30, 49, 73, 99], [0, 19, 43, 69, 89, 107], [17, 35, 39, 58, 82, 108], [18, 38, 56, 60, 79, 103]]
\item 1 \{1=40848, 2=524808, 3=2217336, 4=2345208\} [[8, 39, 40, 41, 42, 70], [0, 1, 2, 3, 31, 80], [30, 31, 32, 33, 61, 110], [6, 7, 8, 9, 37, 86], [23, 54, 55, 56, 57, 85], [15, 16, 17, 18, 46, 95], [27, 28, 29, 30, 58, 107], [24, 25, 26, 27, 55, 104], [26, 57, 58, 59, 60, 88], [29, 60, 61, 62, 63, 91], [25, 74, 105, 106, 107, 108], [41, 72, 73, 74, 75, 103], [7, 56, 87, 88, 89, 90], [10, 59, 90, 91, 92, 93], [13, 62, 93, 94, 95, 96], [22, 71, 102, 103, 104, 105], [9, 10, 11, 12, 40, 89], [44, 75, 76, 77, 78, 106], [18, 19, 20, 21, 49, 98], [11, 42, 43, 44, 45, 73], [3, 4, 5, 6, 34, 83], [16, 65, 96, 97, 98, 99], [19, 68, 99, 100, 101, 102], [17, 48, 49, 50, 51, 79], [47, 78, 79, 80, 81, 109], [20, 51, 52, 53, 54, 82], [12, 13, 14, 15, 43, 92], [0, 28, 77, 108, 109, 110], [4, 53, 84, 85, 86, 87], [1, 50, 81, 82, 83, 84], [32, 63, 64, 65, 66, 94], [2, 33, 34, 35, 36, 64], [21, 22, 23, 24, 52, 101], [38, 69, 70, 71, 72, 100], [14, 45, 46, 47, 48, 76], [35, 66, 67, 68, 69, 97], [5, 36, 37, 38, 39, 67], [1, 8, 19, 27, 62, 75], [8, 54, 61, 78, 92, 100], [0, 11, 25, 53, 96, 100], [10, 18, 53, 66, 103, 110], [1, 9, 44, 57, 94, 101], [45, 52, 69, 83, 91, 110], [24, 31, 48, 62, 70, 89], [11, 54, 58, 69, 80, 94], [5, 13, 32, 78, 85, 102], [0, 7, 24, 38, 46, 65], [22, 29, 40, 48, 83, 96], [27, 31, 42, 53, 67, 95], [9, 46, 53, 64, 72, 107], [13, 20, 31, 39, 74, 87], [8, 21, 58, 65, 76, 84], [23, 36, 73, 80, 91, 99], [7, 35, 78, 82, 93, 104], [24, 28, 39, 50, 64, 92], [28, 35, 46, 54, 89, 102], [11, 57, 64, 81, 95, 103], [2, 10, 29, 75, 82, 99], [7, 26, 72, 79, 96, 110], [6, 41, 54, 91, 98, 109], [0, 37, 44, 55, 63, 98], [18, 25, 42, 56, 64, 83], [3, 14, 28, 56, 99, 103], [4, 21, 35, 43, 62, 108], [6, 10, 21, 32, 46, 74], [14, 27, 64, 71, 82, 90], [36, 40, 51, 62, 76, 104], [15, 19, 30, 41, 55, 83], [11, 24, 61, 68, 79, 87], [9, 16, 33, 47, 55, 74], [4, 23, 69, 76, 93, 107], [5, 51, 58, 75, 89, 97], [19, 26, 37, 45, 80, 93], [21, 25, 36, 47, 61, 89], [5, 18, 55, 62, 73, 81], [1, 18, 32, 40, 59, 105], [1, 12, 23, 37, 65, 108], [39, 46, 63, 77, 85, 104], [0, 35, 48, 85, 92, 103], [6, 20, 28, 47, 93, 100], [30, 37, 54, 68, 76, 95], [4, 11, 22, 30, 65, 78], [10, 17, 28, 36, 71, 84], [3, 38, 51, 88, 95, 106], [2, 15, 52, 59, 70, 78], [18, 22, 33, 44, 58, 86], [29, 42, 79, 86, 97, 105], [30, 34, 45, 56, 70, 98], [39, 43, 54, 65, 79, 107], [33, 40, 57, 71, 79, 98], [1, 29, 72, 76, 87, 98], [3, 17, 25, 44, 90, 97], [7, 14, 25, 33, 68, 81], [23, 66, 70, 81, 92, 106], [11, 19, 38, 84, 91, 108], [26, 39, 76, 83, 94, 102], [12, 49, 56, 67, 75, 110], [12, 16, 27, 38, 52, 80], [17, 30, 67, 74, 85, 93], [8, 16, 35, 81, 88, 105], [17, 63, 70, 87, 101, 109], [10, 38, 81, 85, 96, 107], [4, 12, 47, 60, 97, 104], [16, 23, 34, 42, 77, 90], [12, 26, 34, 53, 99, 106], [5, 19, 47, 90, 94, 105], [0, 14, 22, 41, 87, 94], [17, 60, 64, 75, 86, 100], [5, 48, 52, 63, 74, 88], [20, 33, 70, 77, 88, 96], [20, 63, 67, 78, 89, 103], [8, 51, 55, 66, 77, 91], [2, 13, 21, 56, 69, 106], [5, 16, 24, 59, 72, 109], [9, 23, 31, 50, 96, 103], [31, 38, 49, 57, 92, 105], [3, 7, 18, 29, 43, 71], [12, 19, 36, 50, 58, 77], [27, 34, 51, 65, 73, 92], [14, 57, 61, 72, 83, 97], [34, 41, 52, 60, 95, 108], [9, 20, 34, 62, 105, 109], [25, 32, 43, 51, 86, 99], [3, 40, 47, 58, 66, 101], [36, 43, 60, 74, 82, 101], [7, 15, 50, 63, 100, 107], [14, 60, 67, 84, 98, 106], [2, 45, 49, 60, 71, 85], [6, 43, 50, 61, 69, 104], [42, 46, 57, 68, 82, 110], [13, 41, 84, 88, 99, 110], [2, 16, 44, 87, 91, 102], [33, 37, 48, 59, 73, 101], [2, 48, 55, 72, 86, 94], [21, 28, 45, 59, 67, 86], [15, 29, 37, 56, 102, 109], [0, 4, 15, 26, 40, 68], [4, 32, 75, 79, 90, 101], [15, 22, 39, 53, 61, 80], [26, 69, 73, 84, 95, 109], [1, 20, 66, 73, 90, 104], [6, 13, 30, 44, 52, 71], [3, 10, 27, 41, 49, 68], [32, 45, 82, 89, 100, 108], [42, 49, 66, 80, 88, 107], [9, 13, 24, 35, 49, 77], [8, 22, 50, 93, 97, 108], [6, 17, 31, 59, 102, 106], [12, 61, 73, 98, 102, 107], [10, 22, 47, 51, 56, 72], [9, 58, 70, 95, 99, 104], [11, 15, 20, 36, 85, 97], [7, 19, 44, 48, 53, 69], [3, 8, 24, 73, 85, 110], [22, 34, 59, 63, 68, 84], [0, 49, 61, 86, 90, 95], [25, 37, 62, 66, 71, 87], [5, 9, 14, 30, 79, 91], [4, 16, 41, 45, 50, 66], [6, 55, 67, 92, 96, 101], [1, 13, 38, 42, 47, 63], [20, 24, 29, 45, 94, 106], [40, 52, 77, 81, 86, 102], [34, 46, 71, 75, 80, 96], [2, 18, 67, 79, 104, 108], [19, 31, 56, 60, 65, 81], [10, 35, 39, 44, 60, 109], [8, 12, 17, 33, 82, 94], [16, 28, 53, 57, 62, 78], [46, 58, 83, 87, 92, 108], [3, 52, 64, 89, 93, 98], [7, 32, 36, 41, 57, 106], [31, 43, 68, 72, 77, 93], [43, 55, 80, 84, 89, 105], [14, 18, 23, 39, 88, 100], [17, 21, 26, 42, 91, 103], [28, 40, 65, 69, 74, 90], [0, 5, 21, 70, 82, 107], [15, 64, 76, 101, 105, 110], [4, 29, 33, 38, 54, 103], [37, 49, 74, 78, 83, 99], [2, 6, 11, 27, 76, 88], [23, 27, 32, 48, 97, 109], [13, 25, 50, 54, 59, 75], [1, 26, 30, 35, 51, 100], [7, 17, 39, 45, 58, 62], [19, 29, 51, 57, 70, 74], [16, 26, 48, 54, 67, 71], [28, 38, 60, 66, 79, 83], [49, 59, 81, 87, 100, 104], [15, 21, 34, 38, 94, 104], [6, 12, 25, 29, 85, 95], [12, 18, 31, 35, 91, 101], [31, 41, 63, 69, 82, 86], [43, 53, 75, 81, 94, 98], [0, 13, 17, 73, 83, 105], [13, 23, 45, 51, 64, 68], [2, 24, 30, 43, 47, 103], [8, 30, 36, 49, 53, 109], [1, 11, 33, 39, 52, 56], [21, 27, 40, 44, 100, 110], [25, 35, 57, 63, 76, 80], [40, 50, 72, 78, 91, 95], [7, 11, 67, 77, 99, 105], [52, 62, 84, 90, 103, 107], [55, 65, 87, 93, 106, 110], [10, 14, 70, 80, 102, 108], [4, 14, 36, 42, 55, 59], [10, 20, 42, 48, 61, 65], [37, 47, 69, 75, 88, 92], [2, 58, 68, 90, 96, 109], [0, 6, 19, 23, 79, 89], [5, 27, 33, 46, 50, 106], [18, 24, 37, 41, 97, 107], [22, 32, 54, 60, 73, 77], [34, 44, 66, 72, 85, 89], [1, 5, 61, 71, 93, 99], [9, 15, 28, 32, 88, 98], [4, 8, 64, 74, 96, 102], [3, 16, 20, 76, 86, 108], [46, 56, 78, 84, 97, 101], [3, 9, 22, 26, 82, 92], [2, 7, 12, 22, 66, 74], [3, 11, 50, 55, 60, 70], [30, 38, 77, 82, 87, 97], [39, 47, 86, 91, 96, 106], [12, 20, 59, 64, 69, 79], [0, 8, 47, 52, 57, 67], [18, 26, 65, 70, 75, 85], [17, 22, 27, 37, 81, 89], [5, 44, 49, 54, 64, 108], [1, 45, 53, 92, 97, 102], [27, 35, 74, 79, 84, 94], [33, 41, 80, 85, 90, 100], [21, 29, 68, 73, 78, 88], [32, 37, 42, 52, 96, 104], [26, 31, 36, 46, 90, 98], [42, 50, 89, 94, 99, 109], [36, 44, 83, 88, 93, 103], [2, 41, 46, 51, 61, 105], [7, 51, 59, 98, 103, 108], [3, 13, 57, 65, 104, 109], [8, 13, 18, 28, 72, 80], [0, 10, 54, 62, 101, 106], [20, 25, 30, 40, 84, 92], [35, 40, 45, 55, 99, 107], [4, 9, 19, 63, 71, 110], [11, 16, 21, 31, 75, 83], [1, 6, 16, 60, 68, 107], [5, 10, 15, 25, 69, 77], [15, 23, 62, 67, 72, 82], [4, 48, 56, 95, 100, 105], [38, 43, 48, 58, 102, 110], [24, 32, 71, 76, 81, 91], [29, 34, 39, 49, 93, 101], [6, 14, 53, 58, 63, 73], [9, 17, 56, 61, 66, 76], [14, 19, 24, 34, 78, 86], [23, 28, 33, 43, 87, 95], [23, 26, 38, 44, 61, 74], [11, 71, 74, 86, 92, 109], [9, 18, 38, 45, 78, 90], [18, 27, 47, 54, 87, 99], [22, 25, 31, 45, 88, 109], [12, 24, 54, 63, 83, 90], [25, 46, 70, 73, 79, 93], [0, 33, 45, 75, 84, 104], [10, 23, 83, 86, 98, 104], [19, 22, 28, 42, 85, 106], [17, 20, 32, 38, 55, 68], [24, 33, 53, 60, 93, 105], [2, 8, 25, 38, 98, 101], [13, 16, 22, 36, 79, 100], [9, 52, 73, 97, 100, 106], [5, 12, 45, 57, 87, 96], [13, 37, 40, 46, 60, 103], [6, 18, 48, 57, 77, 84], [10, 34, 37, 43, 57, 100], [7, 28, 52, 55, 61, 75], [53, 56, 68, 74, 91, 104], [4, 18, 61, 82, 106, 109], [19, 43, 46, 52, 66, 109], [3, 36, 48, 78, 87, 107], [3, 23, 30, 63, 75, 105], [50, 53, 65, 71, 88, 101], [13, 34, 58, 61, 67, 81], [56, 59, 71, 77, 94, 107], [10, 31, 55, 58, 64, 78], [2, 62, 65, 77, 83, 100], [8, 11, 23, 29, 46, 59], [20, 23, 35, 41, 58, 71], [3, 12, 32, 39, 72, 84], [9, 39, 48, 68, 75, 108], [5, 8, 20, 26, 43, 56], [3, 33, 42, 62, 69, 102], [8, 68, 71, 83, 89, 106], [21, 30, 50, 57, 90, 102], [5, 22, 35, 95, 98, 110], [7, 10, 16, 30, 73, 94], [38, 41, 53, 59, 76, 89], [14, 21, 54, 66, 96, 105], [1, 25, 28, 34, 48, 91], [27, 39, 69, 78, 98, 105], [2, 14, 20, 37, 50, 110], [4, 28, 31, 37, 51, 94], [0, 30, 39, 59, 66, 99], [3, 15, 45, 54, 74, 81], [1, 15, 58, 79, 103, 106], [6, 36, 45, 65, 72, 105], [5, 65, 68, 80, 86, 103], [4, 7, 13, 27, 70, 91], [11, 17, 34, 47, 107, 110], [11, 18, 51, 63, 93, 102], [6, 49, 70, 94, 97, 103], [47, 50, 62, 68, 85, 98], [7, 31, 34, 40, 54, 97], [28, 49, 73, 76, 82, 96], [1, 7, 21, 64, 85, 109], [2, 5, 17, 23, 40, 53], [37, 58, 82, 85, 91, 105], [15, 27, 57, 66, 86, 93], [4, 17, 77, 80, 92, 98], [14, 17, 29, 35, 52, 65], [29, 32, 44, 50, 67, 80], [9, 21, 51, 60, 80, 87], [7, 20, 80, 83, 95, 101], [1, 22, 46, 49, 55, 69], [0, 12, 42, 51, 71, 78], [18, 30, 60, 69, 89, 96], [16, 40, 43, 49, 63, 106], [4, 25, 49, 52, 58, 72], [13, 26, 86, 89, 101, 107], [2, 19, 32, 92, 95, 107], [11, 14, 26, 32, 49, 62], [27, 36, 56, 63, 96, 108], [19, 40, 64, 67, 73, 87], [5, 11, 28, 41, 101, 104], [12, 21, 41, 48, 81, 93], [59, 62, 74, 80, 97, 110], [41, 44, 56, 62, 79, 92], [6, 26, 33, 66, 78, 108], [35, 38, 50, 56, 73, 86], [22, 43, 67, 70, 76, 90], [8, 14, 31, 44, 104, 107], [26, 29, 41, 47, 64, 77], [16, 37, 61, 64, 70, 84], [0, 20, 27, 60, 72, 102], [12, 55, 76, 100, 103, 109], [44, 47, 59, 65, 82, 95], [3, 46, 67, 91, 94, 100], [17, 24, 57, 69, 99, 108], [31, 52, 76, 79, 85, 99], [16, 19, 25, 39, 82, 103], [15, 24, 44, 51, 84, 96], [16, 29, 89, 92, 104, 110], [40, 61, 85, 88, 94, 108], [0, 43, 64, 88, 91, 97], [30, 42, 72, 81, 101, 108], [21, 33, 63, 72, 92, 99], [8, 15, 48, 60, 90, 99], [6, 15, 35, 42, 75, 87], [2, 9, 42, 54, 84, 93], [32, 35, 47, 53, 70, 83], [0, 9, 29, 36, 69, 81], [24, 36, 66, 75, 95, 102], [1, 4, 10, 24, 67, 88], [1, 14, 74, 77, 89, 95], [34, 55, 79, 82, 88, 102], [6, 39, 51, 81, 90, 110], [10, 13, 19, 33, 76, 97], [6, 24, 40, 56, 80, 82], [33, 51, 67, 83, 107, 109], [12, 30, 46, 62, 86, 88], [14, 38, 40, 75, 93, 109], [6, 22, 38, 62, 64, 99], [9, 25, 41, 65, 67, 102], [12, 28, 44, 68, 70, 105], [11, 35, 37, 72, 90, 106], [3, 21, 37, 53, 77, 79], [5, 7, 42, 60, 76, 92], [13, 29, 53, 55, 90, 108], [21, 39, 55, 71, 95, 97], [5, 29, 31, 66, 84, 100], [2, 4, 39, 57, 73, 89], [10, 26, 50, 52, 87, 105], [18, 36, 52, 68, 92, 94], [1, 36, 54, 70, 86, 110], [17, 19, 54, 72, 88, 104], [3, 19, 35, 59, 61, 96], [2, 26, 28, 63, 81, 97], [1, 17, 41, 43, 78, 96], [15, 31, 47, 71, 73, 108], [23, 25, 60, 78, 94, 110], [4, 20, 44, 46, 81, 99], [8, 10, 45, 63, 79, 95], [15, 33, 49, 65, 89, 91], [24, 42, 58, 74, 98, 100], [7, 23, 47, 49, 84, 102], [0, 18, 34, 50, 74, 76], [0, 16, 32, 56, 58, 93], [30, 48, 64, 80, 104, 106], [27, 45, 61, 77, 101, 103], [11, 13, 48, 66, 82, 98], [9, 27, 43, 59, 83, 85], [20, 22, 57, 75, 91, 107], [8, 32, 34, 69, 87, 103], [14, 16, 51, 69, 85, 101]]
\item 1 \{0=1110, 1=42624, 2=621378, 3=2187588, 4=2275500\} [[8, 39, 40, 41, 42, 70], [0, 1, 2, 3, 31, 80], [30, 31, 32, 33, 61, 110], [6, 7, 8, 9, 37, 86], [23, 54, 55, 56, 57, 85], [15, 16, 17, 18, 46, 95], [27, 28, 29, 30, 58, 107], [24, 25, 26, 27, 55, 104], [26, 57, 58, 59, 60, 88], [29, 60, 61, 62, 63, 91], [25, 74, 105, 106, 107, 108], [41, 72, 73, 74, 75, 103], [7, 56, 87, 88, 89, 90], [10, 59, 90, 91, 92, 93], [13, 62, 93, 94, 95, 96], [22, 71, 102, 103, 104, 105], [9, 10, 11, 12, 40, 89], [44, 75, 76, 77, 78, 106], [18, 19, 20, 21, 49, 98], [11, 42, 43, 44, 45, 73], [3, 4, 5, 6, 34, 83], [16, 65, 96, 97, 98, 99], [19, 68, 99, 100, 101, 102], [17, 48, 49, 50, 51, 79], [47, 78, 79, 80, 81, 109], [20, 51, 52, 53, 54, 82], [12, 13, 14, 15, 43, 92], [0, 28, 77, 108, 109, 110], [4, 53, 84, 85, 86, 87], [1, 50, 81, 82, 83, 84], [32, 63, 64, 65, 66, 94], [2, 33, 34, 35, 36, 64], [21, 22, 23, 24, 52, 101], [38, 69, 70, 71, 72, 100], [14, 45, 46, 47, 48, 76], [35, 66, 67, 68, 69, 97], [5, 36, 37, 38, 39, 67], [41, 43, 55, 63, 78, 86], [5, 71, 73, 85, 93, 108], [51, 57, 70, 83, 92, 109], [3, 9, 22, 35, 44, 61], [3, 7, 17, 24, 58, 74], [1, 13, 21, 36, 44, 110], [7, 15, 30, 38, 104, 106], [30, 36, 49, 62, 71, 88], [33, 37, 47, 54, 88, 104], [14, 54, 58, 68, 75, 109], [0, 34, 50, 90, 94, 104], [56, 58, 70, 78, 93, 101], [25, 41, 81, 85, 95, 102], [27, 33, 46, 59, 68, 85], [48, 54, 67, 80, 89, 106], [0, 13, 26, 35, 52, 105], [50, 52, 64, 72, 87, 95], [12, 16, 26, 33, 67, 83], [1, 9, 24, 32, 98, 100], [23, 25, 37, 45, 60, 68], [5, 14, 31, 84, 90, 103], [0, 4, 14, 21, 55, 71], [2, 42, 46, 56, 63, 97], [29, 31, 43, 51, 66, 74], [4, 20, 60, 64, 74, 81], [1, 17, 57, 61, 71, 78], [36, 42, 55, 68, 77, 94], [9, 15, 28, 41, 50, 67], [19, 35, 75, 79, 89, 96], [30, 34, 44, 51, 85, 101], [9, 13, 23, 30, 64, 80], [4, 12, 27, 35, 101, 103], [22, 38, 78, 82, 92, 99], [11, 20, 37, 90, 96, 109], [5, 7, 19, 27, 42, 50], [21, 25, 35, 42, 76, 92], [18, 22, 32, 39, 73, 89], [12, 18, 31, 44, 53, 70], [26, 28, 40, 48, 63, 71], [1, 11, 18, 52, 68, 108], [5, 22, 75, 81, 94, 107], [6, 21, 29, 95, 97, 109], [8, 10, 22, 30, 45, 53], [7, 60, 66, 79, 92, 101], [8, 48, 52, 62, 69, 103], [31, 47, 87, 91, 101, 108], [6, 10, 20, 27, 61, 77], [33, 39, 52, 65, 74, 91], [0, 6, 19, 32, 41, 58], [38, 40, 52, 60, 75, 83], [3, 37, 53, 93, 97, 107], [45, 51, 64, 77, 86, 103], [1, 14, 23, 40, 93, 99], [10, 18, 33, 41, 107, 109], [2, 4, 16, 24, 39, 47], [8, 17, 34, 87, 93, 106], [10, 26, 66, 70, 80, 87], [24, 28, 38, 45, 79, 95], [2, 68, 70, 82, 90, 105], [65, 67, 79, 87, 102, 110], [8, 15, 49, 65, 105, 109], [20, 22, 34, 42, 57, 65], [6, 40, 56, 96, 100, 110], [10, 23, 32, 49, 102, 108], [39, 43, 53, 60, 94, 110], [18, 24, 37, 50, 59, 76], [9, 17, 83, 85, 97, 105], [35, 37, 49, 57, 72, 80], [53, 55, 67, 75, 90, 98], [14, 16, 28, 36, 51, 59], [4, 57, 63, 76, 89, 98], [15, 19, 29, 36, 70, 86], [17, 19, 31, 39, 54, 62], [11, 51, 55, 65, 72, 106], [32, 34, 46, 54, 69, 77], [47, 49, 61, 69, 84, 92], [7, 20, 29, 46, 99, 105], [16, 32, 72, 76, 86, 93], [59, 61, 73, 81, 96, 104], [3, 11, 77, 79, 91, 99], [21, 27, 40, 53, 62, 79], [2, 9, 43, 59, 99, 103], [2, 19, 72, 78, 91, 104], [3, 18, 26, 92, 94, 106], [24, 30, 43, 56, 65, 82], [39, 45, 58, 71, 80, 97], [4, 17, 26, 43, 96, 102], [36, 40, 50, 57, 91, 107], [6, 14, 80, 82, 94, 102], [44, 46, 58, 66, 81, 89], [2, 11, 28, 81, 87, 100], [8, 25, 78, 84, 97, 110], [6, 12, 25, 38, 47, 64], [28, 44, 84, 88, 98, 105], [7, 23, 63, 67, 77, 84], [62, 64, 76, 84, 99, 107], [11, 13, 25, 33, 48, 56], [15, 21, 34, 47, 56, 73], [16, 69, 75, 88, 101, 110], [0, 8, 74, 76, 88, 96], [13, 66, 72, 85, 98, 107], [12, 20, 86, 88, 100, 108], [10, 63, 69, 82, 95, 104], [27, 31, 41, 48, 82, 98], [3, 16, 29, 38, 55, 108], [5, 12, 46, 62, 102, 106], [13, 29, 69, 73, 83, 90], [0, 15, 23, 89, 91, 103], [5, 45, 49, 59, 66, 100], [42, 48, 61, 74, 83, 100], [1, 54, 60, 73, 86, 95], [2, 7, 51, 73, 102, 107], [27, 49, 78, 83, 89, 94], [24, 46, 75, 80, 86, 91], [4, 33, 38, 44, 49, 93], [6, 28, 57, 62, 68, 73], [9, 14, 20, 25, 69, 91], [33, 55, 84, 89, 95, 100], [27, 32, 38, 43, 87, 109], [21, 26, 32, 37, 81, 103], [3, 8, 14, 19, 63, 85], [21, 43, 72, 77, 83, 88], [24, 29, 35, 40, 84, 106], [1, 45, 67, 96, 101, 107], [18, 23, 29, 34, 78, 100], [19, 48, 53, 59, 64, 108], [36, 58, 87, 92, 98, 103], [15, 37, 66, 71, 77, 82], [10, 39, 44, 50, 55, 99], [2, 8, 13, 57, 79, 108], [3, 25, 54, 59, 65, 70], [15, 20, 26, 31, 75, 97], [9, 31, 60, 65, 71, 76], [42, 64, 93, 98, 104, 109], [7, 36, 41, 47, 52, 96], [12, 34, 63, 68, 74, 79], [0, 22, 51, 56, 62, 67], [16, 45, 50, 56, 61, 105], [39, 61, 90, 95, 101, 106], [5, 10, 54, 76, 105, 110], [0, 5, 11, 16, 60, 82], [4, 48, 70, 99, 104, 110], [30, 52, 81, 86, 92, 97], [13, 42, 47, 53, 58, 102], [6, 11, 17, 22, 66, 88], [18, 40, 69, 74, 80, 85], [12, 17, 23, 28, 72, 94], [1, 30, 35, 41, 46, 90], [14, 39, 49, 56, 77, 107], [14, 24, 61, 64, 70, 88], [18, 28, 35, 56, 86, 104], [19, 22, 28, 46, 83, 93], [18, 30, 54, 83, 96, 103], [9, 46, 49, 55, 73, 110], [0, 24, 53, 66, 73, 99], [2, 12, 49, 52, 58, 76], [4, 22, 59, 69, 106, 109], [15, 44, 57, 64, 90, 102], [3, 27, 56, 69, 76, 102], [7, 10, 16, 34, 71, 81], [2, 27, 37, 44, 65, 95], [34, 37, 43, 61, 98, 108], [9, 16, 42, 54, 78, 107], [5, 15, 52, 55, 61, 79], [6, 16, 23, 44, 74, 92], [21, 33, 57, 86, 99, 106], [8, 26, 51, 61, 68, 89], [3, 32, 45, 52, 78, 90], [7, 14, 35, 65, 83, 108], [14, 27, 34, 60, 72, 96], [6, 43, 46, 52, 70, 107], [31, 34, 40, 58, 95, 105], [20, 50, 68, 93, 103, 110], [24, 34, 41, 62, 92, 110], [1, 8, 29, 59, 77, 102], [18, 47, 60, 67, 93, 105], [4, 7, 13, 31, 68, 78], [2, 20, 45, 55, 62, 83], [14, 44, 62, 87, 97, 104], [4, 41, 51, 88, 91, 97], [23, 41, 66, 76, 83, 104], [20, 38, 63, 73, 80, 101], [23, 36, 43, 69, 81, 105], [10, 47, 57, 94, 97, 103], [5, 18, 25, 51, 63, 87], [26, 39, 46, 72, 84, 108], [3, 15, 39, 68, 81, 88], [11, 36, 46, 53, 74, 104], [0, 29, 42, 49, 75, 87], [21, 50, 63, 70, 96, 108], [28, 31, 37, 55, 92, 102], [21, 31, 38, 59, 89, 107], [8, 21, 28, 54, 66, 90], [5, 35, 53, 78, 88, 95], [1, 19, 56, 66, 103, 106], [5, 23, 48, 58, 65, 86], [1, 7, 25, 62, 72, 109], [1, 27, 39, 63, 92, 105], [1, 38, 48, 85, 88, 94], [6, 30, 59, 72, 79, 105], [11, 29, 54, 64, 71, 92], [8, 38, 56, 81, 91, 98], [9, 38, 51, 58, 84, 96], [12, 41, 54, 61, 87, 99], [6, 18, 42, 71, 84, 91], [35, 45, 82, 85, 91, 109], [7, 44, 54, 91, 94, 100], [17, 30, 37, 63, 75, 99], [17, 47, 65, 90, 100, 107], [13, 16, 22, 40, 77, 87], [26, 44, 69, 79, 86, 107], [9, 21, 45, 74, 87, 94], [5, 26, 56, 74, 99, 109], [13, 50, 60, 97, 100, 106], [10, 13, 19, 37, 74, 84], [16, 53, 63, 100, 103, 109], [20, 33, 40, 66, 78, 102], [25, 28, 34, 52, 89, 99], [4, 30, 42, 66, 95, 108], [5, 30, 40, 47, 68, 98], [2, 32, 50, 75, 85, 92], [0, 7, 33, 45, 69, 98], [12, 19, 45, 57, 81, 110], [23, 33, 70, 73, 79, 97], [8, 18, 55, 58, 64, 82], [2, 15, 22, 48, 60, 84], [17, 35, 60, 70, 77, 98], [26, 36, 73, 76, 82, 100], [9, 19, 26, 47, 77, 95], [11, 21, 58, 61, 67, 85], [3, 13, 20, 41, 71, 89], [0, 10, 17, 38, 68, 86], [4, 11, 32, 62, 80, 105], [29, 47, 72, 82, 89, 110], [3, 10, 36, 48, 72, 101], [6, 35, 48, 55, 81, 93], [16, 19, 25, 43, 80, 90], [15, 27, 51, 80, 93, 100], [2, 23, 53, 71, 96, 106], [11, 24, 31, 57, 69, 93], [29, 39, 76, 79, 85, 103], [12, 24, 48, 77, 90, 97], [8, 33, 43, 50, 71, 101], [32, 42, 79, 82, 88, 106], [0, 37, 40, 46, 64, 101], [17, 42, 52, 59, 80, 110], [14, 32, 57, 67, 74, 95], [11, 41, 59, 84, 94, 101], [0, 12, 36, 65, 78, 85], [1, 4, 10, 28, 65, 75], [24, 36, 60, 89, 102, 109], [17, 27, 64, 67, 73, 91], [15, 25, 32, 53, 83, 101], [20, 30, 67, 70, 76, 94], [9, 33, 62, 75, 82, 108], [12, 22, 29, 50, 80, 98], [22, 25, 31, 49, 86, 96], [6, 13, 39, 51, 75, 104], [3, 40, 43, 49, 67, 104], [23, 26, 38, 42, 62, 90], [27, 71, 74, 86, 90, 110], [13, 34, 59, 67, 82, 86], [5, 13, 28, 32, 70, 91], [1, 12, 37, 42, 51, 69], [17, 25, 40, 44, 82, 103], [11, 19, 34, 38, 76, 97], [13, 38, 46, 61, 65, 103], [6, 24, 67, 78, 103, 108], [1, 6, 15, 33, 76, 87], [22, 43, 68, 76, 91, 95], [19, 30, 55, 60, 69, 87], [10, 21, 46, 51, 60, 78], [0, 9, 27, 70, 81, 106], [21, 65, 68, 80, 84, 104], [2, 5, 17, 21, 41, 69], [13, 17, 55, 76, 101, 109], [7, 12, 21, 39, 82, 93], [3, 21, 64, 75, 100, 105], [34, 45, 70, 75, 84, 102], [9, 53, 56, 68, 72, 92], [0, 43, 54, 79, 84, 93], [3, 46, 57, 82, 87, 96], [19, 40, 65, 73, 88, 92], [40, 51, 76, 81, 90, 108], [34, 55, 80, 88, 103, 107], [10, 14, 52, 73, 98, 106], [2, 30, 74, 77, 89, 93], [7, 11, 49, 70, 95, 103], [37, 48, 73, 78, 87, 105], [22, 33, 58, 63, 72, 90], [12, 56, 59, 71, 75, 95], [11, 14, 26, 30, 50, 78], [2, 10, 25, 29, 67, 88], [12, 55, 66, 91, 96, 105], [25, 46, 71, 79, 94, 98], [8, 11, 23, 27, 47, 75], [5, 8, 20, 24, 44, 72], [0, 18, 61, 72, 97, 102], [4, 29, 37, 52, 56, 94], [38, 41, 53, 57, 77, 105], [28, 39, 64, 69, 78, 96], [4, 15, 40, 45, 54, 72], [1, 26, 34, 49, 53, 91], [32, 35, 47, 51, 71, 99], [1, 22, 47, 55, 70, 74], [0, 25, 30, 39, 57, 100], [6, 31, 36, 45, 63, 106], [3, 47, 50, 62, 66, 86], [3, 28, 33, 42, 60, 103], [14, 22, 37, 41, 79, 100], [2, 40, 61, 86, 94, 109], [17, 45, 89, 92, 104, 108], [2, 14, 18, 38, 66, 110], [35, 38, 50, 54, 74, 102], [16, 41, 49, 64, 68, 106], [18, 62, 65, 77, 81, 101], [7, 28, 53, 61, 76, 80], [28, 49, 74, 82, 97, 101], [7, 22, 26, 64, 85, 110], [19, 44, 52, 67, 71, 109], [8, 16, 31, 35, 73, 94], [37, 58, 83, 91, 106, 110], [14, 42, 86, 89, 101, 105], [11, 15, 35, 63, 107, 110], [3, 12, 30, 73, 84, 109], [8, 36, 80, 83, 95, 99], [22, 27, 36, 54, 97, 108], [13, 18, 27, 45, 88, 99], [25, 36, 61, 66, 75, 93], [4, 9, 18, 36, 79, 90], [31, 52, 77, 85, 100, 104], [5, 9, 29, 57, 101, 104], [41, 44, 56, 60, 80, 108], [13, 24, 49, 54, 63, 81], [23, 31, 46, 50, 88, 109], [31, 42, 67, 72, 81, 99], [3, 23, 51, 95, 98, 110], [26, 29, 41, 45, 65, 93], [10, 35, 43, 58, 62, 100], [20, 23, 35, 39, 59, 87], [10, 31, 56, 64, 79, 83], [15, 59, 62, 74, 78, 98], [11, 39, 83, 86, 98, 102], [4, 25, 50, 58, 73, 77], [4, 19, 23, 61, 82, 107], [4, 8, 46, 67, 92, 100], [29, 32, 44, 48, 68, 96], [6, 49, 60, 85, 90, 99], [14, 17, 29, 33, 53, 81], [24, 68, 71, 83, 87, 107], [16, 27, 52, 57, 66, 84], [7, 18, 43, 48, 57, 75], [16, 37, 62, 70, 85, 89], [9, 34, 39, 48, 66, 109], [16, 21, 30, 48, 91, 102], [5, 33, 77, 80, 92, 96], [2, 6, 26, 54, 98, 101], [20, 28, 43, 47, 85, 106], [7, 32, 40, 55, 59, 97], [1, 5, 43, 64, 89, 97], [8, 12, 32, 60, 104, 107], [15, 58, 69, 94, 99, 108], [10, 15, 24, 42, 85, 96], [17, 20, 32, 36, 56, 84], [0, 20, 48, 92, 95, 107], [1, 16, 20, 58, 79, 104], [9, 52, 63, 88, 93, 102], [0, 44, 47, 59, 63, 83], [19, 24, 33, 51, 94, 105], [6, 50, 53, 65, 69, 89]]
\item 1 \{0=2220, 1=31968, 2=582750, 3=2148516, 4=2362746\} [[8, 39, 40, 41, 42, 70], [0, 1, 2, 3, 31, 80], [30, 31, 32, 33, 61, 110], [6, 7, 8, 9, 37, 86], [23, 54, 55, 56, 57, 85], [15, 16, 17, 18, 46, 95], [27, 28, 29, 30, 58, 107], [24, 25, 26, 27, 55, 104], [26, 57, 58, 59, 60, 88], [29, 60, 61, 62, 63, 91], [25, 74, 105, 106, 107, 108], [41, 72, 73, 74, 75, 103], [7, 56, 87, 88, 89, 90], [10, 59, 90, 91, 92, 93], [13, 62, 93, 94, 95, 96], [22, 71, 102, 103, 104, 105], [9, 10, 11, 12, 40, 89], [44, 75, 76, 77, 78, 106], [18, 19, 20, 21, 49, 98], [11, 42, 43, 44, 45, 73], [3, 4, 5, 6, 34, 83], [16, 65, 96, 97, 98, 99], [19, 68, 99, 100, 101, 102], [17, 48, 49, 50, 51, 79], [47, 78, 79, 80, 81, 109], [20, 51, 52, 53, 54, 82], [12, 13, 14, 15, 43, 92], [0, 28, 77, 108, 109, 110], [4, 53, 84, 85, 86, 87], [1, 50, 81, 82, 83, 84], [32, 63, 64, 65, 66, 94], [2, 33, 34, 35, 36, 64], [21, 22, 23, 24, 52, 101], [38, 69, 70, 71, 72, 100], [14, 45, 46, 47, 48, 76], [35, 66, 67, 68, 69, 97], [5, 36, 37, 38, 39, 67], [28, 33, 59, 76, 89, 98], [11, 75, 79, 82, 94, 102], [3, 7, 10, 22, 30, 50], [15, 50, 72, 86, 93, 109], [39, 43, 46, 58, 66, 86], [63, 67, 70, 82, 90, 110], [1, 9, 29, 93, 97, 100], [1, 13, 21, 41, 105, 109], [12, 26, 33, 49, 66, 101], [2, 43, 48, 74, 91, 104], [26, 48, 62, 69, 85, 102], [2, 19, 32, 41, 82, 87], [14, 31, 44, 53, 94, 99], [23, 40, 53, 62, 103, 108], [30, 34, 37, 49, 57, 77], [54, 58, 61, 73, 81, 101], [21, 25, 28, 40, 48, 68], [33, 37, 40, 52, 60, 80], [15, 29, 36, 52, 69, 104], [32, 54, 68, 75, 91, 108], [0, 16, 33, 68, 90, 104], [14, 36, 50, 57, 73, 90], [20, 42, 56, 63, 79, 96], [17, 34, 47, 56, 97, 102], [11, 20, 61, 66, 92, 109], [36, 40, 43, 55, 63, 83], [29, 51, 65, 72, 88, 105], [8, 49, 54, 80, 97, 110], [2, 9, 25, 42, 77, 99], [48, 52, 55, 67, 75, 95], [13, 18, 44, 61, 74, 83], [23, 45, 59, 66, 82, 99], [60, 64, 67, 79, 87, 107], [9, 23, 30, 46, 63, 98], [10, 15, 41, 58, 71, 80], [7, 20, 29, 70, 75, 101], [0, 4, 7, 19, 27, 47], [12, 16, 19, 31, 39, 59], [0, 35, 57, 71, 78, 94], [4, 21, 56, 78, 92, 99], [21, 35, 42, 58, 75, 110], [4, 9, 35, 52, 65, 74], [6, 10, 13, 25, 33, 53], [1, 18, 53, 75, 89, 96], [5, 46, 51, 77, 94, 107], [17, 39, 53, 60, 76, 93], [1, 6, 32, 49, 62, 71], [13, 30, 65, 87, 101, 108], [10, 23, 32, 73, 78, 104], [0, 20, 84, 88, 91, 103], [4, 12, 32, 96, 100, 103], [40, 45, 71, 88, 101, 110], [0, 14, 21, 37, 54, 89], [11, 33, 47, 54, 70, 87], [6, 20, 27, 43, 60, 95], [5, 22, 35, 44, 85, 90], [3, 38, 60, 74, 81, 97], [51, 55, 58, 70, 78, 98], [25, 30, 56, 73, 86, 95], [6, 41, 63, 77, 84, 100], [5, 27, 41, 48, 64, 81], [7, 24, 59, 81, 95, 102], [4, 17, 26, 67, 72, 98], [6, 26, 90, 94, 97, 109], [3, 17, 24, 40, 57, 92], [11, 18, 34, 51, 86, 108], [16, 21, 47, 64, 77, 86], [8, 25, 38, 47, 88, 93], [9, 13, 16, 28, 36, 56], [7, 12, 38, 55, 68, 77], [37, 42, 68, 85, 98, 107], [18, 22, 25, 37, 45, 65], [18, 32, 39, 55, 72, 107], [2, 66, 70, 73, 85, 93], [42, 46, 49, 61, 69, 89], [1, 14, 23, 64, 69, 95], [22, 27, 53, 70, 83, 92], [16, 29, 38, 79, 84, 110], [17, 81, 85, 88, 100, 108], [2, 11, 52, 57, 83, 100], [24, 28, 31, 43, 51, 71], [3, 23, 87, 91, 94, 106], [34, 39, 65, 82, 95, 104], [19, 24, 50, 67, 80, 89], [13, 26, 35, 76, 81, 107], [10, 18, 38, 102, 106, 109], [11, 28, 41, 50, 91, 96], [5, 69, 73, 76, 88, 96], [6, 22, 39, 74, 96, 110], [8, 30, 44, 51, 67, 84], [9, 44, 66, 80, 87, 103], [8, 72, 76, 79, 91, 99], [7, 15, 35, 99, 103, 106], [1, 4, 16, 24, 44, 108], [3, 29, 46, 59, 68, 109], [31, 36, 62, 79, 92, 101], [3, 19, 36, 71, 93, 107], [57, 61, 64, 76, 84, 104], [10, 27, 62, 84, 98, 105], [20, 37, 50, 59, 100, 105], [0, 26, 43, 56, 65, 106], [5, 14, 55, 60, 86, 103], [5, 12, 28, 45, 80, 102], [14, 78, 82, 85, 97, 105], [45, 49, 52, 64, 72, 92], [15, 19, 22, 34, 42, 62], [8, 15, 31, 48, 83, 105], [12, 47, 69, 83, 90, 106], [2, 24, 38, 45, 61, 78], [27, 31, 34, 46, 54, 74], [8, 17, 58, 63, 89, 106], [9, 21, 27, 32, 76, 90], [12, 42, 54, 60, 65, 109], [7, 14, 17, 41, 66, 107], [10, 24, 54, 66, 72, 77], [3, 15, 21, 26, 70, 84], [24, 32, 37, 46, 70, 106], [2, 27, 68, 79, 86, 89], [19, 48, 56, 61, 70, 94], [1, 8, 11, 35, 60, 101], [0, 12, 18, 23, 67, 81], [25, 54, 62, 67, 76, 100], [8, 19, 26, 29, 53, 78], [1, 15, 45, 57, 63, 68], [0, 30, 42, 48, 53, 97], [11, 22, 29, 32, 56, 81], [5, 30, 71, 82, 89, 92], [18, 30, 36, 41, 85, 99], [21, 29, 34, 43, 67, 103], [21, 62, 73, 80, 83, 107], [15, 23, 28, 37, 61, 97], [9, 50, 61, 68, 71, 95], [6, 36, 48, 54, 59, 103], [0, 8, 13, 22, 46, 82], [26, 37, 44, 47, 71, 96], [3, 44, 55, 62, 65, 89], [34, 48, 78, 90, 96, 101], [29, 40, 47, 50, 74, 99], [2, 26, 51, 92, 103, 110], [31, 45, 75, 87, 93, 98], [43, 57, 87, 99, 105, 110], [0, 6, 11, 55, 69, 99], [7, 31, 67, 96, 104, 109], [19, 55, 84, 92, 97, 106], [1, 30, 38, 43, 52, 76], [3, 8, 52, 66, 96, 108], [23, 48, 89, 100, 107, 110], [16, 52, 81, 89, 94, 103], [35, 46, 53, 56, 80, 105], [13, 49, 78, 86, 91, 100], [37, 51, 81, 93, 99, 104], [5, 10, 19, 43, 79, 108], [21, 33, 39, 44, 88, 102], [10, 39, 47, 52, 61, 85], [7, 43, 72, 80, 85, 94], [22, 36, 66, 78, 84, 89], [4, 33, 41, 46, 55, 79], [14, 39, 80, 91, 98, 101], [6, 14, 19, 28, 52, 88], [20, 31, 38, 41, 65, 90], [7, 36, 44, 49, 58, 82], [4, 11, 14, 38, 63, 104], [2, 46, 60, 90, 102, 108], [22, 51, 59, 64, 73, 97], [18, 59, 70, 77, 80, 104], [6, 18, 24, 29, 73, 87], [32, 43, 50, 53, 77, 102], [28, 57, 65, 70, 79, 103], [1, 25, 61, 90, 98, 103], [27, 35, 40, 49, 73, 109], [24, 65, 76, 83, 86, 110], [22, 58, 87, 95, 100, 109], [6, 47, 58, 65, 68, 92], [34, 63, 71, 76, 85, 109], [24, 36, 42, 47, 91, 105], [2, 7, 16, 40, 76, 105], [17, 42, 83, 94, 101, 104], [12, 53, 64, 71, 74, 98], [4, 18, 48, 60, 66, 71], [18, 26, 31, 40, 64, 100], [11, 36, 77, 88, 95, 98], [9, 39, 51, 57, 62, 106], [7, 21, 51, 63, 69, 74], [5, 16, 23, 26, 50, 75], [12, 24, 30, 35, 79, 93], [2, 13, 20, 23, 47, 72], [9, 17, 22, 31, 55, 91], [6, 12, 17, 61, 75, 105], [27, 39, 45, 50, 94, 108], [13, 42, 50, 55, 64, 88], [9, 15, 20, 64, 78, 108], [1, 10, 34, 70, 99, 107], [5, 8, 32, 57, 98, 109], [3, 11, 16, 25, 49, 85], [19, 33, 63, 75, 81, 86], [17, 28, 35, 38, 62, 87], [38, 49, 56, 59, 83, 108], [16, 45, 53, 58, 67, 91], [12, 20, 25, 34, 58, 94], [14, 25, 32, 35, 59, 84], [23, 34, 41, 44, 68, 93], [0, 41, 52, 59, 62, 86], [10, 17, 20, 44, 69, 110], [28, 42, 72, 84, 90, 95], [4, 28, 64, 93, 101, 106], [40, 54, 84, 96, 102, 107], [15, 27, 33, 38, 82, 96], [10, 46, 75, 83, 88, 97], [31, 60, 68, 73, 82, 106], [4, 40, 69, 77, 82, 91], [4, 13, 37, 73, 102, 110], [2, 5, 29, 54, 95, 106], [8, 33, 74, 85, 92, 95], [1, 37, 66, 74, 79, 88], [25, 39, 69, 81, 87, 92], [13, 27, 57, 69, 75, 80], [16, 30, 60, 72, 78, 83], [15, 56, 67, 74, 77, 101], [20, 45, 86, 97, 104, 107], [3, 9, 14, 58, 72, 102], [3, 33, 45, 51, 56, 100], [0, 5, 49, 63, 93, 105], [2, 10, 14, 49, 65, 67], [9, 33, 67, 73, 94, 105], [14, 22, 26, 61, 77, 79], [9, 26, 41, 45, 54, 83], [31, 47, 49, 95, 103, 107], [2, 6, 15, 44, 81, 98], [5, 20, 24, 33, 62, 99], [15, 32, 47, 51, 60, 89], [1, 22, 33, 48, 72, 106], [8, 12, 21, 50, 87, 104], [2, 17, 21, 30, 59, 96], [10, 21, 36, 60, 94, 100], [10, 16, 37, 48, 63, 87], [5, 42, 59, 74, 78, 87], [36, 53, 68, 72, 81, 110], [14, 16, 62, 70, 74, 109], [0, 29, 66, 83, 98, 102], [2, 39, 56, 71, 75, 84], [34, 50, 52, 98, 106, 110], [11, 15, 24, 53, 90, 107], [44, 52, 56, 91, 107, 109], [22, 28, 49, 60, 75, 99], [1, 17, 19, 65, 73, 77], [23, 31, 35, 70, 86, 88], [4, 25, 36, 51, 75, 109], [19, 25, 46, 57, 72, 96], [1, 7, 28, 39, 54, 78], [23, 60, 77, 92, 96, 105], [9, 43, 49, 70, 81, 96], [14, 18, 27, 56, 93, 110], [7, 23, 25, 71, 79, 83], [7, 11, 46, 62, 64, 110], [15, 49, 55, 76, 87, 102], [11, 13, 59, 67, 71, 106], [3, 20, 35, 39, 48, 77], [11, 26, 30, 39, 68, 105], [8, 10, 56, 64, 68, 103], [2, 4, 50, 58, 62, 97], [14, 51, 68, 83, 87, 96], [5, 9, 18, 47, 84, 101], [33, 50, 65, 69, 78, 107], [2, 37, 53, 55, 101, 109], [0, 15, 39, 73, 79, 100], [1, 47, 55, 59, 94, 110], [4, 20, 22, 68, 76, 80], [8, 23, 27, 36, 65, 102], [14, 29, 33, 42, 71, 108], [17, 54, 71, 86, 90, 99], [18, 35, 50, 54, 63, 92], [6, 35, 72, 89, 104, 108], [16, 32, 34, 80, 88, 92], [4, 15, 30, 54, 88, 94], [17, 25, 29, 64, 80, 82], [6, 40, 46, 67, 78, 93], [38, 46, 50, 85, 101, 103], [6, 30, 64, 70, 91, 102], [21, 38, 53, 57, 66, 95], [10, 26, 28, 74, 82, 86], [22, 38, 40, 86, 94, 98], [6, 21, 45, 79, 85, 106], [3, 32, 69, 86, 101, 105], [19, 35, 37, 83, 91, 95], [3, 37, 43, 64, 75, 90], [19, 30, 45, 69, 103, 109], [12, 29, 44, 48, 57, 86], [0, 9, 38, 75, 92, 107], [29, 37, 41, 76, 92, 94], [21, 55, 61, 82, 93, 108], [7, 13, 34, 45, 60, 84], [18, 52, 58, 79, 90, 105], [4, 8, 43, 59, 61, 107], [4, 10, 31, 42, 57, 81], [11, 19, 23, 58, 74, 76], [6, 23, 38, 42, 51, 80], [3, 27, 61, 67, 88, 99], [3, 18, 42, 76, 82, 103], [3, 12, 41, 78, 95, 110], [5, 13, 17, 52, 68, 70], [13, 24, 39, 63, 97, 103], [12, 46, 52, 73, 84, 99], [0, 24, 58, 64, 85, 96], [13, 19, 40, 51, 66, 90], [28, 34, 55, 66, 81, 105], [8, 45, 62, 77, 81, 90], [26, 63, 80, 95, 99, 108], [16, 27, 42, 66, 100, 106], [30, 47, 62, 66, 75, 104], [8, 16, 20, 55, 71, 73], [31, 37, 58, 69, 84, 108], [5, 7, 53, 61, 65, 100], [28, 44, 46, 92, 100, 104], [35, 43, 47, 82, 98, 100], [0, 34, 40, 61, 72, 87], [41, 49, 53, 88, 104, 106], [1, 5, 40, 56, 58, 104], [32, 40, 44, 79, 95, 97], [27, 44, 59, 63, 72, 101], [25, 41, 43, 89, 97, 101], [9, 24, 48, 82, 88, 109], [24, 41, 56, 60, 69, 98], [13, 29, 31, 77, 85, 89], [1, 12, 27, 51, 85, 91], [20, 28, 32, 67, 83, 85], [0, 17, 32, 36, 45, 74], [25, 31, 52, 63, 78, 102], [20, 57, 74, 89, 93, 102], [12, 36, 70, 76, 97, 108], [11, 48, 65, 80, 84, 93], [7, 18, 33, 57, 91, 97], [26, 34, 38, 73, 89, 91], [16, 22, 43, 54, 69, 93], [6, 16, 57, 82, 101, 107], [22, 41, 47, 57, 67, 108], [1, 42, 67, 86, 92, 102], [13, 32, 38, 48, 58, 99], [4, 23, 29, 39, 49, 90], [7, 26, 32, 42, 52, 93], [1, 20, 26, 36, 46, 87], [11, 17, 27, 37, 78, 103], [24, 49, 68, 74, 84, 94], [2, 8, 18, 28, 69, 94], [12, 37, 56, 62, 72, 82], [5, 15, 25, 66, 91, 110], [21, 46, 65, 71, 81, 91], [18, 43, 62, 68, 78, 88], [17, 23, 33, 43, 84, 109], [27, 52, 71, 77, 87, 97], [14, 20, 30, 40, 81, 106], [19, 38, 44, 54, 64, 105], [9, 34, 53, 59, 69, 79], [10, 29, 35, 45, 55, 96], [5, 11, 21, 31, 72, 97], [39, 64, 83, 89, 99, 109], [8, 14, 24, 34, 75, 100], [30, 55, 74, 80, 90, 100], [3, 28, 47, 53, 63, 73], [9, 19, 60, 85, 104, 110], [7, 48, 73, 92, 98, 108], [0, 25, 44, 50, 60, 70], [15, 40, 59, 65, 75, 85], [16, 35, 41, 51, 61, 102], [4, 45, 70, 89, 95, 105], [36, 61, 80, 86, 96, 106], [33, 58, 77, 83, 93, 103], [0, 10, 51, 76, 95, 101], [3, 13, 54, 79, 98, 104], [6, 31, 50, 56, 66, 76], [2, 12, 22, 63, 88, 107]]
\item 1 \{1=26640, 2=504828, 3=2198688, 4=2398044\} [[8, 39, 40, 41, 42, 70], [0, 1, 2, 3, 31, 80], [30, 31, 32, 33, 61, 110], [6, 7, 8, 9, 37, 86], [23, 54, 55, 56, 57, 85], [15, 16, 17, 18, 46, 95], [27, 28, 29, 30, 58, 107], [24, 25, 26, 27, 55, 104], [26, 57, 58, 59, 60, 88], [29, 60, 61, 62, 63, 91], [25, 74, 105, 106, 107, 108], [41, 72, 73, 74, 75, 103], [7, 56, 87, 88, 89, 90], [10, 59, 90, 91, 92, 93], [13, 62, 93, 94, 95, 96], [22, 71, 102, 103, 104, 105], [9, 10, 11, 12, 40, 89], [44, 75, 76, 77, 78, 106], [18, 19, 20, 21, 49, 98], [11, 42, 43, 44, 45, 73], [3, 4, 5, 6, 34, 83], [16, 65, 96, 97, 98, 99], [19, 68, 99, 100, 101, 102], [17, 48, 49, 50, 51, 79], [47, 78, 79, 80, 81, 109], [20, 51, 52, 53, 54, 82], [12, 13, 14, 15, 43, 92], [0, 28, 77, 108, 109, 110], [4, 53, 84, 85, 86, 87], [1, 50, 81, 82, 83, 84], [32, 63, 64, 65, 66, 94], [2, 33, 34, 35, 36, 64], [21, 22, 23, 24, 52, 101], [38, 69, 70, 71, 72, 100], [14, 45, 46, 47, 48, 76], [35, 66, 67, 68, 69, 97], [5, 36, 37, 38, 39, 67], [2, 18, 22, 26, 93, 107], [23, 64, 67, 86, 99, 109], [13, 44, 54, 72, 97, 105], [1, 26, 67, 70, 89, 102], [22, 25, 44, 57, 67, 92], [6, 31, 39, 58, 89, 99], [25, 28, 47, 60, 70, 95], [12, 22, 47, 88, 91, 110], [7, 38, 48, 66, 91, 99], [0, 19, 50, 60, 78, 103], [54, 68, 74, 90, 94, 98], [3, 22, 53, 63, 81, 106], [3, 21, 46, 54, 73, 104], [1, 20, 33, 43, 68, 109], [14, 27, 37, 62, 103, 106], [0, 14, 20, 36, 40, 44], [9, 27, 52, 60, 79, 110], [16, 24, 43, 74, 84, 102], [11, 52, 55, 74, 87, 97], [8, 21, 31, 56, 97, 100], [0, 18, 43, 51, 70, 101], [2, 43, 46, 65, 78, 88], [20, 61, 64, 83, 96, 106], [19, 22, 41, 54, 64, 89], [51, 65, 71, 87, 91, 95], [3, 7, 11, 78, 92, 98], [34, 37, 56, 69, 79, 104], [1, 32, 42, 60, 85, 93], [22, 30, 49, 80, 90, 108], [15, 29, 35, 51, 55, 59], [17, 27, 45, 70, 78, 97], [63, 77, 83, 99, 103, 107], [5, 21, 25, 29, 96, 110], [66, 80, 86, 102, 106, 110], [29, 39, 57, 82, 90, 109], [15, 19, 23, 90, 104, 110], [27, 41, 47, 63, 67, 71], [2, 69, 83, 89, 105, 109], [14, 55, 58, 77, 90, 100], [21, 35, 41, 57, 61, 65], [7, 32, 73, 76, 95, 108], [10, 13, 32, 45, 55, 80], [9, 23, 29, 45, 49, 53], [57, 71, 77, 93, 97, 101], [2, 8, 24, 28, 32, 99], [3, 28, 36, 55, 86, 96], [12, 16, 20, 87, 101, 107], [30, 44, 50, 66, 70, 74], [9, 34, 42, 61, 92, 102], [17, 58, 61, 80, 93, 103], [8, 49, 52, 71, 84, 94], [31, 34, 53, 66, 76, 101], [42, 56, 62, 78, 82, 86], [4, 7, 26, 39, 49, 74], [3, 13, 38, 79, 82, 101], [24, 38, 44, 60, 64, 68], [6, 25, 56, 66, 84, 109], [9, 13, 17, 84, 98, 104], [18, 32, 38, 54, 58, 62], [1, 9, 28, 59, 69, 87], [1, 5, 72, 86, 92, 108], [11, 21, 39, 64, 72, 91], [37, 40, 59, 72, 82, 107], [12, 37, 45, 64, 95, 105], [4, 35, 45, 63, 88, 96], [39, 53, 59, 75, 79, 83], [6, 10, 14, 81, 95, 101], [8, 18, 36, 61, 69, 88], [13, 21, 40, 71, 81, 99], [4, 12, 31, 62, 72, 90], [7, 15, 34, 65, 75, 93], [7, 10, 29, 42, 52, 77], [13, 16, 35, 48, 58, 83], [2, 15, 25, 50, 91, 94], [2, 12, 30, 55, 63, 82], [28, 31, 50, 63, 73, 98], [19, 27, 46, 77, 87, 105], [6, 24, 49, 57, 76, 107], [20, 30, 48, 73, 81, 100], [14, 24, 42, 67, 75, 94], [6, 16, 41, 82, 85, 104], [36, 50, 56, 72, 76, 80], [5, 15, 33, 58, 66, 85], [11, 24, 34, 59, 100, 103], [10, 18, 37, 68, 78, 96], [0, 4, 8, 75, 89, 95], [60, 74, 80, 96, 100, 104], [5, 11, 27, 31, 35, 102], [10, 41, 51, 69, 94, 102], [9, 19, 44, 85, 88, 107], [0, 10, 35, 76, 79, 98], [16, 19, 38, 51, 61, 86], [11, 17, 33, 37, 41, 108], [40, 43, 62, 75, 85, 110], [15, 40, 48, 67, 98, 108], [0, 25, 33, 52, 83, 93], [45, 59, 65, 81, 85, 89], [12, 26, 32, 48, 52, 56], [33, 47, 53, 69, 73, 77], [4, 29, 70, 73, 92, 105], [26, 36, 54, 79, 87, 106], [6, 20, 26, 42, 46, 50], [5, 46, 49, 68, 81, 91], [3, 17, 23, 39, 43, 47], [23, 33, 51, 76, 84, 103], [1, 4, 23, 36, 46, 71], [8, 14, 30, 34, 38, 105], [17, 30, 40, 65, 106, 109], [48, 62, 68, 84, 88, 92], [16, 47, 57, 75, 100, 108], [5, 18, 28, 53, 94, 97], [21, 32, 87, 92, 103, 109], [33, 38, 49, 55, 78, 89], [9, 14, 25, 31, 54, 65], [30, 35, 46, 52, 75, 86], [15, 26, 81, 86, 97, 103], [45, 50, 61, 67, 90, 101], [1, 7, 30, 41, 96, 101], [0, 5, 16, 22, 45, 56], [27, 32, 43, 49, 72, 83], [5, 60, 65, 76, 82, 105], [3, 8, 19, 25, 48, 59], [9, 20, 75, 80, 91, 97], [51, 56, 67, 73, 96, 107], [6, 11, 22, 28, 51, 62], [6, 17, 72, 77, 88, 94], [18, 23, 34, 40, 63, 74], [36, 41, 52, 58, 81, 92], [21, 26, 37, 43, 66, 77], [42, 47, 58, 64, 87, 98], [54, 59, 70, 76, 99, 110], [39, 44, 55, 61, 84, 95], [10, 16, 39, 50, 105, 110], [12, 23, 78, 83, 94, 100], [18, 29, 84, 89, 100, 106], [15, 20, 31, 37, 60, 71], [0, 11, 66, 71, 82, 88], [2, 13, 19, 42, 53, 108], [7, 13, 36, 47, 102, 107], [24, 29, 40, 46, 69, 80], [8, 63, 68, 79, 85, 108], [12, 17, 28, 34, 57, 68], [4, 10, 33, 44, 99, 104], [4, 27, 38, 93, 98, 109], [1, 24, 35, 90, 95, 106], [2, 57, 62, 73, 79, 102], [48, 53, 64, 70, 93, 104], [3, 14, 69, 74, 85, 91], [10, 22, 65, 70, 83, 86], [3, 20, 29, 65, 72, 102], [1, 15, 21, 79, 88, 105], [15, 22, 36, 42, 100, 109], [43, 52, 69, 76, 90, 96], [4, 17, 20, 55, 67, 110], [14, 19, 32, 35, 70, 82], [11, 18, 48, 60, 77, 86], [8, 17, 53, 60, 90, 102], [26, 33, 63, 75, 92, 101], [2, 37, 49, 92, 97, 110], [5, 10, 23, 26, 61, 73], [14, 21, 51, 63, 80, 89], [7, 50, 55, 68, 71, 106], [9, 16, 30, 36, 94, 103], [19, 28, 45, 52, 66, 72], [1, 44, 49, 62, 65, 100], [6, 23, 32, 68, 75, 105], [11, 20, 56, 63, 93, 105], [4, 18, 24, 82, 91, 108], [15, 27, 44, 53, 89, 96], [34, 46, 89, 94, 107, 110], [1, 14, 17, 52, 64, 107], [2, 9, 39, 51, 68, 77], [0, 17, 26, 62, 69, 99], [19, 31, 74, 79, 92, 95], [9, 26, 35, 71, 78, 108], [8, 44, 51, 81, 93, 110], [6, 12, 70, 79, 96, 103], [49, 58, 75, 82, 96, 102], [14, 23, 59, 66, 96, 108], [46, 55, 72, 79, 93, 99], [18, 30, 47, 56, 92, 99], [23, 28, 41, 44, 79, 91], [23, 30, 60, 72, 89, 98], [7, 19, 62, 67, 80, 83], [10, 19, 36, 43, 57, 63], [4, 21, 28, 42, 48, 106], [24, 36, 53, 62, 98, 105], [16, 28, 71, 76, 89, 92], [0, 7, 21, 27, 85, 94], [21, 33, 50, 59, 95, 102], [27, 39, 56, 65, 101, 108], [28, 37, 54, 61, 75, 81], [17, 22, 35, 38, 73, 85], [4, 13, 30, 37, 51, 57], [6, 36, 48, 65, 74, 110], [20, 25, 38, 41, 76, 88], [3, 33, 45, 62, 71, 107], [9, 21, 38, 47, 83, 90], [3, 61, 70, 87, 94, 108], [12, 19, 33, 39, 97, 106], [40, 49, 66, 73, 87, 93], [2, 38, 45, 75, 87, 104], [12, 24, 41, 50, 86, 93], [31, 43, 86, 91, 104, 107], [1, 18, 25, 39, 45, 103], [38, 43, 56, 59, 94, 106], [55, 64, 81, 88, 102, 108], [22, 34, 77, 82, 95, 98], [32, 39, 69, 81, 98, 107], [6, 18, 35, 44, 80, 87], [11, 14, 49, 61, 104, 109], [32, 37, 50, 53, 88, 100], [2, 7, 20, 23, 58, 70], [13, 22, 39, 46, 60, 66], [4, 16, 59, 64, 77, 80], [13, 25, 68, 73, 86, 89], [29, 34, 47, 50, 85, 97], [12, 18, 76, 85, 102, 109], [7, 16, 33, 40, 54, 60], [10, 53, 58, 71, 74, 109], [0, 6, 64, 73, 90, 97], [0, 12, 29, 38, 74, 81], [5, 8, 43, 55, 98, 103], [1, 10, 27, 34, 48, 54], [5, 14, 50, 57, 87, 99], [29, 36, 66, 78, 95, 104], [3, 10, 24, 30, 88, 97], [16, 25, 42, 49, 63, 69], [4, 47, 52, 65, 68, 103], [7, 24, 31, 45, 51, 109], [22, 31, 48, 55, 69, 75], [5, 41, 48, 78, 90, 107], [17, 24, 54, 66, 83, 92], [25, 34, 51, 58, 72, 78], [8, 13, 26, 29, 64, 76], [3, 15, 32, 41, 77, 84], [34, 43, 60, 67, 81, 87], [11, 16, 29, 32, 67, 79], [35, 42, 72, 84, 101, 110], [41, 46, 59, 62, 97, 109], [6, 13, 27, 33, 91, 100], [3, 9, 67, 76, 93, 100], [20, 27, 57, 69, 86, 95], [31, 40, 57, 64, 78, 84], [26, 31, 44, 47, 82, 94], [35, 40, 53, 56, 91, 103], [52, 61, 78, 85, 99, 105], [9, 15, 73, 82, 99, 106], [8, 15, 45, 57, 74, 83], [2, 5, 40, 52, 95, 100], [2, 11, 47, 54, 84, 96], [5, 12, 42, 54, 71, 80], [25, 37, 80, 85, 98, 101], [28, 40, 83, 88, 101, 104], [8, 11, 46, 58, 101, 106], [0, 30, 42, 59, 68, 104], [37, 46, 63, 70, 84, 90], [0, 58, 67, 84, 91, 105], [1, 13, 56, 61, 74, 77], [11, 15, 68, 70, 80, 107], [12, 21, 44, 58, 69, 108], [1, 58, 63, 76, 97, 104], [16, 21, 34, 55, 62, 70], [23, 38, 42, 95, 97, 107], [1, 12, 51, 66, 75, 98], [10, 31, 38, 46, 103, 108], [28, 33, 46, 67, 74, 82], [22, 27, 40, 61, 68, 76], [4, 25, 32, 40, 97, 102], [10, 17, 25, 82, 87, 100], [44, 46, 56, 83, 98, 102], [0, 13, 34, 41, 49, 106], [0, 23, 37, 48, 87, 102], [2, 4, 14, 41, 56, 60], [21, 36, 45, 68, 82, 93], [1, 22, 29, 37, 94, 99], [37, 42, 55, 76, 83, 91], [41, 43, 53, 80, 95, 99], [27, 42, 51, 74, 88, 99], [10, 21, 60, 75, 84, 107], [1, 6, 19, 40, 47, 55], [9, 18, 41, 55, 66, 105], [23, 25, 35, 62, 77, 81], [9, 24, 33, 56, 70, 81], [3, 42, 57, 66, 89, 103], [7, 18, 57, 72, 81, 104], [8, 22, 33, 72, 87, 96], [6, 29, 43, 54, 93, 108], [3, 26, 40, 51, 90, 105], [13, 24, 63, 78, 87, 110], [25, 30, 43, 64, 71, 79], [8, 23, 27, 80, 82, 92], [14, 29, 33, 86, 88, 98], [20, 34, 45, 84, 99, 108], [7, 28, 35, 43, 100, 105], [4, 11, 19, 76, 81, 94], [34, 39, 52, 73, 80, 88], [20, 22, 32, 59, 74, 78], [5, 20, 24, 77, 79, 89], [13, 20, 28, 85, 90, 103], [29, 31, 41, 68, 83, 87], [19, 24, 37, 58, 65, 73], [16, 23, 31, 88, 93, 106], [1, 11, 38, 53, 57, 110], [2, 16, 27, 66, 81, 90], [5, 7, 17, 44, 59, 63], [55, 60, 73, 94, 101, 109], [2, 17, 21, 74, 76, 86], [2, 10, 67, 72, 85, 106], [47, 49, 59, 86, 101, 105], [0, 39, 54, 63, 86, 100], [17, 19, 29, 56, 71, 75], [18, 33, 42, 65, 79, 90], [4, 15, 54, 69, 78, 101], [2, 6, 59, 61, 71, 98], [0, 53, 55, 65, 92, 107], [17, 31, 42, 81, 96, 105], [40, 45, 58, 79, 86, 94], [8, 35, 50, 54, 107, 109], [3, 16, 37, 44, 52, 109], [6, 21, 30, 53, 67, 78], [14, 28, 39, 78, 93, 102], [4, 9, 22, 43, 50, 58], [10, 15, 28, 49, 56, 64], [26, 28, 38, 65, 80, 84], [35, 37, 47, 74, 89, 93], [7, 14, 22, 79, 84, 97], [32, 34, 44, 71, 86, 90], [12, 27, 36, 59, 73, 84], [6, 45, 60, 69, 92, 106], [13, 18, 31, 52, 59, 67], [7, 12, 25, 46, 53, 61], [50, 52, 62, 89, 104, 108], [31, 36, 49, 70, 77, 85], [36, 51, 60, 83, 97, 108], [49, 54, 67, 88, 95, 103], [11, 25, 36, 75, 90, 99], [3, 18, 27, 50, 64, 75], [5, 9, 62, 64, 74, 101], [8, 10, 20, 47, 62, 66], [30, 45, 54, 77, 91, 102], [4, 61, 66, 79, 100, 107], [19, 26, 34, 91, 96, 109], [5, 32, 47, 51, 104, 106], [3, 56, 58, 68, 95, 110], [8, 12, 65, 67, 77, 104], [5, 19, 30, 69, 84, 93], [20, 35, 39, 92, 94, 104], [0, 9, 32, 46, 57, 96], [14, 18, 71, 73, 83, 110], [46, 51, 64, 85, 92, 100], [15, 30, 39, 62, 76, 87], [17, 32, 36, 89, 91, 101], [26, 41, 45, 98, 100, 110], [14, 16, 26, 53, 68, 72], [9, 48, 63, 72, 95, 109], [33, 48, 57, 80, 94, 105], [1, 8, 16, 73, 78, 91], [0, 15, 24, 47, 61, 72], [43, 48, 61, 82, 89, 97], [11, 13, 23, 50, 65, 69], [5, 13, 70, 75, 88, 109], [7, 64, 69, 82, 103, 110], [11, 26, 30, 83, 85, 95], [52, 57, 70, 91, 98, 106], [3, 12, 35, 49, 60, 99], [24, 39, 48, 71, 85, 96], [2, 29, 44, 48, 101, 103], [6, 15, 38, 52, 63, 102], [38, 40, 50, 77, 92, 96]]
\item 1 \{0=2220, 1=32856, 2=571428, 3=2169384, 4=2352312\} [[8, 39, 40, 41, 42, 70], [0, 1, 2, 3, 31, 80], [30, 31, 32, 33, 61, 110], [6, 7, 8, 9, 37, 86], [23, 54, 55, 56, 57, 85], [15, 16, 17, 18, 46, 95], [27, 28, 29, 30, 58, 107], [24, 25, 26, 27, 55, 104], [26, 57, 58, 59, 60, 88], [29, 60, 61, 62, 63, 91], [25, 74, 105, 106, 107, 108], [41, 72, 73, 74, 75, 103], [7, 56, 87, 88, 89, 90], [10, 59, 90, 91, 92, 93], [13, 62, 93, 94, 95, 96], [22, 71, 102, 103, 104, 105], [9, 10, 11, 12, 40, 89], [44, 75, 76, 77, 78, 106], [18, 19, 20, 21, 49, 98], [11, 42, 43, 44, 45, 73], [3, 4, 5, 6, 34, 83], [16, 65, 96, 97, 98, 99], [19, 68, 99, 100, 101, 102], [17, 48, 49, 50, 51, 79], [47, 78, 79, 80, 81, 109], [20, 51, 52, 53, 54, 82], [12, 13, 14, 15, 43, 92], [0, 28, 77, 108, 109, 110], [4, 53, 84, 85, 86, 87], [1, 50, 81, 82, 83, 84], [32, 63, 64, 65, 66, 94], [2, 33, 34, 35, 36, 64], [21, 22, 23, 24, 52, 101], [38, 69, 70, 71, 72, 100], [14, 45, 46, 47, 48, 76], [35, 66, 67, 68, 69, 97], [5, 36, 37, 38, 39, 67], [11, 54, 58, 61, 87, 92], [24, 29, 59, 102, 106, 109], [33, 37, 40, 66, 71, 101], [42, 46, 49, 75, 80, 110], [2, 45, 49, 52, 78, 83], [6, 10, 13, 39, 44, 74], [12, 17, 47, 90, 94, 97], [23, 66, 70, 73, 99, 104], [18, 23, 53, 96, 100, 103], [15, 19, 22, 48, 53, 83], [6, 11, 41, 84, 88, 91], [0, 5, 35, 78, 82, 85], [1, 27, 32, 62, 105, 109], [26, 69, 73, 76, 102, 107], [15, 20, 50, 93, 97, 100], [3, 7, 10, 36, 41, 71], [21, 26, 56, 99, 103, 106], [30, 34, 37, 63, 68, 98], [18, 22, 25, 51, 56, 86], [27, 31, 34, 60, 65, 95], [24, 28, 31, 57, 62, 92], [29, 72, 76, 79, 105, 110], [14, 57, 61, 64, 90, 95], [9, 14, 44, 87, 91, 94], [8, 51, 55, 58, 84, 89], [39, 43, 46, 72, 77, 107], [36, 40, 43, 69, 74, 104], [12, 16, 19, 45, 50, 80], [2, 32, 75, 79, 82, 108], [20, 63, 67, 70, 96, 101], [3, 8, 38, 81, 85, 88], [0, 4, 7, 33, 38, 68], [21, 25, 28, 54, 59, 89], [1, 4, 30, 35, 65, 108], [5, 48, 52, 55, 81, 86], [17, 60, 64, 67, 93, 98], [9, 13, 16, 42, 47, 77], [32, 38, 53, 77, 89, 102], [8, 13, 22, 61, 67, 82], [5, 29, 41, 54, 95, 101], [11, 24, 65, 71, 86, 110], [14, 38, 50, 63, 104, 110], [4, 19, 56, 61, 70, 109], [8, 20, 33, 74, 80, 95], [22, 27, 66, 78, 84, 93], [3, 15, 21, 30, 70, 75], [10, 15, 54, 66, 72, 81], [28, 33, 72, 84, 90, 99], [3, 42, 54, 60, 69, 109], [30, 42, 48, 57, 97, 102], [0, 9, 49, 54, 93, 105], [5, 17, 30, 71, 77, 92], [22, 28, 43, 80, 85, 94], [2, 7, 16, 55, 61, 76], [23, 35, 48, 89, 95, 110], [0, 6, 15, 55, 60, 99], [23, 28, 37, 76, 82, 97], [12, 53, 59, 74, 98, 110], [14, 26, 39, 80, 86, 101], [8, 14, 29, 53, 65, 78], [18, 30, 36, 45, 85, 90], [1, 38, 43, 52, 91, 97], [20, 26, 41, 65, 77, 90], [28, 34, 49, 86, 91, 100], [5, 20, 44, 56, 69, 110], [0, 40, 45, 84, 96, 102], [14, 19, 28, 67, 73, 88], [8, 21, 62, 68, 83, 107], [19, 24, 63, 75, 81, 90], [11, 23, 36, 77, 83, 98], [11, 16, 25, 64, 70, 85], [9, 50, 56, 71, 95, 107], [19, 25, 40, 77, 82, 91], [14, 20, 35, 59, 71, 84], [6, 46, 51, 90, 102, 108], [1, 7, 22, 59, 64, 73], [13, 50, 55, 64, 103, 109], [4, 13, 52, 58, 73, 110], [3, 12, 52, 57, 96, 108], [4, 41, 46, 55, 94, 100], [24, 36, 42, 51, 91, 96], [2, 14, 27, 68, 74, 89], [2, 15, 56, 62, 77, 101], [9, 21, 27, 36, 76, 81], [23, 29, 44, 68, 80, 93], [12, 24, 30, 39, 79, 84], [5, 10, 19, 58, 64, 79], [3, 43, 48, 87, 99, 105], [2, 17, 41, 53, 66, 107], [25, 30, 69, 81, 87, 96], [6, 12, 21, 61, 66, 105], [15, 27, 33, 42, 82, 87], [13, 18, 57, 69, 75, 84], [31, 36, 75, 87, 93, 102], [26, 32, 47, 71, 83, 96], [29, 35, 50, 74, 86, 99], [6, 47, 53, 68, 92, 104], [16, 22, 37, 74, 79, 88], [21, 33, 39, 48, 88, 93], [1, 10, 49, 55, 70, 107], [7, 44, 49, 58, 97, 103], [1, 40, 46, 61, 98, 103], [3, 9, 18, 58, 63, 102], [35, 40, 49, 88, 94, 109], [7, 13, 28, 65, 70, 79], [17, 29, 42, 83, 89, 104], [27, 39, 45, 54, 94, 99], [38, 44, 59, 83, 95, 108], [0, 12, 18, 27, 67, 72], [7, 12, 51, 63, 69, 78], [0, 39, 51, 57, 66, 106], [3, 44, 50, 65, 89, 101], [13, 19, 34, 71, 76, 85], [4, 9, 48, 60, 66, 75], [0, 41, 47, 62, 86, 98], [37, 42, 81, 93, 99, 108], [11, 17, 32, 56, 68, 81], [37, 43, 58, 95, 100, 109], [2, 8, 23, 47, 59, 72], [35, 41, 56, 80, 92, 105], [29, 34, 43, 82, 88, 103], [4, 10, 25, 62, 67, 76], [1, 6, 45, 57, 63, 72], [7, 46, 52, 67, 104, 109], [16, 21, 60, 72, 78, 87], [34, 39, 78, 90, 96, 105], [17, 23, 38, 62, 74, 87], [34, 40, 55, 92, 97, 106], [6, 18, 24, 33, 73, 78], [11, 35, 47, 60, 101, 107], [25, 31, 46, 83, 88, 97], [20, 25, 34, 73, 79, 94], [10, 47, 52, 61, 100, 106], [31, 37, 52, 89, 94, 103], [20, 32, 45, 86, 92, 107], [5, 11, 26, 50, 62, 75], [5, 18, 59, 65, 80, 104], [8, 32, 44, 57, 98, 104], [32, 37, 46, 85, 91, 106], [4, 43, 49, 64, 101, 106], [9, 15, 24, 64, 69, 108], [36, 48, 54, 63, 103, 108], [26, 31, 40, 79, 85, 100], [10, 16, 31, 68, 73, 82], [33, 45, 51, 60, 100, 105], [1, 16, 53, 58, 67, 106], [2, 26, 38, 51, 92, 98], [17, 22, 31, 70, 76, 91], [30, 40, 54, 62, 73, 80], [0, 42, 53, 56, 76, 94], [1, 18, 60, 71, 74, 94], [11, 13, 21, 37, 80, 102], [33, 43, 57, 65, 76, 83], [3, 11, 22, 29, 90, 100], [26, 48, 68, 70, 78, 94], [5, 16, 23, 84, 94, 108], [1, 9, 25, 68, 90, 110], [21, 41, 43, 51, 67, 110], [12, 54, 65, 68, 88, 106], [1, 19, 36, 78, 89, 92], [9, 19, 33, 41, 52, 59], [27, 37, 51, 59, 70, 77], [18, 38, 40, 48, 64, 107], [19, 37, 54, 96, 107, 110], [39, 50, 53, 73, 91, 108], [32, 54, 74, 76, 84, 100], [10, 53, 75, 95, 97, 105], [6, 14, 25, 32, 93, 103], [3, 13, 27, 35, 46, 53], [16, 34, 51, 93, 104, 107], [5, 7, 15, 31, 74, 96], [20, 42, 62, 64, 72, 88], [1, 15, 23, 34, 41, 102], [8, 10, 18, 34, 77, 99], [17, 39, 59, 61, 69, 85], [10, 17, 78, 88, 102, 110], [21, 31, 45, 53, 64, 71], [15, 25, 39, 47, 58, 65], [1, 44, 66, 86, 88, 96], [23, 45, 65, 67, 75, 91], [12, 23, 26, 46, 64, 81], [3, 23, 25, 33, 49, 92], [4, 21, 63, 74, 77, 97], [7, 24, 66, 77, 80, 100], [7, 25, 42, 84, 95, 98], [0, 20, 22, 30, 46, 89], [24, 34, 48, 56, 67, 74], [9, 51, 62, 65, 85, 103], [2, 5, 25, 43, 60, 102], [3, 14, 17, 37, 55, 72], [15, 57, 68, 71, 91, 109], [10, 28, 45, 87, 98, 101], [45, 55, 69, 77, 88, 95], [10, 27, 69, 80, 83, 103], [7, 14, 75, 85, 99, 107], [33, 44, 47, 67, 85, 102], [24, 35, 38, 58, 76, 93], [36, 47, 50, 70, 88, 105], [12, 20, 31, 38, 99, 109], [12, 32, 34, 42, 58, 101], [6, 26, 28, 36, 52, 95], [27, 38, 41, 61, 79, 96], [39, 49, 63, 71, 82, 89], [9, 17, 28, 35, 96, 106], [38, 60, 80, 82, 90, 106], [0, 11, 14, 34, 52, 69], [2, 22, 40, 57, 99, 110], [36, 46, 60, 68, 79, 86], [13, 30, 72, 83, 86, 106], [15, 35, 37, 45, 61, 104], [60, 70, 84, 92, 103, 110], [15, 26, 29, 49, 67, 84], [0, 8, 19, 26, 87, 97], [6, 16, 30, 38, 49, 56], [11, 33, 53, 55, 63, 79], [16, 33, 75, 86, 89, 109], [0, 16, 59, 81, 101, 103], [5, 66, 76, 90, 98, 109], [5, 8, 28, 46, 63, 105], [2, 4, 12, 28, 71, 93], [8, 30, 50, 52, 60, 76], [7, 21, 29, 40, 47, 108], [6, 17, 20, 40, 58, 75], [9, 20, 23, 43, 61, 78], [3, 19, 62, 84, 104, 106], [6, 48, 59, 62, 82, 100], [2, 24, 44, 46, 54, 70], [17, 19, 27, 43, 86, 108], [4, 11, 72, 82, 96, 104], [21, 32, 35, 55, 73, 90], [29, 51, 71, 73, 81, 97], [35, 57, 77, 79, 87, 103], [0, 10, 24, 32, 43, 50], [2, 63, 73, 87, 95, 106], [14, 36, 56, 58, 66, 82], [3, 45, 56, 59, 79, 97], [13, 31, 48, 90, 101, 104], [6, 22, 65, 87, 107, 109], [54, 64, 78, 86, 97, 104], [9, 29, 31, 39, 55, 98], [2, 13, 20, 81, 91, 105], [18, 29, 32, 52, 70, 87], [30, 41, 44, 64, 82, 99], [18, 28, 42, 50, 61, 68], [13, 56, 78, 98, 100, 108], [42, 52, 66, 74, 85, 92], [4, 47, 69, 89, 91, 99], [57, 67, 81, 89, 100, 107], [51, 61, 75, 83, 94, 101], [48, 58, 72, 80, 91, 98], [41, 63, 83, 85, 93, 109], [5, 27, 47, 49, 57, 73], [4, 18, 26, 37, 44, 105], [14, 16, 24, 40, 83, 105], [12, 22, 36, 44, 55, 62], [4, 22, 39, 81, 92, 95], [8, 11, 31, 49, 66, 108], [7, 50, 72, 92, 94, 102], [1, 8, 69, 79, 93, 101], [5, 9, 22, 32, 72, 97], [2, 42, 67, 86, 90, 103], [10, 29, 33, 46, 56, 96], [21, 46, 65, 69, 82, 92], [17, 21, 34, 44, 84, 109], [15, 40, 59, 63, 76, 86], [24, 49, 68, 72, 85, 95], [11, 15, 28, 38, 78, 103], [22, 41, 45, 58, 68, 108], [13, 32, 36, 49, 59, 99], [9, 34, 53, 57, 70, 80], [1, 20, 24, 37, 47, 87], [7, 26, 30, 43, 53, 93], [3, 16, 26, 66, 91, 110], [36, 61, 80, 84, 97, 107], [7, 17, 57, 82, 101, 105], [12, 37, 56, 60, 73, 83], [6, 31, 50, 54, 67, 77], [27, 52, 71, 75, 88, 98], [18, 43, 62, 66, 79, 89], [4, 14, 54, 79, 98, 102], [3, 28, 47, 51, 64, 74], [0, 25, 44, 48, 61, 71], [14, 18, 31, 41, 81, 106], [8, 12, 25, 35, 75, 100], [4, 23, 27, 40, 50, 90], [10, 20, 60, 85, 104, 108], [39, 64, 83, 87, 100, 110], [8, 48, 73, 92, 96, 109], [30, 55, 74, 78, 91, 101], [5, 45, 70, 89, 93, 106], [0, 13, 23, 63, 88, 107], [2, 6, 19, 29, 69, 94], [16, 35, 39, 52, 62, 102], [19, 38, 42, 55, 65, 105], [33, 58, 77, 81, 94, 104], [1, 11, 51, 76, 95, 99], [6, 23, 42, 71, 79, 106], [12, 41, 49, 76, 87, 104], [8, 27, 56, 64, 91, 102], [22, 33, 50, 69, 98, 106], [8, 16, 43, 54, 71, 90], [10, 21, 38, 57, 86, 94], [23, 31, 58, 69, 86, 105], [2, 10, 37, 48, 65, 84], [4, 31, 42, 59, 78, 107], [25, 36, 53, 72, 101, 109], [3, 20, 39, 68, 76, 103], [9, 38, 46, 73, 84, 101], [7, 34, 45, 62, 81, 110], [26, 34, 61, 72, 89, 108], [1, 12, 29, 48, 77, 85], [6, 35, 43, 70, 81, 98], [15, 44, 52, 79, 90, 107], [11, 19, 46, 57, 74, 93], [0, 29, 37, 64, 75, 92], [2, 21, 50, 58, 85, 96], [19, 30, 47, 66, 95, 103], [14, 33, 62, 70, 97, 108], [14, 22, 49, 60, 77, 96], [13, 24, 41, 60, 89, 97], [5, 24, 53, 61, 88, 99], [9, 26, 45, 74, 82, 109], [11, 30, 59, 67, 94, 105], [5, 13, 40, 51, 68, 87], [18, 47, 55, 82, 93, 110], [4, 15, 32, 51, 80, 88], [0, 17, 36, 65, 73, 100], [3, 32, 40, 67, 78, 95], [17, 25, 52, 63, 80, 99], [16, 27, 44, 63, 92, 100], [7, 18, 35, 54, 83, 91], [20, 28, 55, 66, 83, 102], [1, 28, 39, 56, 75, 104], [0, 21, 79, 91, 95, 104], [16, 28, 32, 41, 48, 69], [6, 64, 76, 80, 89, 96], [6, 27, 85, 97, 101, 110], [12, 70, 82, 86, 95, 102], [13, 25, 29, 38, 45, 66], [1, 5, 14, 21, 42, 100], [7, 19, 23, 32, 39, 60], [18, 76, 88, 92, 101, 108], [28, 40, 44, 53, 60, 81], [0, 58, 70, 74, 83, 90], [40, 52, 56, 65, 72, 93], [19, 31, 35, 44, 51, 72], [3, 24, 82, 94, 98, 107], [8, 15, 36, 94, 106, 110], [1, 13, 17, 26, 33, 54], [10, 14, 23, 30, 51, 109], [46, 58, 62, 71, 78, 99], [2, 11, 18, 39, 97, 109], [34, 46, 50, 59, 66, 87], [49, 61, 65, 74, 81, 102], [55, 67, 71, 80, 87, 108], [9, 67, 79, 83, 92, 99], [2, 9, 30, 88, 100, 104], [4, 16, 20, 29, 36, 57], [7, 11, 20, 27, 48, 106], [52, 64, 68, 77, 84, 105], [43, 55, 59, 68, 75, 96], [5, 12, 33, 91, 103, 107], [25, 37, 41, 50, 57, 78], [37, 49, 53, 62, 69, 90], [31, 43, 47, 56, 63, 84], [15, 73, 85, 89, 98, 105], [22, 34, 38, 47, 54, 75], [4, 8, 17, 24, 45, 103], [10, 22, 26, 35, 42, 63], [3, 61, 73, 77, 86, 93]]
\item 1 \{1=40848, 2=522810, 3=2164500, 4=2400042\} [[8, 39, 40, 41, 42, 70], [0, 1, 2, 3, 31, 80], [30, 31, 32, 33, 61, 110], [6, 7, 8, 9, 37, 86], [23, 54, 55, 56, 57, 85], [15, 16, 17, 18, 46, 95], [27, 28, 29, 30, 58, 107], [24, 25, 26, 27, 55, 104], [26, 57, 58, 59, 60, 88], [29, 60, 61, 62, 63, 91], [25, 74, 105, 106, 107, 108], [41, 72, 73, 74, 75, 103], [7, 56, 87, 88, 89, 90], [10, 59, 90, 91, 92, 93], [13, 62, 93, 94, 95, 96], [22, 71, 102, 103, 104, 105], [9, 10, 11, 12, 40, 89], [44, 75, 76, 77, 78, 106], [18, 19, 20, 21, 49, 98], [11, 42, 43, 44, 45, 73], [3, 4, 5, 6, 34, 83], [16, 65, 96, 97, 98, 99], [19, 68, 99, 100, 101, 102], [17, 48, 49, 50, 51, 79], [47, 78, 79, 80, 81, 109], [20, 51, 52, 53, 54, 82], [12, 13, 14, 15, 43, 92], [0, 28, 77, 108, 109, 110], [4, 53, 84, 85, 86, 87], [1, 50, 81, 82, 83, 84], [32, 63, 64, 65, 66, 94], [2, 33, 34, 35, 36, 64], [21, 22, 23, 24, 52, 101], [38, 69, 70, 71, 72, 100], [14, 45, 46, 47, 48, 76], [35, 66, 67, 68, 69, 97], [5, 36, 37, 38, 39, 67], [11, 41, 54, 87, 91, 94], [24, 28, 31, 59, 89, 102], [15, 19, 22, 50, 80, 93], [9, 42, 46, 49, 77, 107], [8, 38, 51, 84, 88, 91], [6, 39, 43, 46, 74, 104], [3, 7, 10, 38, 68, 81], [1, 4, 32, 62, 75, 108], [0, 33, 37, 40, 68, 98], [20, 33, 66, 70, 73, 101], [6, 10, 13, 41, 71, 84], [5, 18, 51, 55, 58, 86], [2, 15, 48, 52, 55, 83], [17, 47, 60, 93, 97, 100], [30, 34, 37, 65, 95, 108], [26, 56, 69, 102, 106, 109], [23, 53, 66, 99, 103, 106], [17, 30, 63, 67, 70, 98], [27, 31, 34, 62, 92, 105], [20, 50, 63, 96, 100, 103], [26, 39, 72, 76, 79, 107], [11, 24, 57, 61, 64, 92], [29, 42, 75, 79, 82, 110], [3, 36, 40, 43, 71, 101], [0, 4, 7, 35, 65, 78], [14, 44, 57, 90, 94, 97], [14, 27, 60, 64, 67, 95], [23, 36, 69, 73, 76, 104], [5, 35, 48, 81, 85, 88], [8, 21, 54, 58, 61, 89], [9, 13, 16, 44, 74, 87], [18, 22, 25, 53, 83, 96], [12, 16, 19, 47, 77, 90], [21, 25, 28, 56, 86, 99], [2, 32, 45, 78, 82, 85], [12, 45, 49, 52, 80, 110], [1, 29, 59, 72, 105, 109], [0, 10, 22, 44, 79, 85], [13, 21, 32, 42, 60, 105], [27, 46, 54, 65, 75, 93], [2, 13, 17, 77, 99, 104], [0, 11, 21, 39, 84, 103], [9, 19, 31, 53, 88, 94], [11, 46, 52, 78, 88, 100], [18, 37, 45, 56, 66, 84], [5, 27, 32, 41, 52, 56], [53, 75, 80, 89, 100, 104], [6, 51, 70, 78, 89, 99], [3, 14, 24, 42, 87, 106], [20, 55, 61, 87, 97, 109], [21, 31, 43, 65, 100, 106], [14, 36, 41, 50, 61, 65], [32, 54, 59, 68, 79, 83], [4, 26, 61, 67, 93, 103], [24, 43, 51, 62, 72, 90], [6, 16, 28, 50, 85, 91], [29, 51, 56, 65, 76, 80], [8, 30, 35, 44, 55, 59], [2, 37, 43, 69, 79, 91], [2, 11, 22, 26, 86, 108], [26, 48, 53, 62, 73, 77], [12, 31, 39, 50, 60, 78], [18, 28, 40, 62, 97, 103], [4, 30, 40, 52, 74, 109], [8, 19, 23, 83, 105, 110], [34, 40, 66, 76, 88, 110], [44, 66, 71, 80, 91, 95], [20, 42, 47, 56, 67, 71], [56, 78, 83, 92, 103, 107], [2, 62, 84, 89, 98, 109], [15, 60, 79, 87, 98, 108], [50, 72, 77, 86, 97, 101], [3, 22, 30, 41, 51, 69], [36, 55, 63, 74, 84, 102], [1, 27, 37, 49, 71, 106], [4, 16, 38, 73, 79, 105], [17, 39, 44, 53, 64, 68], [33, 52, 60, 71, 81, 99], [28, 34, 60, 70, 82, 104], [2, 12, 30, 75, 94, 102], [35, 57, 62, 71, 82, 86], [6, 25, 33, 44, 54, 72], [0, 18, 63, 82, 90, 101], [12, 17, 26, 37, 41, 101], [7, 29, 64, 70, 96, 106], [1, 5, 65, 87, 92, 101], [31, 37, 63, 73, 85, 107], [59, 81, 86, 95, 106, 110], [14, 49, 55, 81, 91, 103], [6, 11, 20, 31, 35, 95], [7, 13, 39, 49, 61, 83], [4, 12, 23, 33, 51, 96], [9, 14, 23, 34, 38, 98], [10, 16, 42, 52, 64, 86], [1, 23, 58, 64, 90, 100], [16, 24, 35, 45, 63, 108], [8, 43, 49, 75, 85, 97], [3, 48, 67, 75, 86, 96], [21, 40, 48, 59, 69, 87], [6, 24, 69, 88, 96, 107], [38, 60, 65, 74, 85, 89], [5, 40, 46, 72, 82, 94], [3, 13, 25, 47, 82, 88], [11, 33, 38, 47, 58, 62], [13, 19, 45, 55, 67, 89], [0, 19, 27, 38, 48, 66], [7, 19, 41, 76, 82, 108], [3, 21, 66, 85, 93, 104], [9, 28, 36, 47, 57, 75], [15, 25, 37, 59, 94, 100], [8, 18, 36, 81, 100, 108], [15, 20, 29, 40, 44, 104], [1, 7, 33, 43, 55, 77], [4, 8, 68, 90, 95, 104], [16, 22, 48, 58, 70, 92], [19, 25, 51, 61, 73, 95], [9, 54, 73, 81, 92, 102], [6, 17, 27, 45, 90, 109], [7, 11, 71, 93, 98, 107], [39, 58, 66, 77, 87, 105], [22, 28, 54, 64, 76, 98], [0, 5, 14, 25, 29, 89], [4, 10, 36, 46, 58, 80], [9, 27, 72, 91, 99, 110], [24, 34, 46, 68, 103, 109], [15, 34, 42, 53, 63, 81], [41, 63, 68, 77, 88, 92], [12, 57, 76, 84, 95, 105], [21, 26, 35, 46, 50, 110], [10, 32, 67, 73, 99, 109], [42, 61, 69, 80, 90, 108], [2, 24, 29, 38, 49, 53], [3, 8, 17, 28, 32, 92], [25, 31, 57, 67, 79, 101], [0, 45, 64, 72, 83, 93], [10, 14, 74, 96, 101, 110], [18, 23, 32, 43, 47, 107], [30, 49, 57, 68, 78, 96], [5, 15, 33, 78, 97, 105], [1, 9, 20, 30, 48, 93], [5, 16, 20, 80, 102, 107], [47, 69, 74, 83, 94, 98], [1, 13, 35, 70, 76, 102], [7, 15, 26, 36, 54, 99], [10, 18, 29, 39, 57, 102], [23, 45, 50, 59, 70, 74], [12, 22, 34, 56, 91, 97], [17, 52, 58, 84, 94, 106], [1, 16, 40, 53, 56, 61], [0, 9, 51, 67, 100, 105], [21, 37, 70, 75, 81, 90], [29, 35, 47, 54, 84, 101], [17, 23, 35, 42, 72, 89], [0, 16, 49, 54, 60, 69], [26, 32, 44, 51, 81, 98], [13, 28, 52, 65, 68, 73], [5, 11, 23, 30, 60, 77], [19, 32, 35, 40, 91, 106], [10, 43, 48, 54, 63, 105], [11, 14, 19, 70, 85, 109], [0, 42, 58, 91, 96, 102], [5, 8, 13, 64, 79, 103], [21, 38, 77, 83, 95, 102], [2, 5, 10, 61, 76, 100], [10, 15, 21, 30, 72, 88], [7, 20, 23, 28, 79, 94], [10, 34, 47, 50, 55, 106], [10, 25, 49, 62, 65, 70], [1, 34, 39, 45, 54, 96], [43, 58, 82, 95, 98, 103], [27, 43, 76, 81, 87, 96], [34, 49, 73, 86, 89, 94], [19, 34, 58, 71, 74, 79], [10, 23, 26, 31, 82, 97], [31, 36, 42, 51, 93, 109], [4, 19, 43, 56, 59, 64], [9, 26, 65, 71, 83, 90], [2, 41, 47, 59, 66, 96], [16, 21, 27, 36, 78, 94], [31, 46, 70, 83, 86, 91], [40, 55, 79, 92, 95, 100], [3, 19, 52, 57, 63, 72], [37, 52, 76, 89, 92, 97], [35, 41, 53, 60, 90, 107], [32, 38, 50, 57, 87, 104], [1, 6, 12, 21, 63, 79], [23, 29, 41, 48, 78, 95], [11, 18, 48, 65, 104, 110], [16, 29, 32, 37, 88, 103], [28, 33, 39, 48, 90, 106], [22, 37, 61, 74, 77, 82], [15, 32, 71, 77, 89, 96], [13, 26, 29, 34, 85, 100], [4, 55, 70, 94, 107, 110], [20, 26, 38, 45, 75, 92], [7, 40, 45, 51, 60, 102], [13, 18, 24, 33, 75, 91], [4, 28, 41, 44, 49, 100], [38, 44, 56, 63, 93, 110], [24, 41, 80, 86, 98, 105], [4, 37, 42, 48, 57, 99], [3, 33, 50, 89, 95, 107], [7, 22, 46, 59, 62, 67], [6, 23, 62, 68, 80, 87], [12, 29, 68, 74, 86, 93], [2, 14, 21, 51, 68, 107], [6, 36, 53, 92, 98, 110], [4, 17, 20, 25, 76, 91], [0, 17, 56, 62, 74, 81], [36, 52, 85, 90, 96, 105], [4, 9, 15, 24, 66, 82], [5, 12, 42, 59, 98, 104], [14, 53, 59, 71, 78, 108], [11, 50, 56, 68, 75, 105], [0, 6, 15, 57, 73, 106], [0, 30, 47, 86, 92, 104], [25, 30, 36, 45, 87, 103], [5, 44, 50, 62, 69, 99], [14, 20, 32, 39, 69, 86], [22, 35, 38, 43, 94, 109], [19, 24, 30, 39, 81, 97], [18, 34, 67, 72, 78, 87], [3, 20, 59, 65, 77, 84], [3, 12, 54, 70, 103, 108], [18, 35, 74, 80, 92, 99], [2, 7, 58, 73, 97, 110], [28, 43, 67, 80, 83, 88], [7, 31, 44, 47, 52, 103], [13, 37, 50, 53, 58, 109], [46, 61, 85, 98, 101, 106], [3, 45, 61, 94, 99, 105], [3, 9, 18, 60, 76, 109], [6, 22, 55, 60, 66, 75], [13, 46, 51, 57, 66, 108], [12, 28, 61, 66, 72, 81], [39, 55, 88, 93, 99, 108], [1, 25, 38, 41, 46, 97], [2, 8, 20, 27, 57, 74], [16, 31, 55, 68, 71, 76], [24, 40, 73, 78, 84, 93], [30, 46, 79, 84, 90, 99], [6, 48, 64, 97, 102, 108], [49, 64, 88, 101, 104, 109], [2, 9, 39, 56, 95, 101], [22, 27, 33, 42, 84, 100], [5, 17, 24, 54, 71, 110], [9, 25, 58, 63, 69, 78], [25, 40, 64, 77, 80, 85], [27, 44, 83, 89, 101, 108], [8, 47, 53, 65, 72, 102], [8, 11, 16, 67, 82, 106], [8, 15, 45, 62, 101, 107], [11, 17, 29, 36, 66, 83], [1, 52, 67, 91, 104, 107], [8, 14, 26, 33, 63, 80], [1, 14, 17, 22, 73, 88], [15, 31, 64, 69, 75, 84], [33, 49, 82, 87, 93, 102], [7, 12, 18, 27, 69, 85], [2, 18, 44, 70, 88, 105], [22, 40, 57, 65, 81, 107], [10, 28, 45, 53, 69, 95], [23, 49, 67, 84, 92, 108], [12, 38, 64, 82, 99, 107], [19, 37, 54, 62, 78, 104], [6, 14, 30, 56, 82, 100], [5, 31, 49, 66, 74, 90], [0, 8, 24, 50, 76, 94], [11, 37, 55, 72, 80, 96], [9, 17, 33, 59, 85, 103], [8, 34, 52, 69, 77, 93], [7, 25, 42, 50, 66, 92], [15, 41, 67, 85, 102, 110], [2, 28, 46, 63, 71, 87], [9, 35, 61, 79, 96, 104], [20, 46, 64, 81, 89, 105], [3, 29, 55, 73, 90, 98], [12, 20, 36, 62, 88, 106], [25, 43, 60, 68, 84, 110], [6, 32, 58, 76, 93, 101], [10, 27, 35, 51, 77, 103], [13, 30, 38, 54, 80, 106], [14, 40, 58, 75, 83, 99], [16, 34, 51, 59, 75, 101], [15, 23, 39, 65, 91, 109], [5, 21, 47, 73, 91, 108], [4, 21, 29, 45, 71, 97], [0, 26, 52, 70, 87, 95], [1, 19, 36, 44, 60, 86], [16, 33, 41, 57, 83, 109], [3, 11, 27, 53, 79, 97], [13, 31, 48, 56, 72, 98], [7, 24, 32, 48, 74, 100], [1, 18, 26, 42, 68, 94], [4, 22, 39, 47, 63, 89], [17, 43, 61, 78, 86, 102], [42, 54, 74, 88, 95, 97], [21, 33, 53, 67, 74, 76], [8, 22, 29, 31, 87, 99], [3, 15, 35, 49, 56, 58], [15, 27, 47, 61, 68, 70], [0, 20, 34, 41, 43, 99], [2, 16, 23, 25, 81, 93], [11, 25, 32, 34, 90, 102], [14, 28, 35, 37, 93, 105], [1, 57, 69, 89, 103, 110], [5, 7, 63, 75, 95, 109], [9, 21, 41, 55, 62, 64], [0, 12, 32, 46, 53, 55], [5, 19, 26, 28, 84, 96], [9, 29, 43, 50, 52, 108], [2, 4, 60, 72, 92, 106], [6, 18, 38, 52, 59, 61], [4, 11, 13, 69, 81, 101], [17, 31, 38, 40, 96, 108], [48, 60, 80, 94, 101, 103], [18, 30, 50, 64, 71, 73], [51, 63, 83, 97, 104, 106], [27, 39, 59, 73, 80, 82], [33, 45, 65, 79, 86, 88], [12, 24, 44, 58, 65, 67], [6, 26, 40, 47, 49, 105], [3, 23, 37, 44, 46, 102], [45, 57, 77, 91, 98, 100], [54, 66, 86, 100, 107, 109], [39, 51, 71, 85, 92, 94], [36, 48, 68, 82, 89, 91], [30, 42, 62, 76, 83, 85], [24, 36, 56, 70, 77, 79], [13, 20, 22, 78, 90, 110], [7, 14, 16, 72, 84, 104], [10, 17, 19, 75, 87, 107], [1, 8, 10, 66, 78, 98], [5, 43, 53, 57, 70, 93], [7, 30, 53, 91, 101, 105], [13, 23, 27, 40, 63, 86], [2, 6, 19, 42, 65, 103], [4, 27, 50, 88, 98, 102], [6, 29, 67, 77, 81, 94], [34, 44, 48, 61, 84, 107], [17, 55, 65, 69, 82, 105], [28, 38, 42, 55, 78, 101], [7, 17, 21, 34, 57, 80], [21, 44, 82, 92, 96, 109], [37, 47, 51, 64, 87, 110], [1, 24, 47, 85, 95, 99], [22, 32, 36, 49, 72, 95], [9, 32, 70, 80, 84, 97], [3, 26, 64, 74, 78, 91], [18, 41, 79, 89, 93, 106], [12, 35, 73, 83, 87, 100], [20, 58, 68, 72, 85, 108], [8, 46, 56, 60, 73, 96], [31, 41, 45, 58, 81, 104], [25, 35, 39, 52, 75, 98], [15, 38, 76, 86, 90, 103], [0, 13, 36, 59, 97, 107], [3, 16, 39, 62, 100, 110], [10, 20, 24, 37, 60, 83], [16, 26, 30, 43, 66, 89], [8, 12, 25, 48, 71, 109], [19, 29, 33, 46, 69, 92], [0, 23, 61, 71, 75, 88], [4, 14, 18, 31, 54, 77], [5, 9, 22, 45, 68, 106], [11, 49, 59, 63, 76, 99], [1, 11, 15, 28, 51, 74], [2, 40, 50, 54, 67, 90], [14, 52, 62, 66, 79, 102], [10, 33, 56, 94, 104, 108]]
\item 1 \{0=1110, 1=43512, 2=566766, 3=2209788, 4=2307024\} [[8, 39, 40, 41, 42, 70], [0, 1, 2, 3, 31, 80], [30, 31, 32, 33, 61, 110], [6, 7, 8, 9, 37, 86], [23, 54, 55, 56, 57, 85], [15, 16, 17, 18, 46, 95], [27, 28, 29, 30, 58, 107], [24, 25, 26, 27, 55, 104], [26, 57, 58, 59, 60, 88], [29, 60, 61, 62, 63, 91], [25, 74, 105, 106, 107, 108], [41, 72, 73, 74, 75, 103], [7, 56, 87, 88, 89, 90], [10, 59, 90, 91, 92, 93], [13, 62, 93, 94, 95, 96], [22, 71, 102, 103, 104, 105], [9, 10, 11, 12, 40, 89], [44, 75, 76, 77, 78, 106], [18, 19, 20, 21, 49, 98], [11, 42, 43, 44, 45, 73], [3, 4, 5, 6, 34, 83], [16, 65, 96, 97, 98, 99], [19, 68, 99, 100, 101, 102], [17, 48, 49, 50, 51, 79], [47, 78, 79, 80, 81, 109], [20, 51, 52, 53, 54, 82], [12, 13, 14, 15, 43, 92], [0, 28, 77, 108, 109, 110], [4, 53, 84, 85, 86, 87], [1, 50, 81, 82, 83, 84], [32, 63, 64, 65, 66, 94], [2, 33, 34, 35, 36, 64], [21, 22, 23, 24, 52, 101], [38, 69, 70, 71, 72, 100], [14, 45, 46, 47, 48, 76], [35, 66, 67, 68, 69, 97], [5, 36, 37, 38, 39, 67], [11, 54, 58, 61, 87, 92], [24, 29, 59, 102, 106, 109], [33, 37, 40, 66, 71, 101], [42, 46, 49, 75, 80, 110], [2, 45, 49, 52, 78, 83], [6, 10, 13, 39, 44, 74], [12, 17, 47, 90, 94, 97], [23, 66, 70, 73, 99, 104], [18, 23, 53, 96, 100, 103], [15, 19, 22, 48, 53, 83], [6, 11, 41, 84, 88, 91], [0, 5, 35, 78, 82, 85], [1, 27, 32, 62, 105, 109], [26, 69, 73, 76, 102, 107], [15, 20, 50, 93, 97, 100], [3, 7, 10, 36, 41, 71], [21, 26, 56, 99, 103, 106], [30, 34, 37, 63, 68, 98], [18, 22, 25, 51, 56, 86], [27, 31, 34, 60, 65, 95], [24, 28, 31, 57, 62, 92], [29, 72, 76, 79, 105, 110], [14, 57, 61, 64, 90, 95], [9, 14, 44, 87, 91, 94], [8, 51, 55, 58, 84, 89], [39, 43, 46, 72, 77, 107], [36, 40, 43, 69, 74, 104], [12, 16, 19, 45, 50, 80], [2, 32, 75, 79, 82, 108], [20, 63, 67, 70, 96, 101], [3, 8, 38, 81, 85, 88], [0, 4, 7, 33, 38, 68], [21, 25, 28, 54, 59, 89], [1, 4, 30, 35, 65, 108], [5, 48, 52, 55, 81, 86], [17, 60, 64, 67, 93, 98], [9, 13, 16, 42, 47, 77], [8, 43, 53, 80, 101, 108], [10, 57, 70, 97, 106, 110], [0, 12, 69, 75, 91, 98], [7, 14, 27, 39, 96, 102], [30, 43, 70, 79, 83, 94], [7, 54, 67, 94, 103, 107], [54, 60, 76, 83, 96, 108], [7, 17, 44, 65, 72, 83], [14, 21, 32, 67, 77, 104], [1, 28, 37, 41, 52, 99], [26, 47, 54, 65, 100, 110], [48, 54, 70, 77, 90, 102], [9, 21, 78, 84, 100, 107], [2, 15, 27, 84, 90, 106], [1, 8, 21, 33, 90, 96], [10, 37, 46, 50, 61, 108], [24, 30, 46, 53, 66, 78], [15, 21, 37, 44, 57, 69], [5, 12, 23, 58, 68, 95], [12, 25, 52, 61, 65, 76], [39, 52, 79, 88, 92, 103], [51, 57, 73, 80, 93, 105], [3, 19, 26, 39, 51, 108], [4, 51, 64, 91, 100, 104], [3, 14, 49, 59, 86, 107], [1, 48, 61, 88, 97, 101], [9, 22, 49, 58, 62, 73], [1, 5, 16, 63, 76, 103], [13, 23, 50, 71, 78, 89], [15, 28, 55, 64, 68, 79], [42, 55, 82, 91, 95, 106], [10, 19, 23, 34, 81, 94], [2, 13, 60, 73, 100, 109], [18, 24, 40, 47, 60, 72], [8, 15, 26, 61, 71, 98], [33, 46, 73, 82, 86, 97], [4, 8, 19, 66, 79, 106], [7, 16, 20, 31, 78, 91], [16, 25, 29, 40, 87, 100], [2, 9, 20, 55, 65, 92], [0, 16, 23, 36, 48, 105], [2, 37, 47, 74, 95, 102], [39, 45, 61, 68, 81, 93], [14, 35, 42, 53, 88, 98], [6, 12, 28, 35, 48, 60], [5, 18, 30, 87, 93, 109], [27, 33, 49, 56, 69, 81], [42, 48, 64, 71, 84, 96], [13, 22, 26, 37, 84, 97], [4, 31, 40, 44, 55, 102], [20, 27, 38, 73, 83, 110], [21, 34, 61, 70, 74, 85], [22, 31, 35, 46, 93, 106], [12, 18, 34, 41, 54, 66], [33, 39, 55, 62, 75, 87], [6, 63, 69, 85, 92, 105], [9, 15, 31, 38, 51, 63], [10, 20, 47, 68, 75, 86], [19, 28, 32, 43, 90, 103], [21, 27, 43, 50, 63, 75], [11, 18, 29, 64, 74, 101], [36, 49, 76, 85, 89, 100], [0, 6, 22, 29, 42, 54], [5, 32, 53, 60, 71, 106], [8, 29, 36, 47, 82, 92], [24, 37, 64, 73, 77, 88], [7, 11, 22, 69, 82, 109], [22, 32, 59, 80, 87, 98], [0, 57, 63, 79, 86, 99], [7, 34, 43, 47, 58, 105], [6, 17, 52, 62, 89, 110], [30, 36, 52, 59, 72, 84], [2, 29, 50, 57, 68, 103], [36, 42, 58, 65, 78, 90], [8, 35, 56, 63, 74, 109], [45, 51, 67, 74, 87, 99], [23, 44, 51, 62, 97, 107], [6, 19, 46, 55, 59, 70], [25, 34, 38, 49, 96, 109], [16, 26, 53, 74, 81, 92], [1, 11, 38, 59, 66, 77], [0, 13, 40, 49, 53, 64], [28, 38, 65, 86, 93, 104], [9, 66, 72, 88, 95, 108], [12, 24, 81, 87, 103, 110], [17, 24, 35, 70, 80, 107], [34, 44, 71, 92, 99, 110], [27, 40, 67, 76, 80, 91], [18, 31, 58, 67, 71, 82], [4, 11, 24, 36, 93, 99], [19, 29, 56, 77, 84, 95], [20, 41, 48, 59, 94, 104], [3, 16, 43, 52, 56, 67], [13, 20, 33, 45, 102, 108], [25, 35, 62, 83, 90, 101], [31, 41, 68, 89, 96, 107], [3, 15, 72, 78, 94, 101], [1, 10, 14, 25, 72, 85], [4, 14, 41, 62, 69, 80], [5, 40, 50, 77, 98, 105], [3, 9, 25, 32, 45, 57], [5, 26, 33, 44, 79, 89], [17, 38, 45, 56, 91, 101], [45, 58, 85, 94, 98, 109], [0, 11, 46, 56, 83, 104], [4, 13, 17, 28, 75, 88], [10, 17, 30, 42, 99, 105], [3, 60, 66, 82, 89, 102], [2, 23, 30, 41, 76, 86], [11, 32, 39, 50, 85, 95], [6, 18, 75, 81, 97, 104], [9, 54, 69, 79, 93, 101], [23, 38, 74, 80, 82, 90], [9, 19, 33, 41, 60, 105], [3, 13, 27, 35, 54, 99], [11, 26, 62, 68, 70, 78], [15, 30, 40, 54, 62, 81], [16, 64, 70, 75, 89, 109], [6, 16, 30, 38, 57, 102], [16, 34, 82, 88, 93, 107], [16, 22, 27, 41, 61, 79], [9, 24, 34, 48, 56, 75], [3, 48, 63, 73, 87, 95], [33, 48, 58, 72, 80, 99], [5, 7, 15, 59, 74, 110], [1, 15, 23, 42, 87, 102], [43, 49, 54, 68, 88, 106], [4, 9, 23, 43, 61, 109], [4, 22, 70, 76, 81, 95], [0, 45, 60, 70, 84, 92], [27, 42, 52, 66, 74, 93], [6, 14, 33, 78, 93, 103], [11, 17, 19, 27, 71, 86], [0, 10, 24, 32, 51, 96], [34, 40, 45, 59, 79, 97], [26, 32, 34, 42, 86, 101], [19, 25, 30, 44, 64, 82], [20, 26, 28, 36, 80, 95], [35, 41, 43, 51, 95, 110], [1, 49, 55, 60, 74, 94], [14, 50, 56, 58, 66, 110], [2, 8, 10, 18, 62, 77], [17, 23, 25, 33, 77, 92], [39, 54, 64, 78, 86, 105], [4, 52, 58, 63, 77, 97], [2, 17, 53, 59, 61, 69], [35, 50, 86, 92, 94, 102], [5, 11, 13, 21, 65, 80], [20, 35, 71, 77, 79, 87], [14, 29, 65, 71, 73, 81], [12, 20, 39, 84, 99, 109], [8, 23, 59, 65, 67, 75], [1, 6, 20, 40, 58, 106], [5, 20, 56, 62, 64, 72], [7, 25, 73, 79, 84, 98], [12, 27, 37, 51, 59, 78], [40, 46, 51, 65, 85, 103], [28, 34, 39, 53, 73, 91], [0, 8, 27, 72, 87, 97], [38, 53, 89, 95, 97, 105], [0, 15, 25, 39, 47, 66], [10, 28, 76, 82, 87, 101], [30, 45, 55, 69, 77, 96], [10, 58, 64, 69, 83, 103], [12, 22, 36, 44, 63, 108], [14, 20, 22, 30, 74, 89], [6, 51, 66, 76, 90, 98], [13, 19, 24, 38, 58, 76], [11, 31, 49, 97, 103, 108], [18, 63, 78, 88, 102, 110], [29, 44, 80, 86, 88, 96], [0, 44, 59, 95, 101, 103], [6, 21, 31, 45, 53, 72], [2, 22, 40, 88, 94, 99], [1, 7, 12, 26, 46, 64], [13, 61, 67, 72, 86, 106], [5, 41, 47, 49, 57, 101], [18, 33, 43, 57, 65, 84], [25, 31, 36, 50, 70, 88], [37, 43, 48, 62, 82, 100], [29, 35, 37, 45, 89, 104], [0, 14, 34, 52, 100, 106], [31, 37, 42, 56, 76, 94], [2, 21, 66, 81, 91, 105], [21, 36, 46, 60, 68, 87], [3, 17, 37, 55, 103, 109], [2, 4, 12, 56, 71, 107], [8, 44, 50, 52, 60, 104], [11, 47, 53, 55, 63, 107], [3, 47, 62, 98, 104, 106], [46, 52, 57, 71, 91, 109], [19, 37, 85, 91, 96, 110], [7, 21, 29, 48, 93, 108], [17, 32, 68, 74, 76, 84], [32, 38, 40, 48, 92, 107], [10, 16, 21, 35, 55, 73], [12, 57, 72, 82, 96, 104], [2, 38, 44, 46, 54, 98], [3, 11, 30, 75, 90, 100], [42, 57, 67, 81, 89, 108], [23, 29, 31, 39, 83, 98], [4, 18, 26, 45, 90, 105], [7, 13, 18, 32, 52, 70], [5, 24, 69, 84, 94, 108], [7, 55, 61, 66, 80, 100], [41, 56, 92, 98, 100, 108], [8, 28, 46, 94, 100, 105], [3, 18, 28, 42, 50, 69], [6, 50, 65, 101, 107, 109], [5, 25, 43, 91, 97, 102], [13, 31, 79, 85, 90, 104], [15, 60, 75, 85, 99, 107], [22, 28, 33, 47, 67, 85], [26, 41, 77, 83, 85, 93], [1, 9, 53, 68, 104, 110], [36, 51, 61, 75, 83, 102], [9, 17, 36, 81, 96, 106], [4, 10, 15, 29, 49, 67], [1, 19, 67, 73, 78, 92], [32, 47, 83, 89, 91, 99], [24, 39, 49, 63, 71, 90], [8, 14, 16, 24, 68, 83], [8, 22, 45, 64, 107, 110], [23, 26, 35, 49, 72, 91], [33, 42, 59, 63, 83, 100], [11, 28, 72, 81, 98, 102], [3, 22, 65, 68, 77, 91], [14, 19, 36, 54, 97, 109], [6, 23, 27, 47, 64, 108], [10, 54, 63, 80, 84, 104], [6, 25, 68, 71, 80, 94], [9, 52, 64, 80, 85, 102], [7, 50, 53, 62, 76, 99], [0, 19, 62, 65, 74, 88], [40, 52, 68, 73, 90, 108], [18, 37, 80, 83, 92, 106], [11, 14, 23, 37, 60, 79], [42, 51, 68, 72, 92, 109], [3, 23, 40, 84, 93, 110], [41, 44, 53, 67, 90, 109], [1, 44, 47, 56, 70, 93], [0, 20, 37, 81, 90, 107], [4, 47, 50, 59, 73, 96], [0, 17, 21, 41, 58, 102], [10, 22, 38, 43, 60, 78], [12, 55, 67, 83, 88, 105], [25, 37, 53, 58, 75, 93], [30, 39, 56, 60, 80, 97], [1, 24, 43, 86, 89, 98], [9, 18, 35, 39, 59, 76], [3, 20, 24, 44, 61, 105], [9, 28, 71, 74, 83, 97], [13, 57, 66, 83, 87, 107], [15, 58, 70, 86, 91, 108], [27, 36, 53, 57, 77, 94], [13, 36, 55, 98, 101, 110], [5, 19, 42, 61, 104, 107], [3, 21, 64, 76, 92, 97], [37, 49, 65, 70, 87, 105], [5, 8, 17, 31, 54, 73], [1, 13, 29, 34, 51, 69], [8, 12, 32, 49, 93, 102], [8, 13, 30, 48, 91, 103], [29, 32, 41, 55, 78, 97], [12, 21, 38, 42, 62, 79], [2, 11, 25, 48, 67, 110], [5, 22, 66, 75, 92, 96], [18, 27, 44, 48, 68, 85], [4, 16, 32, 37, 54, 72], [1, 17, 22, 39, 57, 100], [26, 29, 38, 52, 75, 94], [36, 45, 62, 66, 86, 103], [20, 23, 32, 46, 69, 88], [39, 48, 65, 69, 89, 106], [38, 41, 50, 64, 87, 106], [16, 59, 62, 71, 85, 108], [2, 6, 26, 43, 87, 96], [7, 30, 49, 92, 95, 104], [22, 34, 50, 55, 72, 90], [9, 27, 70, 82, 98, 103], [17, 34, 78, 87, 104, 108], [14, 17, 26, 40, 63, 82], [14, 18, 38, 55, 99, 108], [17, 20, 29, 43, 66, 85], [34, 46, 62, 67, 84, 102], [1, 18, 36, 79, 91, 107], [12, 31, 74, 77, 86, 100], [0, 9, 26, 30, 50, 67], [19, 31, 47, 52, 69, 87], [6, 24, 67, 79, 95, 100], [7, 51, 60, 77, 81, 101], [0, 18, 61, 73, 89, 94], [6, 15, 32, 36, 56, 73], [2, 7, 24, 42, 85, 97], [5, 9, 29, 46, 90, 99], [2, 5, 14, 28, 51, 70], [21, 30, 47, 51, 71, 88], [4, 48, 57, 74, 78, 98], [16, 28, 44, 49, 66, 84], [15, 34, 77, 80, 89, 103], [3, 46, 58, 74, 79, 96], [5, 10, 27, 45, 88, 100], [11, 15, 35, 52, 96, 105], [1, 45, 54, 71, 75, 95], [16, 60, 69, 86, 90, 110], [12, 30, 73, 85, 101, 106], [4, 20, 25, 42, 60, 103], [7, 19, 35, 40, 57, 75], [7, 23, 28, 45, 63, 106], [28, 40, 56, 61, 78, 96], [2, 16, 39, 58, 101, 104], [24, 33, 50, 54, 74, 91], [2, 19, 63, 72, 89, 93], [31, 43, 59, 64, 81, 99], [6, 49, 61, 77, 82, 99], [8, 11, 20, 34, 57, 76], [13, 25, 41, 46, 63, 81], [15, 24, 41, 45, 65, 82], [4, 21, 39, 82, 94, 110], [3, 12, 29, 33, 53, 70], [13, 56, 59, 68, 82, 105], [4, 27, 46, 89, 92, 101], [10, 26, 31, 48, 66, 109], [0, 43, 55, 71, 76, 93], [10, 53, 56, 65, 79, 102], [35, 38, 47, 61, 84, 103], [11, 16, 33, 51, 94, 106], [8, 25, 69, 78, 95, 99], [21, 40, 83, 86, 95, 109], [32, 35, 44, 58, 81, 100], [15, 33, 76, 88, 104, 109], [10, 33, 52, 95, 98, 107], [14, 31, 75, 84, 101, 105]]
\end{enumerate}

\section{Acknowledgement}

I would like to express my thanks to Taras Banakh, Alex Ravsky and all participants of Lviv offline Geometry meetings who verified my ideas which later led to more and more generalized algorithms for Steiner system discoveries.

\end{document}